\date{}
\documentclass [11pt]{article}
\usepackage{amsfonts,amsmath,amssymb,amsthm,graphicx,makeidx}

\title{A hypergraph blow-up lemma}
\author{
Peter Keevash \thanks{School of Mathematical Sciences,
Queen Mary, University of London, Mile End Road, London E1 4NS, UK.
Email: p.keevash@qmul.ac.uk.
Research supported in part by
ERC grant 239696 and
EPSRC grant EP/G056730/1.} %NSF grant DMS-0555755
}

\theoremstyle{plain}
\newtheorem{theo}{Theorem}[section] % remove section if desired

\newtheorem{lemma}[theo]{Lemma}
\newtheorem{coro}[theo]{Corollary}

\newtheorem{defn}[theo]{Definition}

\theoremstyle{definition}
\newtheorem{rem}[theo]{Remark}
\newtheorem{eg}[theo]{Example}
\newcommand{\mc}[1]{\mathcal{#1}}
\newcommand{\mb}[1]{\mathbb{#1}}
\newcommand{\nib}[1]{\noindent {\bf #1}}

\newcommand{\Lra}{\Leftrightarrow}
\newcommand{\sm}{\setminus}
\newcommand{\ov}{\overline}
\newcommand{\eps}{\epsilon}
\newcommand{\GG}{\Gamma}
\newcommand{\Sa}{\Sigma}

\newcommand{\gens}[1]{\langle #1 \rangle}

\newcommand{\sub}{\subseteq}
\newcommand{\subn}{\subsetneq}
\newcommand{\nsub}{\nsubseteq}

\newcommand{\ul}{\underline}
\newcommand{\wt}{\widetilde}
\newcommand{\es}{\emptyset}
\newcommand{\ca}{\circledast}

\def\COMMENT#1{}

\def\qed{\hfill $\Box$}

\makeindex

\begin{document}
\maketitle

\begin{abstract}
We obtain a hypergraph generalisation of the graph blow-up lemma proved by
Koml\'os, Sark\"ozy and Szemer\'edi, showing that hypergraphs with sufficient
regularity and no atypical vertices behave as if they were complete for
the purpose of embedding bounded degree hypergraphs.%
\COMMENT{Formerly...  sufficiently regular hypergraphs
with no atypical vertices behave like complete partite hypergraphs for the
purpose of embedding bounded degree hypergraphs. In the course of our
arguments we also obtain various useful lemmas concerning hypergraph
regularity that have independent interest, including a characterisation
in terms of the frequency of certain subcomplexes.
There are many potential applications of our theorem to hypergraph generalisations
of results for graphs that were obtained with the blow-up lemma.
We illustrate the method with a hypergraph generalisation of
a result of K\"uhn and Osthus on packing bipartite graphs.
}
\end{abstract}

\section{Introduction}

% algorithmic version ?

Szemer\'edi's regularity lemma \cite{S1} has impressive applications in many areas
of modern graph theory, including extremal graph theory, Ramsey theory and property testing.
Roughly speaking, it says that any graph can be approximated by an average
with respect to a partition of its vertex set into a bounded number of
classes, the number of classes depending only on the accuracy of the
desired approximation, and not on the number of vertices in the
graph. A key property of this approximation is that it leads to a
`counting lemma', allowing an accurate prediction of the number of
copies of any small fixed graph spanned by some specified classes of the
partition. We refer the reader to \cite{KS} for a survey of the
regularity lemma and its applications. An analogous theory for
hypergraphs has only been developed very recently, with independent
and rather different approaches given by R\"odl et
al. (e.g.\ \cite{NRS, RSc1, RSk}) and Gowers \cite{G2},
subsequently reformulated and developed in \cite{Au, ES, I, R, T1, T2}.
In a very short space of time the power of this hypergraph theory has
already been amply demonstrated, e.g.\ by a multidimensional generalisation of
Szemer\'edi's theorem on arithmetic progressions \cite{G2} and a
linear bound for the Ramsey number of hypergraphs with bounded maximum
degree \cite{CFKO,NORS}.

The blow-up lemma is a powerful tool developed by Koml\'os, Sark\"ozy
and Szemer\'edi \cite{KSS} for using the regularity lemma to embed
spanning subgraphs of bounded degree. An informal statement is that graphs
with sufficient regularity and no atypical vertices behave as if they were complete for
the purpose of embedding bounded degree graphs. In \cite{KSS2, KSS3} they used
it to prove Seymour's conjecture on the minimum degree needed to embed
the $k$th power of a Hamilton cycle, and the Alon-Yuster conjecture on
the minimum degree needed for a graph to have an $H$-factor, for some
fixed graph $H$ (a question finally resolved in \cite{KO3}).
There are many other applications to embedding spanning subgraphs,
see the survey \cite{KO4}. There are also several results on embedding
spanning subhypergraphs, such as perfect matchings or (various definitions of)
Hamilton cycles, see the survey \cite{RR10}. For the most part, the proofs of
the known hypergraph results have not needed any analogue of the blow-up lemma.
An exception is an embedding lemma for some special spanning subhypergraphs
proved in \cite{KO1} for loose Hamilton cycles in $3$-graphs, although
the `absorbing' method of R\"odl, Ruci\'nski and Szemer\'edi was
subsequently shown to be a simpler method for Hamilton cycles \cite{RRS,HS,KMO}.
Another partial hypergraph blow-up lemma is to embedding bounded degree
subgraphs of linear size, obtained independently in \cite{CFKO}
and in \cite{NORS} (for $3$-graphs). The application to linear Ramsey numbers
of bounded degree hypergraphs was also subsequently proved by simpler means in \cite{CFS}.
However, one would not necessarily expect that embedding lemmas can always
be avoided by using alternative methods, so a hypergraph blowup lemma
would be a valuable tool.%

In this paper we prove such a result that gives conditions
for embedding any spanning hypergraph of bounded degree.
We will not attempt a formal statement of our result in this introduction,
as it will take us a considerable amount of work to set up the necessary
notation and terminology, particularly for the key notion of super-regularity,
which has some additional subtleties that do not appear in the graph case.
The proof will be by means of a randomised greedy embedding algorithm,
which is very similar that used in \cite{KSS}. However the analysis
is more involved, due both to the additional complications of hypergraph
regularity theory, and the need to work with an approximating hypergraph
rather than the true hypergraph (see Section 3 for further explanation).
There are many potential applications of our theorem to
hypergraph generalisations of results for graphs that were obtained with the
graph blow-up lemma. We will illustrate the method by proving a hypergraph generalisation
of a result of K\"uhn and Osthus \cite[Theorem 2]{KO2} on packing bipartite graphs.

The rest of this paper is organised as follows.
In the next section we prove the blow-up lemma of Koml\'os, Sark\"ozy
and Szemer\'edi \cite{KSS}. This is mostly for expository purposes,
although there are some small differences in our proof,
and it will be useful to refer back to the basic argument
when discussing additional complications that arise for hypergraphs.
We will prove our main result at first in the special case of $3$-graphs
(with some additional simplifications); this case is already sufficiently
complex to illustrate the main ideas of our proofs.
In section 3 we discuss hypergraph regularity theory
(following the approach of R\"odl et al.)
and motivate and define super-regularity for $3$-graphs.
We prove the $3$-graph blow-up lemma in section 4.
In section 5 we develop some additional theory that is needed
for applications of the $3$-graph blow-up lemma,
based on the Regular Approximation Lemma of R\"odl and Schacht \cite{RSc1}.
We also give a `black box' reformulation of the blow-up lemma
that will be more easily accessible for future applications.
We illustrate this by generalising a result in \cite{KO2}
to packing tripartite $3$-graphs.
The final section concerns the general hypergraph blow-up lemma.
As well as generalising from $3$-graphs to $k$-graphs,
we allow additional generalisations that will be needed in future
applications, including restricted positions and complex-indexed
complexes (defined in that section).
The proof is mostly similar to that for $3$-graphs,
so we only give full details for those aspects that are different.
We conclude with some remarks on potential developments
and applications of the blow-up lemma.

\medskip

\nib{Notation.}\
We will introduce a substantial amount of terminology and notation
throughout the paper, which is summarised in the index.
Before starting our discussion we establish the following basic notation.
We write $[n]=\{1,\cdots,n\}$. \index{$[n]$}
If $X$ is a set and $k$ is a number
then $\binom{X}{k} = \{Y \sub X: |Y|=k\}$, \index{$\binom{X}{k}$}
$\binom{X}{\le k} = \cup_{i \le k} \binom{X}{i}$ and \index{$\binom{X}{\le k}$}
$\binom{X}{< k} = \cup_{i < k} \binom{X}{i}$. \index{$\binom{X}{< k}$}
$a \pm b$ denotes an unspecified real number in the interval $[a-b,a+b]$. \index{$a \pm b$}
It is convenient to
regard a finite set $X$ as being equipped with
the uniform probability measure $\mb{P}(\{x\})=1/|X|$,
so that we can express the average of a function $f$ defined on
$X$ as $\mb{E}_{x \in X} f(x)$. \index{$E$@$\mb{E}$}
A {\em $k$-graph} $H$ consists of a \index{k-graph@$k$-graph}
{\em vertex set} \index{vertex set} $V(H)$ \index{$V(H)$}
and an {\em edge set} \index{edge set} $E(H)$, \index{$E(H)$}
each edge being some $k$-tuple of vertices.
We often identify $H$ with $E(H)$, so that $|H|$ is the
number \index{$|H|$} of edges in $H$.
A {\em $k$-complex} $H$ consists of a vertex set $V(H)$ and an edge set $E(H)$,
\index{k-complex@$k$-complex}
where each edge is a subset of the vertex set of size at most $k$,
that is a {\em simplicial complex}, i.e.\ if \index{simplicial complex}
$S \in E(H)$ and $T \subset S$ then $T \in E(H)$.
For $S \subset V(H)$ the {\em neighbourhood} $(k-|S|)$-graph or $(k-|S|)$-complex
\index{neighbourhood}\index{$(\cdot)$|see{neighbourhood}}
is $H(S) = \{A \subset V(H) \sm S: A \cup S \in E(H) \}$, \index{$H(S)$}
and $|H(S)|$ is the {\em degree} \index{degree} of $S$. We also write \index{$H^S$}
$H^S = \{A \sub V(H): S \sub A \in E(H)\}$. A {\em walk} \index{walk} in $H$ is a sequence
of vertices for which each consecutive pair are contained in some
edge of $H$, and the distance between two vertices is the length
of the shortest walk connecting them.
The {\em vertex neighbourhood} $VN_H(x)$ is the set of vertices%
\index{vertex neighbourhood}\index{$VN_H(x)$}
at distance exactly $1$ from $x$ (so $x$ itself is not included).
We will often have to consider
hierarchies involving many real parameters, and it will be useful to use
the notation $0 < \alpha \ll \beta$ to mean that there is an increasing function \index{$\ll$}
$f(x)$ so that the following argument is valid for $0 < \alpha < f(\beta)$.
The parameter $n$ will always be sufficiently large compared to all other
parameters, and we use the phrase {\em with high probability} to refer to an
\index{with high probability} event that has probability $1-o_n(1)$,
i.e.\ the probability tends to $1$ as $n$ tends to infinity.
\COMMENT{Formerly (incompatible later):
If $H$ is a $k$-graph or $k$-complex and $X \subset V(H)$ the
restriction $H[X]$ is a $k$-graph or $k$-complex on $X$ in which
the edges are all edges of $H$ that are contained in $X$.
}

\section{The graph blow-up lemma}

This section is mostly expository. We introduce the basic notions of regularity
and super-regularity for graphs and prove the blow-up lemma of
Koml\'os, Sark\"ozy and Szemer\'edi. This will serve as a warm-up to the hypergraph
blow-up lemma, as our proof even in the graph case differs slightly from the original
in a few details (although the general approach is the same). It will also be helpful
to establish our notation in this simplified setting, and to refer back to the basic
argument when explaining why certain extra complications arise for hypergraphs.
To streamline the proof we focus on a slightly simplified setting, which still contains
all the ideas needed for the general case. We hope that the general reader will find this
section to be an accessible account of a proof that has a reputation for difficulty!

We start with a brief summary of the key notions in graph regularity,
referring the reader to \cite{KS} for more details.
Consider an $r$-partite graph $G$ with vertex set $V$ partitioned as $V = V_1 \cup \cdots \cup V_r$.
Let $G_{ij}$ be the bipartite subgraph of $G$ with parts $V_i$ and $V_j$,
for $1 \le i \ne j \le r$. The {\em density}\index{density}
of $G_{ij}$ is $d(G_{ij}) = \frac{|G_{ij}|}{|V_i||V_j|}$.
Given $\eps>0$, we say that $G_{ij}$ is {\em $\eps$-regular} if%
\index{regular}\index{$\eps$-regular|see{regular}}
for all subsets $V'_i \sub V_i$ and $V'_j \sub V_j$
with $|V'_i|>\eps|V_i|$ and $|V'_j|>\eps|V_j|$,
writing $G'_{ij}$ for the bipartite subgraph of $G$ with parts $V'_i$ and $V'_j$,
we have $|d(G'_{ij})-d(G_{ij})|<\eps$.
Then we say that $G$ is $\eps$-regular if each $G_{ij}$ is $\eps$-regular.
Informally, we may say that each $G_{ij}$ behaves like a random bipartite graph,
up to accuracy $\eps$. This statement is justified by the counting lemma,
which allows one to estimate the number of copies of any fixed graph $F$,
up to accuracy $O(\eps)$, using a suitable product of densities.
For now we just give an example: if we write $T_{123}(G)$ for the
\index{triangles}\index{$T_{123}(G)$}
set of triangles formed by the graphs $G_{12}$, $G_{13}$, $G_{23}$, then
\begin{equation}\label{eq:tri}
d(T_{123}(G)):=\frac{|T_{123}(G)|}{|V_1||V_2||V_3|} = d(G_{12})d(G_{13})d(G_{23}) \pm 8\eps.
\end{equation}
Remarkably, this powerful property can be applied in any graph $G$,
via Szemer\'edi's Regularity Lemma, which can be informally stated as saying that
we can decompose the vertex set of any graph on $n$ vertices
into $m(\eps)$ parts, such that all but at most $\eps n^2$ edges
belong to bipartite subgraphs that are $\eps$-regular.
\COMMENT{Stick to absolute density notation for graphs.}

The blow-up lemma arises from the desire to embed {\em spanning} graphs in $G$,\index{spanning}
meaning that they use every vertex in $V$. Suppose that $|V_i|=n$ for $1 \le i \le r$.
The argument used to prove the counting lemma can be generalised to embed any bounded
degree graph $H$ provided that all components of $H$ have size $o(n)$,
and $\Omega(n)$ vertices of $G$ are allowed to remain uncovered.
However, one cannot guarantee an embedding of a spanning graph:
the definition of $\eps$-regularity does not prevent the existence
of isolated vertices, so we may not even be able to find a perfect matching.
This observation naturally leads us to the stronger notion of super-regularity.
We say that $G_{ij}$ is {\em $(\eps,d_{ij})$-super-regular} if it is $\eps$-regular
and every vertex has degree at least $(d_{ij}-\eps)n$. It is well-known that
one can delete a small number of vertices from a regular pair
to make it super-regular (see Lemma \ref{2del}).
We say that $G$ is {\em $(\eps,d)$-super-regular} if each $G_{ij}$ is either empty or
$(\eps,d_{ij})$-super-regular for some $d_{ij} \ge d$.\index{super-regular}
Now we can state the graph blow-up lemma.

\begin{theo} \label{2blowup} {\bf (Graph blow-up lemma)} \index{graph blow-up lemma}
Suppose $H$ is an $r$-partite graph on $X = X_1 \cup \cdots \cup X_r$
and $G$ is an $r$-partite graph on $V = V_1 \cup \cdots \cup V_r$,
where $|V_i|=|X_i|=n$ for $1 \le i \le r$
and $H_{ij}$ is only non-empty when $G_{ij}$ is non-empty.
If $H$ has maximum degree at most $D$ and $G$ is $(\eps,d)$-super-regular,
where $0 \ll 1/n \ll \eps \ll d, 1/r, 1/D$,
then $G$ contains a copy of $H$, in which for each $1 \le i \le r$
the vertices of $V_i$ correspond to the vertices of $X_i$.
\COMMENT{
comparison with KSS: notation, differences in analysis: e.g.\ distinguish F from C,
no queue in initial phase for N (unlike hyp: can have buffer smaller than d),
assume equal parts (easy mod to constant ratio)
}
\end{theo}

Informally speaking, Theorem \ref{2blowup} embeds any bounded degree graph into
any super-regular graph. Note that arbitrary part sizes are allowed in \cite{KSS},
but for simplicity we start by considering the case when they are all equal.
The proof is via a random greedy algorithm
for embedding $H$ in $G$, which considers the vertices of $X$
in some order and embeds them to $V$ one at a time. We start by
giving an informal description of the algorithm.%
\index{$H$}\index{$G$}\index{$X$}\index{$V$}

\begin{description}

\item[Initialisation.]
List the vertices of $H$ in a certain order, as follows. Some
vertices at mutual distance at least $4$ are identified as buffer
vertices and put at the end of the list. The neighbours
of the buffer vertices are put at the start of the list.
(The rationale for this order is that we hope to embed these neighbours in a nice manner
while there is still plenty of room in the early stages of the algorithm,
and then the buffer vertices still have many suitable places at the conclusion.)
During the algorithm a queue of priority vertices may arise:
it is initially empty.

\item[Iteration.] Choose the next vertex $x$ to be embedded,
either from the queue if this is non-empty, or otherwise from the list.
The image $\phi(x)$ of $x$ is chosen randomly in $V(G)$
among those free spots that do not unduly restrict the free\index{free}
spots for those unembedded neighbours of $x$.
If some unembedded vertex has too few free spots it is added to the
queue. Stop when all non-buffer vertices have been embedded. If the
number of vertices that have ever been in the queue becomes too large
before this point then abort the algorithm as a failure (this is
an unlikely event).

\item[Conclusion.] Choose a system of distinct representatives among the
free slots for the unembedded vertices to complete the embedding.
(This will be possible with high probability.)

\end{description}

Now we will formally describe the random greedy algorithm to construct an
embedding $\phi:V(H) \to V(G)$ such that $\phi(e) \in E(G)$ for every $e \in E(H)$. \index{$\phi$}
First we introduce more parameters with the hierarchy
\[0 \ll 1/n \ll \eps \ll \eps' \ll \eps_*
\ll p_0 \ll \gamma \ll \delta_Q \ll p \ll d_u
\ll \delta'_Q \ll \delta_B \ll d, 1/r, 1/D.\]
To assist the reader we list here the role of each parameter for easy reference.
Parameters $\eps$, $\eps'$ and $\eps_*$ are used to measure graph regularity.
Parameter $\gamma$ plays the role of $\kappa$ in \cite{KSS}: \index{$\gamma$}
it is used to distinguish various cases at the conclusion of the algorithm
when selecting the system of distinct representatives.
Parameter $d_u$ plays the role of $\gamma$ in \cite{KSS}: \index{$d_u$}
it is a universal lower bound on the proportion of vertices in a class of $G$
free to embed any given vertex of $H$.
The {\em queue threshold}\index{queue threshold} parameter $\delta_Q$\index{$\delta_Q$}
corresponds to $\delta'''$ in \cite{KSS}:
the maximum proportional size for the queue before we will abort with failure.
The {\em buffer}\index{buffer parameter} parameter $\delta_B$\index{$\delta_B$}
corresponds to $\delta'$ in \cite{KSS}: the proportional size of the buffer.
We also introduce two probability parameters that are not explicitly named in \cite{KSS},
although they are key to the proof.
Parameter $p_0$ appears in Lemma \ref{2main} in the upper bound\index{$p_0$}
for the probability that any given set $A$ will be significantly under-represented
in the free images for a vertex.
Parameter $p$ appears in Lemma \ref{2-x-to-v} as a lower bound for\index{$p$}
the probability that a given unused vertex will be free as an image
for a given buffer vertex at the conclusion of the algorithm.
Note that the {\em queue admission}\index{queue admission parameter}
parameter $\delta'_Q$\index{$\delta'_Q$} is similar to but slightly different
from the corresponding parameter $\delta''$ in \cite{KSS}:
for any vertex $z$ and time $t$ we will compare $F_z(t)$ to an earlier free set $F_z(t_z)$,
where $t_z$ is the most recent time at which we embedded a neighbour of $z$.

\begin{description}

\item[Initialisation and notation.]

We choose a buffer\index{buffer} set\index{$B$} $B \subset X$
of vertices at mutual distance at least $4$ in $H$
so that $|B \cap X_i| = \delta_B n$ for $1 \le i \le r$. The maximum degree property
of $H$ implies that we can construct $B$ simply by selecting vertices one-by-one greedily.
For any given vertex in $H$ there are fewer than $D^4$ vertices within distance $4$,
so at any point in the construction of $B$ we have excluded at most $D^4r\delta_B n$ vertices
from any given $X_i$. Thus we can construct $B$ if we choose $\delta_B < 1/(rD^4)$.

Let $N = \cup_{x \in B} N_H(x)$ be the vertices with a neighbour in the buffer. \index{$N$}
Since $H$ has maximum degree $D$ we have $|N \cap X_i| \le Dr\delta_Bn < \sqrt{\delta_B} n$
for $1 \le i \le r$, if we choose $\delta_B < 1/(Dr)^{2}$.

We use $t$ to denote time during the algorithm, by which we mean the number of vertices
of $H$ that have been embedded. At time $t$ we denote the queue by $q(t)$ and write
$Q(t) = \cup_{u \le t}\ q(u)$ for the vertices that have ever been in the queue by time $t$.
Initially we set $q(0)=Q(0)=\es$.\index{queue}\index{$q(t)$}\index{$Q(t)$}

We order the vertices in a list\index{list} $L=L(0)$\index{$L(t)$}
that starts with $N$ and ends with $B$.
Within $N$, we arrange that $N_H(x)$ is consecutive for each $x \in B$.
This is possible by the mutual distance property in $B$, which implies that
the neighbourhoods $N_H(x)$, $x \in B$ are mutually disjoint.
We denote the vertex of $H$ selected for embedding at time $t$ by $s(t)$.\index{$s(t)$}
This will be the first vertex of $L(t-1)$,
unless the queue is non-empty, when this takes priority.

We write $F_x(t)$ for the vertices that are {\em free}\index{free}
to embed a given vertex $x$ of $H$. \index{$F_x(t)$}
Initially we set $F_x(0)=V_x$, where we write $V_x$ for that part $V_i$ of $G$ corresponding \index{$V_x$}
to the part $X_i$ of $H$ that contains $x$.
We also write $X_i(t) = X_i \sm \{s(u): u \le t\}$ for the unembedded vertices of $X_i$ \index{$X_i(t)$}
and $V_i(t) = V_i \sm \{\phi(s(u)): u \le t\}$ for the available positions in $V_i$. \index{$V_i(t)$}
We let $X(t) = \cup_{i=1}^r X_i(t)$ and $V(t) = \cup_{i=1}^r V_i(t)$.
\index{$X(t)$}\index{$V(t)$}
Initially we set $X_i(0)=X_i$ and $V_i(0)=V_i$.
\COMMENT{Don't use $\nu_x$ or $\nu_S$}

\item[Iteration.] At time $t$, while there are still some
unembedded non-buffer vertices, we select a vertex to embed $x=s(t)$
according to the following {\em selection rule}. If the queue $q(t)$ is non-empty
then we let $x$ be any member of the queue; otherwise we let $x$ be the first
vertex of the list $L(t-1)$. (A First In First Out rule for the queue is used
in \cite{KSS}, but this is not essential to the proof.)
We choose the image $\phi(x)$ of $x$ uniformly at random among all elements
$y \in F_x(t-1)$ that are `good', a property that can be informally stated \index{good}
as saying that if we set $\phi(x)=y$ then the free sets at time $t$
for the unembedded neighbours of $x$ will have roughly their `expected' size.

To define this formally, we first need to describe the update rule for the free sets
when we embed $x$ to some vertex $y$. First we set $F_x(t)=\{y\}$. We will have
$F_x(t')=\{y\}$  at all subsequent times $t' \ge t$.
Then for any unembedded $z$ that is not a neighbour of $x$ we set $F_z(t)=F_z(t-1) \sm \{y\}$.
Thus the size of $F_z(t)$ decreases by $1$ if $z$ belongs to the same part of $X$ as $x$,
but otherwise is unchanged.
Finally, for any unembedded $z \in N_H(x)$ we set $F_z(t)=F_z(t-1) \cap N_G(y)$.
Now we say that $y \in F_x(t-1)$ belongs to the {\em good} set $OK_x(t-1)$\index{$OK$}
if for every unembedded $z \in N_H(x)$ we have $|F_z(t)| = (1 \pm 2\eps')d_{xz}|F_z(t-1)|$.
Here $d_{xz}$ denotes that density $d(G_{ij})$ for which $x \in V_i$ and $z \in V_j$.
\COMMENT{Don't use $E_x$}

Having chosen the image $\phi(x)$ of $x$ as a random good element $y$,
we conclude the iteration by updating the list $L(t-1)$ and the queue $q(t-1)$.
First we remove $x$ from whichever of these sets it was taken.
Then we add to the queue any unembedded vertex $z$ for which $F_z(t)$ has become `too small'.
To make this precise, suppose $z \in L(t-1) \sm \{x\}$,
and let $t_z$ be the most recent time at which we embedded
a vertex in $N_H(z)$, or $0$ if there is no such time.
(Note that if $z \in N_H(x)$ then $t_z=t$.)
We add $z$ to $q(t)$ if $|F_z(t)| < \delta'_Q|F_z(t_z)|$.
This defines $L(t)$ and $q(t)$.
\COMMENT{
Former mistake had rule  $|F_z(t)| < \delta'_Q|V_z(t)|$.
Lower bound $|F_z(t)|/|V_z(t)|$ has to depend on $\delta'_Q$,
so cannot directly deduce anything on $|A \cap F_z(t)|/|F_z(t)|$.
Can we fix this approach? Don't see how...
}

Repeat this iteration until the only unembedded vertices are buffer vertices,
but abort with failure if at any time we have
$|Q(t) \cap X_i| > \delta_Q|X_i|$ for some $1 \le i \le r$.
Let $T$ denote the time at which the iterative phase terminates\index{$T$}
(whether with success or failure).

\item[Conclusion.] When all non-buffer vertices have been embedded,
we choose a system of distinct representatives among the free slots $F_x(T)$
for the unembedded vertices $x \in X(T)$ to complete the embedding,
ending with success if this is possible, otherwise aborting with failure.
%either ending with success if this is possible or aborting with failure if it is not possible.

\end{description}

Now we analyse the algorithm described above and show that it is successful with high probability.
We start by recording two standard facts concerning graph regularity.
The first fact states that most vertices in a regular pair have `typical' degree,
and the second that regularity is preserved by restriction to induced subgraphs.
We maintain our notation that $G$ is an $r$-partite graph on $V = V_1 \cup \cdots \cup V_r$.
We fix any pair $(i,j)$, write $G_{ij}$ for the bipartite subgraph spanned by $V_i$ and $V_j$
and denote its density by $d_{ij}$. We give the short proofs of these facts here, \index{$d_{ij}$}
both for completeness and as preparation for similar hypergraph arguments later.

\begin{lemma}\label{2neighbour} {\bf (Typical degrees)} \index{typical degrees}
Suppose $G_{ij}$ is $\eps$-regular. Then all but at most $2\eps|V_i|$ vertices in $V_i$
have degree $(d_{ij} \pm \eps)|V_j|$ in $V_j$.
\end{lemma}

\nib{Proof.}
We claim that there is no set $X \sub V_i$ of size $|X|>\eps|V_i|$
such that every $x \in X$ has degree less than $(d_{ij}-\eps)|V_j|$ in $V_j$.
For the pair $(X,V_j)$ would then induce a subgraph of density less than $d_{ij}-\eps$,
contradicting the definition of $\eps$-regularity. Similarly,
there is no set $X \sub V_i$ of size $|X|>\eps|V_i|$
such that every $x \in X$ has degree greater than $(d_{ij}+\eps)|V_j|$ in $V_j$. \qed

\begin{lemma}\label{2restrict} {\bf (Regular restriction)} \index{regular restriction}
Suppose $G_{ij}$ is $\eps$-regular, and we have sets
$V'_i \sub V_i$ and $V'_j \sub V_j$ with
$|V'_i| \ge \sqrt{\eps} |V_i|$ and $|V'_j| \ge \sqrt{\eps}|V_j|$.
Then the bipartite subgraph $G'_{ij}$ of $G_{ij}$ induced by $V'_i$ and $V'_j$
is $\sqrt{\eps}$-regular of density $d'_{ij} = d_{ij} \pm \eps$.
\end{lemma}

\nib{Proof.}
Since $G_{ij}$ is $\eps$-regular we have $d'_{ij} = d_{ij} \pm \eps$.
Now consider any sets $V''_i \sub V'_i$ and $V''_j \sub V'_j$
with $|V''_i| \ge \sqrt{\eps} |V'_i|$ and $|V''_j| \ge \sqrt{\eps} |V'_j|$.
Then $|V''_i| \ge \eps |V_i|$ and $|V''_j| \ge \eps |V_j|$,
so $V''_i$ and $V''_j$ induce a bipartite subgraph of $G_{ij}$
with density $d_{ij} \pm \eps \subset d'_{ij} \pm \sqrt{\eps}$.
Therefore $G'_{ij}$ is $\sqrt{\eps}$-regular. \qed

Our next lemma shows that the definition of good vertices for the algorithm is sensible,
in that most free vertices are good. Before giving the lemma, we observe that
the number of free vertices in any given class does not become too small
at any point during the algorithm. We can quantify this as
$|V_i(t)| \ge |B \cap V_i| - |Q(t) \cap V_i| \ge (\delta_B-\delta_Q)n \ge \delta_B n/2$,
for any $1 \le i \le r$ and time $t$. To see this, note that
we stop the iterative procedure when the only unembedded vertices are buffer vertices,
and during the procedure a buffer vertex is only embedded if it joins the queue.

\begin{lemma}\label{good} {\bf (Good vertices)} \index{good}
Suppose we embed a vertex $x=s(t)$ of $H$ at time $t$.

\noindent Then $|OK_x(t-1)| \ge (1-\eps_*)|F_x(t-1)|$,
and for every unembedded vertex $z$ we have $|F_z(t)| \ge d_u n$.
\end{lemma}

\nib{Proof.} We argue by induction on $t$. At time $t=0$ the first statement
is vacuous, as we do not embed any vertex at time $0$, and the second statement
follows from the fact that $F_z(0)=V_z$ has size $n$ for all $z$. Now suppose $t \ge 1$.
By induction we have $|F_z(t-1)| \ge d_u n$ for every unembedded vertex $z$.
Then by Lemma \ref{2restrict}, for any unembedded $z \in N_H(x)$,
the bipartite subgraph of $G$ induced by $F_x(t-1)$ and $F_z(t-1)$
is $\eps'$-regular of density $(1 \pm \eps)d_{xz}$.
Applying Lemma \ref{2neighbour}, we see that there are at most $2\eps'|F_x(t-1)|$
vertices $y \in F_x(t-1)$ that do not satisfy
$|N_G(y) \cap F_z(t-1)| = (1 \pm 2\eps')d_{xz}|F_z(t-1)|$.
Summing over at most $D$ neighbours of $x$ and applying the definition of good vertices
in the algorithm we obtain the first statement that
$|OK_x(t-1)| \ge (1-2D\eps')|F_x(t-1)| \ge (1-\eps_*)|F_x(t-1)|$.

Next we prove the second statement. Consider any unembedded vertex $z$.
Let $t_z$ be the most recent time at which we embedded a neighbour of $z$,
or $0$ if there is no such time. If $t_z>0$ then we embedded some neighbour $w=s(t_z)$
of $z$ at time $t_z$. Since we chose the image $\phi(w)$ of $w$ to be a good vertex,
by definition we have $|F_z(t_z)|=(1\pm 2\eps')d_{wz}|F_z(t_z-1)| > \frac{1}{2}d|F_z(t_z-1)|$.
If $z$ is not in the queue then the rule for updating the queue in the algorithm
gives $|F_z(t)| \ge \delta'_Q |F_z(t_z)|$. On the other hand,
suppose $z$ is in the queue, and that it joined the queue at some time $t'<t$.
Since $z$ did not join the queue at time $t'-1$ we have
$|F_z(t'-1)| \ge \delta'_Q |F_z(t_z)|$.
Also, between times $t'$ and $t$ we only embed vertices that are in the queue:
the queue cannot become empty during this time, as then we would have embedded $z$ before $x$.
During this time we embed at most $\delta_Q n$ vertices in $V_z$, as we abort
the algorithm if the number of vertices in $X_z$ that have ever been queued exceeds this.

Thus we have catalogued all possible ways
in which the number of vertices free for $z$ can decrease.
It may decrease by a factor no worse than $d/2$ when a neighbour of $z$ is embedded,
and by a factor no worse than $\delta'_Q$ before the next neighbour of $z$ is embedded,
unless $z$ joins the queue. Also, if $z$ joins the queue we may subtract
at most $\delta_Q n$ from the number of vertices free for $z$.
Define $i$ to be the number of neighbours of $z$ that are embedded before $z$ joins the queue
if it does, or let $i=d(z)$ be the degree of $z$ if $z$ does not join the queue.
Now $z$ has at most $D$ neighbours, and $|F_z(0)|=|V_z|=n$, so%
\COMMENT{Factor $\delta'_Q$ before first neighbour.}
$|F_z(t)| \ge (\delta'_Q d/2)^{D-i}((\delta'_Q d/2)^i\delta'_Q n-\delta_Q n)
\ge (\delta'_Q d/2)^D\delta'_Q n - \delta_Q n > d_u n$. \qed

Next we turn our attention to the time period during which we are embedding $N$,
which we will refer to as the {\em initial phase} of the algorithm.\index{initial phase}
We start by observing that the queue remains empty during the initial phase,
and so $N$ is embedded consecutively in the order given by the list $L$.
To see this, we use a similar argument to that used for the second statement
in Lemma \ref{good}. Consider any unembedded vertex $z$ and suppose that the
queue has remained empty up to the current time $t$. Then we have embedded
at most $|N \cap X_z| < \sqrt{\delta_B} n$ vertices in $V_z$.
Also, if we embed a neighbour $w$ of $z$ the algorithm chooses a good image for it,
so by definition of good, the number of free images for $z$ decreases by a factor
no worse than $(1-2\eps')d_{wz} > d/2$ when we embed $w$.
Since $z$ has at most $D$ neighbours
we get $|F_z(t)| \ge (d/2)^D n - \sqrt{\delta_B} n > \delta'_Q n$.
This shows that no unembedded vertex $z$ is added to the queue
during the initial phase.

Now we want to show that for any buffer vertex $x \in B$ there will be
many free positions for $x$ at the end of the algorithm.
This is the point in the argument where super-regularity is essential.
A vertex $v \in V_x$ will be free for $x$ if it is not used for another vertex
and we embed $N_H(x)$ in $N_G(v)$ during the initial phase. Our next lemma
gives a lower bound on this probability, conditional on any embedding of the
previous vertices not using $v$. We fix some $x \in B$ and write\index{$z_j$}
$N_H(x)=\{z_1,\cdots,z_g\}$,\index{$g$}
with vertices listed in the order that they are embedded.\index{$T_j$}
We let $T_j$ be the time at which $z_j$ is embedded. Since $N$ is embedded
consecutively we have $T_{j+1}=T_j+1$ for $1 \le j \le g-1$. We also define $T_0=T_1-1$.
Since vertices in $B$ are at mutual distance at least $4$ in $H$,\index{$T_0$}
at time $T_0$, when we have only embedded vertices from $N$,
no vertices within distance $2$ of $x$ have been embedded.
(This is the only place at which mutual distance $4$ is important.)

\begin{lemma}\label{2-x-to-v} {\bf (Initial phase)}
For any $v \in V_x$, conditional on any embedding of the vertices $\{s(u):u<T_1\}$
that does not use $v$, with probability at least $p$ we have $\phi(N_H(x)) \sub N_G(v)$.
\COMMENT{Formerly [nonsense]: so $v \in OK_x(T_g)$.}
\end{lemma}

\nib{Proof.} We estimate the probability that $\phi(N_H(x)) \sub N_G(v)$
using arguments similar to those we are using to embed $H$ in $G$.
We want to track the free positions within $N_G(v)$
for each unembedded vertex in $N_H(x)$, so we hope to not only
embed each vertex of $N_H(x)$ in $N_G(v)$ but also to do so in a `good' way,
a property that can be informally stated as saying that the free
positions within $N_G(v)$ will have roughly their expected size.
To define this formally, suppose $1 \le j \le g$ and we are considering the embedding of $z_j$.
We interpret quantities at time $T_j$ with the embedding $\phi(z_j)=y$,
for some unspecified $y \in F_{z_j}(T_j-1)$.
We say that $y \in F_{z_j}(T_j-1) \cap N_G(v)$ is {\em good for $v$},\index{good}
and write $y \in OK^v_{z_j}(T_j-1)$, if for every unembedded $z \in N_H(x)$ we have
$|F_z(T_j) \cap N_G(v)| = (1 \pm 2\eps')d_{zz_j}|F_z(T_j-1) \cap N_G(v)|$.
We let $A_j$ denote the event that $y=\phi(z_j)$ is chosen in $OK^v_{z_j}(T_j-1)$.

We claim that conditional on the events $A_{j'}$ for $j'<j$
and the embedding up to time $T_j-1$,
the probability that $A_j$ holds is at least $d_u/2$.
To see this we argue as in Lemma \ref{good}.
First we show that $|F_z(T_j-1) \cap N_G(v)| \ge d_u n$
for every unembedded neighbour $z$ of $x$.
Note that initially $F_z(0) \cap N_G(v) = V_z \cap N_G(v)$
has size at least $(d-\eps)n$ by super-regularity of $G$.
Up to time $T_0$ we embed at most $|N \cap X_z| \le \sqrt{\delta_B}n$
vertices in $X_z$, and do not embed any neighbours of $z$
by the distance property mentioned before the lemma.
At time $T_{j'}$ with $j'<j$, we are conditioning on the event that
the algorithm chooses an image for $z_{j'}$ that is good for $v$, so the number
of free images for $z$ within $N_G(v)$ decreases by a factor
no worse than $(1-2\eps')d_{zz_{j'}} > d/2$. Thus we indeed have
$|F_z(T_j-1) \cap N_G(v)| \ge (d/2)^D|V_z \cap N_G(v)|-|N \cap X_z|
\ge (d/2)^D(d-\eps)n - \sqrt{\delta_B}n \ge d_u n$.

Next, by Lemma \ref{2restrict}, for any unembedded $z \in N_H(x)$,
the bipartite subgraph of $G$ induced by $F_{z_j}(T_j-1) \cap N_G(v)$ and $F_z(T_j-1) \cap N_G(v)$
is $\eps'$-regular of density $(1 \pm \eps)d_{zz_j}$.
Applying Lemma \ref{2neighbour}, we see that there are at most $2\eps'|F_{z_j}(T_j-1) \cap N_G(v)|$
vertices $y \in F_{z_j}(T_j-1) \cap N_G(v)$ that do not satisfy $|F_z(T_j) \cap N_G(v)|
= |N_G(y) \cap F_z(T_j-1) \cap N_G(v)| = (1 \pm 2\eps')d_{zz_j}|F_z(T_j-1) \cap N_G(v)|$.
Summing over at most $D$ neighbours of $z_j$ we see that
$|OK^v_{z_j}(T_j-1)| \ge (1-2D\eps')|F_{z_j}(T_j-1) \cap N_G(v)|$.
Also, by Lemma \ref{good} we have $|OK_{z_j}(T_j-1)| \ge (1-\eps_*)|F_{z_j}(T_j-1)|$, so
\[ |OK^v_{z_j}(T_j-1) \cap OK_{z_j}(T_j-1)| \ge  |F_{z_j}(T_j-1) \cap N_G(v)|
- 2\eps_* n \ge d_u n/2.\]
Since $\phi(z_j)=y$ is chosen uniformly at random from
$OK_{z_j}(T_j-1)$, we see that $A_j$ holds with probability at least $d_u/2$, as claimed.

To finish the proof, note that if all the events $A_j$ hold then we have $\phi(N_H(x)) \sub N_G(v)$.
Multiplying the conditional probabilities, this holds with probability at least $(d_u/2)^D>p$. \qed

Our next lemma is similar to the `main lemma' of \cite{KSS}.\index{main lemma}
An informal statement is as follows.
Consider any subset $Y$ of a given class $X_i$,
and any `not too small' subset $A$ of vertices of $V_i$ that could potentially be used for $Y$.
Then it is very unlikely that no vertices in $A$ will be used
and yet only a small fraction of the free
positions for every vertex in $Y$ will belong to $A$.

\begin{lemma}\label{2main} {\bf (Main lemma)}
Suppose $1 \le i \le r$, $Y \sub X_i$ and $A \sub V_i$ with $|A|>\eps_* n$.
Let $E_{A,Y}$ be the event that (i) no vertices are embedded in $A$ before the conclusion of
the algorithm, and (ii) for every $z \in Y$ there is some time $t_z$ such that
$|A \cap F_z(t_z)|/|F_z(t_z)| < 2^{-D} |A|/n$. Then $\mb{P}(E_{A,Y}) < p_0^{|Y|}$.
\end{lemma}

\nib{Proof.}
We start by choosing $Y' \sub Y$ with $|Y'| > |Y|/D^2$ so that vertices
in $Y'$ are mutually at distance at least $3$ (this can be done greedily,
using the fact that $H$ has maximum degree $D$).
It suffices to bound the probability of $E_{A,Y'}$.
Note that initially we have $|A \cap F_z(0)|/|F_z(0)| = |A|/n$ for any $z \in X_i$.
Also, if no vertices are embedded in $A$, then $|A \cap F_z(t)|/|F_z(t)|$
can only be less than $|A \cap F_z(t-1)|/|F_z(t-1)|$ for some $z$ and $t$
if we embed a neighbour of $z$ at time $t$.
It follows that if $E_{A,Y'}$ occurs, then for every $z \in Y'$
there is a first time $t_z$ when we embed a neighbour $w$ of $z$ and
have $|A \cap F_z(t_z)|/|F_z(t_z)| < |A \cap F_z(t_z-1)|/2|F_z(t_z-1)|$.

By Lemma \ref{good} we have $|F_w(t_z-1)| \ge d_u n$, and by
choice of $t_z$ we have $|A \cap F_z(t_z-1)|/|F_z(t_z-1)| \ge 2^{-D}|A|/n$,
so $|A \cap F_z(t_z-1)| \ge 2^{-D} \eps_* d_u n \ge \eps_*^2 n$.
Then by Lemma \ref{2restrict}, the bipartite subgraph of $G$ induced
by $A \cap F_z(t_z-1)$ and $F_w(t_z-1)$ is $\eps'$-regular of density
$(1 \pm \eps)d_{zw}$. Applying Lemma \ref{2neighbour},
we see that there are at most $2\eps'|F_w(t_z-1)|$ `exceptional' vertices
$y \in F_w(t_z-1)$ that do not satisfy $|A \cap F_z(t_z)| =
|N_G(y) \cap A \cap F_z(t_z-1)| = (1 \pm 2\eps')d_{zw}|A \cap F_z(t_z-1)|$.
On the other hand, the algorithm chooses $\phi(w)=y$ to be good, in that
$|F_z(t_z)| = (1 \pm 2\eps')d_{zw}|F_z(t_z-1)|$, so we can only have
$|A \cap F_z(t_z)|/|F_z(t_z)| < |A \cap F_z(t_z-1)|/2|F_z(t_z-1)|$
by choosing an exceptional vertex $y$. But $y$ is chosen uniformly at
random from $|OK_w(t_z-1)| \ge (1-\eps_*)|F_w(t_z-1)|$ possibilities
(by Lemma \ref{good}). It follows that, conditional on the prior embedding,
the probability of choosing an exceptional vertex for $y$ is at most
$2\eps'|F_w(t_z-1)|/|OK_w(t_z-1)| < 3\eps'$.

Since vertices of $Y'$ have disjoint neighbourhoods,
we can multiply the conditional probabilities over $z \in Y'$
to obtain an upper bound of $(3\eps')^{|Y'|}$.
Recall that this bound is for a subset of $E_{A,Y'}$
in which we have specified a certain neighbour $w$ for every vertex $z \in Y'$.
Taking a union bound over at most $D^{|Y'|}$ choices for these neighbours
gives $\mb{P}(E_{A,Y}) \le \mb{P}(E_{A,Y'}) \le (3\eps'D)^{|Y'|} < p_0^{|Y|}$. \qed

Now we can prove the following theorem, which implies Theorem \ref{2blowup}.

\begin{theo}\label{2hp}
With high probability the algorithm embeds $H$ in $G$.
\end{theo}

\nib{Proof.}
First we estimate the probability of the iteration phase aborting with failure,
which happens when the number of vertices that have ever been queued is too large.
We can take a union bound over all $1 \le i \le r$ and $Y \sub X_i$ with $|Y|=\delta_Q|X_i|$
of $\mb{P}(Y \sub Q(T))$. Suppose that the event $Y \sub Q(T)$ occurs.
Then $A=V_i(T)$ is a large set of unused vertices,
yet it makes up very little of what is free to any vertex in $Y$.
To formalise this, note that by definition, for every $z \in Y$
there is some time $t$ such that $|F_z(t)| < \delta'_Q|F_z(t_z)|$,
where $t_z<t$ is the most recent time at which we embedded a neighbour of $z$.
Since $A$ is unused we have $A \cap F_z(t) = A \cap F_z(t_z)$, so
$|A \cap F_z(t_z)|/|F_z(t_z)| = |A \cap F_z(t)|/|F_z(t_z)|
\le |F_z(t)|/|F_z(t_z)| < \delta'_Q$.
However, we noted before Lemma \ref{good} that $|A| \ge \delta_B n/2$,
so since $\delta'_Q \ll \delta_B$ we have $|A \cap F_z(t_z)|/|F_z(t_z)| < 2^{-D} |A|/n$.
Taking a union bound over all possibilities for $i$, $Y$ and $A$,
Lemma \ref{2main} implies that the failure probability is at most
$r \cdot 4^n \cdot p_0^{\delta_Q n} < o(1)$, since $p_0 \ll \delta_Q$.

Now we estimate the probability of the conclusion of the algorithm aborting with failure.
Recall that buffer vertices have disjoint neighbourhoods,
the iterative phase finishes at time $T$, and $|X_i(T)| \ge \delta_B n/2$.
By Hall's criterion for finding a system of distinct representatives,
the conclusion fails if and only if there is some $1 \le i \le r$
and $S \sub X_i(T)$ such that $|\cup_{z \in S} F_z(T)| < |S|$.
We divide into cases according to the size of $S$.

\begin{description}
\item[$0 \le |S|/|X_i(T)| \le \gamma$.]
By Lemma \ref{good} we have $|F_z(T)| \ge d_un > \gamma n$
for every unembedded $z$, so this case cannot occur.

\item[$\gamma \le |S|/|X_i(T)| \le 1-\gamma$.]
In this case we use the fact that $A := V_i(T) \sm \cup_{z \in S} F_z(T)$
is a large set of unused vertices which is completely unavailable to any vertex $z$ in $S$:
we have  $|A| \ge |V_i(T)|-|S| \ge \gamma|X_i(T)| \ge \gamma \delta_B n/2 \ge \gamma^2 n$,
yet $A \cap F_z(T) = \es$, so $|A \cap F_z(T)|/|F_z(T)|=0 < 2^{-D} |A|/n$.
As above, taking a union bound over all possibilities for $i$, $S$ and $A$,
Lemma \ref{2main} implies that the failure probability is at most
$r \cdot 4^n \cdot p_0^{\gamma^2 n} < o(1)$, since $p_0 \ll \gamma$.

\item[$1-\gamma \le |S|/|X_i(T)| \le 1$.]
In this case we claim that with high probability $\cup_{z \in S} F_z(T) = V_i(T)$,
so in fact Hall's criterion holds. It suffices to consider sets $S \sub X_i(T)$
of size exactly $(1-\gamma)|X_i(T)|$. The claim fails if there is some $v \in V_i(T)$
such that $v \notin F_z(T)$ for every $z \in S$. Since $v$ is unused,
it must be that we failed to embed $N_H(z)$ in $N_G(v)$ for each $z \in S$.
By Lemma \ref{2-x-to-v}, these events have probability at most $1-p$
conditional on the prior embedding. Multiplying the conditional probabilities
and taking a union bound over all $1 \le i \le r$, $v \in V_i$ and
$S \sub X_i(T)$ of size $(1-\gamma)|X_i(T)|$, the failure probability is at most
$rn \binom{n}{(1-\gamma)n} (1-p)^{(1-\gamma)|X_i(T)|} < o(1)$.
This estimate uses the bounds $\binom{n}{(1-\gamma)n} \le 2^{\sqrt{\gamma}n}$,
$(1-p)^{(1-\gamma)|X_i(T)|} < e^{-p\delta_Bn/4} < 2^{-p^2 n}$ and $\gamma \ll p$.
\end{description}

In all cases we see that the failure probability is $o(1)$. \qed

\section{Regularity and super-regularity of $3$-graphs}\label{3reg}

When considering how to generalise regularity to $3$-graphs,
a natural first attempt is to mirror the definitions used for graphs.
Consider an $r$-partite $3$-graph $H$ with vertex set $V$ partitioned as $V = V_1 \cup \cdots \cup V_r$.
Let $H_{ijk}$ be the tripartite sub-$3$-graph of $H$ with parts $V_i$, $V_j$ and $V_k$, for any $i,j,k$.
The {\em density}\index{density} of $H_{ijk}$ is $d(H_{ijk}) = \frac{|H_{ijk}|}{|V_i||V_j||V_k|}$.
Given $\eps>0$, we say that $H_{ijk}$ is {\em $\eps$-vertex-regular} if
for all subsets $V'_i \sub V_i$, $V'_j \sub V_j$ and $V'_k \sub V_k$
with $|V'_i|>\eps|V_i|$, $|V'_j|>\eps|V_j|$ and $|V'_k|>\eps|V_k|$,
writing $H'_{ijk}$ for the tripartite sub-$3$-graph of $H$ with parts $V'_i$, $V'_j$ and $V'_k$,
we have $|d(H'_{ijk})-d(H_{ijk})|<\eps$. There is a decomposition theorem for this definition
analogous to the Szemer\'edi Regularity Lemma. This is often known as the `weak hypergraph regularity lemma',
as although it does have some applications, the property of vertex-regularity
is not strong enough to prove a counting lemma.\index{vertex regular}

To obtain a counting lemma one needs to take account of densities of triples of $H$
with respect to sets of pairs of vertices, 
as well as densities of pairs with respect to sets of vertices.
Thus we are led to define regularity for simplicial complexes.
We make the following definitions.

\begin{defn}\label{def-complex}
Suppose $X$ is a set partitioned as $X = X_1 \cup \cdots \cup X_r$.
We say $S \sub X$ is {\em $r$-partite} if $|S \cap X_i| \le 1$ for $1 \le i \le r$.
Write $K(X)$ for the {\em complete} collection of all $r$-partite subsets of $X$.
We say that $H$ is an {\em $r$-partite $3$-complex} on $X = X_1 \cup \cdots \cup X_r$
if $H$ consists of $r$-partite sets of size at most $3$ and is a simplicial complex,
i.e.\ if $T \sub S \in H$ then $T \in H$.
The {\em index} of an $r$-partite set $S$ is $i(S) = \{1 \le i \le r: S \cap X_i \ne \es\}$.
For $I \sub [r]$ we write $H_I$ for the set of $S \in H$ with index $i(S)=I$, when this set is non-empty.
If there are no sets $S \in H$ with $i(S)=I$ we can choose to either set $H_I = \es$ or let
$H_I$ be {\em undefined}, provided that if $H_I$ is defined then $H_J$ is defined for all $J \sub I$.
When not explicitly stated, the default is that $H_I$ is undefined
when there are no sets of index $I$. For any $S \in H$ we write $H_S=H_{i(S)}$ for the
naturally defined $|S|$-partite $|S|$-graph in $H$ containing $S$.
\COMMENT{check for any `empty' or $\es$ which should be `undefined'}
\index{partite}\index{$r$-partite|see{partite}}
\index{3-complex@$3$-complex}\index{complete}\index{$K(\cdot)$}\index{index}
\index{$i(S)$}\index{defined}\index{undefined}\index{$H_I$}\index{$H_S$}
\end{defn}

To put this definition in concrete terms, whenever the following sets are defined,
$H_{\{i\}}$ is a subset of $X_i$,
$H_{\{i,j\}}$ is a bipartite graph with parts $H_{\{i\}}$ and $H_{\{j\}}$,
and $H_{\{i,j,k\}}$ is a $3$-graph contained in the set of triangles spanned
by $H_{\{i,j\}}$, $H_{\{i,k\}}$ and $H_{\{j,k\}}$. Of course, the interesting part
of this structure is the $3$-graph together with its underlying graphs.
We also have $H_\es$, which is usually equal to $\{\es\}$,
i.e.\ the set of size $1$ whose element is the empty set,
although it could be empty if all other parts are empty.
It is most natural to take $H_{\{i\}} = X_i$ for $1 \le i \le r$.
We often allow $H_{\{i\}}$ to be a strict subset of $X_i$,
but note that if desired we can make such a complex `spanning'
by changing the ground set to $X' = X'_1 \cup \cdots \cup X'_r$, where $X'_i = H_{\{i\}}$, $1 \le i \le r$.
When unspecified, our default notation is that $H$ is an $r$-partite $3$-complex on $X = X_1 \cup \cdots \cup X_r$.
We will see later in Definitions \ref{def-restrict} and \ref{def-*} why we have been so careful
to distinguish the cases $H_I$ empty and $H_I$ undefined in Definition \ref{def-complex}.

\begin{defn}\label{def-more}
To avoid clumsy notation we will henceforth frequently identify a set with a sequence of its vertices,
e.g.\ writing $H_i = H_{\{i\}}$ and $H_{ijk} = H_{\{i,j,k\}}$.
We also use concatenation for set union, e.g.\ if $S=ij=\{i,j\}$ then $Sk=ijk=\{i,j,k\}$.
For any $I \sub [r]$ we write $H_{I^\le} = \cup_{I' \sub I} H_{I'}$
for the subcomplex of $H$ consisting of all defined $H_{I'}$ with $I' \sub I$.
We also write $H_{I^<} = H_{I^\le} \sm H_I = \cup_{I' \subn I} H_{I'}$.
Similarly, for any $S \sub X$ we write $H_{S^\le} = \cup_{S' \sub S} H_{S'}$ and
$H_{S^<} = \cup_{S' \subn S} H_{S'}$. For any set system $A$ we let $A^\le$ be the complex generated by $A$,
which consists of all subsets of all sets in $A$.
\index{concatenation}\index{$H_{I^\le}$}\index{$H_{I^<}$}\index{$A^\le$}
\end{defn}

It is clear that intersections and unions of complexes are complexes.
We clarify exactly what these constructions are with the following definition.

\begin{defn}\label{def-iu}
Suppose $H$ and $H'$ are  $r$-partite $3$-complexes on $X = X_1 \cup \cdots \cup X_r$.
The {\em union} $H\cup H'$ is the $r$-partite $3$-complex
where $(H \cup H')_S = H_S \cup H'_S$ is defined whenever $H_S$ or $H'_S$ is defined.
The {\em intersection} $H\cap H'$ is the $r$-partite $3$-complex
where $(H \cap H')_S = H_S \cap H'_S$ is defined whenever $H_S$ and $H'_S$ are defined.
\index{union}\index{intersection}
\end{defn}

We often specify complexes as a sequence of sets, e.g.\
$(G_{12},(G_i: 1 \le i \le 4),\{\es\})$ has parts $G_{12}$, $G_i$ for $1 \le i \le 4$,
$\{\es\}$ and is otherwise undefined.
Now we come to a key definition.

\begin{defn}\label{def-regular}
We let $T_{ijk}(H)$ be the set of triangles formed by $H_{ij}$, $H_{jk}$ and $H_{ik}$.
We say that a defined triple $H_{ijk}$ is {\em $\eps$-regular}
if for every subgraph $G$ of $H$ with $|T_{ijk}(G)| > \eps|T_{ijk}(H)|$ we have
$$\frac{|H \cap T_{ijk}(G)|}{|T_{ijk}(G)|} = \frac{|H_{ijk}|}{|T_{ijk}(H)|} \pm \eps.$$
We also say that the entire $3$-complex $H$ is {\em $\eps$-regular} if every defined triple $H_{ijk}$
is $\eps$-regular and every defined graph $H_{ij}$ is $\eps$-regular.
\index{regular}
\end{defn}

Thus $H_{ijk}$ is $\eps$-regular if for any subgraph $G$ with `not too few' triangles of index $ijk$,
the proportion of these triangles in $G$ of index $ijk$ that are triples in $H_{ijk}$
is approximately equal to the proportion of triangles in $H$ of index $ijk$ that are triples in $H_{ijk}$.
Note that we never divide by zero in Definition \ref{def-regular}, as $\es \ne T_{ijk}(G) \sub T_{ijk}(H)$.
The definition applies even if every $H_{ijk}$ is undefined, in which case we can think of $H$ as a $2$-complex
with every graph $H_{ij}$ being $\eps$-regular. It also applies even if every $H_{ij}$ is undefined,
in which case we can think of $H$ as a $1$-complex (with no regularity restriction).%
\footnote{
Technically one should say that $H_{ijk}$ is $\eps$-regular {\em in the complex $H$},
as the definition depends on the graphs $H_{ij}$, $H_{jk}$ and $H_{ik}$. 
For the sake of brevity we will omit this qualification,
as it is not hard to see that when sufficiently regular
these graphs are `almost' determined by $H_{ijk}$: if say $H_{ij}$ had many pairs
not contained in triples of $H_{ijk}$ then $T_{ijk}$ would have many triangles
none of which are triples of $H_{ijk}$, contradicting the definition of $\eps$-regularity.
}

This concept of regularity in $3$-complexes is more powerful than
vertex regularity, in that it does admit a counting lemma.
In order to apply it we also need an analogue of the Szemer\'edi Regularity Lemma,
stating that a general $3$-complex can be decomposed into a bounded number of pieces,
most of which are regular. Such a result does hold, but there is an important
technical proviso that one cannot use the same parameter $\eps$ to measure regularity
for both graphs and triples in the complex. In general, one needs to allow the densities
of the graphs $H_{ij}$ to be much smaller than the parameter used to measure the regularity
of triples. This is known as a {\em sparse} setting, as contrasted with a situation
when all densities are much larger than $\eps$, which is known as a {\em dense} setting.%
\index{sparse setting}\index{dense setting}

In the sparse setting, a counting lemma does hold,
but we couldn't find any way to generalise the proof of the blow-up lemma.
To circumvent this difficulty we will instead apply the Regular Approximation Lemma
of R\"odl and Schacht. Informally stated, this allows us to closely approximate any $3$-graph $G^0$
by another $3$-graph $G$, so that the $3$-complex $G^\le$ generated by $G$
can be decomposed (in a certain sense) into $\eps$-regular complexes.
We will come to the formal statement later in the paper, but we mention it here to motivate the form
of the blow-up lemma that we will prove. We will allow ourselves to work in the dense setting of
$\eps$-regular complexes, but we have to take account of the approximation of $G^0$ by $G$
by `marking' the edges $M = G \sm G^0$ as forbidden. If we succeed in embedding a given $3$-graph $H$ in $G$
without using $M$ then we have succeeded in embedding $H$ in $G^0$. 
(A similar set-up is used for the embedding lemma in \cite{NORS}.)
We will refer to the pair $(G^\le,M)$ as a {\em marked complex}.%
\index{marked complex}\index{$M$}

In the remainder of this section we first motivate and then explain the definition
of super-regularity for $3$-graphs. It turns out that this needs to be significantly more
complicated than for ordinary graphs. It is not sufficient to just forbid vertices of small degree.
To see this, consider as an example a $4$-partite $3$-complex $G$
on $X = X_1 \cup X_2 \cup X_3 \cup X_4$, with $|X_i|=n$, $1 \le i \le 4$
that is almost complete: say there are complete bipartite graphs on every pair of classes
and complete $3$-graphs on every triple of classes, except for one vertex $x$ in $X_4$
for which the neighbourhood $G(x)$ is triangle-free.
We can easily choose each $G(x)_{ij}$ to have size $n^2/4$ by
dividing each class $X_i$, $1 \le i \le 3$ into two parts.
Then $H$ is $O(1/n)$-regular with densities $1-O(1/n)$ and
has minimum degrees at least $n^2/4$ in each triple,
but $x$ is not contained in any tetrahedron $K_4^3$,
so we cannot embed a perfect $K_4^3$-packing.

Another complication is that the definition of super-regularity for $3$-graphs is not `local',
in the sense that super-regularity of a graph $G$ was defined by
a condition for each of its bipartite subgraphs $G_{ij}$.
Instead, we need to define super-regularity for the entire structure $(G,M)$,
where $G$ is an $r$-partite $3$-complex and $M$ is a set of marked edges.
To explain this point we need to look ahead to the analysis of the algorithm
that we will use to prove the blow-up lemma. First we need an important definition that
generalises the process of restricting a graph to a subset of its vertex set:
we may also consider restricting a complex to a subcomplex in the following manner.

\begin{defn}\label{def-restrict}
Suppose $H$ and $G$ are $r$-partite $3$-complexes on $X = X_1 \cup \cdots \cup X_r$ and $G$ is a subcomplex of $H$.
The {\em restriction} $H[G]$ is the $r$-partite $3$-complex on $X = X_1 \cup \cdots \cup X_r$,
where $H[G]_I$ is defined if and only if $H_I$ is defined, and
$H[G]_I$ consists of all $S \in H_I$ such that $A \in G$ for every $A \sub S$ such that $G_A$ is defined.
\index{restriction}\index{$H[\cdot]$}
\end{defn}
To put this in words, a given set $S$ in $H$ belongs to the restriction $H[G]$
if every subset $A$ of $S$ belongs to $G$, provided that the part of $G$ corresponding to $A$ is defined.
At this point we will give some examples to illustrate Definition \ref{def-restrict}
and clarify the distinction between parts being empty or undefined in Definition \ref{def-complex}.
Consider a $3$-partite $2$-complex $H$ on $X = X_1 \cup X_2 \cup X_3$
such that $H_{12}$, $H_{13}$ and $H_{23}$ are non-empty graphs.
Suppose $G$ is a $3$-partite $1$-complex on $X = X_1 \cup X_2 \cup X_3$.
If $G_1$, $G_2$ and $G_3$ are defined then $H[G]$ is the $3$-partite $2$-complex on $X = X_1 \cup X_2 \cup X_3$
with $H[G]_i = G_i$ and $H[G]_{ij}$ equal to the bipartite subgraph of $H_{ij}$ spanned by $G_i$ and $G_j$.
This corresponds to the usual notion of restriction for graphs.
Note that all of the sets $H[G]_i$ and graphs $H[G]_{ij}$ are defined and some may be empty.
However, if any of the $G_i$ is undefined then it behaves as if it were equal to $H_i$.
For example, if $G_1$ is undefined and $G_2$ is defined then $H[G]_1 = H_1$
and $H[G]_{12}$ is the bipartite subgraph of $H_{12}$ spanned by $H_1$ and $G_2$.
This highlights the importance of distinguishing between $G_1$ being empty or $G_1$ being undefined.
Note that Definition \ref{def-restrict} is monotone, in the sense that adding sets to any given
defined part $G_I$ of $G$ does not remove any sets from any given part $H[G]_J$ of the restriction $H[G]$.
We record the following obvious property for future reference:
\begin{equation}\label{eq:res}
H[G]_I=G_I \mbox{ when } G_I \mbox{ is defined.}
\end{equation}

Another obvious property used later concerns the {\em empty complex} $(\{\es\})$,\index{empty complex}
which satisfies $H[(\{\es\})]=H$ for any complex $H$.
Next we will reformulate the definition of regularity for $3$-complexes
using the restriction notation. First we make the following definitions.

\begin{defn}\label{def-density}
Suppose $H$ is an $r$-partite $3$-complex on $X = X_1 \cup \cdots \cup X_r$.
For any $I \sub [r]$ we let $H_I^*$\index{$H_I^*$} denote the set of $S \in K(X)_I$
% $r$-partite $|I|$-tuples $S \sub X$ of index $i(S)=I$
such that any strict subset $T \subn S$ belongs to $H_T$ when defined.
When $H_I$ is defined and $H_I^* \ne \es$, we define the
{\em relative density (at $I$)}\index{relative density} of $H$
as $d_I(H)=|H_I|/|H_I^*|$.\index{$d_I(H)$} We also call $d_I(H)$ the
{\em $I$-density}\index{I-density@$I$-density} of $H$.
\index{I-density@$I$-density|see{relative density}}
We define the {\em absolute density}\index{absolute density}\index{density} 
of $H_I$ as\index{absolute density|see{density}}
$d(H_I)=|H_I|/\prod_{i \in I}|X_i|$.\index{$d(H_I)$}
\COMMENT{We caution the reader not to confuse these similar notations.}
\end{defn}

To illustrate Definition \ref{def-density}, note that $H_i^*=X_i$
and $d_i(H)=d(H_i)=|H_i|/|X_i|$ when defined.
If $H_{ij}$ is defined then $H_{ij}^* = H_i \times H_j$
and $d_{ij}(H)=|H_{ij}|/|H_{ij}^*|$ is the density of the
bipartite graph $H_{ij}$ with parts $H_i$ and $H_j$.
We have $d(H_{ij})=d_{ij}(H)d_i(H)d_j(H)$,
so $d(H_{ij})=d_{ij}(H)$ in the case when $H_i=X_i$ and $H_j=X_j$.
Also, if $H_{ij}$, $H_{ik}$ and $H_{jk}$ are defined
then $H^*_{ijk}=T_{ijk}(H)$ is the set of triangles in $H$ of index $ijk$.
We also note that if any of the $H_{ij}$ is undefined it behaves
as if it were equal to $H^*_{ij}$, e.g.\ if $H_{ij}$ is undefined and
$H_{ik}$ and $H_{jk}$ are defined then $H^*_{ijk}$ is the set
of triangles in $(H^*_{ij},H_{ik},H_{jk})$. As an illustration of Definition \ref{def-restrict},
we note that $H^*_I=K(X)[H_{I^<}]_I$, recalling that $K(X)$ is the complex of $r$-partite
sets and $H_{I^<} = \cup_{I' \subn I} H_{I'}$.

Now suppose $H$ is a $3$-partite $3$-complex on $X = X_1 \cup X_2 \cup X_3$
such that $H_{123}$ is defined and $H^*_{123} \ne \es$.
Suppose $G \sub H$ is a $2$-complex such that $G_{12}$, $G_{13}$ and $G_{23}$ are defined.
Then $G^*_{123}=T_{123}(G)$ and $H^*_{123}=T_{123}(H)$.
By Definition \ref{def-restrict} we have $H[G]_{123} = H_{123} \cap G^*_{123}$
and $H[G]^*_{123}=G^*_{123}$, so by Definition \ref{def-density}
$d_{123}(H[G]) = |H \cap G^*_{123}|/|G^*_{123}|$.
Therefore $H_{123}$ is $\eps$-regular if whenever
$|G^*_{123}| > \eps|H^*_{123}|$ we have $d_{123}(H[G]) = d_{123}(H) \pm \eps$.

For the remainder of this section we let
$H$ be an $r$-partite $3$-complex on $X = X_1 \cup \cdots \cup X_r$,
and $G$ be an $r$-partite $3$-complex on $V = V_1 \cup \cdots \cup V_r$,
with $|V_i|=|X_i|$ for $1 \le i \le r$. We want to find an embedding $\phi$ of $H$ in $G$.
Our algorithm will consider the vertices of $X$ in some order and embed them
one at a time. At some time $t$ in the algorithm, for each $S \in H$ there
will be some $|S|$-graph $F_S(t) \sub G_S$\index{$F_S(t)$} consisting of those sets $P \in G_S$
that are {\em free}\index{free} for $S$, in that mapping $S$ to $P$ is
{\em locally consistent}\index{locally consistent} with\index{$(\cdot)$|see{$F_S(t)$}}
the embedding so far. These free sets will be {\em mutually consistent},\index{mutually consistent}
in that
\begin{equation}\label{3free}\index{$F_{S^\le}(t)$}
F_{S^\le}(t)  = \cup_{S' \sub S} F_{S'}(t)
\end{equation}
is a complex for every $S \in H$. We use the convention that $F_S(t)$ is undefined
for any $S \notin H$. Note that (\ref{3free}) applies even for $S \notin H$.
\COMMENT{
This prevents us from saying `we consider any expression referring to an
undefined quantity as being undefined'. Can we make any similar general statement?
}

Initially we define $F_S(0) = G_S$ for all $S \in H$.
Now suppose we have defined $F_S(t-1)$ for all $S \in H$
and then at time $t$ we embed some vertex $x \in X$ to some vertex $y \in F_x(t-1)$.
{\em We will use this notation consistently throughout the paper.}\index{$t$}\index{$x$}\index{$y$}
Then for any $S \in H$ containing $x$ we can only allow sets in $F_S(t)$
that correspond to mapping $x$ to $y$, i.e.\ $F_S(t) = F_S(t-1)^y$,\index{$F_S(t-1)^y$}
which is our notation for $\{P \in F_S(t-1): y \in P\}$.
Also, for any $S$ in the neighbourhood complex $H(x)$,
i.e.\ a set $S$ not containing $x$ such that $Sx = S \cup \{x\} \in H$,
in order for $F_S(t)$ to be mutually consistent with
$F_{Sx}(t) = \{P \in F_{Sx}(t-1): y \in P\}$,
we can only allow sets in $F_S(t)$ that are in the neighbourhood of $y$,
i.e.\ $F_S(t) = F_{Sx}(t-1)(y)$, which is our notation for\index{$F_{Sx}(t-1)(y)$}
$\{P: P \cup \{y\} \in  F_{S \cup \{x\}}(t-1)\}$.
Finally, we need to consider the effect that embedding $x$ has for sets $S$
that do not contain $x$ and are not even in the neighbourhood complex $H(x)$.
Such a set $S$ may contain a set $S'$ in $H(x)$, so that $F_{S'}(t)$
is affected by the embedding of $x$. Then mutual consistency requires for that any set
$P \in F_S(t)$, the subset of $P$ corresponding to $S'$ must belong to $F_{S'}(t)$.
We need to include these restrictions for all subsets $S'$ of $S$.
Also, as we are using the vertex $y$ to embed $x$
we have to remove it from any future free sets.
Thus we are led to the following definition.
(Lemma \ref{consistent} will show that it is well-defined.)

\begin{defn}{\bf (Update rule)}
\label{def-update}\index{update rule}\index{$S.x$}\index{$C_{S^\le}(t)$}
Suppose $x$ is embedded to $y$ at time $t$ and $S \in H$.

\hspace{0.1cm} If $x \in S$ we define $F_S(t) = F_S(t-1)^y = \{P \in F_S(t-1): y \in P\}$.
If $x \notin S$ we define
$$S.x = (S \sm X_x) \cup \{x\}, \quad
C_{S^\le}(t) =  F_{S.x^\le}(t-1)(y) \quad \mbox{and}
\quad F_S(t) = F_{S^\le}(t-1)[C_{S^\le}(t)]_S \sm y.$$
\end{defn}

Note that $S.x$ is either $Sx=S \cup \{x\}$ if $i(x) \notin i(S)$
or obtained from $Sx$ by deleting the element in $S$
of index $i(x)$ if $i(x) \in i(S)$. Thus $S.x \sub X$ is $r$-partite.
Also, the notation `$\sm y$' means that we delete all sets containing $y$;\index{$\sm y$}
this can have an effect only when $i(x) \in i(S)$.
We will show below in Lemma \ref{consistent} that Definition \ref{def-update} makes sense,
but first we will give an example to illustrate how it works.

\begin{eg}\label{eg1}
Suppose that $H$ and $G$ are $4$-partite $3$-complexes,
and that we have $4$ vertices $x_i \in X_i$, $1 \le i \le 4$
that span a tetrahedron $K_4^3$ in $H$,
i.e.\ $H$ contains every subset of $\{x_1,x_2,x_3,x_4\}$.
Suppose also that we have the edges $x'_1 x'_2 x_3$ and $x'_1 x'_3 x'_4$
and all their subsets for some other $4$ vertices $x'_i \in X_i$, $1 \le i \le 4$,
and that there are no other edges of $H$ containing any $x_i$ or $x'_i$, $1 \le i \le 4$.
Initially we have $F_S(0)=G_S$ for every $S \in H$.
Suppose we start the embedding by mapping $x_1$ to some $v_1 \in V_1$ at time $1$.
Applying Definition \ref{def-update} to sets $S$ containing $x_1$ gives
$F_{x_1}(1)=\{v_1\}$, $F_{x_1 x_i}(1) = \{P \in G_{1i}: v_1 \in P\}$ for $2 \le i \le 4$,
and $F_{x_1 x_i x_j}(1) = \{P \in G_{1ij}: v_1 \in P\}$ for $2 \le i<j \le 4$.

Next we consider some examples of Definition \ref{def-update} for sets not containing $x_1$.
We have $C_{x_i^\le}(1) = (F_{x_1x_i}(0)(v_1),\{\es\}) = (G(v_1)_i,\{\es\})$ for $2 \le i \le 4$.
Then $F_{x_i}(1) = F_{x_i^\le}(0)[C_{x_i^\le}(1)]_i \sm v_1 = G(v_1)_i$.
Similarly, we have $F_{x_i x_j}(1)=G(v_1)_{ij}$ for $2 \le i<j \le 4$.
For $x_2 x_3 x_4$ we have
$$C_{x_2x_3x_4^\le}(1) = \bigcup_{S \sub x_2x_3x_4} F_{Sx_1}(0)(v_1)
= G(v_1)_{23^\le} \cup G(v_1)_{24^\le} \cup G(v_1)_{34^\le}.$$
Therefore $F_{x_2 x_3 x_4}(1) = F_{x_2x_3x_4^\le}(0)[C_{x_2x_3x_4^\le}(1)]_{234} \sm v_1$
consists of all triples in $G_{234}$ that also form a triangle in the neighbourhood of $v_1$,
and so complete $v_1$ to a tetrahedron in $G$.

For $x'_2x_3$, we have $C_{x'_2x_3}(1) = G(v_1)_{3^\le}$,
so $F_{x'_2x_3}(1) = F_{x'_2x_3^\le}(0)[C_{x'_2x_3^\le}(1)]_{23} \sm v_1$
consists of all pairs in $G_{23}$ that contain a $G_{13}$-neighbour of $v_1$.
For $x'_2$, $C_{x'_2{}^\le}(1)=(\{\es\})$ is the empty complex,\index{empty complex}
so $F_{x'_2}(1) = F_{x'_2{}^\le}(0)[C_{x'_2{}^\le}(1)]_2 \sm v_1 = F_{x'_2}(0)$ is unaffected.
(Recall that $J[(\{\es\})]=J$ for any complex $J$.)
Finally we give two examples in which the deletion of $v_1$ has some effect.
For $x'_1x'_2x_3$, we have $x'_1x'_2x_3.x_1=x_1x'_2x_3$,
$F_{x_1x'_2x_3^\le}(0) = F_{x_1x_3^\le}(0) \cup F_{x'_2x_3^\le}(0) = G_{13^\le} \cup G_{23^\le}$,
and $C_{x'_1x'_2x_3^\le}(1)=F_{x_1x'_2x_3^\le}(0)(v_1) = (G(v_1)_3,\{\es\})=G(v_1)_{3^\le}$.
Then $F_{x'_1x'_2x_3^\le}(1) = G_{123^\le}[G(v_1)_{3^\le}] \sm v_1$,
so $F_{x'_1x'_2x_3}(1)$ consists of all triples $T$ in $G_{123}$ not containing $v_1$
such that $T \cap V_3$ is a neighbour of $v_1$.
For $x'_1x'_3$, $C_{x'_1x'_3{}^\le}(1)$ is the empty complex,
so $F_{x'_1 x'_3}(1) = F_{x'_1 x'_3}(0) \sm v_1$
consists of all pairs in $G_{13}$ that do not contain $v_1$.
\end{eg}

Now we prove a lemma which justifies Definition \ref{def-update}
and establishes the `mutual consistency' mentioned above, i.e.\ that $F_{S^\le}(t)$ is a complex.
First we need a definition.

\begin{defn}\label{def-set-restrict}
Suppose $S \sub X$ is $r$-partite and $I \sub i(S)$.
We write $S_I = S \cap \cup_{i \in I} X_i$. \index{$S_I$}
We also write $S_T = S_{i(T)}$ for any $r$-partite set $T$ with $i(T)\sub i(S)$.
\end{defn}

\begin{lemma}\label{consistent}
Suppose $S \sub X$ is $r$-partite and $t \ge 1$.
If $x \notin S$ then $C_{S^\le}(t)$ is a subcomplex of $F_{S^\le}(t-1)$
and $F_{S^\le}(t) = F_{S^\le}(t-1)[C_{S^\le}(t)] \sm y$ is a complex.
If $x \in S$ then $F_S(t) = F_S(t-1)^y$, $F_{S \sm x}(t) = F_S(t-1)(y)$
and $F_{S^\le}(t) = F_{S^\le}(t-1)^y \cup F_{S^\le}(t-1)(y)$ is a complex.
\end{lemma}

\nib{Proof.}
Note that $F_{S^\le}(0)=G_{S^\le}$ is a complex. We prove the statement of the lemma
by induction on $t$. The argument uses the simple observation that if $J$ is any complex
and $v$ is a vertex of $J$ then $J(v)$ and $J^v \cup J(v)$ are subcomplexes of $J$.

First suppose that $x \notin S$. Since $F_{S^\le}(t-1)$ is a complex by induction hypothesis,
$C_{S^\le}(t)  = F_{S.x^\le}(t-1)(y)$ is a subcomplex of
$F_{S.x^\le}(t-1)$, and so of $F_{S^\le}(t-1)$. For any $S' \sub S$
write $J(S') = F_{S'{}^\le}(t-1)[C_{S'{}^\le}(t)] \sm y$.
Then $J(S')$ is a complex, since restriction to a complex gives a complex,
and removing all sets containing $y$ preserves the property of being a complex.
Furthermore, $J(S')_{S'} = F_{S'}(t)$ by Definition \ref{def-update}.
We also have $J(S)_{S'}=J(S')_{S'}$, since a set $A'$ of index $S'$ belongs
to $F_{S'{}^\le}(t-1)$ if and only if it belongs to $F_{S^\le}(t-1)$,
and any $B \sub A'$ belongs to $C_{S'{}^\le}(t)$ if and only if it belongs to $C_{S^\le}(t)$.
Therefore $F_{S^\le}(t) = J(S) = F_{S^\le}(t-1)[C_{S^\le}(t)] \sm y$ is a complex.

Now suppose that $x \in S$. Then $F_S(t) = F_S(t-1)^y$ by definition.
Next, note that $C_{S\sm x^\le}(t) = F_{S^\le}(t-1)(y)$ is a complex,
and $F_{S\sm x^\le}(t-1)[C_{S\sm x^\le}(t)]_{S\sm x} = F_S(t-1)(y)$ by equation (\ref{eq:res}).
Since $S$ is $r$-partite, deleting $y$ has no effect,
and we also have $F_{S\sm x}(t) = F_S(t-1)(y)$.
Therefore $F_{S^\le}(t) = \cup_{x \in S' \sub S} (F_{S'}(t-1)^y \cup F_{S'}(t-1)(y))
= F_{S^\le}(t-1)^y \cup F_{S^\le}(t-1)(y)$ is a complex. \qed

We will also use the following lemma to construct $F_S(t)$
iteratively from $\{F_{S'}(t): S' \subn S\}$.

\begin{lemma}\label{build-update}
Suppose that $S \in H$, $|S| \ge 2$, $x \notin S$ and $S \notin H(x)$.
Write $F_{S^<}(t) = \cup_{S' \subn S}F_{S'}(t)$. \index{$F_{S^<}(t)$}
Then $F_{S^\le}(t) = F_{S^\le}(t-1)[F_{S^<}(t)]$.
\end{lemma}

\nib{Proof.} First note that for any $A$ with $|A|<|S|$
we have $A \in F_{S^\le}(t) \Lra A \in F_{S^<}(t) \Lra A \in F_{S^\le}(t-1)[F_{S^<}(t)]$.
Now suppose $A \in F_S(t)$. Then $A \in F_{S^\le}(t-1)$ and $A' \in F_{S^<}(t)$
for every $A' \subn A$, so $A \in F_{S^\le}(t-1)[F_{S^<}(t)]$.
Conversely, suppose that $A \in F_{S^\le}(t-1)[F_{S^<}(t)]_S$.
For any $S' \subn S$ with $S' \in H(x)$ we have $F_{S'}(t)=F_{S'x}(t-1)(y)$,
so $A_{S'} \in F_{S'x}(t-1)(y) \sub C_{S^\le}(t)$.
Also $C_{S^\le}(t)_{S'}$ is undefined for $S' \notin H(x)$;
in particular, our assumption that $S \notin H(x)$
means that $C_{S^\le}(t)_S$ is undefined. Therefore $A_{S'} \in C_{S^\le}(t)$
for every $S' \sub S$ such that $C_{S^\le}(t)_{S'}$ is defined,
i.e. $A \in F_{S^\le}(t-1)[C_{S^\le}(t)]_S$.
Also, if $i(x) \in i(S)$, then writing $z = S_x = S \cap X_x \in S^<$,
we have $A \cap V_x \in F_z(t)=F_z(t-1)\sm y$, so $y \notin A$.
Therefore $A \in F_{S^\le}(t-1)[C_{S^\le}(t)]_S \sm y = F_S(t)$. \qed

We referred to `local consistency' when describing the update rule because \index{locally consistent}
it only incorporates the effect that embedding $x$ has on sets containing at
least one neighbour of $x$. To illustrate this, recall that in Example \ref{eg1} above
we have $F_{x'_2}(1)=F_{x'_2}(0)=G_2$. Now $H$ contains $x_1x_3$ and $x_1$ is embedded to $v_1$,
so $x_3$ must be embedded to a vertex in $G(v_1)_3$. Also, $H$ contains $x'_2x_3$, so any
image of $x'_2$ must have a neighbour in $G(v_1)_3$. This may not be the case for every vertex
in $G_2$, so there is some non-local information regarding the embedding that has not yet
been incorporated into the free sets at time $1$.
In light of this, we should admit that our description of sets in $F_S(t)$
as `free' is a slight misnomer, as there may be a small number of sets in $F_S(t)$
that cannot be images of $S$ under the embedding. This was the case even for the graph
blow-up lemma, in which we described vertices in $F_x(t)$ as `free' images for $x$
but then only allowed the use of $OK_x(t) \sub F_x(t)$. On the other hand, our definition
of free sets is relatively simple, and does contain enough information for the embedding.
To see this, note that by definition of restriction $F_{S^\le}(t)$ is a subcomplex of $G_{S^\le}$
at every time $t$, and when all vertices of $S$ are embedded by $\phi$ we have
$F_S(t)=\{\phi(S)\}$ with $\phi(S) \in G_S$. Furthermore, by removing all sets that contain $y$
in the definition of $F_S(t)$ we ensure that no vertex is used more than once by $\phi$.
Therefore it does suffice to only consider local consistency in constructing the embedding,
provided that the sets $F_S(t)$ remain non-empty throughout. The advantage is that we have
the following simple update rule for sets $S$ that are not local to $x$.\index{local}

\begin{lemma}\label{not-local}
If $S$ does not contain any vertex in $\{x\} \cup VN_H(x)$ then $F_S(t)=F_S(t-1) \sm y$.
\COMMENT{Name lemmas?}
\end{lemma}

\nib{Proof.} Note that $C_{S^\le}(t) =  F_{S.x^\le}(t-1)(y) = (\{\es\})$ is the empty complex. \qed

We will also need to keep track of the marked triples $M$ during the embedding algorithm.
Initially, we just have some triples in $G$ that are marked as forbidden for any triple of $H$.
Then, as the algorithm proceeds, each pair of $H$ is forbidden certain pairs of $G$,
and each vertex of $H$ is forbidden certain vertices of $G$. We adopt the following notation.

\begin{defn}\label{def-mark}
For any triple $E \in H$ we write $E^t$ for the subset of $E$ that is unembedded at time $t$.
We define the marked subset of $F_{E^t}(t)$ corresponding to $E$ as
\index{$E^t$}\index{marked}\index{$M_{E^t,E}(t)$}
$$M_{E^t,E}(t)=\{P \in F_{E^t}(t): P \cup \phi(E \sm E^t) \in M_E\}.$$
\end{defn}

In words, $M_{E^t,E}(t)$ consists of all sets $P$ in $F_{E^t}(t)$ that cannot be used
as images for $E^t$ in the embedding, because when we add the images of the
embedded part $E \sm E^t$ of $E$ we obtain a marked triple.
To illustrate this, suppose that in Example \ref{eg1} we have some marked triples $M$.
At time $t=1$ we map $x_1$ to $v_1$, and then the free set for $x_2x_3$ is $F_{x_2x_3}(1)=G(v_1)_{23}$.
Since the edges $M_{123}$ are marked as forbidden, we will mark $M(v_1)_{23} \sub G(v_1)_{23}$
as forbidden by defining $M_{x_2x_3,x_1x_2x_3}(1)=M(v_1)_{23}$.
As another illustration, recall that $F_{x_2x_3x_4}(1)$ consists of all triples in $G_{234}$
that also form a triangle in the neighbourhood of $v_1$.
Then $M_{x_2x_3x_4,x_2x_3x_4}(1) = M \cap F_{x_2x_3x_4}(1)$ consists of all triples in $M_{234}$
that also form a triangle in the neighbourhood of $v_1$.

For any triple $E$ containing $x$ such that $E^{t-1}=x$
we will choose $y=\phi(x) \notin M_{x,E}(t-1)$.
This will ensure that $\phi(E) \notin M$.
The following lemma will enable us to track marked subsets.
The proof is obvious, so we omit it.
\COMMENT{
$M_{E^t,E}(t) = \{P \in F_{E^t}(t): P \cup \phi(E\sm E^t) \in M_E\}
= \{P \in F_{\ov{E}}(t-1)(y): P \cup \{y\} \cup \phi(E\sm \ov{E}) \in M_E\} = M_{\ov{E},E}(t-1)(y)$
}

\begin{lemma}\label{track-mark}$ $
\begin{itemize}
\item[(i)] If $x\in E$ then $E^t=E^{t-1}\sm x$
and $M_{E^t,E}(t) = M_{E^{t-1},E}(t-1)(y)$.
\item[(ii)] If $x\notin E$ then $E^t=E^{t-1}$
and $M_{E^t,E}(t) = M_{E^{t-1},E}(t-1) \cap F_{E^t}(t)$.
\end{itemize}
\end{lemma}

We need one more definition before we can define super-regularity.
It provides some alternative notation for describing the update rule,
but it has the advantage of not referring explicitly to any embedding
or to another complex $H$.

\begin{defn}\label{def-3preplus} \index{$G^{I_v}$}
Suppose $G$ is an $r$-partite $3$-complex on $V = V_1 \cup \cdots \cup V_r$,
$1 \le i \le r$, $v \in G_i$ and $I$ is a subcomplex of $\binom{[r]}{\le 3}$.
We define $G^{I_v} = G[\cup_{S \in I} G(v)_S]$.%
\COMMENT{
1. Formerly tried to unify with nhoods... bad idea!
2. Make consistent with general case?
}
\end{defn}

We will now explain the meaning of Definition \ref{def-3preplus} and illustrate it using our
running Example \ref{eg1}. To put it in words, $G^{I_v}$ is the restriction of $G$ obtained
by only taking those sets $A \in G$ such that any subset of $A$ indexed by a set $S$ in $I$
belongs to the neighbourhood $G(v)$, provided that the corresponding part $G(v)_S$ is defined.
In Example \ref{eg1} we have $C_{x_2^\le}(1)=(G(v_1)_2,\{\es\})$
and $F_{x_2}(1) = G_{2^\le}[C_{x_2^\le}(1)]_2 = G(v_1)_2$.
Choosing $I = 2^\le = (\{2\},\{\es\})$ we have $C_{x_2^\le}(1) = \cup_{S \in I} G(v_1)_S$
and so $F_{x_2}(1) = G^{I_{v_1}}_2$.
Similarly, choosing $I = 23^\le$ we have $C_{x_2x_3^\le}(1) = \cup_{S \in I} G(v_1)_S$
and $F_{x_2x_3}(1) = G_{23^\le}[C_{x_2^\le}(1)]_{23} = G(v_1)_{23} = G^{I_{v_1}}_{23}$.
Also, choosing $I = 23^\le \cup 24^\le \cup 34^\le$
we see that $F_{x_2x_3x_4}(1) = G^{I_{v_1}}_{234}$ consists of all triples in $G_{234}$
that also form a triangle in the neighbourhood of $v_1$.
In general, we can use this notation to describe the update rule for any complex $H$
if we embed some vertex $x$ of $H$ to some vertex $v$ of $G$ at time $1$.
If $x \in S \in H$ we have $F_S(1)=G^v$ as before.
If $x \notin S$ we let $I = \{i(S'): S' \sub S, S'x \in H\}$
and then $F_S(1)=G^{I_v}_S \sm v$.

Finally, we can give the definition of super-regularity.
\COMMENT{check density formulae here and throughout for `defined'}

\begin{defn}\label{def-3super} (Super-regularity)\index{super-regular}
Suppose $G$ is an $r$-partite $3$-complex on $V = V_1 \cup \cdots \cup V_r$
and $M \sub G_= := \{S \in G: |S|=3\}$.\index{$G_=$}
We say that the marked complex $(G,M)$ is $(\eps,\eps',d_2,\theta,d_3)$-{\em super-regular} if

(i) $G$ is $\eps$-regular,
and $d_S(G) \ge d_{|S|}$ if $|S|=2,3$ and $G_S$ is defined,
\COMMENT{
Formerly $G_i = V_i$ for $1 \le i \le r$:
assume in theorem to allow super-reg restriction.}

(ii) for every $1 \le i \le r$, $v \in G_i$ and $S$ such that $G_{Si}$ is defined,
$|M(v)_S| \le \theta|G(v)_S|$ if $|S|=2$ and
$G(v)$ is an $\eps'$-regular $2$-complex with 
$d_{S}(G(v)) = (1 \pm \eps')d_{S}(G)d_{Si}(G)$ for $S \ne \es$,

(iii)  for every vertex $v$ and subcomplex $I$ of $\binom{[r]}{\le 3}$,
$|(M \cap G^{I_v})_S| \le \theta|G^{I_v}_S|$ if $|S|=3$,
and $G^{I_v}$ is an $\eps'$-regular $3$-complex with densities (when defined)
$$ d_S(G^{I_v}) = \left\{ \begin{array}{ll}
(1 \pm \eps')d_S(G) & \mbox{ if } S \notin I, \\
(1 \pm \eps')d_S(G)d_{Si}(G) & \mbox{ if } \es \ne S \in I \mbox{ and } G_{Si}\mbox{ is defined.}\\
\end{array}
\right.$$
\end{defn}

Just as one can delete a small number of vertices from an $\eps$-regular graph
to make it super-regular, we will see later (Lemma \ref{3del}) that one can delete a small number
of vertices from an $\eps$-regular marked $3$-complex to make it super-regular.
For now we will just explain the meaning of Definition \ref{def-3super} with reference to our
running example. First we remark that the parameters in the definition are listed
according to their order in the hierarchy, in that $\eps \ll \eps' \ll d_2 \ll \theta \ll d_3$.
Thus we consider a dense setting, in which the regularity parameters $\eps$ and $\eps'$
are much smaller than the density parameters $d_2$ and $d_3$. However, one should note
that the marking parameter $\theta$ has to be larger than the density parameter $d_2$,
which is the source of some technical difficulties in our arguments.
We will bound the marked sets by a increasing sequence of parameters
that remain small throughout the embedding. For now we just see what the definition
of super-regularity tells us about the first step.

Condition (i) just says that $G$ is a regular complex and gives lower bounds for the
relative densities of its parts.
Condition (ii) is analogous to the minimum degree condition in
the definition of super-regularity for graphs. The second part of the condition says
that the neighbourhood is a regular complex, and that its relative densities
are approximately what one would expect (we will explain the formulae later).
The first part says that the proportion of marked edges through any vertex is not too great.
We need this to control the proportion of free sets that we have to mark as
forbidden during the embedding algorithm. To illustrate this, suppose that in Example \ref{eg1}
we have some marked triples $M$. At time $t=1$ we map $x_1$ to $v_1$,
and then the free set for $x_2x_3$ is $F_{x_2x_3}(1)=G(v_1)_{23}$.
Since the edges $M_{123}$ are marked as forbidden, we will mark
$M_{x_2x_3,x_1x_2x_3}(1) = M(v_1)_{23} \sub G(v_1)_{23}$
as forbidden. Condition (ii) ensures that not too great a proportion is forbidden.
Note that the density of the neighbourhood complex $G(v_1)$ will be much
smaller than the marking parameter $\theta$, so this does not follow if we
only make the global assumption that $M$ is a small proportion of $G$.

Condition (iii) is the analogue to condition (ii) for the restrictions that embedding
some vertex can place on the embeddings of sets not containing that vertex.
A very simple illustration is the case $I=(\{\es\})$, which gives $|M_S| \le \theta|G_S|$ when defined.
(This could also be obtained from condition (ii) by summing over vertices $v$.)
For a more substantial illustration,
consider Example \ref{eg1} and the subcomplex $I = 23^\le \cup 24^\le \cup 34^\le$.
We noted before the definition that $F_{x_2x_3x_4}(1) = G^{I_{v_1}}_{234}$ consists of all triples in $G_{234}$
that also form a triangle in the neighbourhood of $v_1$. The second part of condition (iii)
ensures that $F_{x_2x_3x_4^\le}(1) = G^{I_{v_1}}_{234^\le}$ is a regular complex, and that its relative densities
are approximately what one would expect (again, we will explain the formulae later).
The first part of condition (iii) again is needed to control the proportion of free sets that
are marked. We will mark
$M_{x_2x_3x_4,x_2x_3x_4}(1) = M \cap F_{x_2x_3x_4}(1) = (M \cap G^{I_{v_1}})_{234}$ as forbidden in $F_{x_2x_3x_4}(1)$,
and the condition says that this is a small proportion.
Again, since the neighbourhood of $v_1$ is sparse relative to $G$,
this does not follow only from a global assumption that $M$ is a small proportion of $G$.

As another illustration of condition (iii), suppose that we modify Example \ref{eg1} by deleting
the edge $x_1x_2x_3$ from $H$. Then $C_{x_2x_3x_4^\le}(1) =  G(v_1)_{24^\le} \cup G(v_1)_{34^\le}$
and $F_{x_2 x_3 x_4}(1) = F_{x_2x_3x_4^\le}(0)[C_{x_2x_3x_4^\le}(1)]_{234} \sm v_1$
consists of all triples $S \in G_{234}$ such that $S_{24}$ and $S_{34}$ are edges
in the neighbourhood of $v_1$. Taking $I = 24^\le \cup 34^\le$ we have
$F_{x_2x_3x_4}(1) = G^{I_{v_1}}_{234}$. So condition (iii) tells us that also with this modified $H$,
after embedding $x_1$ to $v_1$ the complex $F_{x_2x_3x_4^\le}(1)$ is regular and
we do not mark too much of $F_{x_2x_3x_4}(1)$ as forbidden.

\section{The $3$-graph blow-up lemma}

In this section we prove the following blow-up lemma for $3$-graphs.

\begin{theo} \label{3blowup} {\bf ($3$-graph blow-up lemma)}
\index{3-graph@$3$-graph blowup lemma}
Suppose $H$ is an $r$-partite $3$-complex
on $X = X_1 \cup \cdots \cup X_r$
of maximum degree at most $D$,
$(G,M)$ is an $(\eps,\eps',d_2,\theta,d_3)$-super-regular
$r$-partite marked $3$-complex on $V = V_1 \cup \cdots \cup V_r$,
where $n \le |X_i|=|V_i|=|G_i| \le Cn$ for $1 \le i \le r$,
$G_S$ is defined whenever $H_S$ is defined, and
$0 \ll 1/n \ll \eps \ll \eps' \ll d_2 \ll \theta \ll d_3, 1/r, 1/D, 1/C$.
Then $G \sm M$ contains a copy of $H$, in which for each $1 \le i \le r$
the vertices of $V_i$ correspond to the vertices of $X_i$.%
\COMMENT{
1. Assume spanning here (not in super-reg).
2. Make $H$ a complex.
3. Defined or non-empty?
4. G= earlier.
}
\end{theo}

Theorem \ref{3blowup} is similar in spirit to Theorem \ref{2blowup}: informally speaking,
we can embed any bounded degree $3$-graph in any super-regular marked $3$-complex.
(The parallel would perhaps be stronger if we had also introduced marked edges in the graph
statement; this can be done, but there is no need for it, so we preferred the simpler form.)
We remark that we used the assumption $|V_i|=|X_i|=n$ in Theorem \ref{2blowup} for simplicity,
but the assumption $n \le |V_i|=|X_i| \le Cn$ works with essentially the same proof,
and is more useful in applications. (Arbitrary part sizes are permitted in \cite{KSS},
but this adds complications to the proof, and it is not clear why one would
need them, so we will not pursue this option here.)
There are various other generalisations that are
useful in applications, but we will postpone discussion of them until we give the
general hypergraph blow-up lemma. Theorem \ref{3blowup} is already sufficiently
complex to illustrate the main ideas of our approach, so we prove it first
so as not to distract the reader with additional complications.

The section contains six subsections, organised as follows.
In the first subsection we present the algorithm that
we will use to prove Theorem \ref{3blowup},
and also establish some basic properties of the algorithm.
Over the next two subsections we develop some theory:
the second subsection contains some useful properties
of restriction (Definition \ref{def-restrict});
the third contains some properties of regularity for $3$-graphs,
which are similar to but subtly different from known results in the literature.
Then we start on the analysis of the algorithm, following
the template established in the proof of Theorem \ref{2blowup}.
The fourth subsection concerns good vertices,
and is analogous to Lemma \ref{good}.
The fifth subsection concerning the initial phase is the most technical,
containing three lemmas that play the role of Lemma \ref{2-x-to-v}
for $3$-graphs. The final subsection concerns the conclusion of the
algorithm, and is analogous to Lemma \ref{2main} and Theorem \ref{2hp}.

\subsection{The embedding algorithm}\label{3alg}

As for the graph blow-up lemma, we will prove Theorem \ref{3blowup}
via a random greedy algorithm to construct an embedding
$\phi:V(H) \to V(G)$ such that $\phi(e) \in G \sm M$ for every edge $e$ of $H$. \index{$\phi$}
In outline, it is quite similar to the algorithm used for graphs,
but when it comes to details the marked edges create significant complications.
We introduce more parameters with the hierarchy
\begin{gather*}
0 \le 1/n \ll \eps \ll \eps' \ll \eps_{0,0} \ll \cdots \ll \eps_{12D,3} \ll \eps_*
\ll p_0 \ll \gamma \ll \delta_Q \ll p \ll d_u \ll d_2 \\
\ll \theta \ll \theta_0 \ll \theta'_0 \ll \cdots \ll \theta_{12D} \ll \theta'_{12D} \ll \theta_*
\ll \delta'_Q \ll \delta_B \ll d_3, 1/r, 1/D.
\end{gather*}
Most of these parameters do not require any further comment,
as we explained their role in the graph blow-up lemma,
and they will play the same role here.
We need many more `annotated $\eps$' parameters to measure \index{annotated $\eps$}
regularity here, but this is merely a technical inconvenience.
The parameters $\eps_{i,j}$ with $0 \le i \le 12D$ and $0 \le j \le 3$
satisfy $\eps_{i,j} \ll \eps_{i',j'}$ when $i<i'$ or $i=i'$ and $j<j'$.
Because of the marked edges, we also have new `annotated $\theta$' parameters, \index{annotated $\theta$}
which are used to bound the proportion of free sets that are marked.
It is important to note that the buffer parameter $\delta_B$
and queue admission parameter $\delta'_Q$ are larger than the marking
parameter $\theta$, which in turn is larger than the density parameter $d_2$.
The result is that the queue may become non-empty during the initial phase,
and then the set $N$ of neighbours of the buffer $B$ will not be embedded consecutively
in the order given by the original list $L$. To cope with this, we need a
modified selection rule that allows vertices in $N$ to jump the queue.

\begin{description}

\item[Initialisation and notation.]
We choose a buffer set $B \subset X$ of vertices at mutual distance at least $9$ in $H$
so that $|B \cap X_i| = \delta_B |X_i|$ for $1 \le i \le r$.
Since $n \le |X_i| \le Cn$ for $1 \le i \le r$ and $H$ has maximum degree $D$
we can construct $B$ simply by selecting vertices one-by-one greedily.
For any given vertex there are at most $(2D)^8$ vertices at distance less than $9$,
so at any point in the construction we have excluded at most $(2D)^8r\delta_BCn$
vertices from any given $X_i$. Thus we can construct $B$ provided that $(2D)^8r\delta_BC < 1$.

Let $N = \cup_{x \in B} VN_H(x)$ be the set of vertices with a neighbour in the buffer.
Then $|N \cap X_i| \le 2Dr\delta_BCn < \sqrt{\delta_B} n$ for $1 \le i \le r$.

We order the vertices in a list $L=L(0)$ that starts with $N$ and ends with $B$.
Within $N$, we arrange that $VN_H(x)$ is consecutive for each $x \in B$.
We denote the vertex of $H$ selected for embedding at time $t$ by $s(t)$.
This will usually be the first vertex of $L(t-1)$,
but we will describe some exceptions to this principle in the selection rule below.
\COMMENT{Why consecutive? Can't afford all of $N$ jumping the queue!}

We denote the queue by $q(t)$ and write $Q(t) = \cup_{u \le t}\ q(u)$.
We denote the vertices jumping the queue by $j(t)$ and write $J(t) = \cup_{u \le t}\ j(u)$.
Initially we set $q(0)=Q(0)=j(0)=J(0)=\es$.
\index{jump}\index{$j(t)$}\index{$J(t)$}

We write $F_S(t)$ for the sets of $G_S$ that are {\em free} to embed a given set $S$ of $H$.
We also use the convention that $F_S(t)$ is undefined if $S \notin H$.
Initially we set $F_S(0)=G_S$ for $S \in H$.
We also write $X_i(t) = X_i \sm \{s(u): u \le t\}$ for the unembedded vertices of $X_i$
and $V_i(t) = V_i \sm \{\phi(s(u)): u \le t\}$ for the free positions in $V_i$.
We let $X(t) = \cup_{i=1}^r X_i(t)$ and $V(t) = \cup_{i=1}^r V_i(t)$.
\COMMENT{Initially we set $X_i(0)=X_i$ and $V_i(0)=V_i$.}

\item[Iteration.] At time $t$, while there are still some
unembedded non-buffer vertices, we select a vertex to embed $x=s(t)$
according to the following {\em selection rule}. Informally, the \index{selection rule}
rule is that our top priority is to embed any vertex neighbourhood
$VN_H(x)$ with $x \in B$ as a consecutive sequence before embedding $x$ itself
or any other vertex with distance at most $4$ from $x$, and our second
priority is to embed vertices in the queue. Formally, the rule is:
\COMMENT{May as well put VN before x, although it doesn't matter.}

\begin{itemize}
\item
If $j(t-1) \ne \es$ we select $x=s(t)$ to be any element of $j(t-1)$,
\item
If $j(t-1) = \es$ and $q(t-1) \ne \es$ we consider any element $x'$ of $q(t-1)$.
 \begin{itemize}
 \item If the distance from $x'$ to all vertices in the buffer $B$ is at least $5$
 then we select $x=x'=s(t)$.
 \item Otherwise, there is a vertex $x'' \in B$ at distance at most $4$ from $x'$,
 and $x''$ is unique by the mutual distance property of $B$. If there are any
 unembedded elements of $VN_H(x'')$, we choose one of them to be $x=s(t)$,
 choosing $x'$ itself if $x' \in VN_H(x'')$, and put all other unembedded vertices
 of $VN_H(x'')$ in $j(t)$. If all of $VN_H(x'')$ has been embedded we choose $x=x'=s(t)$.
 \end{itemize}
\item
If $j(t-1) = q(t-1) = \es$ we let $x=s(t)$ be the first vertex of $L(t-1)$.
\end{itemize}

We choose the image $\phi(x)$ of $x$ uniformly at random among all elements
$y \in F_x(t-1)$ that are `good', a property that can be informally stated
as saying that if we set $\phi(x)=y$ then the free sets $F_S(t)$
will be regular, have the correct density, and not create too much
danger of using an edge marked as forbidden. Now we will describe the
formal definition. Note that all expressions at time $t$ are
to be understood with the embedding $\phi(x)=y$, for some unspecified vertex $y$.

\nib{Definitions.}

1. For a vertex $x$ we write $\nu_x(t)$ for the number of \index{$\nu_x(t)$}
elements in $VN_H(x)$ that have been embedded at time $t$.
For a set $S$ we write $\nu_S(t) = \sum_{y \in S} \nu_y(t)$. \index{$\nu_S(t)$}
We also define $\nu'_S(t)$ as follows. \index{$\nu'_S(t)$}
When $|S|=3$ we let $\nu'_S(t) = \nu_S(t)$.
When $|S|=1,2$ we let $\nu'_S(t) = \nu_S(t)+K$, where $K$ is the
maximum value of $\nu'_{Sx'}(t')$ over vertices $x'$ embedded at time
$t' \le t$ with $S \in H(x')$; if there is no such vertex $x'$
we let $\nu'_S(t) = \nu_S(t)$.
\COMMENT{
1. Formerly used $\nu_S$, but this is not sufficient.
2. Formerly (changed to fit $k$-graphs):
When $|S|=2$, we let $\nu'_S(t) = \nu_S(t)+K$, where $K$ is the
maximum value of $\nu_x(t')$ over vertices $x$ embedded at time
$t' \le t$ with $S \in H(x)$; if there is no such vertex $x$
we let $\nu'_S(t) = \nu_S(t)$.
When $S=\{z\}$ has size $1$ we let $\nu'_z(t) = \nu_z(t)+K$, where $K$ is the
maximum value of $\nu'_{xz}(t')$ over vertices $x$ embedded at time
$t' \le t$ with $xz \in H$; if there is no such vertex $x$
we let $\nu'_z(t) = \nu_z(t)$.
}

2. As in Definition \ref{def-update}, for any $r$-partite set $S$
we define $F_S(t)=F_S(t-1)^y$ if $x \in S$
or $F_S(t) = F_{S^\le}(t-1)[F_{S.x^\le}(t-1)(y)]_S \sm y$ if $x \notin S$.
For any sets $S' \sub S$ we write $d_{S'}(F(t))=d_{S'}(F_{S^\le}(t))$;
there is no ambiguity, as the density is the same for any $S$ containing $S'$.
\COMMENT{
Former text, cf complex-coloured later...
we introduce some simplified notation for the available sets.
We let $F(t)$ denote the collection of all $F_S(t)$ and write $F(t)_S=F_S(t)$.
Similarly, for any vertex $y$, we let $F(t)^y$ denote the collection of all $F_S(t)^y$,
and $F(t)(y)$ denote the collection of all $F_S(t)(y)$.
This simplified notation allows us to write the update rule as
$F_S(t)=F_S(t-1)^y$ if $x \in S$ or $F_S(t) = F(t-1)[F(t-1)(y)]_S \sm y$ if $x \notin S$.
Here we just intend $F(t-1)[F(t-1)(y)]_S$ as a shorthand for
$F_{S^\le}(t-1)[F_{S.x^\le}(t-1)(y)]_S$
and note that it is unambiguous in its meaning.
}

3. We say that $S$ is {\em unembedded} if every vertex of $S$ is unembedded, \index{unembedded}
i.e.\ $s(u) \notin S$ for $u \le t$.
We define an {\em exceptional} set $E_x(t-1) \sub F_x(t-1)$ by saying\index{exceptional}\index{$E_x(t-1)$}
$y$ is in $F_x(t-1) \sm E_x(t-1)$ if and only if
for every unembedded $\es \ne S \in H(x)$,
\begin{equation}\label{eq:3alg}\tag{$*_{\ref{3alg}}$}
\left.
\begin{aligned}
&d_S(F(t))=(1 \pm \eps_{\nu'_S(t),0})d_S(F(t-1))d_{Sx}(F(t-1)) \nonumber \\
&\mbox{and } F_S(t) \mbox{ is } \eps_{\nu'_S(t),0}\mbox{-regular when } |S|=2. 
\end{aligned}
\qquad \right\}
\end{equation}
%\begin{align}
%&d_S(F(t))=(1 \pm \eps_{\nu'_S(t),0})d_S(F(t-1))d_{Sx}(F(t-1)) \nonumber \\
%&\mbox{and } F_S(t) \mbox{ is } \eps_{\nu'_S(t),0}\mbox{-regular when } |S|=2. \tag*{$(*_{\ref{3alg}})$}
%\end{align}
Lemma \ref{3exceptional} will imply that $E_x(t-1)$ is small compared to $F_x(t-1)$.%
\COMMENT{1. Nhood only here. 2. Check `primes', i.e.\ 2nd index.}
%$\qquad d_S(F(t))=(1 \pm \eps_{\nu'_S(t),0})d_S(F(t-1))d_{Sx}(F(t-1))$,
%$\qquad$ and $F_S(t)$ is $\eps_{\nu'_S(t),0}$-regular when $|S|=2$. \hfill $(*_{\ref{3alg}})$

4. As in Definition \ref{def-mark}, for any triple $E \in H$ we write $E^t$ for the subset of $E$
that is unembedded at time $t$ and $M_{E^t,E}(t)=\{P \in F_{E^t}(t): P \cup \phi(E \sm E^t) \in M_E\}$.
We define $$D_{x,E}(t-1) = \{ y \in F_x(t-1): |M_{E^t,E}(t)| > \theta_{\nu'_{E^t}(t)} |F_{E^t}(t)|\}.$$
Intuitively, these sets consist of vertices $y$ to which it is {\em dangerous} to embed $x$.
Lemma \ref{3marked} will show that only a small proportion of free vertices are dangerous.%
\index{dangerous}\index{$D_{x,E}(t-1)$}

5. Let $U(x)$\index{$U(x)$} be the set of all triples $E \in H$ with
$E^{t-1} \cap (VN_H(x) \cup x) \ne \es$.
We obtain the set of {\em good}\index{good} elements $OK_x(t-1)$ from $F_x(t-1)$
by deleting $E_x(t-1)$ and $D_{x,E}(t-1)$ for every $E \in U(x)$.
\COMMENT{
1. Changed def U from E meets N(x). Need x when $E^{t-1}=x$.
2. Could we mark pairs and vertices?
}

We embed $x$ as $\phi(x)=y$ where $y$ is chosen uniformly at random from the
good elements of $F_x(t-1)$.
We conclude the iteration by updating $L(t-1)$, $q(t-1)$ and $j(t-1)$.
First we remove $x$ from whichever of these sets it was taken.
Then we add to the queue any unembedded vertex $z$ for which $F_z(t)$ has become `too small'.
To make this precise, suppose $z \in L(t-1) \sm \{x\}$,
and let $t_z$ be the most recent time at which we embedded
a vertex in $VN_H(z)$, or $0$ if there is no such time.
(Note that if $z \in VN_H(x)$ then $t_z=t$.)
We add $z$ to $q(t)$ if $|F_z(t)| < \delta'_Q|F_z(t_z)|$.
This defines $L(t)$, $q(t)$ and $j(t)$.

Repeat this iteration until the only unembedded vertices are buffer vertices,
but abort with failure if at any time we have
$|Q(t) \cap X_i| > \delta_Q|X_i|$ for some $1 \le i \le r$.
Let $T$ denote the time at which the iterative phase terminates\index{$T$}
(whether with success or failure).

\item[Conclusion.]
Suppose $x \in X(T)$ is unembedded at time $T$
and we embed the last vertex of $VN_H(x)$ at time $t_x^N$.\index{$t_x^N$}
We define the following {\em available} sets for $x$.\index{available}
We let $A_x$ be obtained from $F_x(t_x^N)$ by removing all
sets $M_{x,E}(t_x^N)$ for triples $E$ containing $x$.\index{$A_x$}
We let $A'_x = A_x \cap V_x(T)$.\index{$A'_x$}
We choose a system of distinct representatives for $\{A'_x:x\in X(T)\}$
to complete the embedding, either ending with success if this is possible,
or aborting with failure if it is not possible.
\end{description}

To justify this algorithm, we need to show that if it does not abort with failure
then it does embed $H$ in $G \sm M$.
We explained in the previous section why the `local consistency' of the update rule implies
that it embeds $H$ in $G$, so we just need to show that no marked edge is used.
This follows from the following lemma.

\begin{lemma}\label{itworks}
Suppose $x$ is the last unembedded vertex of some triple $E$ at time $t-1$.

Then $D_{x,E}(t-1) = M_{x,E}(t-1)$.
If $\phi(x) \notin M_{x,E}(t-1)$ then $\phi(E) \in G \sm M$.
\end{lemma}

\nib{Proof.}
Note that $E^{t-1}=x$, $E^t = \es$ and $F_{\es}(t) = \{\es\}$ is a set of size $1$.
If we were to choose $y \in M_{x,E}(t-1)$ then we would get $M_{\es,E}(t)=\{\es\}$ and so
$|M_{\es,E}(t)| = 1 > \theta_{\nu'_{E^t}(t)} = \theta_{\nu'_{E^t}(t)}|F_{\es}(t)|$.
On the other hand, if we choose $y \notin M_{x,E}(t-1)$ then we get
$M_{\es,E}(t)=\es$ and so $|M_{\es,E}(t)| = 0 < \theta_{\nu'_{E^t}(t)} = \theta_{\nu'_{E^t}(t)}|F_{\es}(t)|$.
Therefore $D_{x,E}(t-1) = M_{x,E}(t-1)$. The second statement is now clear. \qed

It will often be helpful to use the following terminology pertaining to increments of $\nu_x(t)$.
We think of time as being divided into {\em $x$-regimes}, defined by the property that vertices \index{regime}
of $VN_H(x)$ are embedded at the beginning and end of $x$-regimes, but not during $x$-regimes.
Thus the condition for adding a vertex $z$ to the queue is that the free set for $z$
has shrunk by a factor of $\delta'_Q$ during the current $z$-regime.
Note that each vertex $x$ defines its own regimes, and regimes for different vertices can
intersect in a complicated manner.

Note that any vertex neighbourhood contains at most $2D$ vertices.
Thus in the selection rule, any element of the queue can cause at most $2D$
vertices to jump the queue. Note also that when a vertex neighbourhood
jumps the queue, its vertices are immediately embedded at consecutive times
before any other vertices are embedded.

We collect here a few more simple observations on the algorithm.

\begin{lemma}\label{3observe}$ $
\item[(i)] For any $1 \le i \le r$ and time $t$ we have $|V_i(t)| \ge \delta_B n/2$.
\item[(ii)] For any $t$ we have $|J(t)| \le 2D|Q(t)| \le \sqrt{\delta_Q}n$.
\item[(iii)] We have $\nu_x(t) \le 2D$ for any vertex,
$\nu'_S(t) \le 6D$ when $|S|=3$, $\nu'_S(t) \le 10D$ when $|S|=2$,
and $\nu'_S(t) \le 12D$ when $|S|=1$. Thus the $\eps$-subscripts are
always defined in $(*_{\ref{3alg}})$.
\item[(iv)] For any $z \in VN_H(x)$ we have $\nu_z(t)=\nu_z(t-1)+1$,
so for any $S \in H$ that intersects $VN_H(x)$ we have $\nu_S(t)>\nu_S(t-1)$.
\item[(v)] If $\nu_S(t)>\nu_S(t-1)$ then $\nu'_S(t)>\nu'_S(t-1)$.
\item[(vi)] If $z$ is embedded at time $t'\le t$ and $S\in H(z)$ then
$\nu'_S(t) \ge \nu'_{Sz}(t) > \nu'_{Sz}(t'-1)$.
\end{lemma}

\nib{Proof.} As in the graph blow-up lemma,
we stop the iterative procedure when the only unembedded vertices are buffer vertices,
and during the procedure a buffer vertex is only embedded if it joins the queue.
Therefore $|V_i(t)| \ge |B \cap V_i| - |Q(t) \cap V_i| \ge \delta_Bn-\delta_QCn \ge \delta_B n/2$,
so (i) holds. The fact that an element of the queue
can cause at most $2D$ vertices to jump the queue gives (ii).
Statements (iii) and (iv) are clear from the definitions.
Statement (v) follows because $\nu'_S(t)$ and $\nu'_S(t-1)$ are obtained from
$\nu_S(t)$ and $\nu_S(t-1)$ by adding constants that are maxima of certain sets,
and the set at time $t$ includes the set at time $t-1$.
For (vi) note that $\nu'_S(t)=\nu_S(t)+K$, with $K\ge\nu'_{Sz}(t)$
and $\nu'_{Sz}(t) > \nu'_{Sz}(t'-1)$ by (iv) and (v). \qed

\subsection{Restrictions of complexes}

Before analysing the algorithm in the previous subsection, we need to develop
some more theory. In this subsection we prove a lemma that justifies various
manipulations involving restrictions (Definition \ref{def-restrict}).
We often consider situations when several restrictions are placed on a complex,
and then it is useful to rearrange them. 
We define {\em composition} of complexes as follows.%
\COMMENT{
[comments from old definition of restriction]
Consider making restriction more natural? We could distinguish `empty set'
and `undefined' for parts in a complex. Then undefined can play the role of
complete and monotonicity is established. Would we pay for it elsewhere?
Can still consider $F_S(t)$ for any $S$ and call it undefined for $S \notin H$.
Can still have $2$-complex as case of $3$-complex with level 3 undefined.
Need convention in notation that unspecified parts are undefined rather than empty,
e.g.\ $G_{1^\le}$ has $1$-part $G_1$, $\es$-part $\es$ and all other parts undefined.
We can use $A \in^* G_A$ to mean $A \in G_A$ or $G_A$ is undefined.
Union/intersection defined when either/both defined.
Define * as below, no need for well-formed, as if $G_S$ or $G'_S$
is defined then so is $(G*G')_S$.
[Formerly: We say that $G*G'$ is {\em well-formed}
if whenever $G_S \ne \es$ or $G'_S \ne \es$
for some $S$ we have $(G*G')_S \ne \es$...
If $G*G'$ and $G'*G''$ are well-formed
then $(G*G')*G''=G*(G'*G'')$... need not be well-formed!
}

\begin{defn}\label{def-*}\index{$G*G'$}\index{composition}
We write $x \in^* S$ to mean that $x \in S$ or $S$ is undefined.\index{$\in^*$}
Suppose $G$ and $G'$ are $r$-partite $3$-complexes on $X = X_1 \cup \cdots \cup X_r$.
We define $(G*G')_S$ if $G_S$ or $G'_S$ is defined and say that
$S \in (G*G')_S$ if $A \in^* G_A$ and $A \in^* G'_A$ for any $A \sub S$.
We say that $G,G'$ are {\em separate}\index{separate} if there is no $S \ne \es$
with $G_S$ and $G'_S$ both defined.
Given complexes $G^1,\cdots,G^t$ we write $\bigodot_{i=1}^t J_i = J_1 * \cdots * J_t$.
Similarly we write $\bigodot_{i \in I} J_i$ for the composition
of a collection of complexes $\{J^i: i \in I\}$.\index{$\bigodot$}
\index{*|see{$H_I^*$, $\in^*$, $G*G'$, $\bigodot$, $F(T_j)^{Z*v}_{S^\le}$, $i^*(S)$, $P_A^*$}}
\end{defn}

Note that $G \cap G' \sub G*G' \sub G \cup G'$.
To illustrate the relation $\in^*$, we note that if $G$ is a subcomplex of $H$
then $H[G]$ is the set of all $S \in H$ such that $A \in^* G_A$ for all $A \sub S$.
It follows that $H*G=H[G]$, see (vi) in the following lemma.
More generally, the composition $G*G'$ describes `mutual restrictions',
i.e.\ restrictions that $G$ places on $G'$ and restrictions that $G'$ places on $G$.
We record some basic properties of Definition \ref{def-*} in the following lemma.
Since the statement and proof are heavy in notation, we first make a few remarks
to indicate that the properties are intuitive. Property (i) says that mutual
restrictions can be calculated in any order. Property (ii) says that a
neighbourhood in a mutual restriction is given by the mutual restriction
of the neighbourhoods and the original complexes. Property (iii) says
that separate restrictions act independently. The remaining properties
give rules for rearranging repeated restrictions. The most useful cases
are (v) and the second statement in (vi), which convert a repeated
restriction into a single restriction (the other cases are also used,
but their statements are perhaps less intuitive).
One should note that the distinction made earlier between `empty' and `undefined'
is important here; e.g. $(G*G')_S$ undefined implies that $G_S$ and $G'_S$ are
undefined, but this is not true with `undefined' replaced by `empty'.

\begin{lemma}\label{*props}
Suppose $H$ is an $r$-partite $3$-complex and $G,G',G''$ are subcomplexes of $H$.
\begin{itemize}
\item[(i)] * is a commutative and associative operation on complexes,
\item[(ii)] $(G*G')_{S^\le}(v) = G_{S\sm v^\le}*G'_{S\sm v^\le}*G_{S^\le}(v)*G'_{S^\le}(v)$
for any $v \in S \in G \cup G'$,%
\COMMENT{
Formerly [changed for R-indexed case and also manner of use]:
$(G*G')(v)_S = (G*G'*G(v)*G'(v))_S$ for any vertex $v$ and $S \in (G \cup G')(v)$.
}
\item[(iii)] If $G,G'$ are separate then $G*G' = G \cup G'$ and $H[G][G']=H[G'][G]=H[G \cup G']$,
\item[(iv)] $H[G][G*G']=H[G'][G*G']$,
\item[(v)] If $G'$ is a subcomplex of $H[G]$ then $H[G][G']=H[G*G']$.
\item[(vi)] If $G'$ is a subcomplex of $G$ then $G*G'=G[G']$ and
$H[G][G[G']]=H[G'][G[G']]=H[G][G']$. If also $G'_S$ is defined whenever $G_S$ is defined
then $H[G][G']=H[G']$.
\end{itemize}
\end{lemma}

\nib{Proof.}
(i) By symmetry of the definition we have commutativity $G*G'=G'*G$.
Next we show that $G*G'$ is a complex. Suppose $A \sub S' \sub S \in G*G'$.
Since $S' \sub S \in G \cup G'$ we have $S' \in G \cup G'$.
Since $A \sub S \in G*G'$ we have $A \in^* G_A$ and $A \in^* G'_A$.
Therefore $S' \in G*G'$, so $G*G'$ is a complex.
Now we show associativity, i.e.\ $(G*G')*G''=G*(G'*G'')$.
Suppose $S \in (G*G')*G''$. We claim that for any $A \sub S$ we have
$A \in^* G_A$, $A \in^* G'_A$ and $A \in^* G''_A$.
To see this, we apply the definition of $(G*G')*G''$ to get
$A \in^* (G*G')_A$ and $A \in^* G''_A$.
If $A \in (G*G')_A$ we have $A \in^* G_A$ and $A \in^* G'_A$.
Otherwise, $(G*G')_A$ is undefined, so $G_A$ and $G'_A$ are undefined.
This proves the claim. Now $S \in (G*G')*G''$ implies that $S \in G \cup G' \cup G''$,
so $S \in G$ or $S \in G' \cup G''$. If $S \in G' \cup G''$ then by the claim
we have $A \in^* G'_A$ and $A \in^* G''_A$ for $A\sub S$, so $S \in G'*G''$.
Therefore $S \in G \cup (G'*G'')$.
Also, if $A \sub S$ with $(G'*G'')_A$ defined then, for any $A' \sub A$,
since $A' \sub S$, the claim gives $A' \in^* G'_{A'}$ and $A' \in^* G''_{A'}$.
Therefore $A \in (G'*G'')_A$. This shows that $S \in G*(G'*G'')$, so $(G*G')*G''\sub G*(G'*G'')$.
Also, $G*(G'*G'')=(G''*G')*G \sub G''*(G'*G)=(G*G')*G''$, so $(G*G')*G''=G*(G'*G'')$.

(ii) Suppose $v \in S \in G*G'$. Then $A \in^* G_A$ and $A \in^* G'_A$
for any $A \sub S$. Applying this to $A=A'v$ for any $A' \sub S\sm v$
gives $A' \in^* G_A(v)$ and $A' \in^* G'_A(v)$. Thus 
$S \sm v \in G_{S\sm v^\le}*G'_{S\sm v^\le}*G_{S^\le}(v)*G'_{S^\le}(v)$.
Conversely, suppose that $v \in S \in G \cup G'$ and 
$S \sm v \in G_{S\sm v^\le}*G'_{S\sm v^\le}*G_{S^\le}(v)*G'_{S^\le}(v)$.
Then $A \in^* G_A$, $A \in^* G'_A$, $Av \in^* G_{Av}$ and $Av \in^* G'_{Av}$
for any $A \sub S\sm v$, so $S \in G*G'$.%
\COMMENT{
Formerly:
(ii) Suppose $S \in (G*G')(v)$. Then $Sv \in G*G'$, so $A \in^* G_A$ and $A \in^* G'_A$
for any $A \sub Sv$. Applying this to $A=A'v$ for any $A' \sub S$
gives $A' \in^* G(v)_{A'}$ and $A' \in^* G'(v)_{A'}$. Thus $S \in G*G'*G(v)*G'(v)$.
We deduce that $(G*G')(v)_S \sub (G*G'*G(v)*G'(v))_S$.
Conversely, suppose that $S \in (G \cup G')(v)$ and $S \in G*G'*G(v)*G'(v)$.
Then $A \in^* G_A$, $A \in^* G'_A$, $A \in^* G(v)_A$ and $A \in^* G'(v)_A$
for any $A \sub S$, so $Sv \in G*G'$, i.e.\ $S \in (G*G')(v)$.
Therefore $(G*G')(v)_S = (G*G'*G(v)*G'(v))_S$.
}

(iii) Suppose $\es \ne A \sub S \in G$. Then $A \in G_A$ and $G'_A$ is undefined,
since $G,G'$ are separate. Therefore $S \in G*G'$, so $G \sub G*G'$.
Similarly $G' \sub G*G'$, so $G*G' = G \cup G'$.
Next we note that $G' \sub H[G]$:
if $S \in G'$ then $S \in H$ and $G_A$ is undefined for any $\es \ne A \sub S$.
Similarly $G \sub H[G']$, so $H[G][G']$ and $H[G'][G]$ are well-defined.
Now suppose $S \in H[G][G']$.
Since $S \in H[G]$, for any $A \sub S$ we have $A \in^* G_A$.
Since $S \in H[G][G']$, for any $A \sub S$ we have $A \in^* G'_A$.
This shows that $S \in H[G'][G]$ and $S \in H[G' \cup G]$.
Applying the same argument to any $S \in H[G'][G]$
we deduce that $H[G][G']=H[G'][G] \sub H[G \cup G']$.
Conversely, if $S \in  H[G \cup G']$ then
for any $A \sub S$ we have $A \in^* G_A \cup G'_A$.
Since $G,G'$ are separate, we have $A \in^* G_A$ and $A \in^* G'_A$.
Therefore $S \in H[G][G']=H[G'][G]$.
It follows that $H[G][G']=H[G'][G]=H[G \cup G']$.

(iv) We first show that $H[G][G*G'] \sub H[G'][G*G']$.
Suppose $S \in H[G][G*G']$. Then $S \in H[G] \sub H$.
Since $S \in H[G][G*G']$, for any $A \sub S$ we have $A \in^* (G*G')_A$.
Consider any $A \sub S$ such that $G'_A$ is defined.
Then $(G*G')_A$ is defined, so $A \in (G*G')_A \sub G'_A$. Therefore $S \in H[G']$.
Now consider any $A \sub S$ such that $(G*G')_A$ is defined.
Since $S \in H[G][G*G']$ we have $A \in (G*G')_A$. Therefore $S \in H[G'][G*G']$.
Similarly, $H[G'][G*G'] \sub H[G][G*G']$, so equality holds.

(v) First we show that $H[G*G'] \sub H[G][G']$.
Suppose $S \in H[G*G']$. Then $S \in H$. Also, for any $A \sub S$
we have $A \in^* (G*G')_A$, and so $A \in^* G_A$ and $A \in^* G'_A$.
This implies that $S \in H[G]$, and then that $S \in H[G][G']$.
Now we show that $H[G][G']\sub H[G*G']$.
Suppose $S \in H[G][G']$. Then $S \in H[G] \sub H$.
Since $S \in H[G]$, for any $A \sub S$ we have $A \in^* G_A$.
Since $S \in H[G][G']$, for any $A \sub S$ we have $A \in^* G'_A$.
Thus for any $A' \sub A \sub S$ we have $A' \in^* G_{A'}$
and $A' \in^* G'_{A'}$, so $A \in^* (G*G')_A$.
Therefore $S \in H[G*G']$.

(vi) We first note that $G*G'=G[G']$ is immediate from
Definitions \ref{def-*} and \ref{def-restrict}.
Then $H[G][G[G']]=H[G'][G[G']]$ follows from (iv).
Now we show that $H[G][G[G']]=H[G][G']$.
Suppose that $S \in H[G][G[G']]$. Then $S \in H[G]$.
Consider any $A \sub S$ such that $G'_A$ is defined.
Since $S \in H[G][G[G']]$ we have $A \in G[G']_A$, and so $A \in G'_A$.
Therefore $S \in H[G][G']$.
Conversely, suppose that $S \in H[G][G']$.
Consider any $A \sub S$ such that $G[G']_A$ is defined.
Then $G_A$ is defined, so $A \in G_A$, since $S \in H[G]$.
Also, for any $A' \sub A$ with $G'_{A'}$ defined we have $A' \in G'_{A'}$, since $S \in H[G][G']$.
Therefore $A \in G[G']_A$. This shows that $S \in H[G][G[G']]$. 
The second statement is immediate. \qed

\subsection{Hypergraph regularity properties}

In this subsection we record some useful properties of hypergraph regularity,
analogous to the standard facts we mentioned earlier for graph regularity.
Similar results can be found e.g.\ in \cite{DHNR,FR}, but with stronger assumptions
on the hierarchy of parameters. However, with the same proof,
we obtain Lemma \ref{3neighbour1} under weaker assumptions on the parameters,
which will be crucial to the proof of Lemma \ref{3exceptional}.
We start with two analogues to Lemma \ref{2neighbour},
the first concerning graphs that are neighbourhoods of a vertex,
and the second sets that are neighbourhood of a pair of vertices.

\begin{lemma}\label{3neighbour1} {\bf (Vertex neighbourhoods)} \index{vertex neighbourhood}
Suppose $0 < \eps \ll d$ and $0 < \eta \ll \eta' \ll d$
and $G$ is a $3$-partite $3$-complex on $V = V_1 \cup V_2 \cup V_3$
with all densities $d_S(G) > d$.
Suppose also that $G_{13}$, $G_{12}$ are $\eps$-regular,
and $G_{23}$, $G_{123}$ are $\eta$-regular. Then
for all but at most $(4\eps+2\eta')|G_1|$ vertices $v \in G_1$
we have $|G(v)_j| = (1 \pm \eps)d_{1j}(G)|G_i|$ for $j=2,3$
and $G(v)_{23}$ is an $\eta'$-regular graph of relative density
$d_{23}(G(v)) = (1 \pm \eta')d_{123}(G)d_{23}(G)$.
\COMMENT{
Delicate proof: have to carry $|A^{v_i}_{23}|$ rather than estimating,
otherwise $\eps$ would interfere with $\eta$. No counting lemma required:
$T_{123}(A)$ is sum over $v_i$ and $T_{123}(G) \le |G_1||G_2||G_3|$.
}
\end{lemma}

\nib{Proof.}
First we apply Lemma \ref{2neighbour} to see that all but at most $4\eps|G_1|$
vertices in $G_1$ have degree $(d_{1j}(G) \pm \eps)|G_j|$ in $G_{1j}$, for $j=2,3$.
Let $G'_1$ be the set of such vertices. It suffices to show the claim that
all but at most $2\eta'|G_1|$ vertices $v \in G'_1$ have the following property:
if $A^v_2 \sub G(v)_2$ and $A^v_3 \sub G(v)_3$ are sets
with $|A^v_2| > \eta'|G(v)_2|$ and $|A^v_3| > \eta'|G(v)_3|$,
then the bipartite subgraph $A^v_{23} \sub G_{23}$ spanned by $A^v_2$ and $A^v_3$
has $|A^v_{23}| = (d_{23}(G) \pm \eta)|A^v_2||A^v_3|$ edges,
and the bipartite subgraph $A(v)_{23} \sub G(v)_{23}$ spanned by $A^v_2$ and $A^v_3$
has $|A(v)_{23}| = (1 \pm \eta'/2)d_{123}(G)|A^v_{23}|$ edges.

Suppose for a contradiction that this claim is false.
Note that for any $v \in G'_1$ and sets $A^v_2 \sub G(v)_2$, $A^v_3 \sub G(v)_3$
with $|A^v_2| > \eta'|G(v)_2|$, $|A^v_3| > \eta'|G(v)_3|$ we have
$|A^v_j| > \eta'(d_{1j}(G)-\eps)|G_j| > \eta|G_j|$ for $j=2,3$
so $|A^v_{23}| = (d_{23}(G) \pm \eta)|A^v_2||A^v_3|$ since $G_{23}$ is $\eta$-regular.
Then without loss of generality, we can assume that we have vertices
$v_1,\cdots,v_t \in G'_1$ with $t>\eta'|G_1|$,
and sets $A^{v_i}_2 \sub G(v_i)_2$, $A^{v_i}_3 \sub G(v_i)_3$
with $|A^{v_i}_2| > \eta'|G(v_i)_2|$ and $|A^{v_i}_3| > \eta'|G(v_i)_3|$,
such that $|A(v_i)_{23}| < (1 - \eta'/2)d_{123}(G)|A^{v_i}_{23}|$ for $1 \le i \le t$.
Define tripartite graphs
$A^i = A^{v_i}_{23} \cup \{v_ia: a \in A^{v_i}_2 \cup A^{v_i}_3\}$
and $A = \cup_{i=1}^t A_i$.

We can count the number of triangles in these graphs as
$T_{123}(A) = \sum_{i=1}^t T_{123}(A^i) = \sum_{i=1}^t |A^{v_i}_{23}|$.
Now $t>\eta'|G_1|$,
$|A^{v_i}_{23}| > (d_{23}(G) - \eta)|A^{v_i}_2||A^{v_i}_3|$,
$d_{23}(G)>d$, $|A^{v_i}_j|>\eta'|G(v_i)_j|$
and $|G(v_i)_j|>(d-\eps)|G_j|$ for $1 \le i \le t$, $j=2,3$, so
\[ T_{123}(A) > \eta'|G_1| \cdot (d-\eta) \cdot \eta'(d-\eps)|G_2| \cdot \eta'(d-\eps)|G_3|
> \eta |G_1||G_2||G_3| \ge \eta T_{123}(G).\]
Since $G_{123}$ is $\eta$-regular we have $\frac{|G \cap T_{123}(A)|}{|T_{123}(A)|} = d_{123}(G) \pm \eta$.
Therefore
\[|G \cap T_{123}(A)| > (d_{123}(G)-\eta)|T_{123}(A)| = (d_{123}(G)-\eta) \sum_{i=1}^t |A^{v_i}_{23}|.\]
But we also have
\[|G \cap T_{123}(A)| = \sum_{i=1}^t |A(v_i)_{23}|
< \sum_{i=1}^t (1-\eta'/2)d_{123}(G)|A^{v_i}_{23}|
< (d_{123}(G)-\eta) \sum_{i=1}^t |A^{v_i}_{23}|,\]
contradiction. This proves the required claim. \qed

\begin{lemma}\label{3neighbour2} {\bf (Pair neighbourhoods)} \index{pair neighbourhoods}
Suppose $0 < \eps \ll \eps' \ll d$ and $G$ is an $\eps$-regular
$3$-partite $3$-complex on $V = V_1 \cup V_2 \cup V_3$
with all densities $d_S(G) > d$.
Then for all but at most $\eps' |G_{12}|$ pairs $uv \in G_{12}$
we have $|G(uv)_3| = (1 \pm \eps')d_{123}(G)d_{13}(G)d_{23}(G)|G_3|$.
\COMMENT{Separate parameters? No need here.}
\end{lemma}

\nib{Proof.}
Introduce another parameter $\eta$ with $\eps \ll \eta \ll \eps'$.
%By Lemma \ref{2neighbour}, all but at most $2\eps|G_1|$ vertices $v \in G_1$
%have degree $(d_{13}(G) \pm \eps)|G_3|$ in $G_{13}$.
%Let $G'_1$ be the set of such vertices $v \in G_1$.
By Lemma \ref{3neighbour1},
for all but at most $6\eps|G_1|$ vertices $v \in G_1$
we have $|G(v)_i| = (1 \pm \eps)d_{1i}(G)|G_i|$ for $i=2,3$
and $G(v)_{23}$ is an $\eta$-regular graph of relative density
$d_{23}(G(v)) = (1 \pm \eta)d_{123}(G)d_{23}(G)$.
Let $G'_1$ be the set of such vertices $v \in G_1$.
Then for any $v \in G'_1$, applying Lemma \ref{2neighbour} to $G(v)_{23}$,
we see that for all but at most $2\eta|G(v)_2| \le 2\eta|G_2|$
vertices in $u \in G(v)_2$, the degree of $u$ in $G(v)_{23}$ satisfies
\begin{align*}
|G(uv)_3| & = (d_{23}(G(v)) \pm \eta)|G(v)_3|
= ((1 \pm \eta)d_{123}(G)d_{23}(G) \pm \eta)(d_{13}(G) \pm \eps)|G_3| \\
& = (1 \pm \eps')d_{123}(G)d_{13}(G)d_{23}(G)|G_3|.
\end{align*}
Since $|G_{12}| = d_{12}(G)|G_1||G_2| > d|G_1||G_2|$,
this estimate holds for all pairs $uv \in G_{12}$ except for at most
$6\eps|G_1| \cdot |G_2| + |G_1| \cdot 2\eta|G_2| < \eps'|G_{12}|$. \qed

Next we give an analogue of Lemma \ref{2restrict},
showing that regularity is preserved by restriction.

\begin{lemma}\label{3restrict} {\bf (Regular restriction)} \index{regular restriction}
Suppose $0 < \eps \ll d$,
$G$ is a $3$-partite $3$-complex on $V = V_1 \cup V_2 \cup V_3$
with all densities $d_S(G) > d$, $G_{123}$ is $\eps$-regular,
and $J \sub G$ is a $2$-complex with $|J^*_{123}|>\sqrt{\eps}|G^*_{123}|$.
Then $G[J]_{123}$ is $\sqrt{\eps}$-regular and $d_{123}(G[J])=(1\pm\sqrt{\eps})d_{123}(G)$.
\COMMENT{
1. No regularity assumption on $J$.
2. Lower densities come from $J$ or graph restriction $G[J]$.
3. Change $c$ to $d$ here and elsewhere.
4. Formerly had $|G_{123}| > c|V_1||V_2|V_3|$ but we didn't yet relate relative to absolute density.
}
\end{lemma}

\nib{Proof.}
Since $G_{123}$ is $\eps$-regular,
$|G[J]_{123}| = |G \cap J^*_{123}| = (d_{123}(G) \pm \eps)|J^*_{123}|$
and $d_{123}(G[J]) = |G[J]_{123}|/|G[J]^*_{123}|=|G[J]_{123}|/|J^*_{123}|
= d_{123}(G) \pm \eps$.
Now consider any subcomplex $A$ of $J$ with $|A^*_{123}| > \sqrt{\eps} |J^*_{123}|$.
Then $|A^*_{123}| > \eps|G^*_{123}|$, so since $G_{123}$ is $\eps$-regular,
$|G[J] \cap A^*_{123}| = |G \cap A^*_{123}| = (d_{123}(G) \pm \eps)|A^*_{123}|
= (d_{123}(J[G]) \pm \sqrt{\eps})|A^*_{123}|$, i.e.\ $G[J]_{123}$ is $\sqrt{\eps}$-regular. \qed

It is worth noting the special case of Lemma \ref{3restrict} when $J$ is a $1$-complex.
Then $G[J]$ is obtained from $G$ by discarding some vertices,
i.e.\ a restriction according to the traditional definition.
In particular, we see that regularity implies vertex regularity
(the weak property mentioned at the beginning of Section \ref{3reg}).
We also record the following consequence of Lemma \ref{3restrict}.

\begin{coro}\label{3restrict'}
Suppose $0 < \eps \ll d$,
$G$ is a $3$-partite $3$-complex on $V = V_1 \cup V_2 \cup V_3$
with all densities $d_S(G) > d$ and $G_{123}$ is $\eps$-regular.
Suppose also $0 < \eta \ll d$, $J \sub G$ is a $2$-complex
with all densities $d_S(J)>d$ and $J_{12}$, $J_{13}$, $J_{23}$ are $\eta$-regular.

Then $G[J]_{123}$ is $\sqrt{\eps}$-regular and $d_{123}(G[J])=(1\pm\sqrt{\eps})d_{123}(G)$.
\end{coro}

\nib{Proof.} We have %(equation \ref{eq:tri} in the previous section)
$|J^*_{123}|=|T_{123}(J)|=(1 \pm 8\eta)d_{12}(J)d_{13}(J)d_{23}(J)|J_1||J_2||J_3|
> \frac{1}{2} d^6 |V_1||V_2||V_3| > \sqrt{\eps}|G^*_{123}|$
by the triangle counting lemma (\ref{eq:tri}).
The result now follows from Lemma \ref{3restrict}. \qed

Next we note a simple relationship between relative and absolute densities.

\begin{lemma}\label{3absolute}
Suppose $0 < \eps \ll d$,
$G$ is a $3$-partite $3$-complex on $V = V_1 \cup V_2 \cup V_3$
with all densities $d_S(G) > d$ and $G$ is $\eps$-regular.
Then $d(G_{123}) = (1 \pm 8\eps)\prod_{S \sub 123}d_S(G)$.
\end{lemma}

\nib{Proof.} $d(G_{123})=\frac{|G_{123}|}{|V_1||V_2||V_3|}
= \frac{|G_{123}|}{|T_{123}(G)|} \cdot \frac{|T_{123}(G)|}{|V_1||V_2||V_3|}
= d_{123}(G) \cdot (1 \pm 8\eps)\prod_{S \subn 123} d_S(G)$ by (\ref{eq:tri}). \qed
%By definition we have $|G_i|=d_i|V_i|$
%and $|G_{ij}|=d_{ij}|G_i||G_j|$ for $1 \le i,j \le 3$, $i \ne j$.
%The triangle counting lemma (\ref{eq:tri})
%gives $T_{123}(G) = (1 \pm 8\eps)d_{12}(G)d_{13}(G)d_{23}(G)|G_1||G_2||G_3|$,
%so $d(T_{123}(G)) = (1 \pm 8\eps)\prod_{S \subn 123} d_S(G)$.
%Since $d(G_{123})=\frac{|G_{123}|}{|V_1||V_2||V_3|}
%= \frac{|G_{123}|}{|T_{123}(G)|}\frac{|T_{123}(G)|}{|V_1||V_2||V_3|}
%= d_{123}(G) d(T_{123}(G))$ the result follows. \qed

The following more technical lemma will be useful in the proof of Lemma \ref{3marked}.
Later we will give a more general proof that is slicker, but conceptually more difficult,
as it uses the `plus complex' of Definition \ref{def-plus}.
For the convenience of the reader, in the $3$-graph case we will use a
proof that is somewhat pedestrian, but perhaps easier to follow.

\begin{lemma}\label{3technical}
Suppose $0 < \eps \ll \eps' \ll d$,
$G$ is a $4$-partite $3$-complex on $V = V_1 \cup V_2 \cup V_3 \cup V_4$
with all densities $d_S(G) > d$ and $G$ is $\eps$-regular.
\begin{itemize}
\item[(i)]
For any $P \in G_{123}$ and subcomplex $I$ of $123^<$,
let $G_{P,I}$ be the set of vertices $v \in G_4$
such that $P_S \cup v \in G_{S \cup 4}$ for all $S \in I$.%
\COMMENT{Removed $\es \ne$... doesn't matter.}
Let $B_I$ be the set of $P \in G_{123}$ such that we do not have
$|G_{P,I}| = (1 \pm \eps')|V_4| \prod_{\es \ne S \in I} d_{S\cup 4}(G)$.
Then $|B_I| < \eps' |G_{123}|$.
\item[(ii)]
For any $P'\in G_{12}$ and subcomplex $I'$ of $12^<$
let $G'_{P',I'}$ be the set of vertices $v \in G_4$
such that $P'_{S'} \cup v \in G_{S' \cup 4}$ for all $S' \in I'$.
Let $B'_{I'}$ be the set of $P' \in G_{12}$ such that we do not have
$|G'_{P',I'}| = (1 \pm \eps')|V_4| \prod_{\es \ne S' \in I'} d_{S\cup 4}(G)$.
Then $|B'_{I'}| < \eps' |G_{12}|$.
\end{itemize}
\end{lemma}

\nib{Proof.} Introduce auxiliary constants with a hierarchy
$\eps \ll \eps_1 \ll \eps_2 \ll \eps_3 \ll \eps'$.
We consider selecting the vertices $P_1$, $P_2$, $P_3$ of $P$ in turn,
at each step identifying some exceptional sets $P$ for which
the stated estimate on $G_{P,I}$ might fail.
First we choose $P_1$ so that $|G(P_1)_i| = (1 \pm \eps)d_{1i}(G)|G_i|$
and $G(P_1)_{ij}$ is an $\eps_1$-regular graph of relative density
$d_{ij}(G(P_1)) = (1 \pm \eps_1)d_{1ij}(G)d_{ij}(G)$
for distinct $i,j$ in $\{2,3,4\}$. Applying Lemma \ref{3neighbour1}
with $\eta=\eps$ and $\eta'=\eps_1/4$, we see that this holds for all
but at most $\eps_1|G_1|$ vertices $P_1 \in G_1$. Then the number of exceptional
sets $P$ at this stage is at most
$\eps_1|G_1||V_2||V_3|=\eps_1 d_1(G) d(G_{123})^{-1}|G_{123}| < \sqrt{\eps_1}|G_{123}|$.

Now let $J^1 \sub G_{234^<}$ be the $2$-complex defined as follows.
We define the singleton parts by
$J^1_4$ equals $G(P_1)_4$ if $1 \in I$ or $G_4$ if $1 \notin I$,
$J^1_2$ equals $G(P_1)_2$ if $12 \in I$ or $G_2$ if $12 \notin I$, and
$J^1_3$ equals $G(P_1)_3$ if $13 \in I$ or $G_3$ if $13 \notin I$.
We define the graphs by restriction to the singleton parts of the following:
$G(P_1)_{24}$ if $12 \in I$ or $G_{24}$ if $12 \notin I$,
$G(P_1)_{34}$ if $13 \in I$ or $G_{34}$ if $13 \notin I$,
$G(P_1)_{23}$ if $123 \in I$ or $G_{23}$ if $123 \notin I$.
Then $J^1$ is $\sqrt{\eps_1}$-regular by Lemma \ref{2restrict}.
The graph densities $d_{ij}(J^1)$ are either $(1 \pm \eps_1)d_{1ij}(G)d_{ij}(G)$
or $(1 \pm \eps_1)d_{ij}(G)$, according as we restrict $G(P_1)_{ij}$ or $G_{ij}$.

Let $G^1=G_{234^\le}[J^1]$. Then $G^1_{234}$ is $\eps_1$-regular with
$d_{234}(G^1)=(1\pm\eps_1)d_{234}(G)$ by Corollary \ref{3restrict'}.
Next we choose $P_2$ so that $|G^1(P_2)_i| = (1 \pm \eps_2)d_{2i}(G^1)|G^1_i|$ for $i=3,4$
and $G^1(P_2)_{34}$ is an $\eps_2$-regular graph of relative density
$d_{34}(G^1(P_2)) = (1 \pm \eps_2)d_{234}(G^1)d_{34}(G^1)$. By Lemma \ref{3neighbour1}
this holds for all but at most $\eps_2|G^1_2|$ vertices $P_2 \in G^1_2$, so similarly to above,
the number of exceptional sets $P$ at this stage is at most $\sqrt{\eps_2}|G_{123}|$.
Let $J^2 \sub G^1_{34^\le}$ be the $2$-complex in which
$J^2_4$ is $G^1(P_2)_4$ if $2 \in I$ or $G^1_4$ if $2 \notin I$,
$J^2_3$ is $G^1(P_2)_3$ if $23 \in I$ or $G^1_3$ if $23 \notin I$, and
$J^2_{34}$ is $G^1(P_2)_{34}$ if $23 \in I$ or $G^1_{34}$ if $23 \notin I$.
Then $J^2_{34}$ is $\sqrt{\eps_2}$-regular by Lemma \ref{2restrict},
with $d_{34}(J^2)$ either $(1 \pm \eps_2)d_{234}(G^1)d_{34}(G^1)$ or
$(1 \pm \eps_2)d_{34}(G^1)$, according as we restrict $G^1(P_2)_{34}$ or $G^1_{34}$.

Now we choose $P_3$ so that $|J^2(P_3)_4|=(1 \pm \eps_3)d_{34}(J^2)|J^2_4|$.
By Lemma \ref{2neighbour} this holds for all but at most $\eps_3|J^2_3|$ vertices
$P_3 \in J^2_3$, giving at most $\sqrt{\eps_3}|G_{123}|$ exceptional sets $P$ here.
In total, the number of exceptional sets at any stage is fewer than $\eps'|G_{123}|$.
By construction, $G_{P,I}$ equals $J^2(P_3)_4$ if $3 \in I$ or $J^2_4$ if $3 \notin I$.
If $P$ is not exceptional then we can estimate $|G_{P,I}|$ by tracing back through the stages.
At stage 3 we multiply $|J^2_4|$ by $(1 \pm \eps_3)d_{34}(J^2)$ if $3 \in I$,
where $d_{34}(J^2)$ is $(1 \pm \eps_2)d_{234}(G^1)d_{34}(G^1)$ if $23 \in I$
or $(1 \pm \eps_2)d_{34}(G^1)$ if $23 \notin I$,
where $d_{234}(G^1)=(1\pm\eps_1)d_{234}(G)$
and $d_{34}(G^1)$ is $(1 \pm \eps_1)d_{134}(G)d_{34}(G)$ if $13 \in I$
or $(1 \pm \eps_1)d_{34}(G)$ if $13 \notin I$.
Thus we obtain a factor of $d_{S4}(G)$ whenever $3 \in S \in I$.
At stage 2, we obtain $|J^2_4|$ from $|G^1_4|$ by multiplying by $(1 \pm \eps_2)d_{24}(G^1)$
if $2 \in I$, where $d_{24}(G^1)$ is $(1 \pm \eps_1)d_{124}(G)d_{24}(G)$ if $12 \in I$
or $(1 \pm \eps_1)d_{24}(G)$ if $12 \notin I$.
Thus we obtain a factor of $d_{S4}(G)$ whenever $2 \in S \in I$, $3 \notin S$.
Finally, at stage 1, we obtain $|G^1_4|$ from $|G_4|$ by multiplying by
$(1 \pm \eps_1)d_{14}(G^1)$ if $1 \in I$.
Combining all factors we obtain
$|G_{P,I}| = (1 \pm \eps')|V_4| \prod_{\es \ne S \in I} d_{S4}(G)$.

This proves (i). The proof of (ii) is similar and much simpler
(alternatively it could be deduced from (i)).
We consider selecting the vertices $P'_1$ and $P'_2$ of $P'$ in turn.
We choose $P'_1$ so that $|G(P'_1)_4| = (1 \pm \eps)d_{14}(G)|G_4|$.
We let $G'_4$ be $G_4$ if $1 \notin I$ or $G(P'_1)_4$ if $1 \in I$,
and $G'_{24}$ be the restriction of $G_{24}$ to $G_2$ and $G'_4$.
Then $G'_{24}$ is $\eps_1$-regular with $d_{24}(G')=(1\pm \eps)d_{24}(G)$.
We choose $P'_2$ so that $|G'_{24}(P'_2)|=(1 \pm \eps_1)d_{24}(G')|G'_4|$.
Then $G'_{P',I'}$ is $G'_{24}(P'_2)$ if $2 \in I$ or $G'_4$ if $2\notin I$.
Now $|G'_{P',I'}|$ is obtained from $|G_4|$ by multiplying by
$(1 \pm \eps)d_{14}(G)$ if $1 \in I$ and
$(1 \pm \eps_1)(1\pm \eps)d_{24}(G)$ if $2 \in I$,
so $|G'_{P',I'}| = (1 \pm \eps')|V_4| \prod_{\es \ne S' \in I'} d_{S\cup 4}(G)$.
It is clear that there are at most $\eps'|G_{12}|$ exceptional sets $P'$. \qed

Finally we give another formulation of the neighbourhood Lemmas \ref{3neighbour1}
and \ref{3neighbour2}, showing that most vertices and pairs are close to `average'.

\begin{lemma}\label{3average}\index{average}{\bf (Averaging)}
Suppose $0 < \eps \ll \eps' \ll d$,
$G$ is a $3$-partite $3$-complex on $V = V_1 \cup V_2 \cup V_3$
with all densities $d_S(G) > d$ and $G$ is $\eps$-regular. Then
\begin{itemize}
\item[(i)] for all but at most $\eps'|G_1|$ vertices $v \in G_1$
we have $|G(v)_{23}| = (1 \pm \eps')|G_{123}|/|G_1|$,
\item[(ii)] for all but at most $\eps'|G_{12}|$ pairs $uv \in G_{12}$
we have $|G(uv)_3| = (1 \pm \eps')|G_{123}|/|G_{12}|$.
\end{itemize}
\end{lemma}

\nib{Proof.} By Lemma \ref{3neighbour1} with $\eta=\eps$ and $\eta'=\eps'/4$,
for all but at most $\eps'|G_1|$ vertices $v \in G_1$,
$|G(v)_i| = (1 \pm \eps)d_{1i}(G)|G_i|= (1 \pm \eps)d_{1i}(G)d_i(G)|V_i|$
for $i=2,3$ and $d_{23}(G(v)) = (1 \pm \eps'/4)d_{123}(G)d_{23}(G)$.
For such $v$ we have $|G(v)_{23}|=d_{23}(G(v))|G(v)_2||G(v)_3|
=(1\pm \eps'/3)|V_2||V_3|\prod_{S \sub 123,S \ne 1} d_S(G)$.
Also, $|G_{123}| = (1 \pm 8\eps)|V_1||V_2||V_3|\prod_{S \sub 123}d_S(G)$
by Lemma \ref{3absolute}, so $|G(v)_{23}|= (1 \pm \eps')|G_{123}|/|G_1|$, giving (i).
For (ii), Lemma \ref{3absolute} gives
$|G_{123}|/|G_{12}| = (1 \pm 8\eps)|V_3|\prod_{3 \in S \sub 123}d_S(G)$.
Then by Lemma \ref{3neighbour2}, replacing $\eps'$ with $\eps'/2$,
for all but at most $\eps' |G_{12}|$ pairs $uv \in G_{12}$, 
$|G(uv)_3| = (1 \pm \eps'/2)d_{123}(G)d_{13}(G)d_{23}(G)|G_3|
= (1 \pm \eps'/2)|V_3|\prod_{3 \in S \sub 123}d_S(G) = (1 \pm \eps')|G_{123}|/|G_{12}|$. \qed

\subsection{Good vertices}\index{good}

We start the analysis of the algorithm by showing that most free vertices are good.
Our first lemma handles the definitions for regularity and density in the algorithm.

\begin{lemma}$ $ \label{3exceptional}\index{exceptional}\index{$E_x(t-1)$}
The exceptional set $E_x(t-1)$ defined by $(*_{\ref{3alg}})$ satisfies $|E_x(t-1)| < \eps_*|F_x(t-1)|$,
and $F_S(t)$ is $\eps_{\nu'_S(t),1}$-regular with $d_S(F(t)) \ge d_u$ for every $S \in H$.
\COMMENT{
1. No density formula here, but see lemma 3crude.
2. Changed $uv$ to $vw$ to avoid $d_u$ notation conflict.}
\end{lemma}

\nib{Proof.} We argue by induction on $t$. At time $t=0$ the first statement is vacuous.
The second statement at time $0$ follows from the fact that $F_S(0)=G_S$
and our assumption that $(G,M)$ is  $(\eps,\eps',d_2,\theta,d_3)$-super-regular:
condition (i) in Definition \ref{def-3super} tells us that $G_S$ is $\eps$-regular,
with $d_S(G) \ge d_{|S|}$ if $|S|=2,3$. Also $d_S(G) = 1$ if $|S|=0,1$,
as we assumed that $G_i=V_i$ in the hypotheses of Theorem \ref{3blowup}.
Now suppose $t \ge 1$ and $\es \ne S \in H$ is unembedded, so $x \notin S$.
We consider various cases for $S$ to establish the bound on the exceptional set
and the regularity property, postponing the density bound until later in the proof.

We start with the case when $S \in H(x)$. Suppose first that $S=vw$ with $xvw \in H$.
By induction $F_{S'}(t-1)$ is $\eps_{\nu'_{S'}(t-1),1}$-regular
and $d_{S'}(F(t-1)) \ge d_u$ for every $S' \sub xvw$.
Write $\nu=\max\{\nu'_{xv}(t-1),\nu'_{xw}(t-1)\}$
and $\nu^*=\max\{\nu'_{vw}(t-1),\nu'_{xvw}(t-1)\}$.
We claim that $\nu'_{vw}(t)>\nu^*$. This holds by Lemma \ref{3observe}:
(iv) gives $\nu_{vw}(t)>\nu_{vw}(t-1)$,
(v) gives $\nu'_{vw}(t)>\nu'_{vw}(t-1)$, and
(vi) gives $\nu'_{vw}(t)>\nu'_{xvw}(t-1)$.
Now applying Lemma \ref{3neighbour1}, for all but at most
$(4\eps_{\nu,1}+2\eps_{\nu^*,2})|F_x(t-1)|$ vertices $y \in F_x(t-1)$
we have $|F_v(t)|=|F_{xv}(t-1)(y)|=(1\pm\eps_{\nu,1})d_{xv}(F(t-1))|F_v(t-1)|$,
$|F_w(t)|=|F_{xw}(t-1)(y)|=(1\pm\eps_{\nu,1})d_{xw}(F(t-1))|F_w(t-1)|$,
and $F_{vw}(t)=F_{xvw}(t-1)(y)$ is an $\eps_{\nu^*,2}$-regular graph
with $d_{vw}(F(t))=(1\pm\eps_{\nu^*,2})d_{xvw}(F(t-1))d_{vw}(F(t-1))$.
Since $\nu'_{vw}(t)>\nu^*$ we have $(*_{\ref{3alg}})$ when $S=vw$ for such $y$.
Note that it is important for this argument that Lemma \ref{3neighbour1}
makes no assumption of any relationship between $\nu$ and $\nu^*$.
For future reference we also note that the density bounds at time $t-1$
imply that $d_{vw}(F(t)) > d_u^2/2$; we will show a lower bound of $d_u$ later,
but this interim bound will be useful before then.

The argument when $S=\{v\} \in H(x)$ has size $1$ is similar and more straightforward.
By Lemma \ref{2neighbour},
for all but at most $2\eps_{\nu'_{xv}(t-1),1}|F_x(t-1)|$ vertices $y \in F_x(t-1)$
we have $|F_v(t)|=|F_{xv}(t-1)(y)|=(1\pm\eps_{\nu'_{xv}(t-1),1})d_{xv}(F(t-1))|F_v(t-1)|$.
Also, we have $\nu'_v(t) > \nu'_{xv}(t-1)$ by Lemma \ref{3observe}(vi),
so $(*_{\ref{3alg}})$ holds when $S=\{v\}$ for such $y$.
We also note for future reference that $d_v(F(t))>d_u^2/2$.
In the argument so far we have excluded at most $6\eps_{12D,2}|F_x(t-1)|$
vertices $y \in F_x(t-1)$ for each of at most $3D$ sets $S \in H(x)$
with $|S|=1$ or $|S|=2$; this gives the required bound on $E_x(t-1)$.
We also have the required regularity property of $F_S(t)$,
but for now we postpone showing the density bounds.

Next we consider the case when $S \in H$ and $S \notin H(x)$.
If $S=\{v\}$ has size $1$ then $|F_v(t)|=|F_v(t-1)\sm y| \ge |F_v(t-1)|-1$,
so $d_v(F(t)) \ge d_v(F(t-1))-1/n > d_u/2$, say.
Next suppose that $S=vw$ has size $2$.
Recall that Lemma \ref{build-update} gives $F_S(t) = F_{S^\le}(t-1)[F_{S^<}(t)]_S$.
Then $F_{vw}(t)$ is the bipartite subgraph of $F_{vw}(t-1)$ induced by $F_v(t)$ and $F_w(t)$.
We have $F_v(t)=F_{xv}(t-1)(y)$ if $xv \in H$ or $F_v(t-1) \sm y$ if $xv \notin H$.
Similarly, $F_w(t)=F_{xw}(t-1)(y)$ if $xw \in H$ or $F_w(t-1) \sm y$ if $xw \notin H$.
Since we choose $y \notin E_x(t-1)$, if $xv \in H$ then
$|F_v(t)|=(1 \pm \eps_{\nu'_v(t),1})d_v(F(t-1))d_{xv}(F(t-1)) > \frac{1}{2}d_u|F_v(t-1)|$,
and if $xv \notin H$ then $|F_v(t)|=|F_v(t-1) \sm y| \ge |F_v(t)|-1$.
Similarly, if $xw \in H$ then $|F_w(t)| > \frac{1}{2}d_u|F_w(t-1)|$,
and if $xw \notin H$ then $|F_w(t)|\ge |F_w(t)|-1$.
Now $F_{vw}(t-1)$ is $\eps_{\nu'_{vw}(t-1),1}$-regular, so by Lemma \ref{2restrict},
$F_{vw}(t)$ is $\eps_{\nu'_{vw}(t-1),2}$-regular
and $d_{vw}(F(t)) = (1 \pm \eps_{\nu'_{vw}(t-1),2})d_{vw}(F(t-1))$.
This gives the required regularity property for $F_{vw}(t)$ in the case
that $vw$ intersects $VN_H(x)$, when we have $\nu'_{vw}(t)>\nu'_{vw}(t-1)$
by Lemma \ref{3observe}(v). In fact, we are only required to show
that $F_{vw}(t)$ is $\eps_{\nu'_{vw}(t),1}$-regular, but this stronger
regularity property will be useful for the case when $vw$ and $VN_H(x)$ are disjoint.
Now consider the case that $vw$ and $VN_H(x)$ are disjoint.
Let $t'$ be the most recent time at which we embedded a vertex $x'$ with a neighbour in $vw$.
Then by Lemma \ref{not-local}, $F_{vw^\le}(t)$ is obtained from $F_{vw^\le}(t')$ just by deleting
all sets containing any vertices that are embedded between time $t'+1$ and $t$.
Thus $F_{vw}(t)$ is the bipartite subgraph of $F_{vw}(t')$ spanned by $F_v(t)$ and $F_w(t)$.
Now $F_{vw} (t')$ is $\eps_{\nu'_{vw}(t'-1),2}$-regular,
by the stronger regularity property just mentioned above.
Since $\nu'_{vw}(t) \ge \nu'_{vw}(t') > \nu'_{vw}(t'-1)$,
Lemma \ref{2restrict} gives the required regularity property for $F_{vw}(t)$.
For future reference, we note that in either case the bound $d_{vw}(F(t-1))>d_u$
implies that $d_{vw}(F(t))>d_u/2$. %(say).

Continuing with the case when $S \in H$ and $S \notin H(x)$,
we now suppose that $|S|=3$. Again we use $F_S(t) = F_{S^\le}(t-1)[F_{S^<}(t)]_S$.%
\COMMENT{
Formerly... $= F_S(t-1) \cap T_S(F_{S^<}(t))$, where $T_S(F_{S^<}(t))$
denotes the set of triangles formed by the graphs $F_{S'}(t)$
with $S' \sub S$, $|S'|=2$.
}
If $S' \subn S$, whether $S' \in H(x)$ or $S' \notin H(x)$,
we have shown above that $d_{S'}(F(t))>d_u^2/2$, and if $|S'|=2$ that
$F_{S'}(t)$ is $\eps_{\nu'_{S'}(t),1}$-regular.
Since $F_S(t-1)$ is $\eps_{\nu'_S(t-1),1}$-regular,
Corollary \ref{3restrict'} implies that $F_S(t)$ is $\eps_{\nu'_S(t-1),2}$-regular
and $d_S(F(t)) = (1 \pm \eps_{\nu'_S(t-1),2})d_S(F(t-1))$.
This gives the required regularity property for $F_S(t)$ in the case
that $S$ intersects $VN_H(x)$, when we have $\nu'_S(t)>\nu'_S(t-1)$.
As above, although we are only required to show that $F_S(t)$ is $\eps_{\nu'_S(t),1}$-regular,
this stronger regularity property will be useful for the case when $S$ and $VN_H(x)$ are disjoint.
Suppose $S$ and $VN_H(x)$ are disjoint.
Let $t'$ be the most recent time at which we embedded a vertex $x'$ with a neighbour in $S$.
Then by Lemma \ref{not-local}, $F_{S^\le}(t)$ is obtained from $F_{S^\le}(t')$ just by deleting
all sets containing any vertices that are embedded between time $t'+1$ and $t$.
Equivalently, $F_S(t) = F_S(t')[((F_v(t):v\in S),\{\es\})]$.
Now $F_S(t')$ is $\eps_{\nu'_S(t'-1),2}$-regular,
by the stronger regularity property just mentioned above.
Since $\nu'_S(t) \ge \nu'_S(t') > \nu'_S(t'-1)$,
Corollary \ref{3restrict'} gives the required regularity property for $F_S(t)$.

Now we have established the bound on $E_x(t-1)$ and the regularity properties,
so it remains to show the density bounds.
First we consider any unembedded $S \in H$ with $|S|=3$.
We claim that 
\begin{equation}\label{refwants}
F_S(t) = G_{S^\le}[F_{S^<}(t)]_S.
\end{equation}
To see this we use induction. In the base case $t=0$ we have $F_S(0)=G_S$
and $G_{S^\le}[F_{S^<}(0)] = G_{S^\le}[G_{S^<}]=G_{S^\le}$,
so $G_{S^\le}[F_{S^<}(0)]=G_S$. For $t>0$, Lemma \ref{build-update} gives
$F_S(t) = F_{S^\le}(t-1)[F_{S^<}(t)]_S$, i.e.\ $F_S(t)$ consists of all triples
in $F_S(t-1)$ that form triangles in $F_{S^<}(t)$. The induction hypothesis
gives $F_S(t-1) = G_{S^\le}[F_{S^<}(t-1)]_S$, and so we can write
$F_S(t)=G_{S^\le}[F_{S^<}(t-1)][F_{S^<}(t)]_S$. Now $F_{S^<}(t) \sub F_{S^<}(t-1)$,
so Lemma \ref{*props}(vi) gives $G_{S^\le}[F_{S^<}(t-1)][F_{S^<}(t)]
= G_{S^\le}[F_{S^<}(t-1)[F_{S^<}(t)]] = G_{S^\le}[F_{S^<}(t)]$.
This proves (\ref{refwants}).
Also, we showed above that for every $S' \subn S$,
we have $d_{S'}(F(t))>d_u^2/2$, and if $|S'|=2$ then
$F_{S'}(t)$ is $\eps_{\nu'_{S'}(t),1}$-regular.
Since $G_S$ is $\eps$-regular, Corollary \ref{3restrict'} shows that
$F_S(t)$ is $\sqrt{\eps}$-regular and $d_S(F(t))=(1\pm\sqrt{\eps})d_S(G) > d_3/2$.
This gives the required bound $d_S(F(t)) > d_u$, although we will also
use the stronger bound of $d_3/2$ below.

Next consider any unembedded pair $vw \in H$.
Let $t' \le t$ be the most recent time at which we embedded a vertex $x'$
with $x'vw \in H$, or let $t'=0$ if there is no such vertex $x'$.
Note that we have $t'=t$ if $xvw \in H$.
For $t^* \le t$, let $J(t^*)$ be the $1$-complex $(F_v(t^*),F_w(t^*),\{\es\})$.
We claim that $F_{vw}(t^*) = F_{vw^\le}(t')[J(t^*)]_{vw}$ for $t'\le t^*\le t$.
This follows by induction, similarly to the argument when $|S|=3$.
When $t=t'$ the claim follows from $F_{vw^\le}(t')[J(t')]=F_{vw^\le}(t')$.
For $t^*>t'$, we have $vw \notin H(x)$, so Lemma \ref{build-update} gives
$F_{vw}(t^*) = F_{vw^\le}(t^*-1)[J(t^*)]_{vw}$.
Since $F_{vw}(t^*-1)=F_{vw^\le}(t')[J(t^*-1)]_{vw}$ by induction,
Lemma \ref{*props}(vi) gives $F_{vw}(t^*) = F_{vw^\le}(t')[J(t^*-1)][J(t^*)]_{vw}
= F_{vw^\le}(t')[J(t^*)]_{vw}$, as claimed.
Now we claim that $d_{vw}(F(t)) >  (d_3/4)^{i_t} d_2/2$,
where we temporarily use $i_t$ to denote the number of embedded
vertices $x'$ at time $t$ with $x'vw \in H$.
To see this, we argue by induction, noting that initially
$d_{vw}(F(0))=d_{vw}(G)>d_2$. Also, if $i_t=0$ then
$F_{vw}(t)=G_{vw^\le}[J(t)]_{vw}$, so $d_{vw}(F(t))>d_2/2$ by Lemma \ref{2restrict}.
Now suppose that $i_t>0$, so that $t'$ and $x'$ are defined above.
By induction we have $d_{vw}(F(t'-1)) > (d_3/4)^{i_t-1} d_2/2$.
Also $d_{x'vw}(F(t'-1))>d_3/2$ by the lower bound just proved
for relative densities of triples, so $(*_{\ref{3alg}})$ gives
$d_{vw}(F(t'))=(1 \pm \eps_{\nu'_{vw}(t'),0})d_{vw}(F(t'-1))d_{x'vw}(F(t'-1))
> (d_3/4)^{i_t-1} d_2/2 \cdot d_3/2$. Since $F_{vw}(t) = F_{vw^\le}(t')[J(t)]_{vw}$,
Lemma \ref{2restrict} gives $d_{vw}(F(t)) > d_{vw}(F(t'))/2
> (d_3/4)^{i_t} d_2/2$, as claimed. Since $i_t \le D$
we have $d_{vw}(F(t)) > 2d_3^{2D} d_2$, say.
In particular we have the required bound of $d_{vw}(F(t)) > d_u$.

Finally we consider any unembedded vertex $z$.
Let $t_z\le t$ be the most recent time at which we embedded a neighbour $w=s(t_z)$ of $z$,
or $t_z=0$ if there is no such time.
If $t_z>0$ then by $(*_{\ref{3alg}})$ and the above bound for pair densities we have
$d_z(F(t_z))>d_{wz}(F(t_z-1))d_z(F(t_z-1))/2 > d_3^{2D} d_2 d_z(F(t_z-1))$.
Now we consider cases according to whether $z$ is in the list $L(t-1)$,
the queue $q(t-1)$ or the queue jumpers $j(t-1)$.
Suppose first that $z \in L(t-1)$.
Then the rule for updating the queue in the algorithm
gives $|F_z(t)| \ge \delta'_Q |F_z(t_z)|$.
Next suppose that $z \in j(t-1) \cup q(t-1)$,
and $z$ first joined $j(t') \cup q(t')$ at some time $t'<t$.
Since $z$ did not join the queue at time $t'-1$ we have
$|F_z(t'-1)| \ge \delta'_Q |F_z(t_z)|$.
Also, between times $t'$ and $t$ we only embed vertices that are in the queue
or jumping the queue, or otherwise we would have embedded $z$ before $x$.
Now $|Q(t) \cap X_z| \le \delta_Q Cn$, otherwise we abort the algorithm,
and $|J(t) \cap X_z| \le \sqrt{\delta_Q}n$ by Lemma \ref{3observe}(ii),
so we embed at most $2\sqrt{\delta_Q}n$ vertices in $V_z$
between times $t'$ and $t$. Thus we have catalogued all possible ways
in which the number of vertices free for $z$ can decrease.
It may decrease by a factor of $d_3^{2D} d_2$ when a new $z$-regime starts,
and by a factor $\delta'_Q$ during a $z$-regime before $z$ joins the queue.
Also, if $z$ joins the queue or jumps the queue it may decrease by
at most $2\sqrt{\delta_Q}n$ in absolute size.
Since $z$ has at most $2D$ neighbours, and $|F_z(0)|=|V_z|>n$, we have
$|F_z(t)| \ge (\delta'_Qd_3^{2D} d_2)^{2D}\delta'_Q|V_z| - 2\sqrt{\delta_Q}n > d_u|V_z|$. \qed 

In the preceding proof we needed to track the $\eps$-subscripts in great detail
to be sure that they always fall in the range allowed by our hierarchy.
From now on it will often suffice and be more convenient to use a crude upper bound
of $\eps_*$ for any epsilon parameter. We summarise some useful estimates in the following lemma.

\begin{lemma}\label{3crude}$ $\COMMENT{latex: compact with vspace or redef itemize?}
\begin{itemize}
\item[(i)]
If $\es \ne S \in H(x)$ then $d_S(F(t))=(1\pm \eps_*)d_S(F(t-1))d_{Sx}(F(t-1))$

and $|F_S(t)|= (1 \pm \eps_*)|F_{Sx}(t-1)|/|F_x(t-1)|$.
\item[(ii)]
If $S \notin H(x)$ then $d_S(F(t))=(1\pm \eps_*)d_S(F(t-1))$.
\item[(iii)]
If $S \in H$ then $d(F_S(t))=(1\pm\eps_*)\prod_{T\sub S}d_T(F(t))$.
\item[(iv)]
If $S' \sub S \in H$ then
\[\frac{|F_S(t)|}{|F_{S'}(t)||F_{S\sm S'}(t)|}=\frac{d(F_S(t))}{d(F_{S'}(t))d(F_{S\sm S'}(t))}
= (1\pm 4\eps_*)\prod_{T: T\sub S,T \nsub S',T \nsub S\sm S'}d_T(F(t)).\]
\COMMENT{Formerly:
$\frac{|F_S(t)|}{|F_{S'}(t)||F_{S\sm S'}(t)|}=\frac{d(F_S(t))}{d(F_{S'}(t))d(F_{S\sm S'}(t))}>d_u^2/2$.
}
\item[(v)]
If $S' \sub S \in H$ then $|F_S(t)(P)|=(1\pm\eps_*)|F_S(t)|/|F_{S'}(t)|$
for all but at most $\eps_*|F_{S'}(t)|$ sets $P \in F_{S'}(t)$.
\item[(vi)] Statements (iii-v) hold replacing $F_{S^\le}(t)$ by $F_{S^\le}(t)[\GG]$
for any $\eps_{12D,3}$-regular subcomplex $\GG$ of $F_{S^\le}(t)$,
such that $d_T(\GG) \ge \eps_*^2$ when defined.%
\COMMENT{
1. Thought needed more than $\eps_*$ (so also in W-lemma) but ok.
Update: now think W-lemma needs eg $\eps_*^2$ here.
2. Formerly: If also $d_T(\GG) \ge d_u$ when defined
for $T \nsub S'$ and $T \nsub S\sm S'$ then (iv) holds.
}
\end{itemize}
\end{lemma}

\nib{Proof.} The first formula in (i) is a weaker form of $(*_{\ref{3alg}})$.
For the second formula, suppose first that $S=v$ has size $1$.
Then $|F_v(t)|=d_v(F(t))|V_v|=(1\pm \eps_*)d_{xv}(F(t-1))d_v(F(t-1))|V_v|
= (1\pm \eps_*)d_{xv}(F(t-1))|F_v(t-1)|= (1 \pm \eps_*)|F_{xv}(t-1)|/|F_x(t-1)|$.
In the case when $S=uv$ has size $2$ we have
$d_{S'}(F(t))=(1\pm \eps_{12D,1})d_{S'}(F(t-1))d_{S'x}(F(t-1))$
for $S'$ equal to $u$, $v$ or $uv$. Since $\eps_{12D,1} \ll \eps_*$,
the formula follows from the same calculations as in Lemma \ref{3average}(i).
This proves (i). For (ii), note that
if $|S|=1$ then $F_S(t)=F_S(t) \sm y$ has size $|F_S(t-1)|$ or $|F_S(t-1)|-1$.
Also, if $|S|=2,3$ we have $F_S(t) = F_{S^\le}(t-1)[F_{S^<}(t)]_S$ by Lemma \ref{build-update}.
Statement (ii) then follows from Lemma \ref{2restrict} if $|S|=2$ or
Lemma \ref{3restrict} if $|S|=3$.
For (iii) we apply Lemma \ref{3absolute} when $|S|=3$
or the identity $d(F_S(t))=\prod_{T\sub S}d_T(F(t))$ when $|S|\le 2$.
Statement (iv) follows by definition of absolute density and (iii).
For (v) we apply Lemma \ref{3average}.
For (vi) we apply regular restriction to see that $F_{S^\le}(t)[\GG]$
is $\eps_{12D,3}$-regular and then the same proofs. \qed

Our next lemma concerns the definitions for marked edges in the algorithm.

\begin{lemma}\label{3marked}$ $
\begin{itemize}
\item[(i)] For every triple $E \in H$ we have
$|M_{E^t,E}(t)| < \theta'_{\nu'_{E^t}(t)} |F_{E^t}(t)|$, and in fact

$|M_{E^t,E}(t)| \le \theta_{\nu'_{E^t}(t)} |F_{E^t}(t)|$ for $E \in U(x)$.
\item[(ii)] For every $x$ and triple $E \in U(x)$ we have
$|D_{x,E}(t-1)| < \theta_{\nu'_{E^t}(t)} |F_x(t-1)|$.
\end{itemize}
\end{lemma}

\nib{Proof.}
Throughout we use the notation $\ov{E}=E^{t-1}$, $\nu=\nu'_{\ov{E}}(t-1)$, $\nu^*=\nu'_{E^t}(t)$.%
\index{$Ebar$@$\ov{E}$}\index{$\nu$}\index{$\nu^*$}

(i) To verify the bound for $t=0$ we use our assumption that $(G,M)$ is
$(\eps,\eps',d_2,\theta,d_3)$-super-regular. We take $I=(\{\es\})$,
when for any $v$ we have $G^{I_v}=G$ by Definition \ref{def-3preplus}.
Then condition (iii) in Definition \ref{def-3super} gives $|M_E| \le \theta |G_E|$.
Since $E^0=E$, $M_{E,E}(0)=M_E$ and $F_E(0)=G_E$ we have the required bound.
Now suppose $t>0$. When $E \in U(x)$ we have $|M_{E^t,E}(t)| \le \theta_{\nu^*} |F_{E^t}(t)|$
by definition, since the algorithm chooses $y=\phi(x) \notin D_{x,E}(t-1)$.
Now suppose $E \notin U(x)$, and let $t'<t$ be the most recent time
at which we embedded a vertex $x'$ with $E \in U(x')$. Then $E^{t'}=E^t$, $\nu'_{E^t}(t')=\nu^*$,
and $|M_{E^{t},E}(t')| \le \theta_{\nu'_{E^t}(t')} |F_{E^{t}}(t')|$ by the previous case.
For any $z \in E^t$, we can bound $|F_z(t)|$ using the same argument
as that used at the end of the proof of Lemma \ref{3exceptional}.
We do not embed any neighbour of $z$ between time $t'+1$ and $t$,
so the size of the free set for $z$ can only decrease by
a factor of $\delta'_Q$ and an absolute term of $2\sqrt{\delta_Q} n$.
Since $d_z(F(t')) \ge d_u \gg \delta_Q$ we have
$|F_z(t)| \ge \delta'_Q |F_z(t')| - 2\sqrt{\delta_Q} n \ge \frac{1}{2}\delta'_Q |F_z(t')|$.
By Lemma \ref{build-update}, for every $\es \ne S \sub E^t$,
$F_S(t)$ is obtained from $F_S(t')$ by restricting to the $1$-complex $((F_z(t):z \in S),\{\es\})$.
If $|S|=2,3$ then regular restriction (Lemmas \ref{2restrict} and \ref{3restrict})
gives $d_S(F(t)) = (1 \pm \eps_*)d_S(F(t'))$.
Now $d(F_{E^t}(t)) = (1 \pm \eps_*) \prod_{S \sub E^t} d_S(F(t))$,
by Lemma \ref{3crude}(iii), and similarly $d(F_{E^t}(t')) = (1 \pm \eps_*) \prod_{S \sub E^t} d_S(F(t'))$.
This gives
$$\frac{|F_{E^t}(t)|}{|F_{E^t}(t')|} = (1 \pm 3\eps_*) \prod_{S \sub E^t} \frac{d_S(F(t))}{d_S(F(t'))}
= (1 \pm 10\eps_*) \prod_{z \in E^t} \frac{d_z(F(t))}{d_z(F(t'))} > (\delta'_Q/2)^3/2.$$
Therefore $|M_{E^t,E}(t)| \le |M_{E^t,E}(t')| \le \theta_{\nu^*} |F_{E^t}(t')|
< 2 (\delta'_Q/2)^{-3} \theta_{\nu^*} |F_{E^t}(t)| < \theta'_{\nu^*} |F_{E^t}(t)|$.

(ii) First we consider the case $x \in E$.
Then $E^t = \ov{E} \sm x$ and $\nu<\nu^*$ by Lemma \ref{3observe}(vi).
Also $F_{E^t}(t)=F_{\ov{E}}(t-1)(y)$
and $M_{E^t,E}(t) = M_{\ov{E},E}(t-1)(y)$ (see Lemma \ref{track-mark}),
so
$$D_{x,E}(t-1) = \{ y \in F_x(t-1): |M_{\ov{E},E}(t-1)(y)| > \theta_{\nu^*} |F_{\ov{E}}(t-1)(y)| \}.$$
If $\ov{E}=\{x\}$ has size $1$ then $D_{x,E}(t-1)=M_{x,E}(t-1)$ by Lemma \ref{itworks},
so $|D_{x,E}(t-1)| < \theta'_{\nu} |F_x(t-1)|< \theta_{\nu^*} |F_x(t-1)|$ by (i).
If $|\ov{E}| \ge 2$ then
$|F_{E^t}(t)|= |F_{\ov{E}}(t-1)(y)| =  (1 \pm \eps_*)|F_{\ov{E}}(t-1)|/|F_x(t-1)|$
by Lemma \ref{3crude} when $y \notin E_x(t-1)$. Now
\begin{align*}
& \sum_{y \in D_{x,E}(t-1)} |M_{\ov{E},E}(t-1)(y)|
 > \theta_{\nu^*} \sum_{y \in D_{x,E}(t-1) \sm E_x(t-1)} |F_{\ov{E}}(t-1)(y)| \\
& \qquad >  (1-\eps_*)\theta_{\nu^*} (|D_{x,E}(t-1)|-\eps_*|F_x(t-1)|)|F_{\ov{E}}(t-1)|/|F_x(t-1)|.
\end{align*}
We also have an upper bound
$$\sum_{y \in D_{x,E}(t-1)} |M_{\ov{E},E}(t-1)(y)|
\le \sum_{y \in F_x(t-1)} |M_{\ov{E},E}(t-1)(y)|
= |M_{\ov{E},E}(t-1)| < \theta'_{\nu} |F_{\ov{E}}(t-1)|,$$
where the last inequality holds by (i).
Therefore
\[ \frac{|D_{x,E}(t-1)|}{|F_x(t-1)|} < \frac{\theta'_{\nu}}{(1-\eps_*)\theta_{\nu^*}}
+ \eps_* < \theta_{\nu^*}.\]

Now we consider the case when $x \notin E$. Then $E^t=E^{t-1}=\ov{E}$.
Note that when $E \in U(x)$ we have $\ov{E} \cap VN_H(x) \ne \es$, so $\nu^*>\nu$.
Suppose first that $\ov{E} \in H(x)$. Then $F_{\ov{E}}(t)=F_{\ov{E}x}(t-1)(y)$
and $M_{\ov{E},E}(t) = M_{\ov{E},E}(t-1) \cap F_{\ov{E}}(t)$ (see Lemma \ref{track-mark}), so
$$D_{x,E}(t-1) = \{ y \in F_x(t-1): |M_{\ov{E},E}(t-1) \cap F_{\ov{E}x}(t-1)(y)|
> \theta_{\nu^*} |F_{\ov{E}x}(t-1)(y)| \}.$$
Similarly to the previous case, by Lemma \ref{3crude} we have
\begin{align*}
& \Sa := \sum_{y \in D_{x,E}(t-1)} |M_{\ov{E},E}(t-1) \cap F_{\ov{E}x}(t-1)(y)|
 > \theta_{\nu^*} \sum_{y \in D_{x,E}(t-1)\sm E_x(t-1)} |F_{\ov{E}x}(t-1)(y)| \\
& \qquad >  (1-\eps_*)\theta_{\nu^*} (|D_{x,E}(t-1)|-\eps_*|F_x(t-1)|)|F_{\ov{E}x}(t-1)|/|F_x(t-1)|.
\end{align*}
We also have $\Sa \le \sum_{y \in F_x(t-1)} |M_{\ov{E},E}(t-1) \cap F_{\ov{E}x}(t-1)(y)|$.
This sum counts all pairs $(y,P)$ with $P \in M_{\ov{E},E}(t-1)$, $y \in F_x(t-1)$
and $Py \in F_{\ov{E}x}(t-1)$, so we can rewrite it as
\[\Sa \le \sum_{P \in M_{\ov{E},E}(t-1)} |F_{\ov{E}x}(t-1)(P)|.\]
By Lemma \ref{3crude}(v) we have
$|F_{\ov{E}x}(t-1)(P)| = (1 \pm \eps_*)\frac{|F_{\ov{E}x}(t-1)|}{|F_{\ov{E}}(t-1)|}$
for all but at most $\eps_*|F_{\ov{E}}(t-1)|$ sets $P \in F_{\ov{E}}(t-1)$.
Therefore
$$\Sa \le |M_{\ov{E},E}(t-1)|(1+\eps_*)\frac{|F_{\ov{E}x}(t-1)|}{|F_{\ov{E}}(t-1)|}
+ \eps_*|F_{\ov{E}}(t-1)||F_x(t-1)|.$$
Combining this with the lower bound on $\Sa$ we obtain
$$(1-\eps_*)\theta_{\nu^*}\left(\frac{|D_{x,E}(t-1)|}{|F_x(t-1)|}-\eps_*\right)
< (1+\eps_*)\frac{|M_{\ov{E},E}(t-1)|}{|F_{\ov{E}}(t-1)|}
+ \eps_*\frac{|F_{\ov{E}}(t-1)||F_x(t-1)|}{|F_{\ov{E}x}(t-1)|}.$$
Now $|M_{\ov{E},E}(t-1)| < \theta'_{\nu} |F_{\ov{E}}(t-1)|$ by (i), and
$\frac{|F_{\ov{E}}(t-1)||F_x(t-1)|}{|F_{\ov{E}x}(t-1)|} \le 2d_u^{-1} \ll \eps_*^{-1}$
by Lemma \ref{3crude}(iv), so
\[\frac{|D_{x,E}(t-1)|}{|F_x(t-1)|} <
\frac{(1+\eps_*)\theta'_{\nu}+\sqrt{\eps_*}}{(1-\eps_*)\theta_{\nu^*}}
+ \eps_* < \theta_{\nu^*}.\]

It remains to consider the case when $x \notin E$ and $\ov{E}=E^t \notin H(x)$.
Since $\ov{E} \cap VN_H(x) \ne \es$ we have $|\ov{E}| \ge 2$.
\COMMENT{
This is why we changed the definition of $U(x)$ for $E$ to $E^{t-1}$.
[Attempts to make $E$ work: Put another prime on thetas?
No, whatever M-ratio is used in def D appears with neighbour of v
and then degrades with leaks, so will not be satisfied later in def D.
Also, didn't (i) already handle M-ratio? Also, replace $\nu_E(t)$ by $\nu'_{E^t}(t)$:
more natural, less error prone, should be monotone.
}
Then $F_{\ov{E}^\le}(t)=F_{\ov{E}^\le}(t-1)[F_{\ov{E}^<}(t)]$ by Lemma \ref{build-update}
and $M_{\ov{E},E}(t) = M_{\ov{E},E}(t-1) \cap F_{\ov{E}}(t)$, so
$$D_{x,E}(t-1) = \left\{ y \in F_x(t-1): \frac{|M_{\ov{E},E}(t-1) \cap
F_{\ov{E}^\le}(t-1)[F_{\ov{E}^<}(t)]| }{ |F_{\ov{E}^\le}(t-1)[F_{\ov{E}^<}(t)]| }
> \theta_{\nu^*} \right\}.$$
Let $I$ be the subcomplex of $\ov{E}^<$ consisting of all $S \sub \ov{E}$ with $S \in H(x)$.
Then $P \in F_{\ov{E}^\le}(t-1)[F_{\ov{E}^<}(t)]$ if and only if $P \in F_{\ov{E}^\le}(t-1)$
and $P_Sy \in F_{Sx}(t-1)$ for all $S \in I$.
When we choose $y \notin E_x(t-1)$, Lemma \ref{3crude} gives
$d_S(F(t))=(1\pm \eps_*)d_S(F(t-1))d_{Sx}(F(t-1))$ for $\es \ne S \in I$ by (i),
$d_S(F(t))=(1\pm \eps_*)d_S(F(t-1))$ for $S \sub \ov{E}$ with $S \notin I$ by (ii),
and 
\[d(F_{\ov{E}^\le}(t-1)[F_{\ov{E}^<}(t)]) = (1 \pm \eps_*)
\prod_{S \sub \ov{E}} d_S(F(t-1)) \prod_{\es \ne S \in I} d_{Sx}(F(t-1))\] 
by (vi), so 
\[|F_{\ov{E}^\le}(t-1)[F_{\ov{E}^<}(t)]| = (1 \pm 20\eps_*)|F_{\ov{E}^\le}(t-1)|
\prod_{\es \ne S \in I} d_{Sx}(F(t-1)).\]
As in the previous cases we have
\begin{eqnarray*}
\Sa & := & \sum_{y \in D_{x,E}(t-1)} |M_{\ov{E},E}(t-1) \cap F_{\ov{E}^\le}(t-1)[F_{\ov{E}^<}(t)]| \\
& > & \theta_{\nu^*} \sum_{y \in D_{x,E}(t-1)\sm E_x(t-1)} |F_{\ov{E}^\le}(t-1)[F_{\ov{E}^<}(t)]| \\
& > & (1-20\eps_*)\theta_{\nu^*} (|D_{x,E}(t-1)|-\eps_*|F_x(t-1)|) \
|F_{\ov{E}^\le}(t-1)| \prod_{\es \ne S \in I} d_{Sx}(F(t-1)) .
\end{eqnarray*}
For any $P \in F_{\ov{E}}(t-1)$, let $F_{P,I}$ be the set of $y \in F_x(t-1)$
such that $P_S y \in F_{S x}(t-1)$ for all $S \in I$.%
\COMMENT{Removed $\es\ne$... cf lemma.}
Let $B_I$ be the set of $P \in F_{\ov{E}}(t-1)$ such that we do not have
\[|F_{P,I}| = (1 \pm \eps_*)|F_x(t-1)| \prod_{\es \ne S \in I} d_{Sx}(F(t-1)).\]
Then  $|B_I| \le \eps_*|F_{\ov{E}}(t-1)|$ by Lemma \ref{3technical}.
Now $\Sa \le \sum_{y \in F_x(t-1)} |M_{\ov{E},E}(t-1) \cap F_{\ov{E}^\le}(t-1)[F_{\ov{E}^<}(t)]|$,
which counts all pairs $(y,P)$ with $P \in M_{\ov{E},E}(t-1)$ and $y \in F_{P,I}$, so
$$\Sa \le |M_{\ov{E},E}(t-1)|(1 \pm \eps_*)|F_x(t-1)| \prod_{\es \ne S \in I} d_{Sx}(F(t-1))
+ \eps_*|F_{\ov{E}}(t-1)||F_x(t-1)|.$$
Combining this with the lower bound on $\Sa$ we obtain
$$(1-20\eps_*)\theta_{\nu^*}\left(\frac{|D_{x,E}(t-1)|}{|F_x(t-1)|}-\eps_*\right)
< (1+\eps_*)\frac{|M_{\ov{E},E}(t-1)|}{|F_{\ov{E}}(t-1)|}
+ \eps_* \prod_{\es \ne S \in I} d_{Sx}(F(t-1))^{-1}.$$
Now $|M_{\ov{E},E}(t-1)| < \theta'_{\nu} |F_{\ov{E}}(t-1)|$ by (i)
and all densities are at least $d_u \gg \eps_*$,
so again we have
\begin{equation}
\frac{|D_{x,E}(t-1)|}{|F_x(t-1)|} < \frac{(1+\eps_*)\theta'_{\nu}+\sqrt{\eps_*}}{(1-20\eps_*)\theta_{\nu^*}}
+ \eps_* < \theta_{\nu^*}. \tag*{$\Box$}
\end{equation}

The following corollary is now immediate from Lemmas \ref{3exceptional} and \ref{3marked}.
Recall that $OK_x(t-1)$ is obtained from $F_x(t-1)$ by deleting $E_x(t-1)$
and $D_{x,E}(t-1)$ for $E \in U(x)$, and note that since $H$ has maximum degree
at most $D$ we have $|U(x)| \le 2D^2$.

\begin{coro}\label{3ok} $|OK_x(t-1)| > (1-\theta_*)|F_x(t-1)|$. \index{$OK$}
\end{coro}

\subsection{The initial phase}\index{initial phase}

This subsection concerns the initial phase of the algorithm,
during which we embed the neighbourhood $N$ of the buffer $B$.
There are several issues that make the analysis of this phase
significantly more complicated than that of the graph blow-up lemma
(which was given in Lemma \ref{2-x-to-v}).
As we mentioned earlier, the buffer $B$ is larger than before,
so we embed many more vertices during the initial phase, and the queue may open.
One potential problem is that there may be some vertex $v$ and class $V_i$
so that $G(v) \cap V_i$ is used excessively by the embedding
-- this is a concern, since $d_u \ll \delta_B$.
The first lemma in this subsection shows that with high probability
this does not happen.%

Our goal is to show that for any vertex $x \in B$\index{$x$} there will be many 
{\em available}\index{available} vertices $v \in V_x$\index{$v$} such that we embed 
$H(x)$ in $(G\sm M)(v)$ during the initial phase. Then if $v$ is not used before the 
conclusion of the algorithm we will be able to embed $x$ as $\phi(x)=v$.
We need to ensure that for every neighbour $z$ of $x$
and for every triple $E$ of $H$ containing a neighbour of $z$,
the choice of image for $z$ is good, in that the subcomplex of
the free sets for $E^\le$ that is consistent with mapping $x$ to $v$
is suitably regular and does not have too many marked edges.
As in Lemma \ref{2-x-to-v}, we aim to give a lower bound on the probability
of this event, conditional on the previous embedding.
The third lemma in this subsection achieves this.

The marked edges also add a complication to the conclusion of the algorithm,
in which we need to verify Hall's condition for a system of distinct
representatives of the available images for the unembedded buffer vertices.
We need to show that for that any $W \sub V_x$ that is not too small, the
probability that $W$ does not contain a vertex available for $x$ is quite small.
This is achieved by the second lemma in this subsection. We present it before
the lemma on mapping $x$ to $v$ because its proof is similar in spirit
and somewhat simpler.

First we recall the key properties of the selection rule during the initial phase.
Although the queue may become non-empty, jumping ensures that
we embed all vertex neighbourhoods $VN_H(x)$, $x \in B$ at consecutive times,
and before $x$ or any other vertices at distance at most $4$ from $x$.
Then Lemma \ref{not-local} implies that if we start embedding $VN_H(x)$ just after some time $T_0$
then $F_z(T_0) = V_z(T_0)$ consists of all vertices in $V_z$ that have not
yet been used by the embedding, for every $z$ at distance at most $3$ from $x$.
We also recall that $|B \cap V_z| = \delta_B |V_z|$, $|N \cap V_z| < \sqrt{\delta_B}|V_z|$,
$|Q(T_0) \cap V_z| \le \delta_Q|V_z|$ and 
$|J(T_0)  \cap V_z| \le \sqrt{\delta_Q}|V_z|$ by Lemma \ref{3observe}(ii).
Thus for any $z$ at distance at most $3$ from $x$ we have
\begin{equation}\label{eq:3fresh}
|F_z(T_0)|=|V_z(T_0)|>(1-2\sqrt{\delta_B})|V_z|.
\end{equation}
We need the following supermartingale formulation of
the Azuma-Hoeffding inequality, which can be easily derived from the martingale
formulation quoted later as Theorem \ref{azuma}.

\begin{theo}\label{ahsuper} \index{supermartingale}
Suppose $Z_0, \cdots, Z_n$ is a supermartingale, i.e.\ a sequence of random variables
satisfying $\mb{E}(Z_{i+1}|Z_0,\cdots,Z_i) \le Z_i$, and that $|Z_i-Z_{i-1}| \le c_i$,
$1 \le i \le n$, for some constants $c_i$. Then for any $t \ge 0$,
\[\mb{P}(Z_n-Z_0 \ge t) \le 2\exp \left( - \frac{t^2}{2\sum_{i=1}^n c_i^2} \right).\]
\end{theo}

Let $T_I$ be the time at which we embed the last vertex of $N$, ending the initial phase.
Then $T_I \le \sum_{i=1}^r 2\sqrt{\delta_B}|V_i| < \delta_B^{1/3}n$ by (\ref{eq:3fresh}).
Our first lemma shows that there are many available vertices in all neighbourhoods
at time $T_I$. Note that by super-regularity the assumption $|G(v)_j| \ge d_u|V_j|$ in
the lemma holds for every $j$ such that $G_{i(v)j}$ is defined.\index{$T_I$}
\COMMENT{Former comments:
1. Enough to do at $T_I$, as this is present at earlier times,
possibly rather sparsely but $d_u$ is ok for qr restriction.
2. Why not replace $G(v)_{j^*}$ with fixed $A \sub V_{j^*}$?
(Give error prob and take union bound outside proof.)
We need bounds on $G(v)_{k^*} \cap V_{k^*}(t)$ to control $j^*$?
Actually we just need few $j^*$-atypical $k^*$-vertices, but
we'll keep the nhood formulation.
3. Use single martingale bound, not induction.
4. Optional stopping?
5. Class ratios for nhood? Avoid by tracking those that contribute,
bounded by VN-ratio.
6. Remember $N$ ratio is $\delta_B^{1/2}$.
7. Really we only need $G(v)$ formulation, as we can handle further nhoods as they
are embedded, but this is convenient.
}

\begin{lemma} \label{3uniform}
With high probability, for every vertex $v\in G$ and $1 \le j \le r$
with $|G(v)_j| \ge d_u|V_j|$ we have%
\COMMENT{
1. Change to $ij\in G$, $v\in G_i$, density automatic? Keep it:
useful to state density, and anyway R-indexed cannot write ij.
2. Not used in 3initial, used in 3-x-to-v, but keep here.
}
\[|G(v)_j \cap V_j(T_I)| > (1-\delta_B^{1/3})|G(v)_j|.\]
\end{lemma}

\nib{Proof.}
Suppose $|G(v)_j| \ge d_u|V_j|$. The ratio $|G(v)_j \cap V_j(t)|/|G(v)_j|$
only decreases when we embed a vertex in $G(v)_j$. We separate the analysis 
according to two effects. One effect is that we embed a vertex
$z \in N \cap X_j$ to an image $\phi(z) \in G(v)_j$,
where $G(v)_j$ is not too large a fraction of the free images for $z$. 
The other effect is that we we embed a vertex
$z \in N \cap X_j$ to an image $\phi(z) \in G(v)_j$,
and the embedding of some neighbour $w$ of $z$ previously caused
the fraction of $G(v)_j$ in the free images for $z$ to increase significantly.
To analyse these effects we write $T_z$ for the time at which a vertex $z$ is embedded
and define the following sets:%
\COMMENT{
1. New notation: $T_z$ time $z$ embedded, $t_z^N$ time last z-nbr embedded (only used for x).
2. Formerly: The third effect (easily analysed) is embedding vertices from the queue in $G(v)_j$.
... We will only consider this effect for those $z$ where $G(v)_j$ is not too large a fraction
of the free images for $z$. Otherwise we will attribute this embedding of $z$ to a second effect...
3. Jumps in N so covered.
}
\begin{itemize}
\item \index{$\Lambda_j(t)$}
Let $\Lambda_j(t)$ be the set of embedded vertices $z \in N \cap X_j$
such that $\phi(z) \in G(v)_j$ and
\[|G(v)_j \cap F_z(T_z-1)|/|F_z(T_z-1)| < 2^{2D}|G(v)_j|/|V_j|. \]
\item \index{$\Pi_{\ell,j}(t)$}
Let $\Pi_{\ell,j}(t)$ be the set of embedded vertices $w \in N \cap X_\ell$
such that
$$|G(v)_j \cap F_z(T_w)|/|F_z(T_w)| \ge
2 |G(v)_j \cap F_z(T_w-1)|/|F_z(T_w-1)| \ge d_u^2$$
for some $z \in VN_H(w) \cap N \cap X_j$.
\end{itemize}

We claim that any vertex embedded in $G(v)_j$ up to time $T_I$ is either
in the queue, or in $\Lambda_j(T_I)$,
or an $H$-neighbour of some $w \in \Pi_{\ell,j}(T_I)$ for some $\ell$.
To see this, suppose $z$ is embedded in $G(v)_j$ and
is neither in the queue nor in $\Lambda_j(T_I)$.
Since $z \in N$ we have $z \in VN_H(x)$ for some $x \in B$.
Suppose we start embedding $VN_H(x)$ just after time $T_0$.
Then $F_z(T_0)=V_j(T_0)$ has size at least $(1-2\sqrt{\delta_B})|V_j|$ by (\ref{eq:3fresh}),
so $|G(v)_j \cap F_z(T_0)|/|F_z(T_0)| \le |G(v)_j|/(1-2\sqrt{\delta_B})|V_j|$.
We also have $|G(v)_j \cap F_z(T_z-1)|/|F_z(T_z-1)| \ge 2^{2D}|G(v)_j|/|V_j|$,
since $z \notin \Lambda_j(T_I)$. Since $|VN_H(x)\sm z| \le 2D-1$,
there must be some $g$, $1 \le g \le T_z-T_0$
so that at time $T_0+g$ we embed some $w \in VN_H(z)$ and get
$|G(v)_j \cap F_z(T_w)|/|F_z(T_w)| \ge
2 |G(v)_j \cap F_z(T_w-1)|/|F_z(T_w-1)|$.
Here we recall that we embed $VN_H(x)$ consecutively
and so $w \in VN_H(x) \sub N$.
Now for any $g'$, $1 \le g' < T_z-T_0$,
since $|F_z(T_z-1)| \ge d_u|V_j| \ge d_u|F_z(T_0+g')|$ we have
\[\frac{|G(v)_j \cap F_z(T_0+g')|}{|F_z(T_0+g')|}
\ge \frac{|G(v)_j \cap F_z(T_z-1)|}{d_u^{-1}|F_z(T_z-1)|}
\ge 2^{2D}d_u|G(v)_j|/|V_j| > 2^{2D}d_u^2 > d_u^2.\]
Therefore $w \in \Pi_{\ell,j}(T_I)$, where $\ell=i(w)$, which proves the claim.
Since any $w$ has at most $2D$ neighbours we deduce that
\begin{equation}\label{eq:3uniform}
|G(v)_j \cap V_j(T_I)| \ge |G(v)_j| - |Q(T_I) \cap V_j|
- |\Lambda_j(T_I)| - 2D \sum_{\ell=1}^r |\Pi_{\ell,j}(T_I)|.
\end{equation}

Consider $Z_j(t) = |\Lambda_j(t)| - 2^{2D+1}d_j(G(v))|V_j \sm V_j(t)|$.
We claim that $Z_j(0), \cdots, Z_j(T_I)$ is a supermartingale.
To see this, suppose we embed some vertex $z \in X_j$ at time $T_z$.
We can assume that $z\in N$ and $|G(v)_j \cap F_z(T_z-1)|/|F_z(T_z-1)| < 2^{2D}d_j(G(v))$,
or otherwise $z \notin \Lambda_j(T_z)$ by definition,
so $|\Lambda_j(T_z)|=|\Lambda_j(T_z-1)|$ and $Z_j(T_z) < Z_j(T_z-1)$.
Since $\phi(z)$ is chosen randomly in $OK_z(T_z-1) \sub F_z(T_z-1)$
of size at least $|F_z(T_z-1)|/2$ (by Corollary \ref{3ok}), we have
\begin{align*}
& \mb{E}[|\Lambda_j(T_z)|-|\Lambda_j(T_z-1)|] = \mb{P}(\phi(z) \in G(v)_j)
= \frac{|G(v)_j \cap OK_z(T_z-1)|}{|OK_z(T_z-1)|} \\
& < \frac{2|G(v)_j \cap F_z(T_z-1)|}{|F_z(T_z-1)|}
< 2^{2D+1}d_j(G(v)).
\end{align*}
We also have $|V_j \sm V_j(T_z)|=|V_j \sm V_j(T_z-1)|+1$,
so the decrease in the second term of $Z_j(t)$ more than compensates
for the increase in the first, i.e.\ $\mb{E}[Z_j(T_z)-Z_j(T_z-1)] < 0$.
Thus we have a supermartingale. Since $|Z_j(t)-Z_j(t-1)|\le 1$,
$|V_j \sm V_j(T_I)|\le 2\sqrt{\delta_B}|V_j|$ by (\ref{eq:3fresh}),
and $T_I < \delta_B^{1/3}n$, Theorem \ref{ahsuper} gives
\begin{align*}
& \mb{P}\big[ |\Lambda_j(T_I)|>2^{2D+3}d_j(G(v))\sqrt{\delta_B}|V_j| \big]
 \le \mb{P}\big[ Z_j(m)>2^{2D+2}d_j(G(v))\sqrt{\delta_B}|V_j| \big] \\
& < 2\exp \big[-(2^{2D}d_j(G(v))\sqrt{\delta_B}|V_j|)^2/2T_I\big]
< e^{-\sqrt{n}}, \quad \mbox{ }(\mbox{say, for sufficiently large }n).
\end{align*}

Next consider $Y_{\ell,j}(t) = |\Pi_{\ell,j}(t)| - \eps_* |\Pi'_{\ell,j}(t)|$,
where $\Pi'_{\ell,j}(t)$ consists of all vertices in $X_{\ell}$ with at least
one $H$-neighbour in $X_{j}$ that have been embedded at time $t$.
We claim that $Y_{\ell,j}(0), \cdots, Y_{\ell,j}(T_I)$ is a supermartingale.
To see this, suppose we embed some vertex $w \in X_{\ell}$
at time $T_w$. Consider $z \in VN_H(w) \cap N \cap X_{j}$:
we can assume this set is non-empty, otherwise
$w \notin \Pi_{\ell,j}(T_w) \cup \Pi'_{\ell,j}(t)$, so $Y_{\ell,j}(T_w) =Y_{\ell,j}(T_w-1)$.
We can also assume that $w\in N$ and
$|G(v)_j \cap F_z(T_w)|/|F_z(T_w)| \ge
2 |G(v)_j \cap F_z(T_w-1)|/|F_z(T_w-1)| \ge d_u^2$
for some $z \in VN_H(w) \cap N \cap X_j$,
otherwise $w \notin \Pi_{\ell,j}(T_w)$ by definition,
so $|\Pi_{\ell,j}(T_w)|=|\Pi_{\ell,j}(T_w-1)|$ and $Y_{\ell,j}(T_w) < Y_{\ell,j}(T_w-1)$.

By Lemmas \ref{3exceptional} and Lemma \ref{2restrict},
$F_{wz}(T_w-1)[G(v)_j \cap F_z(T_w-1)]$ is $\eps_{12D,2}$-regular.
Applying Lemma \ref{2neighbour}, we see that there are at most
$2\eps_{12D,2}|F_w(T_w-1)|$ `exceptional' vertices $y \in F_w(T_w-1)$
such that embedding $\phi(w)=y$ will not satisfy
\begin{align*}
|G(v)_j \cap F_z(T_w)| & = |F_{zw}(T_w-1)(y) \cap G(v)_j \cap F_z(T_w-1)| \\
& = (1 \pm \eps_*)d_{zw}(F(T_w-1))|G(v)_j \cap F_z(T_w-1)|.
\end{align*}
On the other hand, the algorithm chooses $\phi(w)=y$ to satisfy $(*_{\ref{3alg}})$,
so $|F_z(T_w)| = (1 \pm \eps_*)d_{zw}(F(T_w-1))|F_z(T_w-1)|$. Thus we have
$|G(v)_j \cap F_z(T_w)|/|F_z(T_w)| < 2|G(v)_z \cap F_z(T_w-1)|/|F_z(T_w-1)|$,
unless we choose an exceptional vertex $y$.  But $y$ is chosen uniformly at
random from $|OK_w(T_w-1)| \ge (1-\theta_*)|F_w(T_w-1)|$ possibilities
(by Corollary \ref{3ok}), so $y$ is exceptional with probability at most $3\eps_{12D,2}$.
Therefore $\mb{E}\big[|\Pi_{\ell,j}(T_w)|-|\Pi_{\ell,j}(T_w-1)|\big] < 3\eps_{12D,2}$.

We also have $|\Pi'_{\ell,j}(T_w)|=|\Pi'_{\ell,j}(T_w-1)|+1$,
so the decrease in the second term of $Y_{\ell,j}(t)$ more than compensates
for the increase in the first, i.e.\ $\mb{E}[Y_{\ell,j}(T_w)-Y_{\ell,j}(T_w-1)] < 0$.
Thus we have a supermartingale. We also have $|Y_{\ell,j}(t)-Y_{\ell,j}(t-1)|\le 1$.
Also, $|V_j \sm V_j(T_I)|\le 2\sqrt{\delta_B}|V_j|$ by (\ref{eq:3fresh}),
so $|\Pi'_{\ell,j}(T_i)| \le 2D\sqrt{\delta_B}|V_{j}|$, since $H$ has maximum degree $D$.
Then Theorem \ref{ahsuper} gives
\begin{align*}
& \mb{P}\big[|\Pi_{\ell,j}(T_I)|>2d_u\sqrt{\delta_B}|V_{j}|\big]
< \mb{P}\big[Y_{\ell,j}(t)>d_u\sqrt{\delta_B}|V_{j}|\big] \\
& < 2\exp\big[-(d_u\sqrt{\delta_B}|V_{\ell}|)^2/2T_I\big]
< e^{-\sqrt{n}}, \quad \mbox{ }(\mbox{say, for sufficiently large }n).
\end{align*}

Taking a union bound over $v \in V$ and $1 \le j,\ell \le r$,
with high probability we have
$|\Lambda_j(T_I)| \le 2^{2D+3}d_{j}(G(v))\sqrt{\delta_B}|V_{j}|
=  2^{2D+3}\sqrt{\delta_B}|G(v)_{j}|$
and
$|\Pi_{\ell,j}(T_I)| \le 2d_u\sqrt{\delta_B}|V_{j}| \le 2\sqrt{\delta_B}|G(v)_{j}|$.
Therefore (\ref{eq:3uniform}) gives
\begin{align*}
& |G(v)_{j} \cap V_{j}(T_I)| \ge |G(v)_{j}| - |Q(T_I) \cap V_{j}|
- |\Lambda_j(T_I)| - 2D \sum_{\ell=1}^r |\Pi_{\ell,j}(T_I)|.\\
& \ge |G(v)_{j}| - \delta_Q|V_{j}|
-  2^{2D+3}\sqrt{\delta_B}|G(v)_{j}| - 4Dr\sqrt{\delta_B}|G(v)_{j}|
> (1-\delta_B^{1/3})|G(v)_{j}|,
\end{align*}
since $|G(v)_{j}| \ge d_u|V_{j}|$ and
$\eps_* \ll \delta_Q \ll d_u \ll \delta_B \ll 1/r,1/D$. \qed
%This completes the proof. \qed

For the remainder of this subsection we fix a vertex $x \in B$
and write $VN_H(x) = \{z_1,\cdots,z_g\}$,\index{$g$}
with vertices listed in the order that they are embedded.
Since $H$ has maximum degree $D$ we have $g \le 2D$.
We let $T_j$ be the time at which $z_j$ is embedded.\index{$T_j$}\index{$z_j$}
By the selection rule, $VN_H(x)$ jumps the queue and is embedded
at consecutive times: $T_{j+1}=T_j+1$ for $1 \le j \le g-1$.
For convenience we also define $T_0=T_1-1$.\index{$T_0$}
Note that since $H$ is an $r$-partite complex, no vertex of $VN_H(x)$ lies in $X_x$.
The selection rule also ensures that at time $T_0$
no vertices at distance at most $4$ from $x$ have been embedded,
so for any $z$ with distance at most $3$ from $x$ we have $F_z(T_0)=V_z(T_0)$.
(The need for this property in Lemma \ref{3-x-to-v} explains why we needed to choose
the buffer vertices at mutual distance at least $9$.)

\begin{rem}
Note that the argument of the Lemma \ref{3uniform} can be applied replacing
the sets $G(v)_j$ by any sufficiently large subsets of $G_j$,
provided that they are sufficiently few in number to use a union bound.
For example, we can define a complex $G\gens{S}=\cap_{v \in S} G(v)$
for any $S \sub V$, and show that with high probability,
for every $S \sub V$ with $|S| \le 2D$ and $1 \le j \le r$
with $|G\gens{S}_j| \ge d_u|V_j|$ we have
$|G\gens{S}_j \cap V_j(T_I)| > (1-\delta_B^{1/3})|G\gens{S}_j|$.
It follows that with high probability every $z \in VN_H(x)$ would have many available
vertices in $G(v)_z$ throughout the initial phase, even if we did not use
queue jumping in the selection rule. However, the analysis is simpler
if we do use queue jumping, and then we only need the version
of Lemma \ref{3uniform} stated above.
\COMMENT{
1. Queue jumping causes problems for arbitrary part sizes, so this might fix them.
2. Prefer consecutive so can condition on past without interlacing.
}
\end{rem}

Recall that the available\index{available} set\index{$A_x$} $A_x$
is obtained from $F_x(t_x^N)$ by removing all
sets $M_{x,E}(t_x^N)$ for triples $E$ containing $x$. Here $t_x^N=T_g$.
If $x$ is unembedded at the conclusion of the algorithm at time $T$
we will embed $x$ in $A'_x = A_x \cap V_x(T)$.
Our second lemma shows that for that any $W \sub V_x$ that is not too small, the
probability that $W$ does not contain a vertex available for $x$ is quite small.

\begin{lemma}\label{3initial}\index{$W$}
For any $W \sub V_x$ with $|W|> \eps_*|V_x|$,
conditional on any embedding of the vertices $\{s(u): u < T_1\}$
that does not use any vertex of $W$,
we have $\mb{P}[A_x \cap W = \es] < \theta_*$.
\COMMENT{
1. Thought I needed $|W|/|V_x| \gg \eps_*$ (eg $\gamma^2$ ok) but ok as it cancels.
See also 3crude comment.
2. Okay to assume $W$ unused, as
we're trying to show that $W$ does not survive.
If $W \sub V_x(T_g)$, then for any $S \sub B_x$,
each $x \in S$ had $W \sub V_x(T_1)$ so can apply this lemma.
}
\end{lemma}

\nib{Proof.}
We apply similar arguments to those we are using for the entire embedding,
defining variants of various structures that incorporate restriction to $W$.
Suppose $1 \le j \le g$ and that we are considering the embedding of $z_j$.
We interpret quantities at time $T_j$ with the embedding $\phi(z_j)=y$,\index{$y$}
for some as yet unspecified $y \in F_{z_j}(T_j-1)$. Write $W_j = W \cap F_x(T_j)$.
At time $T_0$ we have $W \sub F_x(T_0) = V_x(T_0)$, so $W_0 = F_x(T_0) \cap W = W$.
Since $z_j \in H(x)$ we have $F_x(T_j)=F_{xz_j}(T_j-1)(y)$, so 
\[W_j = W_{j-1} \cap F_{xz_j}(T_j-1)(y).\] 
For convenient notation we will use $[W_j]$
to denote restriction of hypergraphs and complexes to $W_j$, in that for $x \in S \in H$
we write $F_{S^\le}(T_j)[W_j]$ for $F_{S^\le}(T_j)[(W_j,\{\es\})]$ and
$F_S(T_j)[W_j]$ for $F_{S^\le}(T_j)[(W_j,\{\es\})]_S$. \index{$[W_j]$}
\COMMENT{
1. changed $U(x)$ to triples containing $x$
2. [cut] Also, $E^{T_0}=E$ for any triple $E$ containing $x$, so $M_{x,E}(T_0)$ is undefined.
3. changed second eps index from 2 to 1, cf also 3 to 2 in Claim C.
4. why not give $|W_j|$ using $|W_{j-1}|$? More convenient, but $j=0$ awkward
}

We define {\em exceptional} sets $E^W_{z_j}(T_j-1) \sub F_{z_j}(T_j-1)$
by the property that $y$ is in $F_{z_j}(T_j-1) \sm E^W_{z_j}(T_j-1)$
if and only if\index{exceptional}\index{$E^W_{z_j}(T_j-1)$}

$\qquad |W_j| = (1 \pm \eps_{\nu'_x(T_j),1})|W||F_x(T_j)|/|F_x(T_0)|$. \hfill $(*_{\ref{3initial}})$

Thus if we embed $z_j$ to $y \notin E^W_{z_j}(T_j-1)$, then the intersection of the free
set for $x$ with $W$ is roughly what would be `expected'. Next, to control marked edges,
we define {\em dangerous} vertices similarly to before, except that here we incorporate
the restriction to $W_j$. For any triple $E$ containing $x$ we define
\index{dangerous}\index{$D^W_{z_j,E}(T_j-1)$}
\COMMENT{Always have $x \in E^{T_j}$ and $E \in U(z_j)$ so no slow leak.}
$$D^W_{z_j,E}(T_j-1) = \{y \in F_{z_j}(T_j-1):
|M_{E^{T_j},E}(T_j)[W_j]| > \theta_{\nu'_{E^{T_j}}(T_j)} |F_{E^{T_j}}(T_j)[W_j]|  \}.$$

Now we define events as follows: \index{$A_{i,j}$}
\begin{itemize}
\item
$A_{1,j}$, $j \ge 0$ is the event that property $(*_{\ref{3initial}})$ above holds.
\item
$A_{2,j}$, $j \ge 0$ is the event that for every triple $E$ containing $x$ we have
\COMMENT{
Use prime so that 3marked gives base case.
[$z_j$ is in an edge containing $x$, so affects $E$ even if $z_j \notin E$...
we do have slow leaks in x-to-v, as there we track edges in $U(x)$ and
not of all of these are affected by $z_j$.]
}
$$|M_{E^{T_j},E}(T_j)[W_j]| \le \theta_{\nu'_{E^{T_j}}(T_j)} |F_{E^{T_j}}(T_j)[W_j]|.$$
\item
$A_{3,j}$, $j \ge 1$ is the event that $y=\phi(z_j)$ is chosen in $OK^W_{z_j}(T_j-1)$,
defined to be the subset of $F_{z_j}(T_j-1)$ obtained by deleting the sets
$E^W_{z_j}(T_j-1)$ and $D^W_{z_j,E}(T_j-1)$ for $E$ containing $x$.
For convenient notation we also define $A_{3,0}$ to be the event that
holds with probability $1$.
\index{$OK$}
\end{itemize}

We divide the remainder of the proof into a series of claims.

\nib{Claim A.} The events $A_{1,0}$, $A_{2,0}$ and $A_{3,0}$ hold.

\nib{Proof.} As noted above, we have $W_0=W$, so $A_{1,0}$ holds.
Also, $A_{3,0}$ holds by definition, so it remains to show $A_{2,0}$.
Consider any triple $E=xzz'$ containing $x$. Recall that all vertices at distance at most $4$
from $x$ are unembedded. Then $E^{T_0}=E$, and by Lemma \ref{not-local},
$F_{E^\le}(T_0)$ is the restriction of $G_{E^\le}$ to the $1$-complex $((F_u(T_0):u \in E),\{\es\})$,
so $F_{E^\le}(T_0)[W]$ is the restriction of $G_{E^\le}$ to $(W,F_z(T_0),F_{z'}(T_0),\{\es\})$.
Since $G$ is $\eps$-regular, $F_{E^\le}(T_0)[W]$ is $\eps'$-regular by regular restriction
and $d_S(F_{E^\le}(T_0)[W]) = (1\pm\eps')d_S(G)$ for $S \sub E$, $|S| \ge 2$.
Then
\[|F_E(T_0)[W]|=(1\pm 20\eps')|W||F_z(T_0)||F_{z'}(T_0)|\prod_{S \sub E,|S| \ge 2}d_S(G)\]
by Lemma \ref{3absolute}. Now
\[|M_{E,E}(T_0)[W]| = \sum_{v \in W} |M_{E,E}(T_0)(v)|
\le \sum_{v \in W} |M_E(v)| \le \theta \sum_{v \in W} |G_E(v)|\]
by condition (ii) of super-regularity (Definition \ref{def-3super}).
This condition also says for any $v \in G_x$ that $G_E(v)$ is $\eps'$-regular
and $d_S(G(v))=(1\pm\eps')d_S(G)d_{Sx}(G)$ for $\es \ne S \sub E \sm x$, so
\[|G_E(v)| = d_{zz'}(G(v))|G(v)_z||G(v)_{z'}|
= (1\pm 5\eps')|V_z||V_{z'}|\prod_{\es \ne S \sub zz'}d_S(G)d_{Sx}(G).\]
Here we note that $d_z(G)=d_{z'}(G)=1$, as our hypotheses for Theorem \ref{3blowup}
include the assumption $G_i=V_i$ for $1 \le i \le r$.
Equation (\ref{eq:3fresh}) gives $|F_z(T_0)|>(1-2\sqrt{\delta_B})|V_z(T_0)|$
and $|F_{z'}(T_0)|>(1-2\sqrt{\delta_B})|V_{z'}(T_0)|$, so
\COMMENT{Even $|F_z(T_0)|=|V_z(T_0)|>\delta_Bn/2$ suffices.}
\[\frac{|M_{E,E}(T_0)[W]|}{|F_E(T_0)[W]|} < \frac{(1+ 5\eps')\theta|V_z||V_{z'}|}{(1-20\eps')|F_z(T_0)||F_{z'}(T_0)|}
< (1+30\eps')(1-2\sqrt{\delta_B})^{-2}\theta < \theta_0.\]

\nib{Claim B.} If $A_{3,j}$ holds then $A_{1,j}$ and $A_{2,j}$ hold.

\nib{Proof.} This follows directly from the definitions: if $y \notin E^W_{z_j}(T_j-1)$
then $(*_{\ref{3initial}})$ holds, and if $y \notin D^W_{z_j,E}(T_j-1)$
then $|M_{E^{T_j},E}(T_j)[W_j]| \le \theta_{\nu'_{E^{T_j}}(T_j)} |F_{E^{T_j}}(T_j)[W_j]|$.

\nib{Claim C.} If $A_{1,j-1}$ holds then
$|E^W_{z_j}(T_j-1) \sm E_{z_j}(T_j-1)| < \eps_* |F_{z_j}(T_j-1)|$.

\nib{Proof.}
By Lemma \ref{3exceptional}, $F_{xz_j}(T_j-1)$ is $\eps_{\nu'_{xz_j}(T_j-1),1}$-regular,
so by Lemma \ref{2restrict}, $F_{xz_j}(T_j-1)[W_{j-1}]$ is $\eps_{\nu'_{xz_j}(T_j-1),2}$-regular.
Then by Lemma \ref{2neighbour}, for all but at most $\eps_{\nu'_{xz_j}(T_j-1),2}|F_{z_j}(T_j-1)|$
vertices $y \in F_{z_j}(T_j-1)$ we have $|W_j| = |W_{j-1} \cap F_{xz_j}(T_j-1)(y)|
= (1 \pm \eps_{\nu'_{xz_j}(T_j-1),2})d_{xz_j}(F(T_j-1))|W_{j-1}|$.
Since $A_{1,j-1}$ holds, $|W_{j-1}| = (1 \pm \eps_{\nu'_x(T_j-1),1})|W||F_x(T_j-1)|/|F_x(T_0)|$.
Also, by $(*_{\ref{3alg}})$, for $y \notin E_{z_j}(T_j-1)$ we have
$|F_x(T_j)|=(1 \pm \eps_{\nu'_{x}(T_j-1),0})d_{xz_j}(F(T_j-1))|F_x(T_j-1)|$.
Now $z_j \in VN_H(x)$, so $\nu_x(T_j)=\nu_x(T_j-1)+1$,
and $\nu'_x(T_j) > \max\{\nu'_x(T_j-1),\nu'_{xz}(T_j-1)\}$.
Combining all estimates, for all but at most $\eps_*|F_{z_j}(T_j-1)|$
vertices $y \in F_{z_j}(T_j-1) \sm  E_{z_j}(T_j-1)$ we have
\begin{align*}
|W_j| & = (1 \pm \eps_{\nu'_{xz_j}(T_j-1),2})d_{xz_j}(F(T_j-1))|W_{j-1}| \\
& = (1 \pm \eps_{\nu'_{xz_j}(T_j-1),2})(1 \pm \eps_{\nu'_x(T_j-1),1})d_{xz_j}(F(T_j-1))|W||F_x(T_j-1)|/|F_x(T_0)| \\
& = (1 \pm \eps_{\nu'_x(T_j),1})|W||F_x(T_j)|/|F_x(T_0)|, \qquad \mbox{i.e. } (*_{\ref{3initial}}).
\end{align*}

\nib{Claim D.} If $A_{1,j-1}$ and $A_{2,j-1}$ hold then for any $E$ containing $x$
we have $$|D^W_{z_j,E}(T_j-1)| < \theta_{\nu'_{E^{T_j}}(T_j)} |F_{z_j}(T_j-1)|.$$

\nib{Proof.} Denote $\ov{E}=E^{T_j-1}$,\index{$\ov{E}$}\index{$\nu$}\index{$B_{z_j}$}
$\nu=\nu'_{\ov{E}}(T_j-1)$ and $\nu^*=\nu'_{E^{T_j}}(T_j)$.%
\COMMENT{
Formerly $B_{z_j} = E_{z_j}(T_j-1) \cup E^W_{z_j}(T_j-1)$...
Also, $|B_{z_j}|<2\eps_* |F_{z_j}(T_j-1)|$
by Lemma \ref{3exceptional} and Claim C.
}
Suppose $A_{1,j-1}$ and $A_{2,j-1}$ hold. By $A_{1,j-1}$ we have
$|W_{j-1}| > \frac{1}{2}d_u|W| > \eps_*^2 |V_x|$.
Consider any triple $E$ containing $x$.
We bound $D^W_{z_j,E}(T_j-1)$ with a similar argument to that used in Lemma \ref{3marked}.

\nib{Case D.1.} First consider the case $z_j \in E$. Then $E^{T_j}=\ov{E} \sm z_j$ and $\nu^*>\nu$.
Since $W_j = W_{j-1} \cap F_{xz_j}(T_j-1)(y)$, Lemma \ref{consistent} gives
$F_{E^{T_j}}(T_j)[W_j] = F_{\ov{E}}(T_j-1)(y)[W_j] = F_{\ov{E}}(T_j-1)[W_{j-1}](y)$. Similarly,
$M_{E^{T_j},E}(T_j)[W_j] = M_{\ov{E},E}(T_j-1)(y)[W_j] = M_{\ov{E},E}(T_j-1)[W_{j-1}](y)$, so
$$D^W_{z_j,E}(T_j-1) = \{y \in F_{z_j}(T_j-1):
|M_{\ov{E},E}(T_j-1)[W_{j-1}](y)| > \theta_{\nu^*} |F_{\ov{E}}(T_j-1)[W_{j-1}](y)|  \}.$$
\COMMENT{
Formerly: By Lemma \ref{3exceptional}, $F_S(T_j)$ is $\eps_{\nu'_S(T_j),1}$-regular
for all $S \sub \ov{E}$, $|S| \ge 2$. [Unnecessary?]
By $A_{1,j-1}$, property $(*_{\ref{3initial}})$ holds at time $T_j-1$. [Meaning?]
}
Let $B'_{z_j}$ be the set of $y \in F_{z_j}(T_j-1) \sm E_{z_j}(T_j-1)$ such that we do not have
$$|F_{\ov{E}}(T_j-1)[W_{j-1}](y)|=(1\pm\eps_*)|F_{\ov{E}}(T_j-1)[W_{j-1}]|/|F_{z_j}(T_j-1)|.$$
Applying Lemma \ref{3crude}(vi) with $\GG=(W_{j-1},\{\es\})$
gives $|B'_{z_j}|<\eps_*|F_{z_j}(T_j-1)|$. 
(Note that we need $y \notin E_{z_j}(T_j-1)$ as this is
implicitly assumed to apply Lemma \ref{3crude}.)
Let $B_{z_j} = B'_{z_j} \cup E_{z_j}(T_j-1)$.
Then $|B_{z_j}|<2\eps_* |F_{z_j}(T_j-1)|$ by Lemma \ref{3exceptional}.
Now%
\COMMENT{1. Fixing: Only needed $E^W$ for precise estimates in C;
now forget it and use cruder $B'_{z_j}$. Also need to exclude $E_{z_j}(T_j-1)$
in definition of $B'_{z_j}$ as implicitly assumed by Lemma \ref{3crude}.
Now we are redefining $B_{z_j}$ as $B'_{z_j} \cup E_{z_j}(T_j-1)$...
2. Formerly: We exclude $B_{z_j}$ in the calculation because the bounds in Lemma \ref{3crude}
implicitly assume that $y \notin E_{z_j}(T_j-1)$. There is some overkill in excluding
$B'_{z_j}$ and $E^W_{z_j}(T_j-1)$ (which is contained in $B_{z_j}$).
}
\begin{align*}
& \sum_{y \in D^W_{z_j,E}(T_j-1)} |M_{\ov{E},E}(T_j-1)[W_{j-1}](y)|
 > \theta_{\nu^*} \sum_{y \in D^W_{z_j,E}(T_j-1) \sm B_{z_j}} |F_{\ov{E}}(T_j-1)[W_{j-1}](y)| \\
& \qquad >  (1-\eps_*)\theta_{\nu^*} (|D^W_{z_j,E}(T_j-1)|-2\eps_*|F_{z_j}(T_j-1)|)
|F_{\ov{E}}(T_j-1)[W_{j-1}]|/|F_{z_j}(T_j-1)|.
\end{align*}
We also have an upper bound
\begin{align*}
& \sum_{y \in D^W_{z_j,E}(T_j-1)} |M_{\ov{E},E}(T_j-1)[W_{j-1}](y)|
\le \sum_{y \in F_{z_j}(T_j-1)} |M_{\ov{E},E}(T_j-1)[W_{j-1}](y)| \\
& \qquad = |M_{\ov{E},E}(T_j-1)[W_{j-1}]| < \theta_{\nu} |F_{\ov{E}}(T_j-1)[W_{j-1}]|
\end{align*}
where the last inequality holds by $A_{2,j-1}$. Therefore \[\frac{|D^W_{z_j,E}(T_j-1)|}{|F_{z_j}(T_j-1)|}
< \frac{\theta_{\nu}}{(1-\eps_*)\theta_{\nu^*}} + 2\eps_* < \theta_{\nu^*}.\]

\nib{Case D.2.} Next consider the case $z_j \notin E$. Then $E^{T_j}=E^{T_j-1}=\ov{E}$.
Also $x \in \ov{E} \cap VN_H(z_j)$, so $\nu^*>\nu$.
Suppose first that $\ov{E} \in H(z_j)$. Then $F_{\ov{E}}(T_j)[W_j]=F_{\ov{E} z_j}(T_j-1)[W_{j-1}](y)$
and $M_{\ov{E},E}(T_j)[W_j] = M_{\ov{E},E}(T_j-1) \cap F_{\ov{E}}(T_j)[W_j]$, so%
\COMMENT{Also equals $M_{\ov{E},E}(T_j-1)[W_j] \cap F_{\ov{E}}(T_j)[W_j]$.}
$$D^W_{z_j,E}(T_j-1) = \left\{ y \in F_{z_j}(T_j-1) :
\frac{ |M_{\ov{E},E}(T_j-1) \cap F_{\ov{E} z_j}(T_j-1)[W_{j-1}](y)| }{
 |F_{\ov{E} z_j}(T_j-1)[W_{j-1}](y)| } > \theta_{\nu^*} \right\}.$$
Similarly to the previous case, letting $B'_{z_j}$ be the set of $y \in F_{z_j}(T_j-1) \sm E_{z_j}(T_j-1)$
such that we do not have
$$|F_{\ov{E} z_j}(T_j-1)[W_{j-1}](y)|=(1\pm\eps_*)|F_{\ov{E} z_j}(T_j-1)[W_{j-1}]|/|F_{z_j}(T_j-1)|,$$
we have $|B'_{z_j}|<\eps_*|F_{z_j}(T_j-1)|$. With $B_{z_j} = B'_{z_j} \cup E_{z_j}(T_j-1)$ we have
\begin{eqnarray*}
\Sa & := & \sum_{y \in D^W_{z_j,E}(T_j-1)} |M_{\ov{E},E}(T_j-1) \cap F_{\ov{E} z_j}(T_j-1)[W_{j-1}](y)| \\
& > & \theta_{\nu^*} \sum_{y \in D^W_{z_j,E}(T_j-1) \sm B_{z_j}} |F_{\ov{E} z_j}(T_j-1)[W_{j-1}](y)| \\
& > &  (1-\eps_*)\theta_{\nu^*} (|D^W_{z_j,E}(T_j-1)|-2\eps_*|F_{z_j}(T_j-1)|)
\frac{|F_{\ov{E} z_j}(T_j-1)[W_{j-1}]|}{|F_{z_j}(T_j-1)|}.
\end{eqnarray*}
We also have $\Sa \le \sum_{y \in F_{z_j}(T_j-1)} |M_{\ov{E},E}(T_j-1) \cap F_{\ov{E} z_j}(T_j-1)[W_{j-1}](y)|$.
This last sum counts all pairs $(y,P)$ with $P \in M_{\ov{E},E}(T_j-1)[W_{j-1}]$, $y \in F_{z_j}(T_j-1)$
and $P y \in F_{\ov{E} z_j}[W_{j-1}](T_j-1)$, so we can rewrite it as
$\Sa \le \sum_{P \in M_{\ov{E},E}(T_j-1)[W_{j-1}]} |F_{\ov{E} z_j}(T_j-1)[W_{j-1}](P)|$.
Then Lemma \ref{3crude} gives
$$|F_{\ov{E} z_j}(T_j-1)[W_{j-1}](P)| = (1 \pm \eps_*)
\frac{|F_{\ov{E} z_j}(T_j-1)[W_{j-1}]|}{|F_{\ov{E}}(T_j-1)[W_{j-1}]|}$$
for all but at most $\eps_*|F_{\ov{E}}(T_j-1)[W_{j-1}]|$ sets $P \in F_{\ov{E}}(T_j-1)[W_{j-1}]$.
Therefore
$$\Sa \le |M_{\ov{E},E}(T_j-1)[W_{j-1}]|(1+\eps_*)
\frac{|F_{\ov{E} z_j}(T_j-1)[W_{j-1}]|}{|F_{\ov{E}}(T_j-1)[W_{j-1}]|}
+ \eps_*|F_{\ov{E}}(T_j-1)[W_{j-1}]||F_{z_j}(T_j-1)|.$$
Combining this with the lower bound on $\Sa$ we obtain
\begin{align*}
& (1-\eps_*)\theta_{\nu^*} (|D^W_{z_j,E}(T_j-1)|/|F_{z_j}(T_j-1)|-2\eps_*) \\
& <  (1+\eps_*) \frac{|M_{\ov{E},E}(T_j-1)[W_{j-1}]|}{|F_{\ov{E}}(T_j-1)[W_{j-1}]|}
\ +\ \eps_*\frac{|F_{\ov{E}}(T_j-1)[W_{j-1}]||F_{z_j}(T_j-1)|}{|F_{\ov{E} z_j}(T_j-1)[W_{j-1}]|} .
\end{align*}
Now $|M_{\ov{E},E}(T_j-1)[W_{j-1}]| < \theta_{\nu} |F_{\ov{E}}(T_j-1)[W_{j-1}]|$ by $A_{2,j-1}$.
Also, since $x \in \ov{E}$, and since $F_{z_j}(T_j-1)[W_{j-1}]=F_{z_j}(T_j-1)$,
we can apply Lemma \ref{3crude}(vi) with $\GG=(W_{j-1},\{\es\})$ to get%
\COMMENT{First (? no!) app with restrict so give $\GG$; ok to repeat later}
$\frac{|F_{\ov{E}}(T_j-1)[W_{j-1}]||F_{z_j}(T_j-1)|}{|F_{\ov{E} z_j}(T_j-1)[W_{j-1}]|}
\le 2d_u^{-1} \ll \eps_*^{-1}$, so
\[\frac{|D^W_{z_j,E}(T_j-1)|}{|F_{z_j}(T_j-1)|} <
\frac{(1+\eps_*)\theta_{\nu}+\sqrt{\eps_*}}{(1-\eps_*)\theta_{\nu^*}} + 2\eps_* < \theta_{\nu^*}.\]

\nib{Case D.3.} It remains to consider the case when $z_j \notin E$ and $\ov{E} \notin H(z_j)$.
Since $x \in \ov{E} \cap VN_H(z_j)$ we have $|\ov{E}|\ge 2$
(otherwise we are in Case D.2).
Now $F_{\ov{E}^\le}(T_j)=F_{\ov{E}^\le}(T_j-1)[F_{\ov{E}^<}(T_j)]$ by Lemma \ref{build-update},
so $F_{\ov{E}^\le}(T_j)[W_j]=F_{\ov{E}^\le}(T_j-1)[F_{\ov{E}^<}(T_j)[W_j]]$ by Lemma \ref{*props}.
Also, $M_{\ov{E},E}(T_j)[W_j] = M_{\ov{E},E}(T_j-1) \cap F_{\ov{E}}(T_j)[W_j]$
by Lemma \ref{track-mark}, so
$$D^W_{z_j,E}(T_j-1) = \left\{ y \in F_{z_j}(T_j-1) :
\frac{ |M_{\ov{E},E}(T_j-1) \cap F_{\ov{E}^\le}(T_j-1)[F_{\ov{E}^<}(T_j)[W_j]]_{\ov{E}}| }{
 |F_{\ov{E}^\le}(T_j-1)[F_{\ov{E}^<}(T_j)[W_j]]_{\ov{E}}| } > \theta_{\nu^*} \right\}.$$
Let $I = \{S \subn \ov{E}: S \in H(z_j)\}$.
Then $P \in F_{\ov{E}^\le}(T_j-1)[F_{\ov{E}^<}(T_j)[W_j]]$ if and only if $P \in F_{\ov{E}^\le}(T_j-1)$,
for $x \notin S \in I$ we have $P_S \in F_S(T_j)$, i.e.\ $P_S y \in F_{S z_j}(T_j-1)$,
and for $x\in S \in I$ we have $P_S \in F_S(T_j)[W_j]$, i.e.\ $P_S y \in F_{S z_j}(T_j-1)[W_j]$.
When we choose $y \notin E_{z_j}(T_j-1)$, Lemma \ref{3crude} gives
$d_S(F(T_j))=(1\pm \eps_*)d_S(F(T_j-1))$ for $S \sub \ov{E}$ with $S \notin I$
and $d_S(F(T_j))=(1\pm \eps_*)d_S(F(T_j-1))d_{Sz_j}(F(T_j-1))$ for $\es \ne S \in I$. 
Let $d'_S$ denote $d_S(F(T_j)[W_j]):=d_S(F_{\ov{E}^\le}(T_j)[W_j])$ if $x \in S$ or $d_S(F(T_j))$ if $x \notin S$.
If $y \notin E^W_{z_j}(T_j-1)$ then $|W_j| > \frac{1}{2}d_u|W| > \eps_*^2 |V_x|$,
so regular restriction gives $d'_S=(1\pm \eps_*)d_S(F(T_j))$ for $S \sub \ov{E}$, $S \ne x$
and $d'_x = d(W_j) = |W_j|/|V_x|$. Applying Lemma \ref{3crude} we have
\begin{align*}
& d(F_{\ov{E}^\le}(T_j-1)[F_{\ov{E}^<}(T_j)[W_j]]_{\ov{E}}) = (1 \pm \eps_*) \prod_{S \sub \ov{E}} d'_S \\
& = (1 \pm 30\eps_*) \frac{|W_j|}{|F_x(T_j-1)|} \prod_{S \sub \ov{E}} d_S(F(T_j-1))
\prod_{\es \ne S \in I} d_{Sz_j}(F(T_j-1)).
\end{align*}
Also, $|F_{\ov{E}}(T_j-1)[W_{j-1}]| = (1 \pm \eps_*) \frac{|W_{j-1}|}{|F_x(T_j-1)|}|F_{\ov{E}}(T_j-1)|$
by Lemma \ref{3crude}.\COMMENT{Details?} If $y \notin E^W_{z_j}(T_j-1)$ then
$\frac{|W_j|}{|F_x(T_j)|} = (1 \pm \eps_*)\frac{|W_{j-1}|}{|F_x(T_j-1)|}$,
so $$ |F_{\ov{E}^\le}(T_j-1)[F_{\ov{E}^<}(T_j)[W_j]]_{\ov{E}}|
= (1 \pm 40\eps_*)|F_{\ov{E}}(T_j-1)[W_{j-1}]|\prod_{\es \ne S \in I} d_{Sz_j}(F(T_j-1)). $$
Similarly to the previous cases, now with $B_{z_j} = E^W_{z_j}(T_j-1) \cup E_{z_j}(T_j-1)$,
we have
\begin{eqnarray*}
\Sa & := & \sum_{y \in D^W_{z_j,E}(T_j-1)} |M_{\ov{E},E}(T_j-1) \cap F_{\ov{E}^\le}(T_j-1)[F_{\ov{E}^<}(T_j)[W_j]]_{\ov{E}}| \\
& > & \theta_{\nu^*} \sum_{y \in D^W_{z_j,E}(T_j-1) \sm B_{z_j}} |F_{\ov{E}^\le}(T_j-1)[F_{\ov{E}^<}(T_j)[W_j]]_{\ov{E}}| \\
& > &  (1-40\eps_*)\theta_{\nu^*} (|D^W_{z_j,E}(T_j-1)|-2\eps_*|F_{z_j}(T_j-1)|) \\
& & \times\ |F_{\ov{E}}(T_j-1)[W_{j-1}]| \prod_{\es \ne S \in I} d_{Sz_j}(F(T_j-1)).
\end{eqnarray*}
For any $P \in F_{\ov{E}}(T_j-1)$, let $F_{P,I}$ be the set of $y \in F_{z_j}(T_j-1)$
such that $P_S y \in F_{S z_j}(T_j-1)$ for all $S \in I$.
Let $B_I$ be the set of $P \in F_{\ov{E}}(T_j-1)$ such that we do not have
\[|F_{P,I}| = (1 \pm \eps_*)|F_{z_j}(T_j-1)| \prod_{\es \ne S \in I} d_{Sz_j}(F(T_j-1)).\]
Then Lemma \ref{3technical} gives $|B_I| \le \eps_*|F_{\ov{E}}(T_j-1)|$.
Also, since $F_{\ov{E}^<}(T_j)[W_j]*W_{j-1}=F_{\ov{E}^<}(T_j)[W_j]$,
Lemma \ref{*props}(iv) gives
\[ F_{\ov{E}^\le}(T_j-1)[F_{\ov{E}^<}(T_j)[W_j]]
= F_{\ov{E}^\le}(T_j-1)[W_{j-1}][F_{\ov{E}^<}(T_j)[W_j]].\]
Now $\Sa \le \sum_{y \in F_{z_j}(T_j-1)}
|M_{\ov{E},E}(T_j-1) \cap F_{\ov{E}^\le}(T_j-1)[W_{j-1}][F_{\ov{E}^<}(T_j)[W_j]]_{\ov{E}}|$,
which counts all pairs $(y,P)$ with $P \in M_{\ov{E},E}(T_j-1)[W_{j-1}]$ and $y \in F_{P,I}$, so
\begin{align*}
\Sa \le & |M_{\ov{E},E}(T_j-1)[W_{j-1}]|(1 \pm \eps_*)|F_{z_j}(T_j-1)| \prod_{\es \ne S \in I} d_{Sz_j}(F(T_j-1)) \\
& +\ \eps_*|F_{\ov{E}}(T_j-1)[W_{j-1}]||F_{z_j}(T_j-1)|.
\end{align*}
Combining this with the lower bound on $\Sa$ we obtain
\begin{align*}
& (1-40\eps_*) \theta_{\nu^*} \left( \frac{|D^W_{z_j,E}(T_j-1)|}{|F_{z_j}(T_j-1)|}
- 2\eps_* \right) \\
& < \ (1+\eps_*)\frac{|M_{\ov{E},E}(T_j-1)[W_{j-1}]|}{|F_{\ov{E}}(T_j-1)[W_{j-1}]|}
+ \eps_* \prod_{\es \ne S \in I} d_{Sz_j}(F(T_j-1))^{-1}.
\end{align*}
Now $|M_{\ov{E},E}(T_j-1)[W_{j-1}]| < \theta_{\nu} |F_{\ov{E}}(T_j-1)[W_{j-1}]|$ by $A_{2,j-1}$,
and all densities are at least $d_u \gg \eps_*$,
so again we have
\[\frac{|D^W_{z_j,E}(T_j-1)|}{|F_{z_j}(T_j-1)|} <
\frac{(1+\eps_*)\theta_{\nu}+\sqrt{\eps_*}}{(1-40\eps_*)\theta_{\nu^*}}
+ 2\eps_* < \theta_{\nu^*}.\]
This proves Claim D.

\nib{Claim E.} Conditional on the events $A_{i,j'}$, $1 \le i \le 3$, $0 \le j' < j$
and the embedding up to time $T_j-1$, the probability that $A_{3,j}$ does not
hold is at most $\theta'_{12D}$.

\nib{Proof.} Since $A_{1,j-1}$ and $A_{2,j-1}$ hold,
Claim C gives $|E^W_{z_j}(T_j-1) \sm E_{z_j}(T_j-1)| < \eps_* |F_{z_j}(T_j-1)|$
and Claim D gives $|D^W_{z_j,E}(T_j-1)| < \theta_{12D} |F_{z_j}(T_j-1)|$
for any $E$ containing $x$. We also have $|OK_{z_j}(T_j-1)| > (1-\theta_*)|F_{z_j}(T_j-1)|$
by Corollary \ref{3ok}. Since $y=\phi(z_j)$ is chosen uniformly at random in
$OK_{z_j}(T_j-1)$, the probability that $y \in E^W_{z_j}(T_j-1)$
or $y \in D^W_{z_j,E}(T_j-1)$ for any $E$ containing $x$
is at most $(\eps_* + D\theta_{12D})/(1-\theta_*) < \theta'_{12D}$.
This proves Claim E.%
\COMMENT{
Ref pointed out that $E_{z_j}$ already excluded from OK so no need to worry about it here.
Formerly: ...$|E_{z_j}(T_j-1)| < \eps_* |F_{z_j}(T_j-1)|$ by Lemma \ref{3exceptional}...
}

To finish the proof of the lemma, suppose that all the events $A_{i,j}$, $1 \le i \le 3$,
$1 \le j \le g$ hold. Then $A_{1,g}$ gives
$|F_x(T_g) \cap W| = |W_g| = (1 \pm \eps_*)|W||F_x(T_g)|/|F_x(T_0)|
> \frac{1}{2}d_u|W| > \eps_*^2 |V_x|$.
Also, since all of $VN_H(x)=\{z_1,\cdots,z_g\}$ has been embedded at time $T_g$,
for every triple $E$ containing $x$ we have $E^{T_g}=x$,
and $|M_{x,E}(T_g) \cap W| < \theta'_{12D}|F_x(T_g) \cap W|$ by $A_{2,g}$.
Now $A_x \cap W$ is obtained from $F_x(T_g) \cap W$
by deleting all $M_{x,E}(T_g) \cap W$ for triples $E$ containing $x$,
so $|A_x \cap W| > (1 - D\theta'_{12D})|F_x(T_g) \cap W|$.
In particular, $A_x \cap W$ is non-empty.
If any event $A_{i,j}$ fails then $A_{3,j}$ fails (by Claim B)
and so by Claim E and a union bound over $1 \le j \le g \le 2D$
we can bound the failure probability by $\theta_*$. \qed

Our final lemma in this subsection is similar to the previous one,
but instead of asking for a set $W$ of vertices to contain
an available vertex for $x$, we ask for some particular vertex $v$
to be available for $x$. Recall that $x\in B$ and we start
embedding $VN_H(x)$ at time $T_1$.

\begin{lemma} \label{3-x-to-v}
For any $v \in V_x$, conditional on any embedding of the vertices $\{s(u): u < T_1\}$
that does not use $v$, with probability at least $p$ we have
$\phi(H(x)) \sub (G \sm M)(v)$, so $v \in A_x$.
\COMMENT{
1. $v$ marked forbidden for $x$ with $\mb{P}<\theta_*$.
2. $p \ll d_u$. 
3. Cut then restored `allocated'.
}
\end{lemma}

\nib{Proof.}
We estimate the probability that $\phi(H(x)) \sub (G \sm M)(v)$
using arguments similar to those we are using to embed $H$ in $G\sm M$.
The structure of the proof is very similar to that of Lemma \ref{3initial}.
Here we will see the purpose of properties (ii) and (iii) in the definition
of super-regularity, which ensure that every $v \in V_x$ is a potential image of $x$.
For $z \in VN_H(x)$ we write\index{$\alpha_z$}
\COMMENT{Formerly: The following notation will be convenient.}
\[\alpha_z = \frac{|F_{xz}(T_0)(v)|}{d_{xz}(F(T_0))|F_z(T_0)|}
= \frac{d_z(F(T_0)(v))}{d_{xz}(F(T_0))d_z(F(T_0))} .\]
We consider a vertex $z$ to be {\em allocated}\index{allocated}
if $z$ is embedded or $z=x$. For $\es\ne S\in H$ unembedded
we define $\nu''_S(t)$ as follows.\index{$\nu''_S(t)$}
When $|S|=3$ we let $\nu''_S(t) = \nu_S(t)$.
When $|S|=1,2$ we let $\nu''_S(t) = \nu_S(t)+K$, where $K$ is the
maximum value of $\nu''_{Sx'}(t')$ over allocated vertices $x'$
with $S \in H(x')$; if there is no such vertex $x'$ we let $\nu''_S(t) = \nu_S(t)$.
Thus $\nu''_S(t)$ is defined similarly to $\nu'_S(t)$,
replacing `embedded' with `allocated'. We have $\nu''_S(t)\ge \nu'_S(t)$
and Lemma \ref{3observe}(iii-vi) hold replacing $\nu'$ with $\nu''$.
\COMMENT{Could say replacing `embedded' with `allocated' in (vi)...}

Suppose $1 \le j \le g$ and that we are considering the embedding of $z_j$.
We interpret quantities at time $T_j$ with the embedding $\phi(z_j)=y$,
for some as yet unspecified $y \in F_{z_j}(T_j-1)$.

We define {\em exceptional} sets $E^v_{z_j}(T_j-1) \sub F_{xz_j}(T_j-1)(v)$
by $y \in F_{z_j}(T_j-1) \sm E^v_{z_j}(T_j-1)$ if and only if
for every unembedded $\es \ne S \in H(x) \cap H(z_j)$,
\index{exceptional}\index{$E^v_{z_j}(T_j-1)$}
\begin{equation}\label{eq:*3xtov}\tag{$*_{\ref{3-x-to-v}}$}
\left.
\begin{aligned}
& F_{Sx}(T_j)(v) \mbox{ is } \eps_{\nu''_S(T_j),1}\mbox{-regular if } |S|=2, \\
& d_S(F(T_j)(v)) = (1 \pm \eps_{\nu''_S(T_j),1})d_S(F(T_j))d_{Sx}(F(T_j)) \mbox{ if } |S|=2, \\
& d_S(F(T_j)(v)) =  (1 \pm \eps_{\nu''_S(T_j),1})d_S(F(T_j))d_{Sx}(F(T_j))\alpha_S \mbox{ if } |S|=1.
\end{aligned}
\qquad \right\}
\end{equation}

Here we use the notation $d_S(F(T_j)(v))=d_S(F_{Sx^\le}(T_j)(v))$.
Let $Y$ be the set of vertices at distance at most $3$ from $x$ in $H$ \index{$Y$}
and let $H' = \{S \in H: S \sub Y\}$. \index{$H'$}
For any $Z \sub Y$ and unembedded $S \in H$ we define\index{$Z$}\index{$F(T_j)^{Z*v}_{S^\le}$}
\[F(T_j)^{Z*v}_{S^\le} = F_{S^\le}(T_j)\left[
\bigcup_{S' \sub Z \cap S, S' \in H(x)} F_{S' x}(T_j)(v)\right].\]
Thus $F(T_j)^{Z*v}_S$ consists of all sets $P \in F_S(T_j)$
such that $P_{S'} v \in F_{S' x}(T_j)$ for all $S' \sub Z \cap S$ with $S' \in H(x)$.
For any triple $E$, we use the notation $\ov{E}=E^{T_j-1}$,\index{$\ov{E}$}
$\nu=\nu''_{\ov{E}}(T_j-1)$\index{$\nu$}
and $\nu^*=\nu''_{E^{T_j}}(T_j)$\index{$\nu^*$}
similarly to the previous lemma, replacing $\nu'$ with $\nu''$.
For $Z \sub Y$ and $E \in U(z_j)$
we define sets of {\em dangerous} vertices by%
\index{dangerous}\index{$D^{Z*v}_{z_j,E}(T_j-1)$}
\[ D^{Z*v}_{z_j,E}(T_j-1) = \left\{y \in F(T_j-1)^{Z*v}_{z_j}:
\left|M_{E^{T_j},E}(T_j) \cap F(T_j)^{Z*v}_{E^{T_j}}\right|
> \theta_{\nu^*}\left|F(T_j)^{Z*v}_{E^{T_j}}\right| \right\}. \]

The strategy of the proof is to analyse the event that
all the complexes $F(T_j)^{Z*v}_{E^{T_j\le}}$ are well-behaved,
meaning informally that they are regular, have roughly expected
densities and do not have too many marked edges. The regularity and density
properties will hold if we choose $y=\phi(z_j) \notin E^v_{z_j}(T_j-1)$,
and the marked edges will be controlled if we choose $y=\phi(z_j) \notin D^{Z*v}_{z_j,E}(T_j-1)$.

We think of $Z$ as the {\em sphere of influence}, \index{sphere of influence}
as it defines the sets which we restrict to be in the neighbourhood of $v$.
Our eventual goal is that all sets in $H(x)$ should be embedded in $G(v)$,
but to achieve this we need to consider arbitrary choices of $Z \sub Y$.
For later use in the proof we record here some properties of $Z$
that follow directly from the definitions.
\begin{itemize}
\item[(i)] $F(T_j)^{Z*v}_S = F(T_j)^{(Z\cap S)*v}_S$.
\item[(ii)] $F(T_j-1)^{Z*v}_{z_j}$ is $F_{z_j}(T_j-1)$ if $z_j \notin Z$
or $F_{xz_j}(T_j-1)(v)$ if $z_j \in Z$.
\item[(iii)] $F(T_j)^{\es*v}_S = F(T_j)_S$, so $D^{\es*v}_{z_j,E}(T_j-1) = D_{z_j,E}(T_j-1)
= \{y \in F_{z_j}(T_j-1): |M_{E^{T_j},E}(T_j)| > \theta_{\nu^*}|F(T_j)_{E^{T_j}}| \}$,
as defined in the description of the algorithm.
\item[(iv)] If $Z'=Z\cup z_j$ then $F(T_j)^{Z'*v}_{E^{T_j}}=F(T_j)^{Z*v}_{E^{T_j}}$
by (i), since $z_j \notin E^{T_j}$, so

$D^{Z'*v}_{z_j,E}(T_j-1) = D^{Z*v}_{z_j,E}(T_j-1) \cap F_{xz_j}(T_j-1)(v)$.
\item[(v)] If $S \sub Z$ and $S \in H(x)$ then $F(T_j)^{Z*v}_S = F_{Sx}(T_j)(v)$.
\end{itemize}

Write $B_{z_j} = E^v_{z_j}(T_j-1) \cup E_{z_j}(T_j-1)$. \index{$B_{z_j}$}
Recall that $U(z_j)$ is the set of triples $E$ with $\ov{E} \cap VN_H(z_j)z_j \ne \es$.
We consider the following events: \index{$A_{i,j}$}

\begin{itemize}
\item
$A_{1,j}$, $j \ge 1$ is the event that property $(*_{\ref{3-x-to-v}})$ above holds.
We also define $A_{1,0}$ to be the event that
$|G(v)_z \cap V_z(T_0)| > (1-\delta_B^{1/3})|G(v)_z|$
for every $z \in VN_H(x)$.
\item
$A_{2,j}$, $j \ge 0$ is the event that for every triple $E \in H'$ and $Z \sub E$ we have
\[\left|M_{E^{T_j},E}(T_j) \cap F(T_j)^{Z*v}_{E^{T_j}}\right|
\le \theta'_{\nu^*}\left|F(T_j)^{Z*v}_{E^{T_j}}\right|.  \]
\item
$A_{3,j}$, $j \ge 0$ is the event that for every triple $E \in H'$, $Z \sub E$
and $\es \ne S \sub E^{T_j}$,

$\qquad F(T_j)^{Z*v}_S$ is $\eps_{\nu''_S(T_j),2}$-regular for $|S| \ge 2$, with
\[d_S(F(T_j)^{Z*v}) = \left\{ \begin{array}{l}
(1 \pm \eps_{\nu''_S(T_j),2})d_S(F(T_j))d_{Sx}(F(T_j)) \mbox{ if } S \sub Z, S \in H(x), |S|=2 \\
(1 \pm \eps_{\nu''_S(T_j),2})d_S(F(T_j))d_{Sx}(F(T_j))\alpha_S \mbox{ if } S \sub Z, S \in H(x), |S|=1 \\
(1 \pm \eps_{\nu''_S(T_j),2})d_S(F(T_j)) \mbox{ otherwise.}
\end{array} \right. \]
\item
$A_{4,j}$, $j \ge 1$ is the event that $y=\phi(z_j)$ is chosen in $OK^v_{z_j}(T_j-1)$,
defined to be the subset of $F_{xz_j}(T_j-1)(v)$ obtained by deleting the sets
$B_{z_j}$ and $D^{Z'*v}_{z_j,E}(T_j-1)$ for all $E \in U(z_j)$, $Z \sub E$, $Z'=Z \cup z_j$.
We also define $A_{4,0}$ to be the event that holds with probability $1$.
\index{$OK$}
\end{itemize}

By property (iv) above, an equivalent definition of $OK^v_{z_j}(T_j-1)$
is the subset of $F_{xz_j}(T_j-1)(v)$ obtained by deleting the sets $B_{z_j}$
and $D^{Z'*v}_{z_j,E}(T_j-1)$ for all $E \in U(z_j)$ and $Z' \sub E \cup z_j$.
Also, since $OK_{z_j}(T_j-1)$ is obtained from $F_{z_j}(T_j-1)$ by deleting $E_{z_j}(T_j-1)$
and $D_{z_j,E}(T_j-1) = D^{\es*v}_{z_j,E}(T_j-1)$ for $E \in U(z_j)$
we have $OK^v_{z_j}(T_j-1) \sub OK_{z_j}(T_j-1)$.

We will use the following notation throughout: $Z$ is a subset of $Y$,
$Z'=Z\cup z_j$, $I = \{S \sub Z: S \in H(x)\}$, $I' = \{S \sub Z': S \in H(x)\}$.
We divide the remainder of the proof into a series of claims.%
\index{$Z$}\index{$Z'$}\index{$I$}\index{$I'$}

\nib{Claim A.} The events $A_{1,0}$, $A_{2,0}$, $A_{3,0}$ and $A_{4,0}$ hold with high probability.

\nib{Proof.} $A_{4,0}$ holds by definition. For any $z \in VN_H(x)$,
$d_z(G(v))=(1\pm\eps')d_{xz}(G)d_z(G) > d_u$ by condition (ii) of super-regularity,
so $A_{1,0}$ holds with high probability by Lemma \ref{3uniform}.
Next recall that no vertex at distance within $4$ of $x$ has been embedded at time $T_0$,
so $F_z(T_0)=V_z(T_0)$ for any $z$ within distance $3$ of $x$.
(This is why we choose the buffer vertices to be at mutual distance at least $9$.)
We have $|F_z(T_0)|>(1-2\sqrt{\delta_B})|V_z|$ by (\ref{eq:3fresh})
and $|F_z(T_0) \cap G(v)_z| > (1-\delta_B^{1/3})|G(v)_z|$ by $A_{1,0}$.
For any $S \in H'$ with $|S| \ge 2$, $F_S(T_0)$ is the restriction of $G_S$
to $((F_z(T_0):z \in S),\{\es\})$. Since $G_S$ is $\eps$-regular,
$F_S(T_0)$ is $\eps'$-regular with $d_S(F(T_0))=(1\pm \eps')d_S(G)$.
It also follows that $\alpha_z > 1-2\delta_B^{1/3}$.%
\COMMENT{Technically want whp for all x and v...}

Now we show that $A_{3,0}$ holds. Consider any triple $E \in H'$ and $Z \sub E$.
Suppose $\es \ne S \sub E$. There are two cases according to whether $S \in I$.
Suppose first that $S \in I$. By property (v) above, $F(T_0)^{Z*v}_{S^\le}=F_{Sx^\le}(T_0)(v)$
is the restriction of $G_{Sx^\le}(v)$ to $((F_z(T_0) \cap G(v)_z:z \in S),\{\es\})$.
If $|S|=2$ then $G_{Sx^\le}(v)$ is $\eps'$-regular
and $d_S(G(v))=(1\pm\eps')d_S(G)d_{Sx}(G)$ by condition (ii) of super-regularity.
Then by regular restriction $F_{S x}(T_0)(v)$ is $\eps_{0,0}$-regular
and $d_S(F(T_0)(v))=(1\pm\eps_{0,0})d_S(G)d_{Sx}(G)
=(1\pm 2\eps_{0,0})d_S(F(T_0))d_{Sx}(F(T_0))$.
Also, if $|S|=1$ then $d_S(F(T_0))d_{Sx}(F(T_0))\alpha_S = d_S(F(T_0)(v))$ by definition.
This gives the properties required by $A_{3,0}$ when $S \in I$.
In fact, we have the stronger statements in which $\eps_{\nu''_S(T_0),2}$ is replaced by $\eps_{0,1}$, say.
On the other hand, if $S \notin I$ then $F(T_0)^{Z*v}_{S^\le}$ is the restriction of $F_{S^\le}(T_0)$
to $\cup_{S' \sub S\cap Z,S' \in H(x)} F_{S' x}(T_0)(v)$. Since $F_{S^\le}(T_0)$
is $\eps'$-regular and $F_{S' x}(T_0)(v)$ is $\eps_{0,1}$-regular for $S' \in I$,
by regular restriction $F(T_0)^{Z*v}_S$ is $\eps_{0,0}$-regular
with $d_S(F(T_0)^{Z*v}) = (1 \pm \eps')d_S(F(T_0))$. Thus $A_{3,0}$ holds.

It remains to show that $A_{2,0}$ holds.
We will abuse notation and let $I$ also denote the subcomplex
$\{i(S):S \in I\}$ of $\binom{[r]}{\le 3}$.
Then $F(0)^{Z*v}_E = G^{I_v}_E$ as defined in Definition \ref{def-3preplus}.
By property (iii) of super-regularity $|M_E \cap G^{I_v}_E| \le \theta|G^{I_v}_E|$
and $G^{I_v}_{E^\le}$ is $\eps'$-regular with $S$-density (for $S\sub E$)
equal to $(1 \pm \eps')d_S(G)d_{S x}(G)$ if $S \in I$
or $(1 \pm \eps')d_S(G)$ otherwise.
Now $F(T_0)^{Z*v}_z$ is $F_z(T_0) \cap G(v)_z$ if $z \in Z \cap H(x)$ or $F_z(T_0)$ otherwise,
and similarly $G^{I_v}_z$ is $G(v)_z$ if $z \in Z \cap H(x)$ or $G_z$ otherwise.
Either way we have $|F(T_0)^{Z*v}_z| > (1-\delta_B^{1/3})|G^{I_v}_z|$
by the estimates recalled above. For $|S| \ge 2$ we showed above that $d_S(F(T_0)^{Z*v})$
is $(1\pm 2\eps_{0,0})d_S(F(T_0))d_{Sx}(F(T_0))$ if $S \in I$
or $d_S(F(T_0)^{Z*v}) = (1 \pm \eps')d_S(F(T_0)$ if $S \notin I$.
Recalling that $d_{S'}(F(T_0))=(1\pm \eps')d_{S'}(G)$ for $S' \in H'$ with $|S'|\ge 2$,
Lemma \ref{3absolute} gives
\[\frac{|F(T_0)^{Z*v}_E|}{|G^{I_v}_E|}=\frac{d(F(T_0)^{Z*v}_E)}{d(G^{I_v}_E)}
= \frac{(1 \pm 9\eps_{0,2})\prod_{z\in E}|F(T_0)^{Z*v}_z|}{(1 \pm 8\eps')\prod_{z\in E}|G_v|}
> (1-10\eps_{0,2})(1-\delta_B^{1/3})^3 > 1/2.\]
Now $|M_{E,E}(T_0) \cap F(T_0)^{Z*v}_E| \le |M_E \cap G^{I_v}_E|
\le \theta|G^{I_v}_E| \le 2\theta|F(T_0)^{Z*v}_E|$,
giving even a stronger bound on the marked edges than is needed.
Thus $A_{2,0}$ holds.
\COMMENT{Cut: By Lemma \ref{3absolute} we have
$d_E(G^{I_v})=(1 \pm 8\eps')\prod_{S\sub E}d_S(G)\prod_{S \in I}d_{S x}(G)$.
}

\nib{Claim B.} Suppose $A_{3,j-1}$ holds.
If $y \notin E^v_{z_j}(T_j-1)$ then $A_{1,j}$ and $A_{3,j}$ hold,
and if $y \notin D^{Z*v}_{z_j,E}(T_j-1)$ then $A_{2,j}$ holds.
Thus $A_{3,j-1}$ and $A_{4,j}$ imply $A_{1,j}$, $A_{2,j}$ and $A_{3,j}$.

\nib{Proof.} Suppose $y \notin E^v_{z_j}(T_j-1)$. Then $A_{1,j}$ holds by definition.
$A_{3,j}$ follows from $A_{1,j}$ similarly to the case $j=0$ considered in Claim A.
To see this, consider any triple $E \in H'$ and $Z \sub E$.
Suppose $\es \ne S \sub E$ is unembedded. If $S \sub Z$ and $S \in H(x) \cap H(z_j)$
then $F(T_j)^{Z*v}_S=F_{Sx}(T_j)(v)$ satisfies $(*_{\ref{3-x-to-v}})$
by $A_{1,j}$, so we have the properties required by $A_{3,j}$.
In fact, we have the stronger statements in which $\eps_{\nu''_S(T_j),2}$ is replaced by $\eps_{\nu''_S(T_j),1}$.
For any other $S$ we use the definition of $F(T_j)^{Z*v}_{S^\le}$ as the restriction of $F_{S^\le}(T_j)$
to $\cup_{S' \sub S\cap Z,S' \in H(x)} F_{S' x}(T_j)(v)$. Note that by Lemmas \ref{build-update} and \ref{*props}
we get the same result if we replace $F_{S' x}(T_j)(v)$ by $F_{S' x}(T_j-1)(v)$ for those $S' \notin H(z_j)$.
Now $F_{S'}(T_j)$ is $\eps_{\nu''_{S'}(T_j),1}$-regular for $S' \sub S$ by Lemma \ref{3exceptional},
$F_{S' x}(T_j)(v)$ is $\eps_{\nu''_{S'}(T_j),1}$-regular for $S' \sub S \cap Z$, $S' \in H(x) \cap H(z_j)$
and $F_{S' x}(T_j-1)(v)$ is $\eps_{\nu''_{S'}(T_j-1),2}$-regular for $S' \in I$ by $A_{3,j-1}$.
So by regular restriction $F(T_j)^{Z*v}_S$ is $\eps_{\nu''_S(T_j),2}$-regular
with $d_S(F(T_j)^{Z*v}) = (1 \pm \eps_{\nu''_S(T_j),2})d_S(F(T_j))$. Thus $A_{3,j}$ holds.%
\COMMENT{Should give more details for app of build-update and *-props?}

Next consider any triple $E \in H'$ and $Z \sub E$.
Suppose $y \notin D^{Z*v}_{z_j,E}(T_j-1)$.
If $E \in U(z_j)$ then by definition we have
$|M_{E^{T_j},E}(T_j) \cap F(T_j)^{Z*v}_{E^{T_j}}|
\le \theta_{\nu^*}|F(T_j)^{Z*v}_{E^{T_j}}|$,
which is a stronger bound than required.
On the other hand, if $E \notin U(z_j)$ then consider the most recent time
$T_{j'} < T_j$ when we embedded $z_{j'}$ with $E \in U(z_{j'})$,
setting $j'=0$ if there is no such time. Then $E^{T_{j'}} = E^{T_j}$,
and by the stronger bound at time $j'$
we have $|M_{E^{T_j},E}(T_{j'}) \cap F(T_{j'})^{Z*v}_{E^{T_j}}|
\le \theta_{\nu^*}|F(T_{j'})^{Z*v}_{E^{T_j}}|$.
(Recall that we also obtained a stronger bound for $A_{2,0}$ in Claim A.)
Now $F(T_j)^{Z*v}_{E^{T_j}}$ is obtained from $F(T_{j'})^{Z*v}_{E^{T_j}}$
by deleting at most $2D$ vertices $\phi(z_{j^*})$, $j'+1 \le j^* \le j$
and the sets containing them, so
$|F(T_j)^{Z*v}_{E^{T_j}}|\ge (1-\eps_*)|F(T_{j'})^{Z*v}_{E^{T_j}}|$;
this can be seen by regular restriction and Lemma \ref{3absolute},
or simply from the fact that a trivial bound for the number of deleted sets
has a lower order of magnitude for large $n$. Therefore
$|M_{E^{T_j},E}(T_j) \cap F(T_j)^{Z*v}_{E^{T_j}}|
\le |M_{E^{T_j},E}(T_{j'}) \cap F(T_{j'})^{Z*v}_{E^{T_j}}|
\le \theta_{\nu^*}|F(T_{j'})^{Z*v}_{E^{T_j}}|
< \theta'_{\nu^*}|F(T_j)^{Z*v}_{E^{T_j}}|$, so $A_{2,j}$ holds.

\nib{Claim C.} If $A_{1,j-1}$ and $A_{3,j-1}$ hold then
$|E^v_{z_j}(T_j-1) \sm E_{z_j}(T_j-1)| < \eps_* |F_{xz_j}(T_j-1)(v)|$.
\COMMENT{$\eps_*|F_x(T_j-1)|$ suffices, but this is clearer.}

\nib{Proof.} For any $S \in H$ we write $\nu''_S=\nu''_S(T_j-1)$ and $\nu^*_S=\nu''_S(T_j)$.
Consider any unembedded $\es \ne S \in H(x) \cap H(z_j)$.%
\index{$\nu''_S$}\index{$\nu^*_S$}

\nib{Case C.1.} Suppose first that $S=z$ has size $1$.
Note that $\nu^*_z > \max\{\nu''_z,\nu''_{z_jz},\nu''_{xz}\}$
by the analogue of Lemma \ref{3observe} for $\nu''$.
We consider two cases according to whether $xz_jz \in H$.

\nib{Case C.1.i.} Suppose that $xz_jz \in H$.
Then  $F_{xz}(T_j)(v)=F_{xz_jz}(T_j-1)(yv)$. % $F_{xz}(T_j)=F_{xz_jz}(T_j-1)(y)$, so
Since $F_{xz_jz}(T_j-1)(v)=F(T_j-1)^{z_jz*v}_{z_jz}$,%
\COMMENT{confused what happens if $S$ not in $Z$, until realised can choose any $Z$!}
using $A_{3,j-1}$ we have $d_{z_jz}(F(T_j-1)(v))=(1\pm\eps_{\nu''_{z_jz},2})
d_{xz_jz}(F(T_j-1))d_{z_jz}(F(T_j-1)) > d_u^2/2$ and
$F_{xz_jz}(T_j-1)(v)$ is $\eps_{\nu''_{z_jz},2}$-regular.
We also have $d_z(F(T_j-1)(v))=(1\pm\eps_{\nu''_z,2})d_{xz}(F(T_j-1))d_z(F(T_j-1))\alpha_z$.
By Lemma \ref{2neighbour}, for all but at most
$\eps_{\nu''_{z_jz},3}|F_{xz_j}(T_j-1)(v)|$ vertices $y \in F_{xz_j}(T_j-1)(v)$
we have $|F_{xz}(T_j)(v)|=|F_{xz_jz}(T_j-1)(v)(y)|
= (1\pm\eps_{\nu''_{z_jz},3})d_{z_jz}(F(T_j-1)(v))|F_{xz}(T_j-1)(v)|$, so
\begin{align*}
&d_z(F(T_j)(v))  =  (1\pm\eps_{\nu''_{z_jz},3})d_{z_jz}(F(T_j-1)(v))d_z(F(T_j-1)(v)) \\
& =  (1\pm\eps_{\nu''_{z_jz},3})(1\pm\eps_{\nu''_{z_jz},2})(1\pm\eps_{\nu''_z,2}) \times \\
& \qquad d_{xz_jz}(F(T_j-1))d_{z_jz}(F(T_j-1)) d_{xz}(F(T_j-1))d_z(F(T_j-1))\alpha_z.
\end{align*}
Also if $y \notin E_{z_j}(T_j-1)$ then $(*_{\ref{3alg}})$ gives
\begin{align*}
&d_{xz}(F(T_j))=(1 \pm \eps_{\nu^*_{xz},0})d_{xz}(F(T_j-1))d_{xz_jz}(F(T_j-1)), \mbox{ and }\\
&d_z(F(T_j))=(1 \pm \eps_{\nu^*_z,0})d_z(F(T_j-1))d_{z_jz}(F(T_j-1)).
\end{align*}
Thus for such $y$ we have the required estimate
$d_z(F(T_j)(v))=(1\pm\eps_{\nu^*_z,1})d_{xz}(F(T_j))d_z(F(T_j))\alpha_z$.
\COMMENT{Also using $\nu^*_z \ge \nu^*_{xz}$...}

\nib{Case C.1.ii.} Suppose that $xz_jz \notin H$. Since $z_jz \in H$ and $xz_j \in H$
we have
\[F_{xz^\le}(T_j)=F_{xz^\le}(T_j-1)[(F_{xz_j}(T_j-1)(y),F_{zz_j}(T_j-1)(y),\{\es\})],\]
i.e.\ $F_{xz}(T_j)$ is the bipartite subgraph of $F_{xz}(T_j-1)$ induced
by $F_{xz_j}(T_j-1)(y)$ and $F_{zz_j}(T_j-1)(y)$. Then we have
$F_{xz}(T_j)(v) = F_{xz}(T_j-1)(v) \cap F_{zz_j}(T_j-1)(y)$.
Now $F_{zz_j}(T_j-1)$ is $\eps_{\nu''_{z_jz},1}$-regular by Lemma \ref{3exceptional}
and $d_z(F(T_j-1)(v))=(1\pm\eps_{\nu''_z,2})d_{xz}(F(T_j-1))d_z(F(T_j-1))\alpha_z > d_u^2/2$ by $A_{3,j-1}$.
Then by Lemmas \ref{2restrict} and \ref{2neighbour}, for all but at most
$\eps_{\nu''_{z_jz},2}|F_{xz_j}(T_j-1)(v)|$ vertices $y \in F_{xz_j}(T_j-1)(v)$
we have $|F_{xz}(T_j)(v)| = |F_{xz}(T_j-1)(v) \cap F_{zz_j}(T_j-1)(y)|
= (1\pm\eps_{\nu''_{z_jz},2})d_{zz_j}(F(T_j-1))|F_{xz}(T_j-1)(v)|$.
This gives
\begin{eqnarray*}
d_z(F(T_j)(v)) & = & (1\pm\eps_{\nu''_{z_jz},2})d_{zz_j}(F(T_j-1))d_z(F(T_j-1)(v))\\
& = & (1\pm\eps_{\nu''_{z_jz},2})(1\pm\eps_{\nu''_z,2})
d_{zz_j}(F(T_j-1))d_{xz}(F(T_j-1))d_z(F(T_j-1))\alpha_z.
\end{eqnarray*}
We also have $d_z(F(T_j))=(1 \pm \eps_{\nu^*_z,0})d_z(F(T_j-1))d_{zz_j}(F(T_j-1))$
if $y \notin E_{z_j}(T_j-1)$ by $(*_{\ref{3alg}})$
and $d_{xz}(F(T_j))=(1 \pm \eps_{\nu^*_{xz},0}^{1/2})d_{xz}(F(T_j-1))$ by Lemma \ref{2restrict},
and using $xz_jz \notin H$.
Thus for such $y$ we have the required estimate
$d_z(F(T_j)(v))=(1\pm\eps_{\nu^*_z,1})d_{xz}(F(T_j))d_z(F(T_j))\alpha_z$.

\nib{Case C.2.} The remaining case is when $S=z'z$ has size $2$.
Note that $xz_jz'z$ is $r$-partite, as $z'z \in H(x) \cap H(z_j)$ and $xz_j \in H$.
Note also that $\nu^*_{z'z} > \max\{\nu''_{z'z},\nu''_{z_jz'z},\nu''_{xz'z}\}$.
Consider the complex $$ J = F_{z_jz'z^\le}(T_j-1)\left[
\bigcup_{S' \sub z_jz'z, S' \in H(x)} F_{S' x}(T_j-1)(v)\right]. $$
We claim that $F_{xz'z}(T_j)(v)=J(y)$. To see this, note first that
$F_{xz'z^\le}(T_j)=F_{xz'z^\le}(T_j-1)[F_{xz'z^<}(T_j)]$ by Lemma \ref{build-update}.
Then by Lemma \ref{*props} we can write
$F_{xz'z^\le}(T_j)=F_{xz'z^\le}(T_j-1)*F_{xz'z^<}(T_j)$ and so
\begin{eqnarray*}
F_{xz'z^\le}(T_j)(v) &=& F_{z'z^\le}(T_j-1)*F_{z'z^\le}(T_j)
*F_{xz'z^\le}(T_j-1)(v)*F_{xz'z^<}(T_j)(v) \\
&=& F_{z'z^\le}(T_j)*F_{xz'z^\le}(T_j-1)(v)*F_{xz'{}^\le}(T_j)(v)*F_{xz^\le}(T_j)(v).
\end{eqnarray*}
Here we used
$F_{z'z^\le}(T_j-1)*F_{z'z^\le}(T_j)=F_{z'z^\le}(T_j-1)[F_{z'z^\le}(T_j)]=F_{z'z^\le}(T_j)$
and $F_{xz'z^<}(T_j)(v) = F_{xz'{}^\le}(T_j)(v) \cup F_{xz^\le}(T_j)(v)
= F_{xz'{}^\le}(T_j)(v)*F_{xz^\le}(T_j)(v)$.
To put the above identity in words: $F_{xz'z}(T_j)(v)$ is the bipartite subgraph of
$F_{z'z^\le}(T_j) \cap F_{xz'z}(T_j-1)(v)$ induced by $F_{xz'}(T_j)(v)$ and $F_{xz}(T_j)(v)$.
Also,
\begin{align*}
J & = F_{z_jz'z^\le}(T_j-1) * \bigodot_{S' \sub z_jz'z, S' \in H(x)} F_{S'x^\le}(T_j-1)(v), 
\mbox{ so by Lemma \ref{*props}(ii)} \\
J(y) & = F_{z'z^\le}(T_j-1) * \bigodot_{S' \sub z'z, S' \in H(x)} F_{S' x^\le}(T_j-1)(v) \\
&\quad * F_{z_jz'z^\le}(T_j-1)(y) * \bigodot_{S' \sub z_jz'z, S' \in H(x)} F_{S' x^\le}(T_j-1)(vy) \\
& = F_{z'z^\le}(T_j) * F_{xz'z^\le}(T_j-1)(v) * \bigodot_{S' \sub z_jz'z, S' \in H(x)} F_{S'x^\le}(T_j-1)(vy).
\end{align*}
Here we used $F_{z_jz'z^\le}(T_j-1)(y)= F_{z'z^\le}(T_j)$
and $F_{z'z^\le}(T_j-1)*F_{z'z^\le}(T_j)=F_{z'z^\le}(T_j-1)[F_{z'z^\le}(T_j)]=F_{z'z^\le}(T_j)$.
We also recall that $S=z'z \in H(x)$, so $S' \in H(x)$ for any $S' \sub z'z$.
Note that $\bigodot_{S' \sub z_jz'z, S' \in H(x)} F_{S'x^\le}(T_j-1)(vy)$ is a $1$-complex
containing $\{\es\}$, $F_{xz_jz}(T_j-1)(vy)$ if $xz_jz \in H$
and $F_{xz_jz'}(T_j-1)(vy)$ if $xz_jz' \in H$.
As in Case C.1, we have $F_{xz}(T_j)(v)=F_{xz_jz}(T_j-1)(vy)$ if $xz_jz \in H$
or $F_{xz}(T_j)(v)=F_{xz}(T_j-1)(v)\cap F_z(T_j)$ if $xz_jz \notin H$,
and similarly  $F_{xz'}(T_j)(v)=F_{xz_jz'}(T_j-1)(vy)$ if $xz_jz' \in H$
or $F_{xz'}(T_j)(v)=F_{xz'}(T_j-1)(v)\cap F_{z'}(T_j)$ if $xz_jz' \notin H$.
Thus $J(y)$ is also equal to
$F_{z'z^\le}(T_j)*F_{xz'z^\le}(T_j-1)(v)*F_{xz'{}^\le}(T_j)(v)*F_{xz^\le}(T_j)(v)$,
which proves that $F_{xz'z^\le}(T_j)(v)=J(y)$.
\COMMENT{Would back-to-basics proof be better?}

For any $S' \sub z_jz'z$,
$F_{S'}(T_j-1)$ is $\eps_{\nu'_{S'},1}$-regular by Lemma \ref{3exceptional} if $S' \in H$
and, using property (v) above,
$F_{S' x}(T_j-1)(v)$ is $\eps_{\nu''_{S'},2}$-regular by $A_{3,j-1}$ if $S' \in H(x)$.
By Lemma \ref{3restrict}, $J_{z_jz'z}$ is $\eps_{\nu''_{z_jz'z},2}$-regular with
\[d_{z_jz'z}(J) = (1 \pm \eps_{\nu''_{z_jz'z},2})d_{z_jz'z}(F(T_j-1)).\]
Similarly, by Lemma \ref{2restrict}, if $S' \sub z_jz'z$, $|S'|=2$, $S' \notin H(x)$
then $J_{S'}$ is  $\eps_{\nu''_{S'},2}$-regular
with $d_{S'}(J) = (1 \pm \eps_{\nu''_{S'},2})d_{S'}(F(T_j-1))$.
On the other hand, if $S' \sub z_jz'z$, $|S'|=2$, $S' \in H(x)$
then $J_{S'}=F_{S' x}(T_j-1)(v)$ is $\eps_{\nu''_S,2}$-regular
with $d_{S'}(J)=(1\pm\eps_{\nu''_S,2})d_{S' x}(F(T_j-1))d_{S'}(F(T_j-1))$ by $A_{3,j-1}$.
In particular,
\[d_{z'z}(J)=(1\pm\eps_{\nu''_{z'z},2})d_{xz'z}(F(T_j-1))d_{z'z}(F(T_j-1)).\]
Now by Lemma \ref{3neighbour1}, for all but at most $6\eps_{12D,3}|F_{xz_j}(T_j-1)(v)|$
vertices $y \in F_{xz_j}(T_j-1)(v)$, $F_{xz'z}(T_j)(v)=J(y)$ is $\eps_{\nu^*_{z'z},0}$-regular
and $d_{z'z}(F(T_j)(v)) = d_{z'z}(J(y)) = (1 \pm \eps_{\nu^*_{z'z},0})d_{z_jz'z}(J)d_{z'z}(J)$.
Also, $d_{z'z}(F(T_j))=(1 \pm \eps_{\nu^*_{z'z},0})d_{z_jz'z}(F(T_j-1))d_{z'z}(F(T_j-1))$
if $y \notin E_{z_j}(T_j-1)$ by $(*_{\ref{3alg}})$,
and $d_{xz'z}(F(T_j))=(1 \pm \eps_{\nu''_{xz'z},2})d_{xz'z}(F(T_j-1))$
by Lemmas \ref{3exceptional} and \ref{3restrict}. Thus
\begin{align*}
& d_{z'z}(F(T_j)(v)) = (1 \pm \eps_{\nu^*_{z'z},0})d_{z_jz'z}(J)d_{z'z}(J) \\
& =  (1 \pm \eps_{\nu^*_{z'z},0})(1 \pm \eps_{\nu''_{z_jz'z},2})(1\pm\eps_{\nu''_{z'z},2})
d_{z_jz'z}(F(T_j-1))d_{xz'z}(F(T_j-1))d_{z'z}(F(T_j-1)) \\
& = (1 \pm \eps_{\nu^*_{z'z},1})d_{z'z}(F(T_j))d_{xz'z}(F(T_j)),
\mbox{ i.e. } (*_{\ref{3-x-to-v}}) \mbox{ holds for } S=zz'.
\end{align*}
Combining the estimates for all cases we have at most $\eps_*|F_{xz_j}(T_j-1)(v)|$
exceptional vertices $y$, so this proves Claim C.

\nib{Claim D.} If $A_{1,j-1}$, $A_{2,j-1}$ and $A_{3,j-1}$ hold then
for any $E \in U(z_j)$, $Z \sub E$, $Z' = Z \cup z_j$ we have
$$|D^{Z'*v}_{z_j,E}(T_j-1)| < \theta_{\nu^*} |F_{xz_j}(T_j-1)(v)|.$$

\nib{Proof.} Note that $\nu^*>\nu$ by Lemma \ref{3observe}, since $E \in U(z_j)$.

\nib{Case D.1.} First consider the case $z_j \in E$. Then $E^{T_j}=\ov{E} \sm z_j$
and $F_{E^{T_j\le}}(T_j)=F_{\ov{E}^\le}(T_j-1)(y)$.
We will show that
\begin{equation}\label{eq:d1}\tag{$+_{\ref{3-x-to-v}}$}
F(T_j)^{Z*v}_{E^{T_j\le}}=F(T_j-1)^{Z'*v}_{\ov{E}^\le}(y).
\end{equation}
Before proving this in general we will illustrate a few cases of this statement.
Suppose that $x \in E$, say $E=xz_jz$ for some $z$.
If $z$ is embedded then $\ov{E}=xz_j$ and $E^{T_j}=x$,
so $F(T_j)^{Z*v}_x = F_x(T_j)=F_{xz_j}(T_j-1)(y)=F(T_j-1)^{Z'*v}_{xz_j}(y)$.%
\COMMENT{Formerly [seems incorrect!]
(In this case Claim D follows directly from Lemma \ref{3marked},
although our argument for the general case will also prove it.)
}
If $z$ is not embedded then $\ov{E}=E$ and $E^{T_j}=xz$,
$F(T_j)^{Z*v}_z$ is $F(T_j)_{xz}(v) = F(T_j-1)_E(vy)$ if $z \in I$
or $F(T_j)_z=F_{z_jz}(T_j-1)(y)$ otherwise,
and $F(T_j)^{Z*v}_{xz}$ is the bipartite subgraph of $F(T_j)_{xz}=F(T_j-1)_E(y)$
spanned by $F(T_j)^{Z*v}_x = F_x(T_j) = F_{xz_j}(T_j-1)(y)$ and $F(T_j)^{Z*v}_z$.
Also, we have $P \in F(T_j-1)^{Z'*v}_E(y)$ if $P y \in F(T_j-1)^{Z'*v}_E$,
i.e.\ $P y \in F(T_j-1)_E$ and $(P y)_S v \in F_{S x}(T_j-1)$
for all $S \sub Z'$, $S \in H(x)$. Equivalently, (i) $P \in F(T_j-1)_E(y) = F(T_j)_{xz}$,
(ii) for $S \sub Z\sm z_j$, $S \in H(x)$ we have $P_S v \in F_{S x}(T_j-1)$,
i.e.\ $P_S \in F_{S x}(T_j-1)(v)$, and (iii) if $Sz_j \in H(x)$
we have $P_S yv \in F_{S xz_j}(T_j-1)$,
i.e.\ $P_S \in F_{S xz_j}(T_j-1)(y)(v)=F_{S x}(T_j)(v)$.
Thus $P \in F(T_j)^{Z*v}_{xz}$.

On the other hand, suppose that $x \notin E$ and
consider the case that $E$ is unembedded, i.e.\ $\ov{E}=E=z_jz'z$ say.
Then $F(T_j)^{Z*v}_{z'z}$ is the bipartite subgraph of $F(T_j)'_{z'z}$
spanned by $F(T_j)'_z$ and $F(T_j)'_{z'}$, where we write
$F(T_j)'_{z'z}$ for $F(T_j)_{xz'z}(v)$ if $z'z \in I$
or $F(T_j)_{z'z}$ otherwise,
$F(T_j)'_z$ for $F(T_j)_{xz}(v)$ if $z \in I$ or $F(T_j)_z$ otherwise, and
$F(T_j)'_{z'}$ for $F(T_j)_{xz'}(v)$ if $z' \in I$ or $F(T_j)_{z'z}$ otherwise.
Recall that $F(T_j)_{xz}(v)$ is $F_{xz_jz}(T_j-1)(vy)$ if $xz_jz \in H$ (see Case C.1.i)
or $F_{xz}(T_j-1)(v) \cap F_{z_jz}(T_j-1)(y)$ if $xz_jz \notin H$ (see Case C.1.ii).
Similar statements hold for $F_{xz'}(T_j)(v)$.
Also, if $z'z \sub Z$ then $F(T_j-1)^{Z'*v}_{E^\le}$ is the complex $J$ defined in Case C.2,
so as shown there $F(T_j)_{xz'z^\le}(v) = F(T_j-1)^{Z'*v}_{E^\le}(y)$.

We deduce $(+_{\ref{3-x-to-v}})$ from the case $A=E^{T_j}$ of the following more general statement,
which will also be used in Cases D.2 and D.3:
\begin{equation}\label{eq:d2}\tag{$\dagger_{\ref{3-x-to-v}}$}
F(T_j)^{Z*v}_{A^\le}=F(T_j-1)^{Z'*v}_{Az_j^\le}(y) \mbox{ for } A \in H(z_j). 
\end{equation}
To see this, note that
$F_{Sx^\le}(T_j)=F_{Sx^\le}(T_j-1)[F_{Sxz_j^\le}(T_j-1)(y)]$ for any $z_j \notin S \in I$
by Definition \ref{def-update} (deleting $y$ has no effect), so by Lemma \ref{*props} we have
\begin{eqnarray*}
F_{Sx^\le}(T_j)(v) &=& (F_{Sx^\le}(T_j-1)*F_{Sxz_j^\le}(T_j-1)(y))(v) \\
&=& F_{S^\le}(T_j-1)*F_{Sz_j^\le}(T_j-1)(y) \\
& & \quad * F_{Sx^\le}(T_j-1)(v) * F_{Sxz_j^\le}(T_j-1)(yv)\\
&=& F_{Sx^\le}(T_j-1)(v) * F_{Sz_j^\le}(T_j-1)(y) * F_{Sxz_j^\le}(T_j-1)(yv),
\end{eqnarray*}
since $F_{Sz_j^\le}(T_j-1)(y)*F_{S^\le}(T_j-1)=F_{Sz_j^\le}(T_j-1)(y)$.
Now by definition and Lemma \ref{*props} we have
\begin{align*}
& F(T_j-1)^{Z'*v}_{Az_j^\le}
= F_{Az_j^\le}(T_j-1)[\cup_{S \in I'} F_{S x}(T_j-1)(v)]
= F_{Az_j^\le}(T_j-1) * \bigodot_{S \in I'} F_{Sx^\le}(T_j-1)(v), \mbox{ so} \\
& F(T_j-1)^{Z'*v}(y)_{A^\le} = F_{A^\le}(T_j-1)*\bigodot_{S \in I} F_{Sx^\le}(T_j-1)(v) \\
& \hspace{0.3\textwidth} *F_{Az_j^\le}(T_j-1)(y)*\bigodot_{S \in I'} F_{Sx^\le}(T_j-1)(vy) \\
& = F_{Az_j^\le}(T_j-1)(y)*\bigodot_{S \in I} F_{Sx^\le}(T_j-1)(v)
*\bigodot_{S \in I'} F_{Sx^\le}(T_j-1)(vy).
\end{align*}
Here we used $F_{A^\le}(T_j-1)*F_{Az_j^\le}(T_j-1)(y)=F_{Az_j^\le}(T_j-1)(y)$.
On the other hand,
\begin{align*}
& F^{Z*v}_{A^\le}(T_j) = F_{A^\le}(T_j)[\cup_{S \in I}F_{Sx}(T_j)(v)]
= F_{A^\le}(T_j)*\bigodot_{S \in I}F_{Sx^\le}(T_j)(v) \\
& = F_{A^\le}(T_j)* \bigodot_{S \in I}\left(
F_{Sx^\le}(T_j-1)(v) * F_{Sz_j^\le}(T_j-1)(y) * F_{Sxz_j^\le}(T_j-1)(yv) \right)\\
& = F_{Az_j^\le}(T_j-1)(y)*\bigodot_{S \in I}
F_{Sx^\le}(T_j-1)(v)*\bigodot_{S' \in I'} F_{S'x^\le}(T_j-1)(vy)
= F(T_j-1)^{Z'*v}(y)_{A^\le}.
\end{align*}
In the second equality above we substituted for % the previous expression
$F_{Sx^\le}(T_j)(v)$, and in the third we set $S'=Sz_j$ and used
$F_{Sz_j^\le}(T_j-1)(y)*F_{Az_j^\le}(T_j-1)(y)=F_{Az_j^\le}(T_j-1)(y)$.
This proves $(\dagger_{\ref{3-x-to-v}})$, and so $(+_{\ref{3-x-to-v}})$.%

Now $M_{E^{T_j},E}(T_j)=M_{\ov{E},E}(T_j-1)(y)$ by Lemma \ref{track-mark},
so we have $M_{E^{T_j},E}(T_j) \cap F(T_j)^{Z*v}_{E^{T_j}}
= M_{\ov{E},E}(T_j-1)(y) \cap F(T_j-1)^{Z'*v}_{\ov{E}^\le}(y)$.
Recalling that $F(T_j-1)^{Z'*v}_{z_j}=F_{xz_j}(T_j-1)(v)$, we have
\[ D^{Z'*v}_{z_j,E}(T_j-1) = \left\{y \in F_{xz_j}(T_j-1)(v):
\frac{\left|\left(M_{\ov{E},E}(T_j-1) \cap F(T_j-1)^{Z'*v}_{\ov{E}^\le}\right)(y)\right|}{
\left|F(T_j-1)^{Z'*v}_{\ov{E}^\le}(y)\right|} > \theta_{\nu^*} \right\}. \]
Also, $F(T_j-1)^{Z'*v}_{\ov{E}^\le}$ is $\eps_{12D,2}$-regular by $A_{3,j-1}$,
so writing $B'_{z_j}$ for the set of vertices
$y \in F(T_j-1)^{Z'*v}_{z_j}$ for which we do not have
$|F(T_j-1)^{Z'*v}_{\ov{E}}(y)|=(1\pm\eps_*)|F(T_j-1)^{Z'*v}_{\ov{E}}|/|F(T_j-1)^{Z'*v}_{z_j}|$,
we have $|B'_{z_j}| < \eps_*|F(T_j-1)^{Z'*v}_{z_j}|$ by Lemma \ref{3average}.%
\COMMENT{
Formerly: Note that $|B_{z_j}| < 2\eps_*|F_{z_j}(T_j-1)| < \sqrt{\eps_*}|F_{xz_j}(T_j-1)(v)|$
by $A_{3,j-1}$, since all densities are at least $d_u \gg \eps_*$.
Earlier: (One can even show that $B'_{z_j} \sub B_{z_j}$, although we do not need this fact.)
}
Now
\COMMENT{
Prefer simpler averaging argument, even though $B'$ should be unnecessary.
Let's check direct calculation... [how I found the previous mistake!]
By $A_{3,j-1}$ and Lemma \ref{3absolute} we have
$d(F(T_j-1)^{Z'*v}_{\ov{E}}) = (1\pm\eps_*)\prod_{S \sub \ov{E}}d_S(F(T_j-1))
\prod_{S\in I'}d_{Sx}(F(T_j-1))$.
If $y \notin E^v_{z_j}(T_j-1)$ we have
$d(F(T_j)^{Z*v}_{E^{T_j}}) = (1\pm\eps_*)\prod_{S \sub E^{T_j}}d_S(F(T_j))
\prod_{S\in I}d_{Sx}(F(T_j))$.
We have $d_S(F(T_j))=(1\pm\eps_*)d_S(F(T_j-1))d_{Sz_j}(F(T_j-1))$ if $\es \ne S \in H(z_j)$
or $d_S(F(T_j))=(1\pm\eps_*)d_S(F(T_j-1))$ if $S \notin H(z_j)$ by Lemma \ref{3crude}.
So $d(F(T_j-1)^{Z'*v}_{\ov{E}})/d(F(T_j)^{Z*v}_{E^{T_j}}) = (1\pm 30\eps_*)
\prod_{z_j \in S \sub \ov{E}} d_S(F(T_j-1))\prod_{z_j \in S\in I'}d_{Sx}(F(T_j-1))
\prod_{\es \ne S' \in H(z_j)} d_{S'z_j}(F(T_j-1))^{-1}
\prod_{S' \sub Z,S' \in H(x), S'x \in H(z_j)} d_{S'xz_j}(F(T_j-1))^{-1}$.
The second and fourth products cancel completely and the
first and third cancel all but $d_z(F(T_j-1))$.
}
\begin{align*}
\Sa & := \sum_{y \in D^{Z'*v}_{z_j,E}(T_j-1)}
\left|\left(M_{\ov{E},E}(T_j-1) \cap F(T_j-1)^{Z'*v}_{\ov{E}^\le}\right)(y)\right| \\
& > \theta_{\nu^*} \sum_{y \in D^{Z'*v}_{z_j,E}(T_j-1) \sm B'_{z_j}}
\left|F(T_j-1)^{Z'*v}_{\ov{E}^\le}(y)\right| \\
& >  (1-\eps_*)\theta_{\nu^*} (|D^{Z'*v}_{z_j,E}(T_j-1)|-\eps_*|F_{xz_j}(T_j-1)(v)|)\
|F(T_j-1)^{Z'*v}_{\ov{E}}|/|F_{xz_j}(T_j-1)(v)|.
\end{align*}
We also have an upper bound
\begin{align*}
& \Sa \le \sum_{y \in F_{xz_j}(T_j-1)(v)}
\left|\left(M_{\ov{E},E}(T_j-1) \cap F(T_j-1)^{Z'*v}_{\ov{E}^\le}\right)(y)\right| \\
& = |M_{\ov{E},E}(T_j-1) \cap F(T_j-1)^{Z'*v}_{\ov{E}^\le}|
< \theta'_{\nu} |F(T_j-1)^{Z'*v}_{\ov{E}^\le}|,
\end{align*}
where the last inequality holds by $A_{2,j-1}$.
Therefore \[\frac{|D^{Z'*v}_{z_j,E}(T_j-1)|}{|F_{xz_j}(T_j-1)(v)|}
< \frac{\theta'_{\nu}}{(1-\eps_*)\theta_{\nu^*}} + \eps_* < \theta_{\nu^*}.\]

\nib{Case D.2.} Next consider the case $z_j \notin E$ and $\ov{E} \in H(z_j)$.
Then $E^{T_j}=E^{T_j-1}=\ov{E}$ and $F_{E^{T_j\le}}(T_j)=F_{\ov{E}z_j^\le}(T_j-1)(y)$.
Also, by ($\dagger_{\ref{3-x-to-v}}$) we have
$F(T_j)^{Z*v}_{\ov{E}^\le}=F(T_j-1)^{Z'*v}_{\ov{E}z_j^\le}(y)$.
Now $M_{\ov{E},E}(T_j)=M_{\ov{E},E}(T_j-1) \cap F_{\ov{E}}(T_j)$ by Lemma \ref{track-mark}
and $F(T_j)^{Z*v}_{\ov{E}} \sub F_{\ov{E}}(T_j)$,
so $M_{\ov{E},E}(T_j) \cap F(T_j)^{Z*v}_{\ov{E}}
= M_{\ov{E},E}(T_j-1) \cap F(T_j-1)^{Z'*v}_{\ov{E}z_j}(y)$.
Since $F(T_j-1)^{Z'*v}_{z_j}=F_{xz_j}(T_j-1)(v)$, we have
\[ D^{Z'*v}_{z_j,E}(T_j-1) = \left\{y \in F_{xz_j}(T_j-1)(v):
\frac{\left|M_{\ov{E},E}(T_j-1) \cap F(T_j-1)^{Z'*v}_{\ov{E}z_j}(y)\right|}{
\left|F(T_j-1)^{Z'*v}_{\ov{E}z_j}(y)\right|} > \theta_{\nu^*}\right\}. \]
Also, since $\ov{E} \in H(z_j)$ we have $\ov{E}z_j \sub E_0^{T_j-1}$
for some triple $E_0 \in H'$, so applying $A_{3,j-1}$ to $E_0$
we see that $F(T_j-1)^{Z'*v}_{\ov{E}z_j^\le}$ is $\eps_{12D,2}$-regular.
Then writing $B'_{z_j}$ for the set of vertices 
$y \in F(T_j-1)^{Z'*v}_{z_j}$ for which we do not have
$|F(T_j-1)^{Z'*v}_{\ov{E}z_j}(y)|=(1\pm\eps_*)|F(T_j-1)^{Z'*v}_{\ov{E}z_j}|/|F(T_j-1)^{Z'*v}_{z_j}|$,
we have $|B'_{z_j}| < \eps_*|F(T_j-1)^{Z'*v}_{z_j}|$ by Lemma \ref{3average}.
Now
\begin{align*}
& \Sa := \sum_{y \in D^{Z'*v}_{z_j,E}(T_j-1)}
\left|M_{\ov{E},E}(T_j-1) \cap F(T_j-1)^{Z'*v}_{\ov{E}z_j}(y)\right|
 > \theta_{\nu^*} \sum_{y \in D^{Z'*v}_{z_j,E}(T_j-1) \sm B'_{z_j}}
\left|F(T_j-1)^{Z'*v}_{\ov{E}z_j}(y)\right| \\
& \qquad >  (1-\eps_*)\theta_{\nu^*} (|D^{Z'*v}_{z_j,E}(T_j-1)|-\eps_*|F_{xz_j}(T_j-1)(v)|)\
|F(T_j-1)^{Z'*v}_{\ov{E}z_j}|/|F_{xz_j}(T_j-1)(v)|.
\end{align*}
We also have
$\Sa \le \sum_{y \in F_{xz_j}(T_j-1)(v)}
\left|M_{\ov{E},E}(T_j-1) \cap F(T_j-1)^{Z'*v}_{\ov{E}z_j}(y)\right|$.
This sum counts all pairs $(y,P)$ with $P \in M_{\ov{E},E}(T_j-1) \cap F(T_j-1)^{Z'*v}_{\ov{E}}$,
$y \in F_{xz_j}(T_j-1)(v)$ and $Py \in F(T_j-1)^{Z'*v}_{\ov{E}z_j}$,
so $\Sa \le \sum_{P \in M_{\ov{E},E}(T_j-1) \cap F(T_j-1)^{Z'*v}_{\ov{E}}} |F(T_j-1)^{Z'*v}_{\ov{E}z_j}(P)|$.
By Lemma \ref{3crude}(vi) we have
\[|F(T_j-1)^{Z'*v}_{\ov{E}z_j}(P)|=(1\pm\eps_*)
\frac{|F(T_j-1)^{Z'*v}_{\ov{E}z_j}|}{|F(T_j-1)^{Z'*v}_{\ov{E}}|}\]
for all but at most $\eps_*|F(T_j-1)^{Z'*v}_{\ov{E}}|$
sets $P \in F(T_j-1)^{Z'*v}_{\ov{E}}$. Therefore
\[\Sa \le |M_{\ov{E},E}(T_j-1) \cap F(T_j-1)^{Z'*v}_{\ov{E}}|(1+\eps_*)
\frac{|F(T_j-1)^{Z'*v}_{\ov{E}z_j}|}{|F(T_j-1)^{Z'*v}_{\ov{E}}|}
+ \eps_*|F(T_j-1)^{Z'*v}_{\ov{E}}||F_{xz_j}(T_j-1)(v)|.\]
Combining this with the lower bound on $\Sa$ gives
\begin{align*}
& (1-\eps_*)\theta_{\nu^*} (|D^{Z'*v}_{z_j,E}(T_j-1)|/|F_{xz_j}(T_j-1)(v)|-\eps_*) \\
& < (1+\eps_*)\frac{|M_{\ov{E},E}(T_j-1) \cap F(T_j-1)^{Z'*v}_{\ov{E}}|}{|F(T_j-1)^{Z'*v}_{\ov{E}}|}
\ +\ \eps_*\frac{|F(T_j-1)^{Z'*v}_{\ov{E}}||F_{xz_j}(T_j-1)(v)|}{|F(T_j-1)^{Z'*v}_{\ov{E}z_j}|} .
\end{align*}
Now $F(T_j-1)^{Z'*v}_{\ov{E}}=F(T_j-1)^{Z*v}_{\ov{E}}$,
$|M_{\ov{E},E}(T_j-1) \cap F(T_j-1)^{Z*v}_{\ov{E}}| < \theta'_{\nu} |F(T_j-1)^{Z*v}_{\ov{E}}|$
by $A_{2,j-1}$, and
$\frac{|F(T_j-1)^{Z'*v}_{\ov{E}}||F_{xz_j}(T_j-1)(v)|}{|F(T_j-1)^{Z'*v}_{\ov{E}z_j}|}
\le 2d_u^{-4} \ll \eps_*^{-1}$ by Lemma \ref{3crude} and $A_{3,j-1}$, so
\[\frac{|D^{Z'*v}_{z_j,E}(T_j-1)|}{|F_{xz_j}(T_j-1)(v)|} < \
\frac{(1+\eps_*)\theta'_{\nu}+\sqrt{\eps_*}}{(1-\eps_*)\theta_{\nu^*}}
+ \eps_* < \theta_{\nu^*}.\]

\nib{Case D.3.} It remains to consider the case when $z_j \notin E$ and $\ov{E} \notin H(z_j)$.
Since $E \in U(z_j)$ and $\ov{E} \notin H(z_j)$ we have $|\ov{E}|\ge 2$.
Then $F_{\ov{E}^\le}(T_j)=F_{\ov{E}^\le}(T_j-1)[F_{\ov{E}^<}(T_j)]$ by Lemma \ref{build-update}.
Also, we claim that\COMMENT{stick to $Z$, although doesn't matter}
\begin{equation}\label{eq:d3}\tag{$\diamondsuit_{\ref{3-x-to-v}}$}
F(T_j)^{Z*v}_{\ov{E}^\le}=F(T_j-1)^{Z*v}_{\ov{E}^\le}[F(T_j)^{Z*v}_{\ov{E}^<}].
\end{equation}
To see this, consider $S \sub \ov{E}$ and $P \in F(T_j)^{Z*v}_S$.
Then by definition $P \in F_S(T_j)$ and
$P_{S'} \in F_{S'x}(T_j)(v)$ for all $S' \sub S \cap Z$ with $S' \in H(x)$.
Since $F_S(T_j) \sub F_S(T_j-1)$ it follows that $P \in F(T_j-1)^{Z*v}_S$.
Also, for any $S'' \sub S' \sub S$ with $S' \subn E$, $S''\sub Z$ and $S'' \in H(x)$
we have $P_{S'} \in F_{S'}(T_j)$ and $P_{S''} \in F_{S''x}(T_j)(v)$,
so $P_{S'} \in F(T_j)^{Z*v}_{S'}$. Thus
$P \in F(T_j-1)^{Z*v}_{\ov{E}^\le}[F(T_j)^{Z*v}_{\ov{E}^<}]$.
Conversely, suppose that $S \sub \ov{E}$ and
$P \in F(T_j-1)^{Z*v}_{\ov{E}^\le}[F(T_j)^{Z*v}_{\ov{E}^<}]_S$.
If $S \ne \ov{E}$ then $P \in F(T_j)^{Z*v}_{\ov{E}^<} \sub F(T_j)^{Z*v}_{\ov{E}^\le}$.
Now suppose $S=\ov{E}$. Then $P \in F_{\ov{E}}(T_j-1)$
and $P_{S'} \in F_{S'}(T_j)$ for $S' \subn S$, so $P \in F_{\ov{E}}(T_j)$.
Also, for $S' \subn S$ we have $P_{S'} \in F(T_j)^{Z*v}_{S'}$,
so $P_{S''} \in F_{S''x}(T_j)(v)$ for all $S'' \sub S' \cap Z$ with $S'' \in H(x)$.
Therefore $P \in F(T_j)^{Z*v}_S$. This proves $(\diamondsuit_{\ref{3-x-to-v}})$.
Note also that since $z\notin\ov{E}$ we can replace $Z$ by $Z'$ 
in $(\diamondsuit_{\ref{3-x-to-v}})$. %without altering either side.

Since $M_{\ov{E},E}(T_j)=M_{\ov{E},E}(T_j-1) \cap F_{\ov{E}}(T_j)$
and $F(T_j)^{Z'*v}_{\ov{E}} \sub F_{\ov{E}}(T_j)$
we have
$M_{\ov{E},E}(T_j) \cap F(T_j)^{Z'*v}_{\ov{E}}
= M_{\ov{E},E}(T_j-1) \cap F(T_j-1)^{Z'*v}_{\ov{E}^\le}[F(T_j)^{Z'*v}_{\ov{E}^<}]$.
Then\COMMENT{
Decided to keep $Z'$. Formerly:
Here we have replaced $Z'$ by $Z$, since $z_j \notin \ov{E}$.
}
\[ D^{Z'*v}_{z_j,E}(T_j-1) = \left\{y \in F_{xz_j}(T_j-1)(v):
\frac{\left|M_{\ov{E},E}(T_j-1) \cap F(T_j-1)^{Z'*v}_{\ov{E}^\le}[F(T_j)^{Z'*v}_{\ov{E}^<}]_{\ov{E}}\right|}{
\left| F(T_j-1)^{Z'*v}_{\ov{E}^\le}[F(T_j)^{Z'*v}_{\ov{E}^<}]_{\ov{E}}\right|}
> \theta_{\nu^*} \right\}. \]
Next note that by $A_{3,j-1}$ and Lemma \ref{3absolute} we have
\[d(F(T_j-1)^{Z'*v}_{\ov{E}}) = (1\pm 8\eps_*)
\prod_{S\sub\ov{E}}d_S(F(T_j-1))\prod_{S\in I}d_{Sx}(F(T_j-1))\prod_{S\in I,|S|=1}\alpha_S.\]
Now consider $y \notin E^v_{z_j}(T_j-1)$. By Claim B we can apply
the properties in $A_{3,j}$ with this choice of $y$, so we also have
\[d(F(T_j)^{Z'*v}_{\ov{E}}) = (1\pm 8\eps_*)
\prod_{S\sub\ov{E}}d_S(F(T_j))\prod_{S\in I}d_{Sx}(F(T_j))\prod_{S\in I,|S|=1}\alpha_S.\]
Now $d_S(F(T_j))$ is $(1\pm \eps_*)d_S(F(T_j-1))d_{Sz_j}(F(T_j-1)$ for $S \in H(z_j)$
by $(*_{\ref{3alg}})$ or $(1\pm \eps_*)d_S(F(T_j-1))$ for $S \notin H(z_j)$ by Lemma \ref{3crude}.
Therefore $d(F(T_j)^{Z'*v}_{\ov{E}}) = (1\pm 40\eps_*)d(F(T_j-1)^{Z'*v}_{\ov{E}}) \times d^*$, where
\begin{align*}
d^* & = \prod_{S\sub\ov{E},S\in H(z_j)}d_{Sz_j}(F(T_j-1))
\prod_{S\in I,Sx\in H(z_j)} d_{Sxz_j}(F(T_j-1)) \\
& = (1\pm 8\eps_*)\prod_{S\sub\ov{E},S\in H(z_j)} d_{Sz_j}(F(T_j-1)^{Z'*v}).
\end{align*}
To see the second equality above, note that for any $S\in I$ with $Sx\in H(z_j)$,
applying $A_{3,j-1}$ to any triple $E_0 \in H$ with $Sz_j \sub E_0^{T_j-1}$,
since $Sz_j \sub Z' = Zz_j$ we have
$d_{Sz_j}(F(T_j-1)^{Z'*v}) = (1\pm \eps_*)d_{Sxz_j}(F(T_j-1))d_{Sz_j}(F(T_j-1))$.
We deduce that
\[\left|F(T_j-1)^{Z'*v}_{\ov{E}^\le}[F(T_j)^{Z'*v}_{\ov{E}^<}]_{\ov{E}}\right|
= (1\pm 50\eps_*)\left|F(T_j-1)^{Z'*v}_{\ov{E}}\right| d^*\]
for such $y \notin E^v_{z_j}(T_j-1)$. Now%
\COMMENT{Here we do mean $B_{z_j}$ rather than $B'_{z_j}$.}
\begin{align*}
\Sa &:= \sum_{y \in D^{Z'*v}_{z_j,E}(T_j-1)}
\left|M_{\ov{E},E}(T_j-1) \cap F(T_j-1)^{Z'*v}_{\ov{E}^\le}[F(T_j)^{Z'*v}_{\ov{E}^<}]_{\ov{E}} \right| \\
& > \theta_{\nu^*} \sum_{y \in D^{Z'*v}_{z_j,E}(T_j-1)\sm B_{z_j}}
\left|  F(T_j-1)^{Z'*v}_{\ov{E}^\le}[F(T_j)^{Z'*v}_{\ov{E}^<}]_{\ov{E}} \right| \\
& >  (1-50\eps_*)\theta_{\nu^*} (|D^{Z'*v}_{z_j,E}(T_j-1)|-\sqrt{\eps_*}|F_{xz_j}(T_j-1)(v)|)\
\left|F(T_j-1)^{Z'*v}_{\ov{E}}\right| d^*.
\end{align*}
(Note that we used the estimate $|B_{z_j}| < 2\eps|F_{z_j}(T_j-1)| < \sqrt{\eps_*}|F_{xz_j}(T_j-1)(v)|$,
by Lemma \ref{3exceptional} and Claim C.)
We also have
$\Sa \le \sum_{y \in F_{xz_j}(T_j-1)(v)}
\left|M_{\ov{E},E}(T_j-1) \cap F(T_j-1)^{Z'*v}_{\ov{E}^\le}[F(T_j)^{Z'*v}_{\ov{E}^<}]_{\ov{E}} \right|$.
This sum counts all pairs $(y,P)$ with 
$P \in M_{\ov{E},E}(T_j-1) \cap F(T_j-1)^{Z'*v}_{\ov{E}}$,
$y \in F_{xz_j}(T_j-1)(v)$ and $P_S \in F(T_j)^{Z'*v}_S$
for all $S \subn \ov{E}$, $S \in H(z_j)$; there is no need
to consider $S \notin H(z_j)$, as by Lemma \ref{build-update}
we get the same expression if we replace the restriction to
$F(T_j)^{Z'*v}_S$ by $F(T_j-1)^{Z'*v}_S$ for such $S$.
Note that $F(T_j)^{Z'*v}_S=F(T_j)^{Z*v}_S$ and
$P_S \in F(T_j)^{Z*v}_S$ $\Lra$ $P_Sy \in F(T_j-1)^{Z'*v}_{Sz_j}$ by ($\dagger_{\ref{3-x-to-v}}$).
Given $P$, let $F^v_{P,Z}$ be the set of $y \in F_{xz_j}(T_j-1)(v)$ satisfying this condition,
and let $B^v_Z$ be the set of $P \in F(T_j-1)^{Z'*v}_{\ov{E}}$ such that we do not have%
\COMMENT{Could eg write $B^v_{Z'}$... a bit fussy...}
\[|F^v_{P,Z}| = (1\pm \eps_*)|F_{xz_j}(T_j-1)(v)|d^* = (1\pm \eps_*)|F(T_j-1)^{Z'*v}_{z_j}|d^*.\]
Lemma \ref{3technical} applied with $G=F(T_j-1)^{Z'*v}_{\ov{E}^\le}$ 
and $I=\{i(S): S \subn \ov{E}, S \in H(z_j)\}$
gives $|B^v_Z| < \eps_*|F(T_j-1)^{Z'*v}_{\ov{E}}|$.
Then
\[\Sa \le |M_{\ov{E},E}(T_j-1) \cap F(T_j-1)^{Z'*v}_{\ov{E}}|(1+\eps_*)|F_{xz_j}(T_j-1)(v)|d^*
+ \eps_*|F(T_j-1)^{Z'*v}_{\ov{E}}||F_{xz_j}(T_j-1)(v)|.\]
Combining this with the lower bound on $\Sa$ we obtain
\[(1-50\eps_*)\theta_{\nu^*} \left(\frac{|D^{Z'*v}_{z_j,E}(T_j-1)|}{|F_{xz_j}(T_j-1)(v)|}-\sqrt{\eps_*}\right)
< (1+\eps_*)\frac{|M_{\ov{E},E}(T_j-1) \cap F(T_j-1)^{Z'*v}_{\ov{E}}|}{\left|F(T_j-1)^{Z'*v}_{\ov{E}}\right| }
+ \frac{\eps_*}{d^*}.\]
Now $d^* \gg \eps_*$, and $|M_{\ov{E},E}(T_j-1) \cap F(T_j-1)^{Z'*v}_{\ov{E}}|<
\theta'_{\nu}|F(T_j-1)^{Z'*v}_{\ov{E}}|$ by $A_{2,j-1}$
(since $F(T_j-1)^{Z'*v}_{\ov{E}} = F(T_j-1)^{Z*v}_{\ov{E}}$).
Then
\[\frac{|D^{Z'*v}_{z_j,E}(T_j-1)|}{|F_{xz_j}(T_j-1)(v)|} < \
\frac{(1+ \eps_*)\theta'_{\nu}+\sqrt{\eps_*}}{(1-50\eps_*)\theta_{\nu^*}}
+ \sqrt{\eps_*} < \theta_{\nu^*}.\]
This completes the proof of Claim D.

\nib{Claim E.} Conditional on the events $A_{i,j'}$, $1 \le i \le 4$, $0 \le j' < j$
and the embedding up to time $T_j-1$ we have $\mb{P}(A_{4,j}) > d_u/2$.

\nib{Proof.} Suppose $A_{1,j-1}$, $A_{2,j-1}$ and $A_{3,j-1}$ hold.
Then Claim D gives $|D^{Z'*v}_{z_j,E}(T_j-1)| < \theta_{12D} |F_{xz_j}(T_j-1)(v)|$
for any $E \in U(z_j)$, $Z \sub E$, $Z' = Z \cup z_j$.
Also $F_{xz_j}(T_j-1)(v) > (1-\delta_B^{1/4})d_u|F_{z_j}(T_j-1)|$ 
by $A_{3,j-1}$ and since $\alpha_{z_j} > 1-2\delta_B^{1/3}$.
As $|E_{z_j}(T_j-1)| < \eps_* |F_{z_j}(T_j-1)|$ by Lemma \ref{3exceptional},
$B_{z_j} = E_{z_j}(T_j-1) \cup E^v_{z_j}(T_j-1)$ has size
$|B_{z_j}| < \sqrt{\eps_*} |F_{xz_j}(T_j-1)(v)|$ by Claim C.
Since $H$ has maximum degree at most $D$ we have at most $2D^2$ choices for $E \in U(z_j)$
then $8$ choices for $Z \sub E$, so
$|OK^v_{z_j}(T_j-1)|/|F_{xz_j}(T_j-1)(v)|> 1-\sqrt{\eps_*}-16D^2\theta_{12D} > 1-\theta_*$.
Now $y=\phi(z_j)$ is chosen uniformly at random in $OK_{z_j}(T_j-1) \sub F_{z_j}(T_j-1)$,
and $OK^v_{z_j}(T_j-1) \sub OK_{z_j}(T_j-1)$. Then Claim E follows from
\[\mb{P}(y \in OK^v_{z_j}(T_j-1))  = \frac{OK^v_{z_j}(T_j-1)}{OK_{z_j}(T_j-1)}
> \frac{(1-\theta_*)|F_{xz_j}(T_j-1)(v)|}{|F_{z_j}(T_j-1)|} 
> (1-\theta_*)(1-\delta_B^{1/4})d_u > d_u/2.\]

To finish the proof of the lemma, note that
if all the events $A_{i,j}$, $1 \le i \le 4$, $1 \le j \le g$ hold,
then the events $A_{4,j}$ imply that $\phi(H(x)) \sub (G \sm M)(v)$.
Multiplying the conditional probabilities given by Claim E
over at most $2D$ vertices of $VN_H(x)$ gives
probability at least $(d_u/2)^{2D} > p$. \qed

\subsection{The conclusion of the algorithm} \index{conclusion}

The algorithm will be successful if the following two conditions hold.
Firstly, it must not abort during the iterative phase because of the
queue becoming too large. Secondly, there must be a system of distinct
representatives for the available images of the unembedded vertices,
which all belong to the buffer $B$.
Recall that $A_x = F_x(t_x^N) \sm \cup_{E \ni x} M_{x,E}(t_x^N)$
is the available set for $x \in B$ at the time $t_x^N$
when the last vertex of $VN_H(x)$ is embedded.
Since $VN_H(x)$ has been embedded, until the conclusion of the algorithm at time $T$,
no further vertices will be marked as forbidden for the image $x$,
and no further neighbourhood conditions will be imposed on the image of $x$,
although some vertices in $A_x$ may be used to embed other vertices in $X_x$.
Thus the set of vertices available to embed $x$ at time $T$
is $A'_x = A_x \cap V_x(T) = F_x(T) \sm \cup_{E \ni x} M_{x,E}(T)$.
Therefore we seek a system of distinct representatives for
$\{A'_x: x \in X(T)\}$.

We start with the `main lemma', which is almost identical to Lemma \ref{2main}.
For completeness we repeat the proof, giving the necessary modifications.

\begin{lemma}\label{3main} \index{main lemma}
Suppose $1 \le i \le r$, $Y \sub X_i$ and $A \sub V_i$ with $|A|>\eps_* n$.
Let $E_{A,Y}$ be the event that (i) no vertices are embedded in $A$ before the conclusion of
the algorithm, and (ii) for every $z \in Y$ there is some time $t_z$ such that
$|A \cap F_z(t_z)|/|F_z(t_z)| < 2^{-2D} |A|/|V_i|$. Then $\mb{P}(E_{A,Y}) < p_0^{|Y|}$.
\end{lemma}

\nib{Proof.}
We start by choosing $Y' \sub Y$ with $|Y'| > |Y|/(2D)^2$ so that vertices
in $Y'$ are mutually at distance at least $3$ (this can be done greedily,
using the fact that $H$ has maximum degree $D$).
It suffices to bound the probability of $E_{A,Y'}$.
Note that initially we have $|A \cap F_z(0)|/|F_z(0)| = |A|/|V_i|$ for all $z\in X_i$.
Also, if no vertices are embedded in $A$, then $|A \cap F_z(t)|/|F_z(t)|$
can only be less than $|A \cap F_z(t-1)|/|F_z(t-1)|$ for some $z$ and $t$
if we embed a neighbour of $z$ at time $t$.
It follows that if $E_{A,Y'}$ occurs, then for every $z \in Y'$
there is a first time $t_z$ when we embed a neighbour $w$ of $z$ and
have $|A \cap F_z(t_z)|/|F_z(t_z)| < |A \cap F_z(t_z-1)|/2|F_z(t_z-1)|$.

By Lemma \ref{3exceptional}, the densities
$d_z(F(t_z-1))$, $d_z(F(t_z-1))$ and $d_{zw}(F(t_z-1))$ are all at least $d_u$
and $F_{zw}(t_z-1)$ is  $\eps_*$-regular. Applying Lemma \ref{2neighbour},
we see that there are at most $\eps_*|F_w(t_z-1)|$ `exceptional' vertices
$y \in F_w(t_z-1)$ that do not satisfy
$|A \cap F_z(t_z)| = |F_{zw}(t_z-1)(y) \cap A \cap F_z(t_z-1)|
= (1 \pm \eps_*)d_{zw}(F(t_z-1))|A \cap F_z(t_z-1)|$.
On the other hand, the algorithm chooses $\phi(w)=y$ to satisfy $(*_{\ref{3alg}})$,
so $|F_z(t_z)| = (1 \pm \eps_*)d_{zw}(F(t_z-1))|F_z(t_z-1)|$. Thus we can only have
$|A \cap F_z(t_z)|/|F_z(t_z)| < |A \cap F_z(t_z-1)|/2|F_z(t_z-1)|$
by choosing an exceptional vertex $y$.  But $y$ is chosen uniformly at
random from $|OK_w(t_z-1)| \ge (1-\theta_*)|F_w(t_z-1)|$ possibilities
(by Corollary \ref{3ok}). It follows that, conditional on the prior embedding,
the probability of choosing an exceptional vertex for $y$ is at most
$\eps_*|F_w(t_z-1)|/|OK_w(t_z-1)| < 2\eps_*$.

Since vertices of $Y'$ have disjoint neighbourhoods,
we can multiply the conditional probabilities over $z \in Y'$
to obtain an upper bound of $(2\eps_*)^{|Y'|}$.
Recall that this bound is for a subset of $E_{A,Y'}$
in which we have specified a certain neighbour $w$ for every vertex $z \in Y'$.
Taking a union bound over at most $(2D)^{|Y'|}$ choices for these neighbours
gives $\mb{P}(E_{A,Y}) \le \mb{P}(E_{A,Y'}) < (4\eps_*D)^{|Y'|} < p_0^{|Y|}$. \qed

Now we can prove the following theorem, which implies Theorem \ref{3blowup}.
The proof is quite similar to the graph case, except that the marked edges
create an additional case when verifying Hall's criterion,
which is covered by Lemma \ref{3initial}.

\begin{theo}\label{3final}
With high probability the algorithm embeds $H$ in $G\sm M$.
\end{theo}

\nib{Proof.}
First we estimate the probability of the iteration phase aborting with failure,
which happens when the number of vertices that have ever been queued is too large.
We can take a union bound over all $1 \le i \le r$ and $Y \sub X_i$ with $|Y|=\delta_Q|X_i|$
of $\mb{P}(Y \sub Q(T))$. Suppose that the event $Y \sub Q(T)$ occurs. Then for every $z \in Y$
there is some time $t$ such that $|F_z(t)| < \delta'_Q|F_z(t_z)|$,
where $t_z<t$ is the most recent time at which we embedded a neighbour of $z$.
Since $A=V_i(T)$ is unused we have $A \cap F_z(t) = A \cap F_z(t_z)$, so
$|A \cap F_z(t_z)|/|F_z(t_z)| = |A \cap F_z(t)|/|F_z(t_z)| \le |F_z(t)|/|F_z(t_z)| < \delta'_Q$.
However, we have $|A| \ge \delta_B n/2$ by Lemma \ref{3observe}(i),
so since $\delta'_Q \ll \delta_B$ we have $|A \cap F_z(t_z)|/|F_z(t_z)| < 2^{-2D} |A|/|V_i|$.
Taking a union bound over all possibilities for $i$, $Y$ and $A$,
Lemma \ref{3main} implies that the failure probability is at most
$r \cdot 4^{Cn} \cdot p_0^{\delta_Q n} < o(1)$, since $p_0 \ll \delta_Q$.

Now we estimate the probability of the conclusion of the algorithm aborting with failure.
By Hall's criterion for finding a system of distinct representatives,
the conclusion fails if and only if there is $1 \le i \le r$
and $S \sub X_i(T)$ such that $|\cup_{z \in S} A'_z| < |S|$.
Recall that $|X_i(T)| \ge \delta_B n/2$ by Lemma \ref{3observe}(i)
and buffer vertices have disjoint neighbourhoods.
We divide into cases according to the size of $S$.

\begin{description}
\item[$0 \le |S|/|X_i(T)| \le \gamma$.]
For every unembedded $z$ and triple $E$ containing $z$
we have $|F_z(T)| \ge d_un $ by Lemma \ref{3exceptional} and
$|M_{z,E}(T)|\le \theta_*|F_z(T)|$ by Lemma \ref{3marked}.
Since $z$ has degree at most $D$ we have
$|A'_z| \ge (1-D\theta_*)d_u n > \gamma n$, so this case cannot occur.

\item[$\gamma \le |S|/|X_i(T)| \le 1/2$.]
We use the fact that $A := V_i(T) \sm \cup_{z \in S} A'_z$
is a large set of unused vertices which cannot be used by any vertex $z$ in $S$:
we have $|A| \ge |V_i(T)|-|S| \ge |X_i(T)|/2 \ge \delta_B n/4$,
yet $A \cap F_z(T) \sub \cup_{E \ni z} M_{z,E}(T)$ has size at most $D\theta_*|F_z(T)|$
by Lemma \ref{3marked}, so $|A \cap F_z(T)|/|F_z(T)| \le D\theta_* < 2^{-2D} |A|/|V_i|$.
As above, taking a union bound over all possibilities for $i$, $S$ and $A$,
Lemma \ref{3main} implies that the failure probability is at most
$r \cdot 4^{Cn} \cdot p_0^{\gamma \delta_B n/2} < o(1)$, since $p_0 \ll \gamma, \delta_B$.

\item[$1/2 \le |S|/|X_i(T)| \le 1-\gamma$.]
We use the fact that $W := V_i(T) \sm \cup_{z \in S} A'_z$
satisfies $W \cap A_z = W \cap A'_z = \es$ for every $z \in S$.
Now $|W| \ge |V_i(T)|-|S| \ge \gamma |X_i(T)| \ge \gamma\delta_B n/2$,
so by Lemma \ref{3initial}, for each $z$ the event $W \cap A_z = \es$ 
has probability at most $\theta_*$ when we embed $VN_H(z)$, conditional on the prior embedding.
Multiplying the conditional probabilities and taking a union bound
over all possibilities for $i$, $S$ and $W$, the failure probability is at most
$r \cdot 4^{Cn} \cdot \theta_*^{\delta_B n/4} < o(1)$, since $\theta_* \ll \delta_B$.

\item[$1-\gamma \le |S|/|X_i(T)| \le 1$.]
We claim that with high probability $\cup_{z \in S} A'_z = V_i(T)$,
so in fact Hall's criterion holds. It suffices to consider sets $S \sub X_i(T)$
of size exactly $(1-\gamma)|X_i(T)|$. The claim fails if there is some $v \in V_i(T)$
such that $v \notin A'_z$ for every $z \in S$. Since $v$ is unused we have $v\notin A_z$,
and by Lemma \ref{3-x-to-v}, for each $z$ the event $v\notin A_z$ has probability at most $1-p$
when we embed $VN_H(z)$, conditional on the prior embedding.
Multiplying the conditional probabilities
and taking a union bound over all $1 \le i \le r$, $v \in V_i$ and
$S \sub X_i(T)$ of size $(1-\gamma)|X_i(T)|$, the failure probability is at most
$rCn \binom{Cn}{(1-\gamma)Cn} (1-p)^{(1-\gamma)|X_i(T)|} < o(1)$.
This estimate uses the bounds $\binom{Cn}{(1-\gamma)Cn} \le 2^{\sqrt{\gamma}n}$,
$(1-p)^{(1-\gamma)|X_i(T)|} < e^{-p\delta_Bn/4} < 2^{-p^2 n}$ and $\gamma \ll p$.
\end{description}

In all cases we see that the failure probability is $o(1)$. \qed

\section{Applying the blow-up lemma}\label{apps}

To demonstrate the utility of the blow-up lemma  we will work through an application
in this section. To warm up, we sketch the proof of K\"uhn and Osthus
\cite[Theorem 2]{KO2} on packing bipartite graphs using the graph blow-up lemma.
Then we generalise this result to packing tripartite $3$-graphs.
We divide this section into four subsections, organised as follows.
In the first subsection we illustrate the use of the graph blow-up lemma,
which is based on a decomposition obtained from Szemer\'edi's Regularity Lemma
and a simple lemma that one can delete a small number of vertices from a regular
pair to make it super-regular.
The second subsection describes some more hypergraph regularity theory for $3$-graphs:
the Regular Approximation Lemma and Counting Lemma of R\"odl and Schacht.
In the third subsection we give the $3$-graph analogue of the super-regular deletion lemma,
which requires rather more work than the graph case. We also give a `black box' reformulation
of the blow-up lemma that will be more accessible for future applications.
We then apply this in the fourth subsection to packing tripartite $3$-graphs.

\subsection{Applying the graph blow-up lemma}

In this subsection we sketch a proof of the following result of K\"uhn and Osthus.
First we give some definitions. For any graph $F$, an {\em $F$-packing} is a collection \index{packing}
of vertex-disjoint copies of $F$. We say that a graph $G$ is {\em $(a\pm b)$-regular}
if the degree of every vertex in $G$ lies between $a-b$ and $a+b$.

\begin{theo}\label{2pack}
For any bipartite graph $F$ with different part sizes and $0<c\le 1$
there is a real $\eps>0$ and positive integers $C,n_0$
such that any $(1\pm \eps)cn$-regular graph $G$ on $n>n_0$ vertices
contains an $F$-packing covering all but at most $C$ vertices.
\end{theo}

Note that the assumption that $F$ has different part sizes is essential.
For example, if $F=C_4$ is a $4$-cycle and $G$ is a complete bipartite graph
with parts of size $(1+\eps)n/2$ and $(1-\eps)n/2$ then any $F$-packing
leaves at least $\eps n$ vertices uncovered.
Also, we cannot expect to cover all vertices even when
the number $f$ of vertices of $F$ divides $n$,
as $G$ may be disconnected and have a component
in which the number of vertices is not divisible by $f$.
\COMMENT{Even a single edge is an example!}

Without loss of generality we can assume $F$ is
a complete bipartite graph $K_{r,s}$ for some $r \ne s$.\index{$K_{r,s}$}
It is convenient to assume that $G$ is bipartite,
having parts $A$ and $B$ of sizes $n/2$.
This can be achieved by choosing $A$ and $B$ randomly:
if $G$ is $(1\pm \eps)cn$-regular then with high probability
the induced bipartite graph is $(1\pm 2\eps)cn/2$-regular.
Then we refine the partition $(A,B)$ using the following
`degree form' of Szemer\'edi's Regularity Lemma.
\COMMENT{Formerly: $n/2 \pm n^{2/3}$...and applying Chernoff bounds}

\begin{lemma}\label{srl-deg}
Suppose $0 < 1/T \ll \eps \ll d < 1$ and $G=(A,B)$ is a bipartite graph with $|A|=|B|=n/2$.
Then there are partitions $A = A_0 \cup A_1 \cup \cdots \cup A_t$
and $B = B_0 \cup B_1 \cup \cdots \cup B_t$ for some $t \le T$
such that $|A_i|=|B_i|=m$ for $1 \le i \le t$ for some $m$
and $|A_0 \cup B_0| \le \eps n$,
and a spanning subgraph $G'$ of $G$ such that
$d_{G'}(x) > d_G(x)-(d+\eps)n$ for every vertex $x$ and
every pair $(A_i,B_j)$ with $1 \le i,j \le k$ induces a bipartite
subgraph of $G'$ that is either empty or $\eps$-regular of density at least $d$.
\end{lemma}

Lemma \ref{srl-deg} can be easily derived from the usual statement of
Szemer\'edi's Regularity Lemma (see e.g.\ \cite[Lemma 41]{KO4} for a non-bipartite version).
We refer to the parts $A_i$ and $B_i$ with $1 \le i \le t$ as {\em clusters} and \index{cluster}
the parts $A_0$ and $B_0$ as {\em exceptional} sets.
We write $(A_i,B_j)_{G'}$ for the bipartite subgraph of $G'$ induced by $A_i$ and $B_j$.
There is a naturally associated {\em reduced graph}, \index{reduced graph}
in which vertices correspond to clusters and edges to dense regular pairs:
$R$ is a weighted bipartite graph on $([t],[t])$ with an edge $(i,j)$ of weight $d_{ij}$
whenever $(A_i,B_j)_{G'}$ is $\eps$-regular of density $d_{ij} \ge d$.
We choose the parameter $d$ to satisfy $\eps \ll d \ll c$.

The next step of the proof is to select a nearly-perfect matching in $R$, \index{matching}
using the defect form of Hall's matching theorem. Using the fact that $G$
is $(1\pm \eps)cn$-regular one can show that for any $I \sub [t]$
we have $|N_R(I)| \ge (1-2(d+2\eps)/c)|I| > (1-\sqrt{d})|I|$ (say)
so $R$ has a matching of size $(1-\sqrt{d})t$.
The details are given in Lemma 11 of \cite{KO2}.

Now it is straightforward to find an $F$-packing covering all but at most
$3\sqrt{d}n$ vertices. For each edge $(i,j)$ of the matching in $R$ we greedily
remove copies of $K_{r,s}$ while possible, alternating which of $A_i$ and $B_j$
contains the part of size $s$ to maintain parts of roughly equal size.
While we still have at least $\eps m$ vertices remaining in each of $A_i$ and $B_j$,
the definition of $\eps$-regularity implies that the remaining
subgraph $(A_i,B_j)_{G'}$ has density at least $d-\eps$. Then the K\"ovari-S\'os-Tur\'an
theorem \cite{KST} implies that we can choose the next copy of $K_{r,s}$.
We can estimate the number of uncovered vertices
by $2\sqrt{d}t \cdot m < 2\sqrt{d}n$ in clusters not covered by the matching,
$2t \cdot \eps m < 2\eps n$ in clusters covered by the matching,
and $\eps n$ in $A_0 \cup B_0$,
so at most $3\sqrt{d}n$ vertices are uncovered.

However, we want to prove the stronger result that there is an $F$-packing
covering all but at most $Cn$ vertices. To do this we first move a small
number of vertices from each cluster to the exceptional sets so as to make
the matching pairs super-regular. This is a standard property of graph regularity;
we include the short proof of the next lemma for comparison with the analogous
statement later for hypergraphs.

\begin{lemma}\label{2del}
Suppose $G=(A,B)$ is an $\eps$-regular bipartite graph of density $d$
with $|A|=|B|=m$. Then there are $A^* \sub A$ and $B^* \sub B$
with $|A^*|=|B^*|=(1-\eps)m$ such that the restriction $G^*$ of $G$ to
$(A^*,B^*)$ is $(2\eps,d)$-super-regular.
\end{lemma}

\nib{Proof.}
Let $A_0 = \{x \in A: d(x) < (d-\eps)m\}$.
Then $|A_0|<\eps m$, otherwise $(A_0,B)$ would contradict
the definition of $\eps$-regularity for $(A,B)$.
Similarly $B_0 = \{x \in B: d(x) < (d-\eps)m\}$ has $|B_0|<\eps m$.
Let $A^*$ be obtained from $A$ by deleting a set of size $\eps m$ containing $A_0$.
Define $B^*$ similarly. Then $|A^*|=|B^*|=(1-\eps)m$.
For any $A' \sub A^*$, $B' \sub B^*$ with $|A'|>2\eps|A^*|$, $|B'|>2\eps|B^*|$
we have $|A'|>\eps m$, $|B'|>2\eps m$, so $(A',B')_G$ has density
$d \pm \eps$ by $\eps$-regularity of $G$. Thus $G^*$ is $2\eps$-regular.
Also, for any vertex $x$ of $G^*$ we have $d_{G^*}(x) \ge d_G(x)-\eps m
\ge (d-\eps)m - \eps m \ge (d-2\eps)(1-\eps)m$,
so $G^*$ is $(2\eps,d)$-super-regular. \qed

We make the matching pairs $(\eps,2d)$-super-regular
and move all discarded vertices and unmatched clusters into
the exceptional sets. For convenient notation, we redefine
$A_1,\cdots,A_{t'}$ and $B_1,\cdots,B_{t'}$, where $t'=(1-\sqrt{d})t$,
to be the parts of the super-regular matched pairs,
and $A_0$, $B_0$ to be the new exceptional sets.
Thus $|A_i|=|B_i|=(1-\eps)m$ for $1 \le i \le t'$
and $|A_0 \cup B_0| \le \eps n + 2t' \eps m + 2\sqrt{d}tm < 3\sqrt{d} n$.

We will select vertex-disjoint copies of $K_{r,s}$ to cover $A_0 \cup B_0$.
We want to do this in such a way that the matching pairs remain super-regular
(with slightly weaker parameters) and the uncovered parts of the clusters
all have roughly equal sizes. Then we will be able to apply the graph blow-up
lemma (Theorem \ref{2blowup}) to pack the remaining vertices in each matching pair
almost perfectly with copies of $K_{r,s}$. (We assumed $|V_i|=|X_i|=n$ in
Theorem \ref{2blowup} for simplicity, but it is easy to replace this
assumption by $n\le|V_i|=|X_i|\le 2n$, say.) The sizes of the uncovered parts
in a pair may not permit a perfect $K_{r,s}$-packing, but it is easy to see
that one can cover all but at most $r+s$ vertices in each pair.
Taking $C=T(r+s)$ we will thus cover all but at most $C$ vertices.

It remains to show how to cover $A_0 \cup B_0$ with vertex-disjoint copies of $K_{r,s}$.
First we set aside some vertices that we will not use so as
to preserve the super-regularity of the matching pairs.
For each matching pair $(i,j)$ we randomly partition
$A_i$ as $A'_i \cup A''_i$ and $B_j$ as $B'_j \cup B''_j$.
We will only use vertices from $A' = \cup_i A'_i$ and $B' = \cup_j B'_j$
when covering $A_0 \cup B_0$. By Chernoff bounds, with high probability\index{Chernoff}
these partitions have the following properties:
\begin{enumerate}
\item
all parts $A'_i$, $A''_i$, $B'_j$, $B''_j$ have sizes $(1-\eps)m/2 \pm m^{2/3}$.
\item
every vertex $x$ has at least $d(x)/2 - 4\sqrt{d} n$ neighbours in $A' \cup B'$.
\item
whenever a vertex $x$ has $K \ge \eps m$ neighbours in a cluster $A_i$
it has $K/2 \pm m^{2/3}$ neighbours in each of $A'_i$, $A''_i$;
a similar statement holds for clusters $B_j$.
\item
each of $N_{G'}(x) \cap B''_j$ and $N_{G'}(y) \cap A''_i$
have size at least $dm/2 - m^{2/3}$
for every matching edge $(i,j)$ and $x \in A_i$, $y \in B_j$.
\end{enumerate}

Now we cover $A_0 \cup B_0$ by the following greedy procedure.
Suppose we are about to cover a vertex $x \in A_0 \cup B_0$, say $x \in A_0$.
We consider a cluster to be {\em heavy} if we have covered more than $d^{1/4}m$ \index{heavy}
of its vertices. Since $|A_0 \cup B_0|<3\sqrt{d} n$ we have covered at most
$3(r+s)\sqrt{d} n$ vertices by copies of $K_{r,s}$, so there are at most
$4(r+s)d^{1/4}t$ heavy clusters. Since $G$ is $(1\pm \eps)cn$-regular,
$x$ has at least $cn/3$ neighbours in $B'$ by property 2.
At most $4(r+s)d^{1/4}n$ of these neighbours lie in heavy clusters
and at most $cn/4$ of them lie in clusters $B_j$ where $x$ has at most $cm/4$
neighbours. Thus we can choose a matching pair $(i,j)$
such that $x$ has at least $cm/4$ neighbours in $B_j$,
so at least $cm/10$ neighbours in $B'_j$ by property 3.
Since $(A_i,B_j)$ is $2\eps$-regular, $(A'_i,N(x) \cap B'_j)$
has density at least $d/2$, so by K\"ovari-S\'os-Tur\'an
we can choose a copy of $K_{r-1,s}$ with $r-1$ vertices in $A'_i$
and $s$ vertices in $N(x) \cap B'_j$. Adding $x$ we obtain a copy
of $K_{r,s}$ covering $x$.

Thus we can cover $A_0 \cup B_0$, only using vertices from $A' \cup B'$,
and by avoiding heavy clusters we never cover more than $d^{1/4}m + \max\{r,s\}$
vertices in any cluster. For each matching edge $(i,j)$,
the uncovered part of $(A_i,B_j)_{G'}$ is $3\eps$-regular,
and has minimum degree at least $dm/3$ by property 4.
Thus it is super-regular, and as described above
we can complete the proof via the blow-up lemma.

\subsection{The regular approximation lemma and dense counting lemma}

In order to apply regularity methods to $3$-graphs we need a result analogous
to the Szemer\'edi Regularity Lemma, decomposing an arbitrary $3$-graph
into a bounded number of $3$-complexes, most of which are regular.
This was achieved by Frankl and R\"odl \cite{FR}, but as we mentioned in Section \ref{3reg}, 
in this sparse setting the parameters are not suitable for our blow-up lemma. We will instead use
the regular approximation lemma, which provides
a dense setting for an approximation of the original $3$-graph. For $3$-graphs
this result is due to Nagle, R\"odl and Schacht \cite{NRS} and in general to R\"odl and Schacht \cite{RSc1}.
A similar result was proved by Tao \cite{T1}. For simplicity
we will just discuss the lemma for $3$-graphs, although the statement for $k$-graphs
is very similar. Note also that we will formulate the results using the notation
established in this paper. We start with a general definition.
\COMMENT{
1. Only need $k=2,3$.
2. Former footnote: R\"odl and Schacht \cite{RSc1} also obtain a decomposition theorem
without approximation which has a less restrictive parameter hierarchy
than that of Gowers, but with a more complicated notion of regularity
that we will not discuss here.
}

\begin{defn}\label{def-partition}
Suppose that $V=V_1 \cup \cdots \cup V_r$ is an $r$-partite set.
A {\em partition $k$-system} $P$ on $V$ is a collection of partitions \index{partition $k$-system}
$P_A$ of $K(V)_A$ for every $A \in \binom{[r]}{\le k}$.
We say that $P$ is a {\em partition $k$-complex}\index{partition $k$-complex}\index{$P$}\index{$P_A$}
if every two sets $S,S'$ in the same cell of $P_A$ \index{cell}
are {\em strongly equivalent}, defined as in \cite{G2}\index{strongly equivalent}
to mean that $S_B$ and $S'_B$ belong to the same cell of $P_B$
for every $B \sub A$.
Given $S \in K(V)$ we write $C^P_S$ for the cell in $P_A$ containing $S$.
We write $C_S=C^P_S$ when there is no danger of ambiguity.
For any $S' \sub S$ we write $C_{S'} \le C_S$ and say that\index{$\le$}
$C_{S'}$ {\em lies under} or is {\em consistent} with $C_S$.\index{lies under}\index{consistent}
We define the {\em cell complex} $C_{S^\le} = \cup_{S' \sub S} C_{S'}$.\index{cell complex}
\COMMENT{
Formerly very confused, here and in HRT!
Gowers partitions don't do this, but this is a device to control the number
of parts; we define qr on pieces created by strong equivalence, which restores
the RS notion. I was also worried about e.g.\ $(V_1,V_2)$ and
each of $V_1$, $V_2$ being further partitioned by the parts of the bipartite
graphs joining them, but these graphs are quasirandom with respect to $V_1$ and
$V_2$, so cannot restrict to subsets and have zero density elsewhere.
}
\end{defn}

We can use a partition $2$-complex to decompose a $3$-graph as follows.

\begin{defn}\label{def-partition2}
Suppose $G$ is an $r$-partite $3$-graph on $V$
and $P$ is a partition $2$-complex on $V$.
We define $G[P]$ to be the coarsest partition $3$-complex \index{$G[P]$}
refining $P$ and the partitions $\{G_S,K(V)_S \sm G_S\}$ for $S \in G$.
\COMMENT{
I was confused for a while about whether the lower levels
are a $(k-1)$-complex or a partition $(k-1)$-complex:
it depends whether we look at all of $G$ or a single polyad.
}
\end{defn}

We make a few remarks here to explain the structures defined in Definitions
\ref{def-partition} and \ref{def-partition2}. The partition $2$-complex $P$
has vertex partitions and graph partitions. The vertex partitions $P_i$
are of the form $V_i = V_i^1 \cup \cdots \cup V_i^{a_i}$ for some $a_i$, for $1 \le i \le r$.
The graph partitions $P_{ij}$ are of the form $K(V)_{ij} = J_{ij}^1 \cup \cdots \cup J_{ij}^{a_{ij}}$
for some $a_{ij}$, for $1 \le i<j \le r$. By strong equivalence, any bipartite graph $J_{ij}^{b_{ij}}$
is spanned by some pair $(V_i^{b_i},V_j^{b_j})$: it cannot cut across several such pairs.
We also say that $V_i^{b_i}$ and $V_j^{b_j}$ lie under or are consistent with $J_{ij}^{b_{ij}}$.
A choice of $i<j<k$, singleton parts $V_i^{b_i}$, $V_j^{b_j}$, $V_k^{b_k}$ and
graph parts $J_{ij}^{b_{ij}}$, $J_{jk}^{b_{jk}}$, $J_{ik}^{b_{ik}}$ such that the
singleton parts are consistent with the graph parts is called a {\em triad}.\index{triad}
(This terminology is used by R\"odl et al.)
If we consider the set of triangles in a triad, then we obtain a partition of
the $r$-partite triples of $V$ as we range over all triads. Another way to describe this
is to say as in \cite{G2} that $S,S'$ are {\em weakly equivalent}\index{weakly equivalent}
when $S_B,S'_B$ are in the same cell of $P_B$ for every strict subset $B \subn A$.
Let $P^*_A$ denote the partition of $K(V)_A$ into weak equivalence classes. \index{$P^*_A$}
Then $P^*_{ijk}$ is the partition of $K(V)_{ijk}$ by triads as described above.
The partition $3$-complex $G[P]$ has two cells for each triad:
for each cell $C$ of $P^*_{ijk}$ we have cells $G_{ijk} \cap C$
and $(K(V) \sm G)_{ijk} \cap C$ of $G[P]_{ijk}$.
For embeddings in $G$ only the cells $G_{ijk} \cap C$ are of interest,
but we include both for symmetry in the definition.
We make the following further definitions.

\begin{defn}\label{def-partition3}
Suppose $P$ is a partition $k$-complex.
We say that $P$ is {\em equitable} if for every $k' \le k$ the $k'$-cells \index{equitable}
all have equal size, i.e.\ $|T|=|T'|$ for every $T \in P_A$, $T' \in P_{A'}$,
$A, A' \in \binom{[r]}{k'}$.
We say that $P$ is {\em $a$-bounded} if $|P_A| \le a$ for every $A$.
\index{bounded}\index{$a$-bounded|see{bounded}}
We say that $P$ is {\em $\eps$-regular} if every \index{regular}
cell complex $C_{S^\le}$ is $\eps$-regular.
We say that $r$-partite $3$-graphs $G^0$ and $G$ on $V$
are {\em $\nu$-close} if $|G^0_A \Delta G_A| < \nu |K(V)_A|$
\index{close}\index{$\nu$-close|see{close}}
for every $A \in \binom{[r]}{3}$. Here $\Delta$ denotes symmetric difference. \index{$\Delta$}
\COMMENT{
1. The application will average to get $\nu^{1/2}$-close
wrt $P$ on all but $\nu^{1/2}$-prop polyads.
2. We say that $G$ is {\em perfectly $\eps$-regular with respect to $P$}
if $C_{S^\le}$ is $\eps$-regular for every $S \in G[P]$.
[Overkill, but good to stress it:
every intersection with polyad clique is $\eps$-regular w.r.t. polyads
for every induced complex $Q = (P \cup G)(x)$.]
}
\end{defn}

The following is a slightly modified version of the Regular Approximation
Lemma \cite[Theorem 14]{RSc1}.%
\footnote{The differences are:
(i) we are starting with an initial partition of $V$, so technically
we are using a simplified form of \cite[Lemma 25]{RSc1},
(ii) a weaker definition of `equitable' is given in \cite{RSc1},
that the singleton cells have equal sizes, but in fact they prove their
result with the definition used here, and
(iii) we omit the parameter $\eta$ in our statement, as by increasing $r$
we can ensure that all but at most $\eta n^3$ edges are $r$-partite.
}
The reader should note the key point of the constant hierarchy:
although the closeness of approximation $\nu$ may be quite large
(it will satisfy $d_2 \ll \nu \ll d_3)$, the regularity parameter
$\eps$ will be much smaller.
\COMMENT{Former comments:
1. I want $r$-partite to use my notation, although it is a
bit more awkward for RS setup. Note that it replaces $\eta$.
2. $\nu$ clashes with blowup, but no confusion will arise.
3. Equitable assumed in L4.2 (RS?), and no harm in keeping it.
4. Don't use densities to keep some similarities with original.
5. For a while I didn't notice the $(k-1)$-complex instead of $k$ and
wondered why there were no exceptional polyads.
}

\begin{theo} (R\"odl-Schacht \cite{RSc1}) \label{ral}
\index{regular approximation lemma}
Suppose integers $n,a,r$ and reals $\eps, \nu$
satisfy $0 < 1/n \ll \eps \ll 1/a \ll \nu, 1/r$
and that $G^0$ is an $r$-partite $3$-graph
on an equitable $r$-partite set of $n$ vertices
$V = V_1 \cup \cdots \cup V_r$, where $a!|n$.
Then there is an $a$-bounded % $\eps$-regular
equitable $r$-partite partition $2$-complex $P$ on $V$
and an $r$-partite $3$-graph $G$ on $V$ that is
$\nu$-close to $G^0$ such that $G[P]$ is $\eps$-regular.
\COMMENT{$G[P]$ is $a^3$-bounded.}
\end{theo}

As we mentioned earlier, $\eps$-regular $3$-complexes are useful because of
a counting lemma that allows one to estimate the number of copies of any fixed
complex $J$, using a suitable product of densities. First we state a counting
lemma for tetrahedra, analogous to the triangle counting lemma in \ref{eq:tri}.
Suppose $0 < \eps \ll d, \gamma$ and $G$ is an $\eps$-regular $4$-partite $3$-complex
on $V = V_1 \cup \cdots \cup V_4$ with all relative densities $d_S(G) \ge d$.
Then $G^*_{[4]}$ is the set of tetrahedra in $G$. We have the estimate
\begin{equation}\label{eq:tet}
d(G^*_{[4]}) := \frac{|G^*_{[4]}|}{|V_1||V_2||V_3||V_4|} = (1\pm\gamma)\prod_{S\sub [4]}d_S(G).
\end{equation}
This follows from a result of Kohayakawa, R\"odl and Skokan \cite[Theorem 6.5]{KRS}:
they proved a counting lemma for cliques in regular $k$-complexes.

More generally, suppose $J$ is an $r$-partite $3$-complex on $Y = Y_1 \cup \cdots \cup Y_r$
and $G$ is an $r$-partite $3$-complex on $V = V_1 \cup \cdots \cup V_r$.
We let $\Phi(Y,V)$ denote the set of all $r$-partite maps \index{$\Phi(Y,V)$}
from $Y$ to $V$: these are maps $\phi: Y \to V$ such that $\phi(Y_i) \sub V_i$ for each $i$.
We say that $\phi$ is a {\em homomorphism} if $\phi(J) \sub G$. \index{homomorphism}
For $I \sub [r]$ we let $G_I : K(V)_I \to \{0,1\}$ also denote the characteristic \index{$G_I$}
function of $G_I$, i.e.\ $G_I(S)$ is $1$ if $S \in G_I$ and $0$ otherwise.
The following general dense counting lemma \index{dense counting lemma}
from \cite{RSc2} gives an estimate for {\em partite homomorphism density}
\index{partite homomorphism density}
$d_J(G)$ of $J$ in $G$, by which we mean the probability that a random $r$-partite map
from $Y$ to $V$ is a homomorphism from $J$ to $G$. \index{$d_J(G)$}
We use the language of homomorphisms for convenient notation,
but note that we can apply the same estimate to the density
of embedded copies of $J$ in $G$, as most maps are injective.
\COMMENT{Intuitively it expresses the fact, in a regular complex, sets of size
$k$ are almost uncorrelated with sets of size less than $k$...}

\begin{theo}(R\"odl-Schacht \cite{RSc2}, see Theorem 13)\label{3count}%
\footnote{
We have rephrased their statement and slightly generalised it by
allowing the sets $Y_i$ to have more than one vertex: this version can easily
be deduced from the case $|Y_i|=1$, $1 \le i \le r$ by defining an auxiliary
complex with the appropriate number of copies of each $V_i$ (see \cite{CFKO}).}
Suppose $0 < \eps \ll d, \gamma, 1/r, 1/j$, that
$J$ and $G$ are $r$-partite $3$-complexes with vertex sets
$Y = Y_1 \cup \cdots \cup Y_r$ and $V = V_1 \cup \cdots \cup V_r$ respectively,
that $|J|=j$, and $G$ is $\eps$-regular with all densities $d_S(G) \ge d$.
Then
$$d_J(G) = \mb{E}_{\phi \in \Phi(Y,V)} \left[ \prod_{A \in J} G_A(\phi(A)) \right]
= (1 \pm \gamma) \prod_{A \in J} d_A(G).$$
\end{theo}

\subsection{Obtaining super-regularity}

Suppose that we want to embed some bounded degree $3$-graph $H$
in another $3$-graph $G^0$ on a set $V$ of $n$ vertices, where $n$ is large. \index{$G^0$}
We fix constants with hierarchy $0 < 1/n \ll \eps \ll 1/a \ll \nu, 1/r \ll 1$. \index{$a$}
We delete at most $a!$ vertices so that the number remaining is divisible by $a!$,
take an equitable $r$-partition $V = V_1 \cup \cdots \cup V_r$,
and apply Theorem \ref{ral} to obtain an $a$-bounded
equitable $r$-partite partition $2$-complex $P$ on $V$
and an $r$-partite $3$-graph $G$ on $V$ that is
$\nu$-close to $G^0$ such that $G[P]$ is $\eps$-regular.
Since $G$ is so regular, our strategy for embedding $H$ in $G^0$
will be to think about embedding it in $G$, subject
to the rule that the edges $M = G \sm G^0$ are marked as `forbidden'. \index{$M$}
Recall that we refer to the pair $(G,M)$ as a marked complex.
To apply the $3$-graph blow-up lemma (Theorem \ref{3blowup}) we need
the following analogue of Lemma \ref{2del}, showing that we can enforce
super-regularity by deleting a small number of vertices.
\COMMENT{Stick to $r$-partite, later $R$-indexed.}

\begin{lemma} \label{3del} Suppose that
$0 < \eps_0 \ll \eps \ll \eps' \ll d_2 \ll \theta \ll d_3, 1/r$, % $d_3 \ll 1/r$?
and $(G,M)$ is a marked $r$-partite $3$-complex on $V= V_1 \cup \cdots \cup V_r$
such that when defined $G_S$ is $\eps_0$-regular, $|M_S| \le \theta|G_S|$
and $d_S(G) \ge d_{|S|}$ if $|S|=2,3$.
Then we can delete at most $2\theta^{1/3}|G_i|$ vertices from each $G_i$, $1 \le i \le r$
to obtain an $(\eps,\eps',d_2/2,2\sqrt{\theta},d_3/2)$-super-regular marked complex $(G^\sharp,M^\sharp)$.
\COMMENT{
1. Can we delete exactly $2\theta^{1/3}|G_i|$? Should be okay by randomly
making up the difference: whp hit nhoods correctly...
2. Formerly: used $d'_i$ and $\theta'$, need $d_3 \ll 1/r$? (No)
3. Remove 2s? Too much effort...
}
\end{lemma}

\nib{Proof.}
The idea is to delete vertices which cause failure of the regularity,
density or marking conditions in Definition \ref{def-3super} (super-regularity).
However, some care must be taken to ensure that this process terminates.
There are three steps in the proof: firstly, we identify
sets $Y_i$ of vertices in $G_i$ that cause the conditions
on marked edges to fail; secondly, we identify sets $Z_i$
of vertices in $G_i$ that either cause the regularity and density conditions
to fail or have atypical neighbourhood in some $Y_j$;
thirdly, we delete the sets $Y_i$ and $Z_i$ and show
that what remains is a super-regular pair.

\nib{Step 1.} Fix $1 \le i \le r$. We will identify a set $Y_i$ of vertices in $G_i$
that are bad with respect to the conditions on marked edges in the definition of super-regularity.
For any $j,k$ such that $G_{ijk}$ is defined we let $Y_{i,jk}$ be the set of vertices $v \in G_i$
for which $|M(v)_{jk}| > \sqrt{\theta}|G(v)_{jk}|$. For any triple $S$ such that $G_S$ is defined
and subcomplex $I$ of $S^\le$ such that $G_{S'i}$ is defined for all $S'\in I$
we let $Y^I_{i,S}$ be the set of vertices $v \in G_i$
for which $|(M \cap G^{I_v})_S| > \sqrt{\theta}|G^{I_v}_S|$.
Let $Y_i$ be the union of all such sets $Y_{i,jk}$ and $Y^I_{i,S}$.
We will show that $|Y_i| < \theta^{1/3}|G_i|$.

First we bound the sets $Y_{i,jk}$.  Let $Z_{i,jk}$ be the set of vertices $v \in G_i$
such that we do not have $|G(v)_{jk}|=(1\pm\eps)|G_{ijk}|/|G_i|$ and $G(v)_S$
is $\eps$-regular with $d_S(G(v))= (1\pm \eps)d_S(G)d_{Si}(G)$ for $\es\ne S \sub jk$.
Since $G_{ijk}$ is $\eps_0$-regular we have $|Z_{i,jk}|<\eps|G_i|$
by Lemma \ref{3average} and Lemma \ref{3neighbour1}.
Therefore $\sum_{v \in Y_{i,jk}} |M(v)_{jk}| > \sqrt{\theta}\sum_{v \in Y_{i,jk}\sm Z_{i,jk}} |G(v)_{jk}|
> \sqrt{\theta}(|Y_{i,jk}|-\eps|G_i|)(1-\eps)|G_{ijk}|/|G_i|$.
We also have an upper bound $\sum_{v \in Y_{i,jk}} |M(v)_{jk}| \le |M_{ijk}| \le \theta|G_{ijk}|$
by the hypotheses of the lemma. This gives $|Y_{i,jk}|/|G_i| < \sqrt{\theta}/(1-\eps)+\eps < 2\sqrt{\theta}$.

Now we bound $Y^I_{i,S}$. Define $Z_{i,S}$ to equal $Z_{i,jk}$ if $S=ijk$
or $Z_{i,ab} \cup Z_{i,bc} \cup Z_{i,ac}$ if $i \notin S = abc$.
If $v \in G_i \sm Z_{i,S}$ then by regular restriction $G^{I_v}_{S^\le}$ is $\sqrt{\eps}$-regular
and $d_{S'}(G^{I_v})$ is $(1\pm\eps)d_{S'}(G)d_{S'i}(G)$ if $\es \ne S' \in I$
or $(1\pm\sqrt{\eps})d_{S'}(G)$ otherwise.
By Lemma \ref{3absolute} we have $d(G_S)=(1\pm 8\eps)\prod_{S'\sub S}d_{S'}(G)$
and $d(G^{I_v}_S)=(1\pm 8\sqrt{\eps})\prod_{S'\sub S}d_{S'}(G^{I_v})
= (1\pm 20\sqrt{\eps})d(G_S)\prod_{\es \ne S' \in I}d_{S'i}(G)$.

Write $\Sa = \sum_{v \in Y^I_{i,S}} |(M \cap G^{I_v})_S|$.
Then\COMMENT{lemma for $G^{I_v}$ calcs?}
\[\Sa > \sqrt{\theta} \sum_{v \in Y^I_{i,S}\sm Z_{i,S}} |G^{I_v}_S|
> \sqrt{\theta}(|Y^I_{i,S}|-3\eps|G_i|)(1-20\sqrt{\eps})|G_S|\prod_{\es \ne S' \in I}d_{S'i}(G).\]
For any $P \in G_S$, let $G_{P,I}$ be the set of $v \in G_i$
such that $P_{S'} v \in G_{S'i}$ for all $\es \ne S' \in I$.
Let $B_I$ be the set of $P \in G_S$ such that we do not have
$|G_{P,I}| = (1 \pm \eps')|G_i| \prod_{\es \ne S' \in I} d_{S'i}(G)$.
Then $|B_I| \le \eps'|G_S|$ by Lemma \ref{3technical}.
Now $\Sa \le \sum_{v \in G_i} |(M \cap G^{I_v})_S|$, which counts all pairs
$(v,P)$ with $P \in M_S$ and $v \in G_{P,I}$, so
\[\Sa \le |M_S|(1+\eps')|G_i|\prod_{\es \ne S' \in I} d_{S'i}(G) + \eps'|G_S||G_i|.\]
Combining this with the lower bound on $\Sa$ we obtain
\[ \sqrt{\theta}(|Y^I_{i,S}|/|G_i|-3\eps)(1-20\sqrt{\eps})
< \frac{|M_S|}{|G_S|} (1+\eps') + \eps'\prod_{\es\ne S' \in I}d_{S'i}(G)^{-1}. \]
Since $|M_S|<\theta|G_S|$ and $\eps' \ll d_2$ we deduce that
$|Y_{i,jk}|/|G_i| < 2\sqrt{\theta}$.

In total we deduce that $|Y_i|/|G_i| < \left( \binom{r-1}{2} + 2^8 \binom{r}{3} \right)2\sqrt{\theta} < \theta^{1/3}$.

\nib{Step 2.} Next we come to the regularity and density conditions.
Recall that $G(v)$ is $\eps$-regular with $d_S(G(v))= (1\pm \eps)d_S(G)d_{Si}(G)$
for $\es\ne S \sub jk$ when $v \notin Z_{i,jk}$, where $|Z_{i,jk}|<\eps|G_i|$.
Now suppose $G_{ij}$ is defined and let $Z'_{i,j}$ be the set of $v \in G_i$ such that
$|G(v)_j \cap Y_j| \ne d_{ij}(G)|Y_j| \pm \eps|G_j|$.
We claim that $|Z'_{i,j}|<2\eps|G_i|$.
To prove this, we can assume that $|Y_j|>\eps|G_j|$, otherwise $Z'_{i,j}$ is empty.
Since $G_{ij}$ is $\eps_0$-regular, $G_{ij}[Y_j]$ is $\eps$-regular by Lemma \ref{2restrict}.
Then the bound on $Z'_{i,j}$ follows from Lemma \ref{2neighbour}.
Let $Z_i$ be the union of all the sets $Z_{i,jk}$ and $Z'_{i,j}$.
Then $|Z_i|/|G_i| < \binom{r}{2}\eps + 2(r-1)2\eps < \sqrt{\eps}$, say.

\nib{Step 3.}
Now we show that deleting $Y_i \cup Z_i$ from $G_i$ for every $1 \le i \le r$
gives an $(\eps,\eps',d_2/2,2\sqrt{\theta},d_3/2)$-super-regular marked complex $(G^\sharp,M^\sharp)$.
To see this, note first that the above upper bounds on $Y_i$ and $Z_i$ show that
we have deleted at most $2\theta^{1/3}$-proportion of each $G_i$.
Since $G$ is $\eps_0$-regular with $d_S(G) \ge d_{|S|}$ if $|S|=2,3$,
regular restriction implies that $G^\sharp$ is $\eps$-regular
with $d_S(G^\sharp)=(1\pm\eps)d_S(G)$ if $|S|=2,3$
and $d_i(G^\sharp)>(1-2\theta^{1/3})d_i(G)$ for $1\le i \le r$.
This gives property (i) of super-regularity.

Now suppose $G_{ij}$ is defined and $v \in G^\sharp_i$.
Then $|G(v)_j| =(1\pm\eps)d_{ij}(G)|G_j|$ and
$|G(v)_j \cap Y_j| =  d_{ij}(G)|Y_j| \pm \eps|G_j|$ since $v\notin Z'_{i,j}$.
Since $G^\sharp_j = G_j \sm (Y_j \cup Z_j)$ and $|Z_j|<\sqrt{\eps}|G_j|$ we have
\[|G^\sharp(v)_j| = |G(v)_j \sm (Y_j\cup Z_j)|
= |G(v)_j| - d_{ij}(G)|Y_j| \pm 2\sqrt{\eps}|G_j|
= (1\pm\eps')d_{ij}(G^\sharp)|G^\sharp_j|.\]
Next suppose that $G_{ijk}$ is defined and $v \in G^\sharp_i$.
Then $d_{jk}(G(v))=(1\pm\eps)d_{jk}(G)d_{ijk}(G)$ and $G(v)_{jk}$ is $\eps$-regular,
so $d_{jk}(G^\sharp(v))=(1\pm\eps')d_{jk}(G^\sharp)d_{ijk}(G^\sharp)$ 
and $G^\sharp(v)_{jk}$ is $\eps'$-regular by regular restriction.
Also, since $v \notin Y_{i,jk}$ we have $|M(v)_{jk}| \le \sqrt{\theta}|G(v)_{jk}|$.
Since
\[|G^\sharp(v)_{jk}| = d_{jk}(G^\sharp(v))|G^\sharp(v)_j||G^\sharp(v)_k|
> (1-\eps')(1-2\theta^{1/3})^2 d_{jk}(G(v))|G(v)_j||G(v)_k|
> \frac{1}{2}|G(v)_{jk}|\]
we have $|M^\sharp(v)_{jk}| \le |M(v)_{jk}| \le 2\sqrt{\theta}|G^\sharp(v)|$.
This gives property (ii) of super-regularity.

Finally, consider any triple $S$ such that $G_S$ is defined
and subcomplex $I$ of $S^\le$ such that $G_{S'i}$ is defined for all $S'\in I$.
Since $v \notin Z_{i,S}$, $G^{I_v}_{S^\le}$ is $\sqrt{\eps}$-regular
and $d_{S'}(G^{I_v})$ is $(1\pm\eps)d_{S'}(G)d_{S'i}(G)$ if $\es \ne S' \in I$
or $(1\pm\sqrt{\eps})d_{S'}(G)$ otherwise. For $j \in S$,
we have $d_j(G^{\sharp I_v})=d_j(G^\sharp)>(1-2\theta^{1/3})d_j(G)$ if $j \notin I$ or
$d_j(G^{\sharp I_v})=d_j(G^\sharp(v)) = (1\pm\eps')d_{ij}(G^\sharp)d_j(G^\sharp) > \frac{1}{2}d_2 d_j(G)$ if $j \in I$.
Then by regular restriction, for $|S'|\ge 2$ with $S' \sub S$, $G^{\sharp I_v}_{S'}$ is $\eps'$-regular
with $d_{S'}(G^{\sharp I_v})=(1\pm\eps^{1/4})d_{S'}(G^{I_v})$ equal to
$(1\pm\eps')d_{S'}(G^\sharp)d_{S'i}(G^\sharp)$ if $\es \ne S' \in I$
or $(1\pm\eps')d_{S'}(G^\sharp)$ otherwise.
Also, by Lemma \ref{3absolute} we have
\[ \frac{|G^{\sharp I_v}_S|}{|G^{I_v}_S|} = \frac{d(G^{\sharp I_v}_S)}{d(G^{I_v}_S)}
= (1\pm 10\eps')\prod_{S'\sub S}\frac{d_{S'}(G^{\sharp I_v})}{d_{S'}(G^{I_v})}
> (1-20\eps')(1-2\theta^{1/3})^3 > 1/2.\]
Since $v \notin Y^I_{i,S}$ we have $|(M \cap G^{I_v})_S| \le \sqrt{\theta}|G^{I_v}_S|$,
so $|(M^\sharp \cap G^{\sharp I_v})_S| \le 2\sqrt{\theta}|G^{\sharp I_v}_S|$.
This gives property (iii) of super-regularity, so the proof is complete. \qed

\begin{rem}
For some applications it may be important to preserve exact equality of the
part sizes, i.e.\ start with $|G_i|=n$ for $1\le i \le r$ and delete exactly
$2\theta^{1/3}n$ vertices from each $G_i$ to obtain super-regularity.
This can be achieved by deleting the sets $Y_i$ and $Z_i$ of Lemma \ref{3del},
together with some randomly chosen additional  vertices to equalise the numbers.
With high probability these additional vertices intersect all vertex neighbourhoods
in approximately the correct proportion, and then the same proof shows that
the resulting marked complex is super-regular. We omit the details.
\end{rem}

For applications it is also useful to know that super-regularity
is preserved when one restricts to subsets that are both large
and have large intersection with every vertex neighbourhood.

\begin{lemma}\label{3super-restrict}{\bf (Super-regular restriction)}
\index{super-regular restriction}
Suppose $0< \eps \ll \eps' \ll \eps'' \ll d_2  \ll \theta \ll d_3, d'$
and $(G,M)$ is a $(\eps,\eps',d_2,\theta,d_3)$-super-regular
marked $r$-partite $3$-complex on $V=V_1\cup\cdots\cup V_r$ with $G_i=V_i$ for $1\le i \le r$.
Suppose also that we have $V'_i \sub V_i$ for $1 \le i \le r$,
write $V'=V'_1\cup\cdots\cup V'_r$, $G'=G[V']$, $M'=M[V']$,
and that $|V'_i|\ge d'|V_i|$ and $|G(v)_i \cap V'_i| \ge d'|G(v)_i|$
whenever $1 \le i,j \le r$, $v \in V'_j$ and $G_{ij}$ is defined.
Then $(G',M')$ is $(\eps',\eps'',d_2/2,\sqrt{\theta},d_3/2)$-super-regular.
\end{lemma}

\nib{Proof.} The argument is similar to Step 3 of the previous lemma.
Suppose $|S|=3$, $G_S$ is defined, $i\in S$, $v\in G_i$.
By Definition \ref{def-3super}(i) for $(G,M)$,
$G_{S^\le}$ is $\eps$-regular with $d_{S'}(G) \ge d_{|S'|}$ for $S'\sub S$, $|S'|=2,3$.
By assumption we have $d_j(G') \ge d'd_j(G)$ for $j \in S$, so regular restriction implies that
$G'_{S^\le}$ is $\eps'$-regular with $d_{S'}(G') \ge d_{|S'|}/2$ for $S'\sub S$, $|S'|=2,3$.
This gives Definition \ref{def-3super}(i) for $(G',M')$.
Similarly, by Definition \ref{def-3super}(ii) for $(G,M)$, writing $S=ijk$,
$G(v)_{jk^\le}$ is $\eps'$-regular with $d_{S'}(G(v))=(1\pm\eps')d_{S'}(G)d_{S'i}(G)$
for $\es\ne S'\sub jk$ and $|M(v)_{jk}| \le \theta|G(v)_{jk}|$.
By assumption we have $d_j(G'(v)) \ge d'd_j(G(v))$
and $d_k(G'(v)) \ge d'd_k(G(v))$, so regular restriction implies that
$G'(v)_{jk^\le}$ is $\eps''$-regular with $d_{S'}(G'(v))=(1\pm\eps'')d_{S'}(G')d_{S'i}(G')$
for $\es\ne S'\sub jk$. Also $|G'(v)_{jk}|/|G(v)_{jk}|=d(G'(v)_{jk})/d(G(v)_{jk})
> (d')^2/2$ so $|M'(v)_{jk}|/|G'(v)_{jk}| \le 2\theta/(d')^2 < \sqrt{\theta}$.
This gives Definition \ref{def-3super}(ii) for $(G',M')$.

Finally, consider any triple $S$ such that $G_S$ is defined
and subcomplex $I$ of $S^\le$ such that $G_{S'i}$ is defined for all $S'\in I$.
By Definition \ref{def-3super}(iii) for $(G,M)$,
$|(M \cap G^{I_v})_S| \le \theta|G^{I_v}_S|$, $G^{I_v}_{S^\le}$ is $\eps'$-regular
and $d_{S'}(G^{I_v})$ is $(1\pm\eps')d_{S'}(G)d_{S'i}(G)$ if $\es \ne S' \in I$
or $(1\pm\eps')d_{S'}(G)$ otherwise.
By assumption $d_j(G'{}^{I_v})$ is $d_j(G')>d'd_j(G)$ if $j \notin I$
or $d_j(G'(v)) \ge d'd_j(G(v))$ if $j \in I$.
Then by regular restriction, for $|S'|\ge 2$ with $S' \sub S$, $G'{}^{I_v}_{S'}$ is $\eps''$-regular
with $d_{S'}(G'{}^{I_v})=(1\pm\sqrt{\eps'})d_{S'}(G^{I_v})$ equal to
$(1\pm\eps'')d_{S'}(G')d_{S'i}(G')$ if $\es \ne S' \in I$
or $(1\pm\eps'')d_{S'}(G')$ otherwise. Also, by Lemma \ref{3absolute} we have
$|G'{}^{I_v}_S|/|G^{I_v}_S| = d(G'{}^{I_v}_S)/d(G^{I_v}_S)
= (1\pm 10\eps')\prod_{S'\sub S}d_{S'}(G'{}^{I_v})/d_{S'}(G^{I_v}) > d'{}^3/2$.
Therefore $|(M \cap G'{}^{I_v})_S|/|G'{}^{I_v}_S| \le 2\theta/d'{}^3 < \sqrt{\theta}$.
This gives Definition \ref{def-3super}(iii) for $(G',M')$. \qed

Now we will present a `black box' reformulation that \index{black box}
goes straight from regularity to embedding, bypassing super-regularity
and the blow-up lemma. This more accessible form of our results will be
more convenient for future applications of the method.
First we give a definition.
\COMMENT{
1. Need $c \gg \theta$ as written. Can we go lower? Could have marking condition
on res pos but not black box. Black box could be intersection of constant number
of nhood complexes with marking condition and allow dense subset of this.
2. Split different uses of $c$? No, too fussy.
3. Need to split $c$ and $c_0$ later, as $c_0 \ll \delta_Q$ and $c \gg \theta$, $c \gg d_u$.
4. Eliminated M.
}

\begin{defn}\label{def-3robust}{\bf (Robustly universal)}\index{robustly universal}
Suppose $J$ is an $r$-partite $3$-complex on $Y=Y_1\cup\cdots\cup Y_r$ with $J_i=Y_i$ for $1\le i \le r$.
We say that $J$ is {\em $c^\sharp$-robustly $D$-universal} if whenever
\begin{itemize}
\item[(i)] $Y'_i \sub Y_i$ with $|Y'_i| \ge c^\sharp|Y_i|$ such that $Y'=\cup_{i=1}^r Y'_i$,
$J'=J[Y']$ satisfy $|J'_S(v)| \ge c^\sharp|J_S(v)|$ whenever
$|S|=3$, $J_S$ is defined, $i\in S$, $v\in J'_i$,
\item[(ii)] $H'$ is an $r$-partite $3$-complex on $X'=X'_1\cup\cdots\cup X'_r$
of maximum degree at most $D$ with $|X'_i|=|Y'_i|$ for $1 \le i \le r$,
\end{itemize}
then $J'$ contains a copy of $H'$, in which vertices of $X'_i$ correspond to vertices of $Y'_i$
for $1 \le i \le r$.\index{$c^\sharp$}
\COMMENT{
Formerly: 1. $(c^*,d^*)$-robustly $D$-universal def
included: $d(J_S)>d^*$ and $|J_S(v)|>d^*|J_S|/|J_i|$
whenever $|S|=3$, $J_S$ is defined, $i\in S$, $v\in J_i$.
2. * below conflicts with earlier, use sharp!
}
\end{defn}

Now we show that one can delete a small number of vertices from a regular complex
with a small number of marked triples to obtain a robustly universal pair.

\begin{theo}\label{3robust}
Suppose $0 \le 1/n \ll \eps \ll d^\sharp \ll d_2 \ll \theta \ll d_3, c^\sharp, 1/D, 1/C$,\index{$d^\sharp$}
that $G$ is an $\eps$-regular $r$-partite $3$-complex on $V=V_1\cup\cdots\cup V_r$
with $d_S(G) \ge d_{|S|}$ for $|S|=2,3$ when defined and $n\le |G_i|=|V_i|\le Cn$
for $1 \le i \le r$, and $M \sub G_=$ with $|M_S| \le \theta|G_S|$ when defined.
Then we can delete at most $2\theta^{1/3}|G_i|$ vertices from $G_i$ for $1 \le i \le r$
to obtain $G^\sharp$ and $M^\sharp$ so that
\begin{itemize}
\item[(i)] $d(G^\sharp_S)>d^\sharp$
and $|G^\sharp_S(v)|>|G^\sharp_S|/2|G^\sharp_i|$
whenever $|S|=3$, $G_S$ is defined, $i\in S$, $v\in G_i$, and
\item[(ii)] $G^\sharp\sm M^\sharp$ is $c^\sharp$-robustly $D$-universal.
\end{itemize}
\end{theo}

\nib{Proof.}
Define additional constants with $\eps \ll \eps_1 \ll \eps_2 \ll \eps_3 \ll d^\sharp$.
Applying Lemma \ref{3del}, we can delete  at most $2\theta^{1/3}|G_i|$ vertices from each $G_i$
to obtain an $(\eps_1,\eps_2,d_2/2,2\sqrt{\theta},d_3/2)$-super-regular marked complex
$(G^\sharp,M^\sharp)$ on some $V^\sharp=V^\sharp_1\cup\cdots\cup V^\sharp_r$.
We will show that $J=G^\sharp\sm M^\sharp$ is $(c^\sharp,d^\sharp)$-robustly $D$-universal.
To see this suppose $|S|=3$,  $G_S$ is defined, $i\in S$, $v\in G'_i$.
By Definition \ref{def-3super}(i) and Lemma \ref{3absolute} we have
$d(G^\sharp_S)=(1\pm 8\eps_1)\prod_{S'\sub S}d_{S'}(G^\sharp)
> (1-8\eps_1)(1-2\theta^{1/3})^3(d_2/2)^3(d_3/2) > d^\sharp$.
Writing $S=ijk$, by Definition \ref{def-3super}(ii) we have
$d(G^\sharp_S(v)) = d_{jk}(G^\sharp(v))d_j(G^\sharp(v))d_k(G^\sharp(v))
= \prod_{\es\ne S' \sub jk} (1\pm\eps_2)d_{S'}(G^\sharp)d_{S'i}(G^\sharp)$
so $|G^\sharp_i||G^\sharp_S(v)|/|G^\sharp_S| = d_i(G^\sharp)d(G^\sharp_S(v))/d(G^\sharp_S)
= (1\pm 8\eps_1)/(1\pm\eps_2)^3 > 1/2$.
Now suppose that $V'_i \sub V^\sharp_i$ and $H'$ are given
as in Definition \ref{def-3robust} applied to $J$.
Then $(G^\sharp[V'],M^\sharp[V'])$ is $(\eps_2,\eps_3,d_2/4,2\theta^{1/4},d_3/4)$-super-regular
by Lemma \ref{3super-restrict}. Applying Theorem \ref{3blowup}, $J' = J[V']$ contains a copy of $H'$,
in which vertices of $X'_i$ correspond to vertices of $V'_i$ for $1 \le i \le r$. \qed

\subsection{Applying the $3$-graph blow-up lemma}

In this subsection we illustrate the $3$-graph blow-up lemma by proving the following
generalisation of Theorem \ref{2pack} to packings of tripartite $3$-graphs.

\begin{theo}\label{3pack} \index{packing}
For any $3$-partite $3$-graph $F$ in which not all part sizes are equal and $0<c_1,c_2<1$
there is a real $\eps>0$ and positive integers $C,n_0$ such that if
$G$ is a $3$-graph on $n > n_0$ vertices $V$ such that
every vertex $v$ has degree $|G(v)| = (1 \pm \eps)c_1n^2$ % was $\binom{n-1}{2}$
and every pair of vertices $u,v$ has degree $|G(uv)| > c_2n$
then $G$ contains an $F$-packing that covers all but at most $C$ vertices.
\end{theo}

We start by recording some auxiliary results that will be used in the proof.
First we need a result of Erd\H{o}s on the number of copies of a $k$-partite $k$-graph.

\begin{theo} \label{erdos} (Erd\H{o}s \cite{E})
For any $a>0$ and $k$-partite $k$-graph $F$ on $f$ vertices there is $b>0$ so that if
$H$ is a $k$-graph on $n$ vertices with at least $a n^k$ edges then
$H$ contains at least $b n^f$ copies of $F$.
\end{theo}

Next we need Azuma's inequality on martingale deviations.

\begin{theo} \label{azuma} (Azuma \cite{Az})\index{martingale}\index{Azuma}
Suppose $Z_0, \cdots, Z_n$ is a martingale, i.e.\ a sequence of random variables
satisfying $\mb{E}(Z_{i+1}|Z_0,\cdots,Z_i)=Z_i$, and that $|Z_i-Z_{i-1}| \le c_i$,
$1 \le i \le n$, for some constants $c_i$. Then for any $t \ge 0$,
$$\mb{P}(|Z_n-Z_0| \ge t) \le 2\exp \left( - \frac{t^2}{2\sum_{i=1}^n c_i^2}
\right).$$
\end{theo}

We also need the following theorem of Kahn, which is a
linear programming generalisation of a theorem of Pippenger on
matchings in regular hypergraphs with small codegrees.
(A {\em matching}\index{matching} is a set of vertex-disjoint edges.)
Suppose $H$ is a $k$-graph and $t:E(H) \to \mb{R}^+$. Write
$t(H) = \sum_{A \in E(H)} t(A)$ and $t'(S) = \sum_{A \in E(H), S \sub A} t(A)$
for $S \sub V(H)$. Let $\mbox{co}(t) = \max_{S \in \binom{V(H)}{2}} t'(S)$.
Say that $t$ is a {\em fractional matching} if $t'(x) = \sum_{A \in E(H): x \in A} t(A) \le 1$
for every $x \in V(H)$.%
\index{$t(H)$}\index{$t'(S)$}\index{$\mbox{co}(t)$}\index{fractional matching}\index{matching}

\begin{theo} \label{kahn} (Kahn \cite{Ka})
For any $\eps>0$ there is $\delta>0$ so that if $H$ is a $k$-graph on $n$ vertices
and $t$ is a fractional matching of $H$ with $\mbox{co}(t) < \delta$
then $H$ has a matching of size at least $t(H)-\eps n$.
\COMMENT{Formerly $(1-\eps)t(H)$, but what if $t(H) \ll \delta$?!}
\end{theo}

Now we will prove Theorem \ref{3pack}.
We may suppose that $F$ is complete $3$-partite, say $F=K(Y)_{123}$ on $Y = Y_1 \cup Y_2 \cup Y_3$.
We introduce parameters with the hierarchy\COMMENT{$\theta$, $\theta'$ unused}
\[0 \ll 1/n_0 \ll \eps \ll \eps' \ll \eps'' \ll \ll d^\sharp \ll d_2 \ll 1/a
\ll \nu \ll 1/r \ll d_3
\ll \delta \ll \gamma \ll \beta \ll \alpha \ll c_1, c_2, 1/|Y|.\]
\index{$\delta$}\index{$\gamma$}\index{$\beta$}\index{$\alpha$}
We delete at most $a!$ vertices of $G$ so that the number remaining is divisible by $a!$,
take an equitable $r$-partition $V = V_1 \cup \cdots \cup V_r$,
and apply Theorem \ref{ral} to obtain an $a$-bounded
equitable $r$-partite partition $2$-complex $P$ on $V$
and an $r$-partite $3$-graph $G'$ on $V$ that is
$\nu$-close to $G$ such that $G'[P]$ is $\eps$-regular.
Since $P$ is $a$-bounded, for every graph $J_{ij} \in P_{ij}$ with $1\le i<j \le r$,
the densities $d(J_i)$, $d(J_j)$ and $d(J_{ij})$ are all at least $1/a$.
We refer to singleton parts in $P$ as {\em clusters}\index{cluster}.
Let $n_1$ be the size of each cluster.\index{$n_1$}
By means of an initial partition we may also assume that $n_1 < \nu n$.
Let $M = G' \sm G$ be the edges marked as `forbidden'.
\COMMENT{
1. Note position of $r$ in hierarchy.
2. $n_1<\nu n$: anything between $1/a$ and $d$ will do;
differs from blow-up notation where part sizes $>n$.
3. Change to $a_1$ clusters per part, size $n_1$.
[4. Cut $d^\sharp$. Replaced at refs request.]
}

Next we define the {\em reduced $3$-graph}\index{reduced $3$-graph} $R$\index{$R$},
a weighted $3$-graph in which vertices
correspond to clusters and triples correspond to cells of $G'[P]$
that are useful for embedding, in that they have large density and few marked edges.
Let $Z = Z_1 \cup \cdots \cup Z_r$ be an $r$-partite set with\index{$Z$}
$|Z_i| = a_1 := |P_i|$ for $1 \le i \le r$,\index{$a_1$}
where $a_1 \le a$ and $n-a! < rn_1a_1 \le n$.
We identify $Z_i$ with $[a_1]$, although it is to be understood that
$Z_i$ and $Z_j$ are disjoint for $i \ne j$, and label the cells of
$P_i$ as $C_{i,1}, \cdots, C_{i,a_1}$.
We identity an $r$-partite triple $S \in K(Z)$
with $S' = \cup_{i \in S} C_{i,S_i}$, where $S_i=S \cap Z_i$.
Write $N=n^2n_1$\index{$N$} and $K(S')=K(V)[S']$.\index{$K(S')$}
We say that $S$ is an edge of $R$ with
{\em weight}\index{weight}  $w(S) = |G'[S']_S|/N$\index{$w(S)$}
if $|G'[S']_S| > d_3|K(S')_S|$ and $|M[S']_S| < \sqrt{\nu}|K(S')_S|$.

We define the {\em weighted degree} $d_w(j)$ of a vertex $j$ in $R$ to be the sum of $w(S)$
\index{weighted degree}\index{$d_w(j)$}
over all $j \in S \in R$. We will delete a small number of vertices from $R$
to obtain a $3$-graph $R'$ that is almost regular with respect to weighted degrees.
For any $i \in A \in \binom{[r]}{3}$ and $j \in Z_i$
we define\index{$B^j_{i,A}$}\index{$Z_{i,A}$}
\[B^j_{i,A} = \{S: j \in S \in K(Z)_A, |M[S']_S| > \sqrt{\nu}|K(S')_S|\}
\mbox{ and } Z_{i,A} = \{j \in Z_i: |B^j_{i,A}| > \nu^{1/4}a_1^2\}.\]
Then
\[ |M_A| = \sum_{j \in Z_i} \sum_{S: j \in S \in K(Z)_A} |M[S']|
> \sum_{j \in Z_{i,A}} \sum_{S \in B^j_{i,A}}  \sqrt{\nu}|K(S')_S|
> |Z_{i,A}| \cdot \nu^{1/4}a_1^2 \cdot \sqrt{\nu}n_1^3.\]
Since $G'$ is $\nu$-close to $G$ (see Definition \ref{def-partition3})
we have $|M_A|<\nu|K(V)_A| < \nu(n/r)^3$,
so $|Z_{i,A}| < \nu^{1/4} (n/rn_1)^3/a_1^2 < 2\nu^{1/4}a_1$.
Let $Z'= Z'_1 \cup \cdots \cup Z'_r$\index{$Z'$} be obtained by deleting all sets $Z_{i,A}$
from $Z$ and let $R'=R[Z']$.\index{$R'$}
Since $\nu \ll 1/r$ we have $|Z'_i|>(1-\nu^{1/5})a_1$ for $1\le i \le r$.

Now we estimate the weighted degrees $d'_w(j)$ in $R'$. Suppose $j \in Z'$.
We have $d'_w(j) = N^{-1} \sum_{S: j \in S \in R} |G'[S']_S|$.\index{$d'_w(j)$}
There are at most $a_1^2$ triples $S\in K(V)$ containing $j$,
so at most $d_3 a_1^2n_1^3 \le d_3 N$ triples in $G'[S']_S$
for such $S$ with $|G'[S']_S| < d_3|K(S')_S|$.%
\COMMENT{Ref says $\sqrt{d_3}$ because use $G'$-density $>d_3$; I don't see why...}
There are at most $\nu^{1/5}a_1^2$ triples $S \nsub Z'$ with $j\in S\in K(V)$,
so at most $\nu^{1/5}N$ triples in $G'[S']_S$ for such $S$.
Since $j \in Z'$ we have $|B^j_{i,A}| \le \nu^{1/4}a_1^2$,
so there are at most $\nu^{1/4}N$ triples in $G'[S']_S$
for triples $S$ with $j\in S\in K(V)$ and $|M[S']_S| > \sqrt{\nu}|K(S')_S|$.
Finally, there are at most $r(n/r)^2n_1=N/r$ triples that are not $r$-partite.
Altogether, at most $\frac{3}{2}d_3N$ triples contributing to $d'_w(j)$
do not belong to $G'[S']_S$ with $j \in S \in R$.
Since $|G(v)| = (1 \pm \eps)c_1n^2$ for all $v \in G$
and $\eps \ll 1/a \ll \nu \ll 1/r \ll d_3$ we have
\[d'_w(j) = N^{-1} \sum_{S: j \in S \in R} |G'[S']_S|
=  N^{-1} \sum_{v \in C_{i,j}} |G(v)| \pm 3d_3/2 = c_1 \pm 2d_3.\]

Define $t:E(R')\to\mb{R}^2$ by $t(S) = w(S)/(c_1+2d_3)$.
Then $t'(j) = d'_w(j)/(c_1+2d_3)\le 1$ for $j \in R'$, so $t$ is a fractional matching.
We have 
\[t(R') = \sum_{j \in R'} \sum_{S:j \in S} t'(j)/3 = \sum_{j \in R'} d'_w(j)/3(c_1+2d_3) > (1-\sqrt{d_3})|Z'|/3.\]
Also, the trivial bound $w(S) < n_1^3/N$ gives
$\mbox{co}(t) < n/n_1 \cdot n_1^3/(c_1+2d_3)N < n_1/c_1n < d_3$.
Applying Theorem \ref{kahn}, there is a matching in $R'$
of size at least $|Z'|/3 - \frac{1}{6}\delta |Z'|$,
i.e.\ at most $\frac{1}{2}\delta|Z'|$ vertices are not covered by the matching.
\COMMENT{$M$ clashes with marked}

We can use an edge $S$ of the matching as follows.
Consider the partition $P^*_S$ of $K(V)_S$ by weak equivalence,\index{weakly equivalent}
i.e.\ we have a cell of $P^*_S$ corresponding to each triad\index{triad}
of consistent bipartite graphs indexed by $S$.
The cells lying over $S'= \cup_{i \in S} C_{i,S_i}$ give a partition of $K(S')$,
which we denote by $C^{S,1} \cup \cdots \cup C^{S,a_S}$.\index{$C^{S,i}$}
Since $P$ is $a$-bounded we have $a_S \le a^3$.
Furthermore, since $P$ is equitable, the triangle counting lemma (\ref{eq:tri})
gives $|C^{S,i}|=(1\pm\eps')|C^{S,j}|$ for any $1\le i,j \le a_S$.
Now at most $2\nu^{1/4}|K(S')_S|$ triples of $G'[S']_S$
can lie in cells $C^{S,i}$ with $|M \cap C^{S,i}| > \nu^{1/4}|C^{S,i}|$;
otherwise we would have at least $(1-\eps')2\nu^{1/4}a_S$ such cells,
giving $|M[S']| \ge (1-\eps')2\nu^{1/4}a_S \cdot \nu^{1/4}(1-\eps')|K(S')_S|/a_S
> (1-3\eps')2\sqrt{\nu}|K(S')_S|$, contradicting $S \in R$.
Since $|G'[S']_S| > d_3|K(S')_S|$ and $\nu \ll d_3$,
more than $\frac{1}{2}d_3|K(S')_S|$ triples of $G'[S']_S$
lie in cells $C^{S,i}$ with $|M \cap C^{S,i}| < \nu^{1/4}|C^{S,i}|$,
so we can choose such a cell $C^{S,i}$ with
$|G' \cap C^{S,i}| > \frac{1}{2}d_3|C^{S,i}|$.

Fix such a cell $C^{S,i}$ for each matching edge $S$
and let $G^S$ be the associated cell complex of $G'[P]$,
i.e.\ $G^S_S = G' \cap C^{S,i}$ and for $S' \subn S$,
$G^S_{S'}$ is the cell of $P_{S'}$ underlying $C^{S,i}$.
Then $G^S$ is $\eps$-regular. Also, $|G^S_S|>\frac{1}{2}d_3|C^{S,i}|$,
so writing $M^S = M \cap G^S$ we have
$|M^S| < \nu^{1/4}|C^{S,i}| < \nu^{1/5}|G^S_S|$.
At this stage, if we were satisfied with an $F$-packing covering all but
$o(n)$ vertices, we could just repeatedly remove copies of $F$ from each $G^S \sm M^S$.%
\footnote{
This could be achieved using the counting lemma to count copies of $F$ in $G^S$
and the `extension lemma' to bound the number of copies 
of $F$ using an edge in $M^S$. Alternatively one could start the proof with
the `regular decomposition lemma' instead of the regular approximation lemma,
then find $F$ using the sparse counting lemma.
}
However, we want to cover all but at most $C$ vertices,
so we will apply the blow-up lemma, using the black box form in the
previous subsection. By Lemma \ref{3robust} we can delete at most $2\nu^{1/15} n_1$
vertices from each cluster so that each $(G^S,M^S)$ becomes a pair $(G^{\sharp S},M^{\sharp S})$
such that $J^S = G^{\sharp S} \sm M^{\sharp S}$ is $1/2$-robustly $|Y|^2$-universal.

Now we gather together all the removed vertices into an exceptional set $A_0$.
This includes at most $a!$ vertices removed at the start of the proof,
at most $\nu^{1/5}n$ vertices in parts indexed by $Z \sm Z'$,
at most $\frac{1}{2}\delta n$ vertices in parts not covered by the matching,
and at most $2\nu^{1/15} n$ vertices deleted in making the pairs robustly universal.
Therefore $|A_0| < \delta n$. For convenient notation we denote the
clusters of $G^{\sharp S}$ by $A_{S,1}$, $A_{S,2}$, $A_{S,3}$.
Thus $A_0$ and $A_{S,1}$, $A_{S,2}$, $A_{S,3}$ for matching edges $S$
partition the vertex set of $G$. To cover the vertices of $A_0$ by
copies of $F$ we need the following claim.

\nib{Claim.} Any vertex $v$ belongs to at least $\beta n^{|Y|-1}$ copies of $F$ in $G$.

\nib{Proof.}
Let $\Phi$ be the set of all pairs $(ab,T)$ such that
$ab \in G(v)$ and $T \in \binom{G(ab) \sm v}{|Y_3|-1}$.
There are $|G(v)| > (1-\eps)c_1 n^2$ choices for $ab$,
and for each $ab$ the minimum degree property gives
at least $\binom{c_2n-1}{|Y_3|-1}$ choices for $T$.
Therefore $|\Phi| \ge (1-\eps)c_1 n^2 \binom{c_2n-1}{|Y_3|-1}$.
Let $\Psi$ be the set of all $T \in \binom{V(G) \sm v}{|Y_3|-1}$
such that there are at least $\frac{1}{3}c_1c_2^{|Y_3|-1} n^2$
pairs $ab \in G(v)$ with $(ab,T) \in \Phi$. Then
$|\Phi| < |\Psi|n^2 + \binom{n}{|Y_3|-1} \cdot \frac{1}{3}c_1c_2^{|Y_3|-1} n^2$,
so
$$|\Psi| > (1-\eps)c_1 \binom{c_2n-1}{|Y_3|-1}
- \binom{n}{|Y_3|-1} \cdot \frac{1}{3}c_1c_2^{|Y_3|-1} >
\frac{1}{3}c_1c_2^{|Y_3|-1} \binom{n}{|Y_3|-1}.$$
For each $T$ in $\Psi$, since $\alpha \ll c_1,c_2$, Theorem \ref{erdos} implies that
the sets $ab \in G(v)$ with $(ab,T) \in \Phi$ span at least
$\alpha n^{|Y_1|+|Y_2|}$ copies of $K_{|Y_1|,|Y_2|}$.
Each of these gives a copy of $F$ containing $v$ when we add $T \cup v$.
Summing over $T$ in $\Psi$ and dividing by $|Y|!$ (a crude estimate for overcounting)
we obtain (since $\beta \ll \alpha$)
at least $\beta n^{|Y|-1}$ copies of $F$ containing $v$. \qed

Next we randomly partition each set $A_{S,j}$ as $A'_{S,j} \cup A''_{S,j}$,
each vertex being placed independently into either class with probability $1/2$.
The reason for this partition is that, as in the proof of Theorem \ref{2pack},
we will be able to use the sets $A'_{S,j}$ when covering the vertices in $A_0$,
whilst preserving the vertices in $A''_{S,j}$ so as to maintain super-regularity.
Theorem \ref{azuma} gives the following properties with high probability:
\begin{enumerate}
\item
$|A'_{S,j}|$ and $|A''_{S,j}|$ are $|A_{S,j}|/2 \pm n^{2/3}$ for every $S$ and $j$,
\item
for every $S$ and $j$ and each vertex $v \in A_{S,j}$,
letting $\{ T_i \}_{i \ne j}$ denote the singleton classes of $G^{\sharp S}(v)$,
$|A'_{S,j} \cap T_i|$ and  $|A''_{S,j} \cap T_i|$ are $|T_i|/2 \pm n^{2/3}$ for $i\ne j$, and
\item
for any vertex $v$ of $G$, there are at least $\gamma n^{|Y|-1}$ copies of $F$
in which all vertices, except possibly $v$, are in $\cup_{S,j} A'_{S,j}$.
\end{enumerate}
In fact, the first two properties are simple applications of Chernoff bounds
(in which the martingale is just a sum of independent variables).
For the third property we use a vertex exposure martingale.
\COMMENT{just say Lipschitz?}
Fix $v$ and let $Z$ be the random variable which is the number of copies of $F$
in which all vertices, except possibly $v$, are in $\cup_{S,j} A'_{S,j}$.
Since $|A_0| < \delta n$, the Claim gives $\mb{E}Z > (1/2)^{|Y|-1}(\beta - \delta)n^{|Y|-1}$.%
\COMMENT{
By the claim there are more than $(\beta - \delta)n^{|Y|-1}$ such copies
(since $|A_0| < \delta n$), so $\mb{E}Z > (1/2)^{|Y|-1}(\beta - \delta)n^{|Y|-1}$.
}
Order the vertices of $\cup_{S,j} A_{S,j}$ as $v_1,\cdots,v_{n'}$,
where $n' > (1-\delta)n$, and define
the random variable $Z_i$ as the conditional expectation of $Z$
given whether $v_{i'}$ is in $A'_{S,j}$ or $A''_{S,j}$ for $i' \le i$.
Then $Z_0 = \mb{E}Z$ and $Z_{n'}=Z$. Also $|Z_i - Z_{i-1}| < n^{|Y|-2}$,
using a crude upper bound on the number of copies of $F$ containing $v$ and
some other vertex $v_j$. Now by Theorem \ref{azuma} we have
$\mb{P}(Z < \gamma n^{|Y|-1}) < \mb{P}(|Z_{n'}-Z_0| > 2^{-|Y|}\beta n^{|Y|-1})
< e^{-\beta^3 n}$.

Now we cover $A_0$ by the following greedy procedure.
Suppose we are about to cover a vertex $v \in A_0$.
We consider a cluster to be {\em heavy} if we have covered more than $\gamma n_1$
of its vertices. Since $|A_0|< \delta n$ we have covered at most \index{heavy}
$|Y|\delta n$ vertices by copies of $F$, so there are at most
$|Y|\delta n / \gamma n_1 < \frac{1}{2}\gamma ra_1$ heavy clusters.
As shown above, there at least $\gamma n^{|Y|-1}$ copies of $F$
that we can use to cover $v$. At most $\frac{1}{2}\gamma ra_1 n_1 n^{|Y|-2}
< \frac{1}{2}\gamma n^{|Y|-1}$ of these use a heavy cluster,
so we can cover $A_0$ while avoiding heavy clusters.

Next we restrict to the vertices not already covered by the copies of $F$
covering $A_0$, where we will use robust universality to finish the packing.
Recall that each $J^S$ is $1/2$-robustly $|Y|^2$-universal.
By properties (1) and (2) of the partitions $A_{S,j} = A'_{S,j} \cup A''_{S,j}$,
on restricting to the uncovered vertices we obtain $J'{}^S$ that
satisfies conditions (i) in Definition \ref{def-3robust}.
(Property (i) of Theorem \ref{3robust} and $d^\sharp \gg 1/n$
shows that the $\pm n^{2/3}$ errors are negligible.)
Also, any $F$-packing has maximum degree less than $|Y|^2$,
so satisfies condition (ii) in Definition \ref{def-3robust}.
Thus we can assume that $J'{}^S$ is complete,
in that we can embed any $F$-packing in $J'{}^S$,
subject only to the constraints given by the sizes
of the uncovered parts of each cluster.
\COMMENT{Komlos:
1. limited refs, make sure I know the proofs!
2. Komlos says at most $k|Y|$ total uncovered, but I don't see it.
}

Now it is not hard to finish the proof with a slightly messy ad hoc
argument, but we prefer to use the elegant argument of K\'omlos \cite[Lemma 12]{Ko}.
Denote the classes of $J'{}^S$ by $J_1, J_2, J_3$.
Since we avoided heavy clusters we have
$(1-2\gamma)n_1 \le |J_i| \le n_1$ for $1 \le i \le 3$.
Let $P^3 = \{\alpha \in [0,1]^3: \alpha_1 \le \alpha_2 \le \alpha_3,
\sum_{i=1}^3 \alpha_i = 1\}$. We can associate a `class vector'
$\alpha(X) \in P^3$ to a $3$-partite set $X = X_1 \cup X_2 \cup X_3$
by $\alpha(X)_i = |X_{\sigma(i)}|/|X|$, for some permutation $\sigma \in S_3$
chosen to put the classes in increasing order by size.
For $\alpha, \beta \in P^3$ write $\alpha \prec \beta$ if
$\alpha_1 \le \beta_1$ and $\alpha_1+\alpha_2 \le \beta_1 + \beta_2$.
Since the classes of $F$ are not all of equal size and
$\gamma \ll 1/|Y|$ we have $\alpha(F) \prec \alpha(J)$.
By a theorem of Hardy, Littlewood and P\'olya, this implies that
there is a doubly stochastic matrix $M$ such that
$\alpha(J) = M\alpha(F)$. By Birkhoff's theorem $M$ is
a convex combination of permutation matrices $M = \sum_i \lambda_i P_i$,
$\sum \lambda_i = 1$. Thus we can write the class vector of $J$ as
$\alpha(J) = \sum_i \lambda_i P_i\alpha(F)$, which is
a convex combination of the permutations of the class vector of $F$.
In fact, although the constant is not important,
since $P^3$ has dimension $3$, we can apply
Carath\'eodory's theorem
to write $\alpha(J) = \sum_{i=1}^3 \mu_i P_i\alpha(F)$,
$\sum_{i=1}^3 \mu_i=1$ as a convex combination using only $3$
permutations of $\alpha(F)$.%
\footnote{This last remark is attributed to Endre Boros in \cite{Ko}.}
Finally, to pack copies of $F$ in $J'{}^S$ we can use
$\lfloor \mu_i|J|/|Y| \rfloor$ copies of $F$ permuted
according to $P_3$, for $1 \le i \le 3$. At most $3|Y|$ vertices
of any $J_i$ are left uncovered because of the rounding,
so in total at most $C=3|Y|ra_1$ vertices will remain uncovered.
This completes the proof. \qed

As for Theorem \ref{2pack}, we needed to assume that not all part sizes of $F$ are equal
and we could not expect to cover all vertices. Some assumption on the degrees of pairs was
convenient, as without it a nearly regular $3$-graph can have some vertices
that do not belong to any copies of $F$. For example, let $G_0$ be a tripartite $3$-graph
on $V = V_1 \cup V_2 \cup V_3$ with $|V_1|=|V_2|=|V_3|=n_0/3$ such that
every vertex $v$ has degree $|G_0(v)| = (1 \pm \eps)c_1n_0^2$. Form $G$ from $G_0$ by adding
new vertices $v_1,\cdots,v_t$ where $t \approx (3c_1)^{-1}$ and edges so that $G(v_i)$
are pairwise disjoint graphs of size $c_1n_0^2$ contained in
$\binom{V_1}{2} \cup \binom{V_2}{2} \cup \binom{V_3}{2}$.
Then $G$ has $n=n_0+t$ vertices and $|G(v)|=(1\pm 2\eps)c_1n^2$ for every vertex $v$.
However, for every new vertex $v_i$, every pair $ab \in G(v_i)$ is only contained
in the edge $v_iab$, so $v_i$ is not contained in any $K_{2,2,3}$ (say).
This example does not show that the assumption on pairs is necessary, as we can still cover
all but $t$ vertices, but it at least indicates that it may not be so easy to remove the assumption.
For simplicity we assumed that every pair has many neighbours, but it is
clear from the proof that this assumption can be relaxed somewhat.
The bottleneck is the Claim, which can be established under the weaker assumption
that for every vertex $v$ there are at least $n^{2-\theta}$ pairs $ab$ in $G(v)$
with $|G(ab)|>n^{1-\theta}$, for some $\theta>0$ depending on $F$.

Finally, we remark that one can apply the general hypergraph blow-up lemma in
the next section and the same proof to obtain the following result (we omit the details).
For any $k$-partite $k$-graph $F$ in which not all part sizes are equal and $0<c_1,c_2<1$
there is a real $\eps>0$ and positive integers $C,n_0$ such that if
$G$ is an $k$-graph on $n > n_0$ vertices $V$ such that
every vertex $v$ has degree $|G(v)| = (1 \pm \eps)c_1n^{k-1}$ % was $\binom{n-1}{2}$
and every $(k-1)$-tuple $S$ of vertices has degree $|G(S)| > c_2n$
then $G$ contains an $F$-packing that covers all but at most $C$ vertices.

\section{General hypergraphs}

In this section we present the general blow-up lemma.
Besides working with $k$-graphs for any $k \ge 3$,
we will introduce the following further generalisations:
\begin{itemize}
\item[(i)] Restricted positions: a small number of sets in $H$ may be \index{restricted positions}
constrained to use a certain subset of their potential images in $G$
(provided that these constraints are regular and not too sparse).
\item[(ii)] Complex-indexed complexes: a structure that provides greater\index{complex-indexed}
flexibility, in particular the possibility of embedding spanning
hypergraphs (such as Hamilton cycles).
\end{itemize}

We divide this section into five subsections organised as follows.
The first subsection contains various definitions needed for the general case,
some of which are similar to those already given for $3$-complexes
and some of which are new. In the second subsection we state the general
blow-up lemma and the algorithm that we use to prove it.
The third subsection contains some properties of hypergraph regularity,
analogous to those proved earlier for $3$-graphs.
We give the analysis of the algorithm in the fourth subsection,
thus proving the general blow-up lemma. Since much of the analysis
is similar to that for $3$-graphs we only give full details for those
aspects of the general case that are different.
The last subsection contains the general cases of the lemmas
to be used in applications of the blow-up lemma, namely Lemmas \ref{3del},
\ref{3super-restrict} (super-regular restriction) and \ref{3robust} (robust universality).

\subsection{Definitions}

We start by defining complex-indexed complexes.

\begin{defn}\label{def-indexed}\index{multicomplex}\index{multi-$k$-complex}\index{$\sub$}\index{copy}
We say $R$ is a {\em multicomplex} on $[r]$ if it consists some
number of copies (possibly 0) of every $I \sub [r]$
which are partially ordered by some relation, which we denote by $\sub$,
such that whenever $I^* \in R$ is a copy of some subset $I$ of $[r]$
and $J$ is a subset of $I$, there is a unique copy $J^*$ of $J$ with $J^* \sub I^*$.
We say that $R$ is a {\em multi-$k$-complex} if $|I| \le k$ for all $I \in R$.
\COMMENT{
1. change $\le$ to $\sub$
2. allow many copies of $\es$?
}

Suppose $V$ is a set partitioned as $V = V_1 \cup \cdots \cup V_r$.
Suppose each $V_i$, $1 \le i \le r$ is further partitioned as
$V_i = \cup_{i^*} V_{i^*}$, where $i^*$ ranges over all copies of $i$ in $R$.
We say $G$ is an {\em $R$-indexed complex} on $V$
\index{R-indexed@$R$-indexed|see{complexed-indexed}}
if it consists of disjoint parts $G_I$ for $I \in R$ (possibly undefined),
such that $G_i \sub V_i$ for singletons $i \in R$,
and $G_{I^\le}:= \cup_{I' \sub I} G_{I'}$ is a complex whenever $G_I$ is defined.
We say that an $R$-indexed complex $J$ on $V$ is an {\em $R$-indexed subcomplex} of $G$
if $J_I \sub G_I$ when defined. %\index{R-indexed subcomplex@$R$-indexed subcomplex}
We say $S \sub V$ is {\em $r$-partite} if $|S \cap V_i| \le 1$ for $1 \le i \le r$.
The {\em multi-index} $i^*(S)$ of $S$ is that $I \in R$ with $S \in G_I$.
\index{multi-index}\index{$i^*(S)$}
If $S \in G$ we write $G_S = G_{i^*(S)}$ for the part of $G$ containing $S$.
\COMMENT{Formerly: If $I' \sub I$ we write $G_{I'} \sub G_I$. Or $\le$?}
\end{defn}

We emphasise that $\sub$ is more restrictive than the inclusion relation
(also denoted $\sub$) between copies considered merely as subsets of $[r]$.
To avoid confusion we never `mix' subsets with copies of subsets.
Thus, if $I \in R$ then $J \sub I$ means $J \in R$ \index{$\in$}
and $(J,I)$ is in the relation $\sub$.
Also, if $I \in R$ then $i \in I$ means $\{i\} \in R$
and $(\{i\},I)$ is in the relation $\sub$.
We also write $i \in R$ to mean $\{i\} \in R$
when the meaning is clear from the context.
As in Definition \ref{def-more}, we henceforth simplify notation
by writing $i$ instead of $\{i\}$. 

\begin{rem}
We are adopting similar notation for complex-indexed complexes as
for usual complexes for ease of discussing analogies between the two situations.
Thus we typically denote a singleton multi-index by $i$ and a set multi-index by $I$.
If we need to distinguish a multi-index from the index of which it is a copy,
we typically use the notation that $i^*$ is a copy of $i$ and $I^*$ is a copy of $I$.
\end{rem}

We illustrate Definition \ref{def-indexed} with the following example.
\COMMENT{
1. Formerly: and using concatenation to denote union....
Problem! Union of multi-indices not well-defined. If $A$, $B$ are copies of $I$, $J$
then maybe many copies of $IJ$ contain $A$ and $B$...
Mostly we use subscripts from $H$ so there is no problem.
Super-regularity needs clarification...
2. Formerly:
As in our notation for $r$-partite complexes, subscripts are to be
understood as their multi-index where appropriate; for example $d_S(G)=d_{i^*(S)}(G)$.
We will distinguish singleton multi-indices by using the
notation $i^* \in V(R)$ (rather than $i \in [r]$).
We will often use $i^* \in A$ rather than $i^* \le A$ as a more natural notation
when $i^*$ is a singleton multi-index.
The degree $d_R(i^*)$ of $i^*$ is the number of multi-indices
that contain it. As for (normal) $k$-complexes, if $R$ is a multi-$k$-complex
we write $R_= = \{A \in R: |A|=k \}$.
}

\begin{figure}
\begin{center}
\includegraphics[height=5cm]{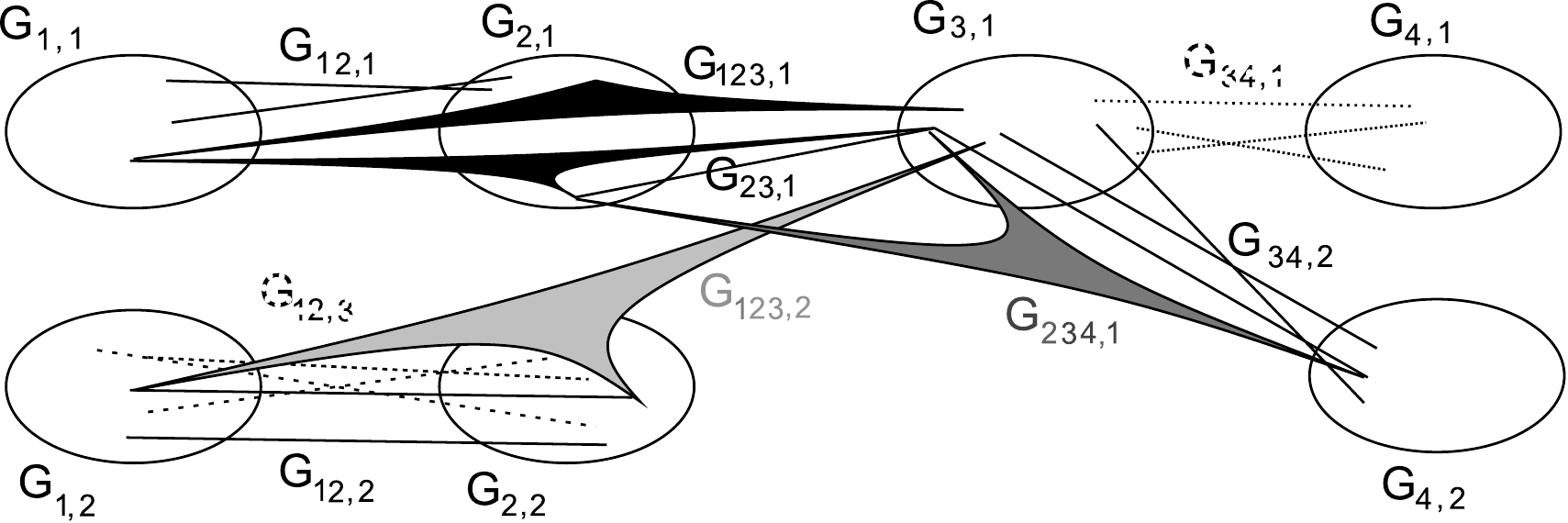}
\end{center}
\caption{A complex-indexed complex}
\label{pic-indexed}
\end{figure}

\begin{eg}\label{eg3} % before eg2 for some reason...
Figure \ref{pic-indexed} depicts an example of an $R$-indexed complex
$G$ in which $R$ is a multi-$3$-complex on $[4]$ (not all parts have
been labelled). The multi-indices have
been represented as ordered pairs $(A,t)$, where $A \in \binom{[4]}{\le 3}$
and $t$ is a number (arbitrarily chosen) to distinguish different copies of $A$.
An example of the inclusion structure is $(34,2) \sub (234,1)$,
since the intended interpretation of our picture
is that for every triple in $G_{234,1}$ its restriction to index $34$ lies in
$G_{34,2}$. Other examples are $(1,2) \sub (123,2)$ and $(2,2) \sub (123,2)$,
but $(1,1) \not\sub (123,2)$ and $(12,3) \not\sub (123,2)$, since
the intended interpretation of our picture is that there are triples in $G_{123,2}$
such that their restriction to index $12$ is a pair in $G_{1,2} \times G_{2,2}$
not belonging to $G_{12,3}$.
\end{eg}

Complex-indexed complexes arise naturally from the partition complexes
needed for regular decompositions of hypergraphs, as in Theorem \ref{ral}.
Suppose $P$ is an $r$-partite partition $k$-complex on $V$.
Recall that $C_S$ denotes the cell containing a set $S$.
We can index the cells of $P$ by a multicomplex $R$ on $[r]$, where for
each cell $C_S$ we form a copy $i^*(S)$ of its (usual) index $i(S)$.
The elements of $R$ are partially ordered by a relation $\sub$
which corresponds to the consistency relation $\le$ discussed above,
\index{consistent}\index{$\le$}
i.e.\ $i^*(S') \sub i^*(S)$ exactly when $C_{S'} \le C_S$.
One could think of $P$ as a `complete $R$-indexed $k$-complex',
in that it contains every $r$-partite set of size at most $k$,
although we remark that many partition $k$-complexes $P$ will give
rise to the same multicomplex $R$, so the phrase is ambiguous.

Note that if we do not allow sets in $R$ to have multiplicity more than $1$
then an $R$-indexed complex is precisely an $r$-partite complex.
The reason for working in the more general context is that $R$-indexed complexes
are the structures that naturally arise from an application of Theorem
\ref{ral}, so a general theory of hypergraph embedding will need to
take this into account. In particular, in order to use every vertex of $V$
we need to consider every part of the partition $P_i$ of $V_i$ for $1 \le i \le r$,
i.e.\ we need multi-index copies of each index $i$.
The multi-index copies of larger index sets are useful because
it may not be possible to choose mutually consistent cells.
To illustrate this point, it may be
helpful to consider an example of a $4$-partite $3$-complex where
each triple has constant density (say $1/10$) and there is no
tetrahedron $K_4^3$. A well-known example of R\"odl
is obtained by independently orienting each pair of vertices
at random and taking the edges to be all triples that form cyclic triangles.
(A similar example is described in \cite{NR}.)
Then we cannot make a consistent choice of cells
in any $4$-partite subcomplex so that each cell has good density.
However, by working with indexed complexes one can embed
using cells from each of the four triples,
provided that the choice of cells is locally consistent.

Much of the notation we set up for $r$-partite $3$-complexes extends
in a straightforward manner to $R$-indexed complexes
for some multi-$k$-complex $R$ on $[r]$.
Throughout we make the following replacements:
replace `$3$' by `$k$', `$r$-partite' by `$R$-indexed',
`$I \sub [r]$' by `$I \in R$', `index' by `multi-index',
$i(S)$ by $i^*(S)$, and understand $\sub$ as the partial order of multi-indices.
Thus we define $G_{I^\le} = \cup_{I' \sub I} G_{I'}$,
$G_{I^<} = \cup_{I' \subn I} G_{I'}$, etc.\ as in Definition \ref{def-more}.
We define restriction of $R$-indexed complexes as in Definition \ref{def-restrict},
and more generally composition of $R$-indexed complexes $G$ and $G'$
as in Definition \ref{def-*}: we define $(G*G')_S$ if $(G \cup G')_S$ is defined
and say that $S \in (G*G')_S$ if $A \in^* G_A$ and $A \in^* G'_A$ for any $A \sub S$.
Lemma \ref{*props} applies to $R$-indexed complexes, with the same proof.%
\index{$G*G'$}\index{restriction}\index{composition}

As in Definition \ref{def-set-restrict}, if $S \sub X$ is $r$-partite
and $I \sub i^*(S)$ we write $S_I = S \cap \cup_{i \in I} X_i$.
We also write $S_T=S_{i(T)}$ for any $r$-partite set $T$ with $i^*(T) \sub i^*(S)$.
As in Definition \ref{def-density} we write $G_I^*$
for the set of $r$-partite $|I|$-tuples $S$ with $i^*(S)=I$
such that $T \in G_T$ when defined for all $T \subn S$.
When defined we have relative densities $d_I(G)=|G_I|/|G_I^*|$
and absolute densities $d(G_I) = |G_I|/\prod_{i \in I} |V_i|$.
(Recall that $i \in I$ means $\{i\} \sub I$ according to $R$.)%
\index{relative density}\index{density}

To define regularity we adopt the reformulation using restriction notation
already discussed for $3$-complexes. For any $I \in R$
such that $G_I$ is defined we say that $G_I$ is {\em $\eps$-regular}
if for every subcomplex $J$ of $G$ with $|J^*_I| > \eps|G^*_I|$
and $J_I$ undefined we have $d_I(G[J])=d_I(G) \pm \eps$.
Note that if $|I|=0,1$ we have $G[J]_I=G_I$,
so $G_I$ is automatically $\eps$-regular for any $\eps$.
We say that $G$ is {\em $\eps$-regular}
if whenever any $G_I$ is defined it is $\eps$-regular.\index{regular}

\begin{rem}\label{rem:nhood}
Some care is needed when forming neighbourhoods in complex-indexed complexes.
The usual definition defines $G_I(v)$ and $G_{I^\le}(v)$.
However the expression $G(v)_I$ may be ambiguous, as there may be several
multi-indices $I'$ such that $I = I' \sm i$, where $i=i^*(v)$.
We avoid this ambiguous expression unless the meaning is clear from the context.
The expression $P \in G(v)$ is unambiguous: it means $Pv \in G$.
\end{rem}

Henceforth we suppose $R$ is a multi-$k$-complex on $[r]$, \index{$R$}
$H$ is an $R$-indexed complex on $X = \cup_{i \in R} X_i$\index{$X$}\index{$H$}
and $G$ is an $R$-indexed complex on $V = \cup_{i \in R} V_i$,\index{$G$}\index{$V$}
with $|V_i|=|X_i|$ for $i \in R$.
As before we want to find an embedding $\phi$ of $H$ in $G$,\index{$\phi$}
via an algorithm that considers the vertices of $X$ in some order
and embeds them one at a time. At some time $t$ in the algorithm, for each $S \in H$ there
will be some $|S|$-graph $F_S(t) \sub G_S$ consisting of those sets $P \in G_S$ \index{$F_S(t)$}
that are `free' for $S$, in that mapping $S$ to $P$ is `locally consistent' with
\index{free}\index{locally consistent}\index{mutually consistent}
the embedding so far. These free sets will be `mutually consistent',
in that $F_{S^\le}(t)  = \cup_{S' \sub S} F_{S'}(t)$ is a complex. \index{$F_{S^\le}(t)$}
With the modifications already mentioned,
we can apply the following definitions and lemmas (and their proofs)
to the general case: Definition \ref{def-update} (the update rule),\index{update rule}
Lemma \ref{consistent} (consistency),
Lemma \ref{build-update} (iterative construction) and
Lemma \ref{not-local} (localisation).\index{local}
Just as some triples were marked as forbidden for $3$-graphs,
in the general case some $k$-tuples $M$ will be marked as forbidden.\index{marked}\index{$M$}
Definition \ref{def-mark} ($M_{E^t,E}(t)$) applies in general\index{$E^t$}\index{$M_{E^t,E}(t)$}
when $E \in H$ is a $k$-tuple. 

Definition \ref{def-3preplus} is potentially ambiguous for $R$-indexed complexes
(see Remark \ref{rem:nhood}), so we adopt the following modified definition.%
\COMMENT{
** Change earlier def for consistency?
Mostly unused, just eg and claims $A_{i,0}$...?
Or maybe this is clearer as it emphasises the ambiguity...
}

\begin{defn}\label{def-preplus}\index{$G^{I_v}$}
Suppose $G$ is an $R$-indexed complex on $V = \cup_{i \in R} V_i$,
$i \in R$, $v \in G_i$ and $I$ is a submulticomplex of $R$.
We define $G^{I_v} = G[\cup_{S \in I} G_S(v)]$.
\end{defn}

Now we give the general definition of super-regularity,
which is very similar to that used for $3$-graphs.

\begin{defn}\label{def-super} (Super-regularity)\index{super-regular}
Suppose $R$ is a multi-$k$-complex on $[r]$,
$G$ is an $R$-indexed complex on $V = \cup_{i \in R} V_i$
and $M \sub G_= := \{P \in G: |P|=k\}$.\index{$G_=$}
We say that $(G,M)$ is $(\eps,\eps',d_a,\theta,d)$-{\em super-regular} if
\COMMENT{
** Rewritten to avoid union of multi-indices, which are not well-defined.
Specified subscripts to avoid ambiguity, e.g. $G(v)_{S'\sm i}$ is not well-defined
as $S'\sm i$ could be covered by several $S''$ with the same index as $S'$.
Are we okay later using $H$-subscripts for densities but omitting subscripts `inside'??
}

(i) $G$ is $\eps$-regular, and if $G_S$ is defined then
$d_S(G) \ge d$ if $|S|=k$ or $d_S(G) \ge d_a$ if $|S|<k$,

(ii) if $G_S$ is defined, $i \in S$ and $v \in G_i$ then
$|M_S(v)| \le \theta|G_S(v)|$ if $|S|=k$ and
$G_{S^\le}(v)$ is $\eps'$-regular with
$d_{S'\sm i}(G_{S^\le}(v)) = (1 \pm \eps')d_{S'\sm i}(G)d_{S'}(G)$
for $i \subn S' \sub S$,

(iii) for every submulticomplex $I$ of $R$,
if $G_S$ is defined, $i \in R$ and $v \in G_i$, then
$|(M \cap G^{I_v})_S| \le \theta|G^{I_v}_S|$ if $|S|=k$,
and $G^{I_v}$ is $\eps'$-regular with densities (when defined)
$$ d_S(G^{I_v}) = \left\{ \begin{array}{ll}
(1 \pm \eps')d_{S}(G)d_{T}(G) & \mbox{ if } \es \ne S = T\sm i, T\in I, G_T \mbox{ defined,}\\
(1 \pm \eps')d_{S}(G) & \mbox{ otherwise.} 
\end{array}
\right.$$
\end{defn}

We remark that the parameters in Definition \ref{def-super}
will satisfy the hierarchy  $\eps \ll \eps' \ll d_a \ll \theta \ll d$.
Thus we work in a dense setting with regularity parameters much smaller than density parameters.
Note that we have two density thresholds $d_a$ and $d$,
where $d_a$ is a lower bound on the $S$-densities when $|S|<k$
and $d$ is a lower bound on the $S$-densities when $|S|=k$.
The marking parameter $\theta$ again lies between these thresholds.
Next we will formalise the notation $d_S(F(t))=d_S(F_{S^\le}(t))$
used for $3$-complexes by defining $F(t)$ as an object in its own right:
a {\em complex-coloured complex}.\index{complex-coloured}

\begin{defn}\label{def-coloured}
Suppose $R$ is a multi-$k$-complex on $[r]$,
$H$ is an $R$-indexed complex on $X = \cup_{i \in R} X_i$
and $G$ is an $R$-indexed complex on $V = \cup_{i \in R} V_i$.
An {\em $H$-coloured complex $F$ (in $G$)}
\index{H-coloured@$H$-coloured|see{complex-coloured}}
consists of $|S|$-graphs $F_S \sub G_S$
such that $F_{S^\le} = \cup_{S' \sub S} F_{S'}$ is a complex for $S \in H$.
We allow $F_S$ to be undefined for some $S \in H$.
We say that an $H$-coloured complex $J$ is an
{\em $H$-coloured subcomplex} of $F$ if $J_S \sub F_S$ when defined.
The {\em restriction} $F[J]$ is the $H$-coloured complex \index{restriction}
with $F[J]_S = F_{S^\le}[J_{S^\le}]_S$.
When $F(t)$ is an $H$-coloured complex at time $t$
we use $F(t)_S$ and $F_S(t)$ interchangeably.\index{$F(t)_S$}\index{$F_S(t)$}
We let $G$ also denote the $H$-coloured complex $F$ in $G$
such that $F_S=G_S$ for all $S \in H$.
When $I$ is a subcomplex of $H$ we let $F_I$ denote the\index{$F_I$}
$H$-coloured complex that consists of $F_S$ for every $S \in I$
and is otherwise undefined.
\COMMENT{
1. This is a more general structure than a complex, as for any $I^* \in R$
that is a copy of some $I \sub [r]$, each $S \in H_{I^*}$ gives rise to a
sub-$|I|$-graph of the same part $K(V)_{I^*} = K(V)_I = \{P \sub V: i(P)=I\}$.
2. Omit(?) spanning, $\le$, $<$, $=$, $F(t)(v)$,
$F(t)^v$ (needs to include $F(t)(v)$ to be complex)
3. Is $F_I$ used?
4. Formerly:
The {\em complete $H$-coloured complex} $K_H(V)$ on $V$ %\index{complete}
has parts $K_H(V)_S = K(V)_{i(S)}$ for all $S \in H$.
(Note that we do mean the usual index $i(S)$ here.)
[No multi-indices defined for $|S|\ge 2$. Nicer to use $G$!]
5. Formerly: For $P \in G$ the {\em neighbourhood} $H$-coloured complex%\index{neighbourhood}
$F(P)$ has $F(P)_S = F_S(P)$ defined when $S\in H(P)$.%\index{$F(P)$}
Problem: should be $F_{IP}(P)$ for some $I$ with the same multi-index as $S$... not wdf!
}
\end{defn}

We use the terminology `coloured' in analogy with various combinatorial questions
involving hypergraphs in which each set can be assigned several colours.
In our case a set $E \in G$ is assigned as colours
all those $S \in H$ for which $E \in F_S$.
Note also that if we had the additional property that the parts $F_S$, $S \in H$
were mutually disjoint, i.e.\ any set in $G$ has at most one colour from $H$,
then $F$ would be an $H$-indexed complex.
(We make this comment just to illustrate the definition:
we will not have cause to consider any $H$-indexed complexes.)

In the proof of Theorem \ref{3blowup} there were several places
where we divided the argument into separate cases.
This will not be feasible for general $k$, so we will introduce
some more notation, which may at first appear somewhat awkward,
but will repay us by unifying cases into a single argument.

\begin{defn}\label{def-plus}
Suppose $R$ is a multi-$k$-complex on $[r]$,
$H$ is an $R$-indexed complex on $X = \cup_{i \in R} X_i$
and $G$ is an $R$-indexed complex on $V = \cup_{i \in R} V_i$.
Fix $x \in X$.\index{$(\cdot)^c$}\index{$(\cdot)^+$}%
\begin{itemize}
\item[(i)]
We define a multi-$(k+1)$-complex $R^+$ on $[r+1]$ as follows.
There is a single multi-index copy of $r+1$, also called $r+1$.
Suppose $x\in X_{i^*}$, where $i^* \in R$ is a copy of some $i \in [r]$.
Consider $I^* \in R$ such that $I^*$ is a copy of $I \sub [r]$.
If $i^* \in I^*$ we let $I^{*c}$ be a copy of $(I\sm i)\cup\{r+1\}$.
If $i^* \notin I^*$ and $I^* \ne J\sm i^*$ for any $J \in R$
we let $I^{*+}$ be a copy of $I\cup\{r+1\}$.
We extend $\sub$ by the rules (when defined)
$J \sub I^{*c}$ for $i^* \notin J \sub I^*$;
$J^c \sub I^{*c}$ for $i^* \in J \sub I^*$; 
$J^c \sub I^{*+}$ for $i^* \in J$ and $J\sm i^* \sub I^*$;
$I^* \sub I^{*+}$; and $J^+\sub I^{*+}$ for $J \sub I^*$.
\COMMENT{
1. Need uniqueness. [Formerly: Only use $I^{*+}$ when $I^*i^* \notin R$, 
but might as well define it generally... no!]
2. Redundant: $\{r+1\} \sub I^{*c}$, $\{r+1\} \sub I^{*+}$
}
\item[(ii)]
Let $X_{r+1} = \{x^c\}$ consist of a single new vertex
that we consider to be a copy of $x$.\index{copy}%
\footnote{
To avoid confusion we should point out that the use of `copy' here
is different to the sense in which multi-indices are copies of normal indices.
}
Let $H^+$ be the $R^+$-indexed complex
$H \cup \{Sx^c: S \in H\}$ on $X^+ = X \cup X_{r+1}$,
where $i^*(Sx^c)$ is $(i^*(Sx))^c$ if $Sx \in H$ or $i^*(S)^+$ if $Sx\notin H$.
If $Sx \in H$ we write $(Sx)^c=Sx^c$. If $x\notin S\in H$ we write $S^c=S$.
If $S\in H\sm H(x)$ we write $S^+=Sx^c$. Note that $i^*((Sx)^c)=(i^*(Sx))^c$
and $i^*(S^+)=(i^*(S))^+$.
\item[(iii)]
Let $V_{r+1}$ be a new set
of vertices disjoint from $V$ having the same size as $V_x=V_{i^*}$.
We think of $V_{r+1}$ as a copy of $V_x$, in that for each $v \in V_x$
there is a copy $v^c \in V_{r+1}$. Let $G^+$
be the $R^+$-indexed complex
$G \cup \{Pv^c: P \in G, v \in V_x\}$ on $V^+ = V \cup V_{r+1}$,
where $i^*(Pv^c)$ is $i^*(Pv)^c$ if $Pv \in G$ or $i^*(P)^+$ if $Pv\notin G$.
If $Pv \in G$ we write $(Pv)^c=Pv^c$. If $v\notin P\in G$ we write $P^c=P$.
If $P \in G\sm G(v)$ we write $P^+=Pv^c$. Note that $i^*((Pv)^c)=(i^*(Pv))^c$
and $i^*(P^+)=(i^*(P))^+$.
\item[(iv)]
If $I \sub H$ or $I \sub G$ we write $I^c = \{A^c: A \in I\}$.
\item[(v)]
Suppose $F$ is an $H$-coloured complex in $G$.
We define $H^+$-coloured complexes $F^c = \bigcup_{S\in H} F_S^c$ and
$F^+ = \bigcup_{S\in H} (F_S \cup F_S^c) = F \cup F^c$.
Suppose $I$ is an $R$-indexed subcomplex of $H$.
The {\em plus complex} is $F^{I+x}=G^+[F \cup F^+_{I^c}]$.%
\index{plus complex}\index{$F^{I+x}$}
\end{itemize}
\end{defn}

We give the following example to illustrate Definition \ref{def-plus}.

\begin{eg}\label{eg2}
As in Example \ref{eg1},
suppose that $H$ and $G$ are $4$-partite $3$-complexes,
that we have $4$ vertices $x_i \in X_i$, $1 \le i \le 4$
that span a tetrahedron $K_4^3$ in $H$,
that we have the edges $x'_1 x'_2 x_3$ and $x'_1 x'_3 x'_4$
and all their subsets for some other $4$ vertices $x'_i \in X_i$, $1 \le i \le 4$,
and that there are no other edges of $H$ containing any $x_i$ or $x'_i$, $1 \le i \le 4$.
We can think of $H$ and $G$ as $R$-indexed complexes with $R = \binom{[4]}{\le 3}$,
i.e.\ we have one copy of each subset of $[4]$ of size at most $3$.

We work through Definition \ref{def-plus}, setting $x=x_1$.
$R^+$ is the subcomplex of $[5]$ that contains $R = \binom{[4]}{\le 3}$
and all sets $S \cup 5$ with $S \in \binom{[4]}{\le 3}$.
$H^+$ is the $R^+$-indexed complex on $X^+ = X \cup X_5$
where $X_5 = \{x_1^c\}$ and $H^+$ consists of all sets $S$ and $Sx_1^c$ with $S \in H$.
$G^+$ is the $R^+$-indexed complex on $V^+ = V \cup V_5$
where $V_5 = \{v^c:v\in V_1\}$ and $G^+$ consists of all sets $S$ and $Sv^c$
with $S \in G$ and $v\in V_1$. Note that $H^+$ and $G^+$
are $5$-partite $4$-complexes.

Let $F(0)$ be the $H$-coloured $3$-complex in which $F(0)_S=G_S$ for all $S \in H$.
We will describe the plus complex $F(0)^{H+x_1}$.
For any $S \in H$ we have $F(0)^{H+x_1}_S = F(0)_S = G_S$ by equation (\ref{eq:res}). 
Similarly, for any $S \in H(x_1)$ we have $F(0)^{H+x_1}_{Sx_1^c} = (F(0)^+_{H^c})_{Sx_1^c}=G_{Sx_1}^c$.
(Recall that $S \in H(x_1)$ iff $Sx_1 \in H$.)
For example, $F(0)^{H+x_1}_{x_2x_1^c} = G_{x_2x_1}^c = \{v_2v_1^c:v_2v_1 \in G_{x_2x_1}\}$.
If $S \in H \sm H(x_1)$ then $F(0)^{H+x_1}_{Sx_1^c}$ consists of all
$P \in G^+_{Sx_1^c}$ such that for all $S' \sub S$ we have
$P_{S'} \in G_{S'}$ and $P_{S'x_1^c} \in G_{S'x_1}^c$ if $S' \in H(x_1)$.
For example, $F(0)^{H+x_1}_{x_1x_2x_3x_1^c}$ consists of all $4$-tuples
$v'_1 v_2 v_3 v_1^c$ where $v_1v_2v_3$ and $v'_1v_2v_3$ are in $G_{123}$.

Now suppose, as in Example \ref{eg1}, that we start the embedding by mapping
$x_1$ to some $v_1 \in V_1$. Let $F(1)$ be the $H$-coloured $3$-complex given
by the update rule: this is worked out in Example \ref{eg1} and
formally defined in Definition \ref{def-update}. We can also describe it
using the plus complex. For example, we saw that $F_{x_2x_3x_4}(1)=F(1)_{x_2x_3x_4}$
consists of all triples in $G_{234}$ that form a triangle in the neighbourhood of $v_1$,
i.e.\ form a tetrahedron with $v_1$. We can write this as
$F(1)_{x_2x_3x_4} = F(0)^{H+x_1}(v_1^c)_{x_2x_3x_4}$,
as by definition $v_2v_3v_4v_1^c \in F(0)^{H+x_1}_{x_2x_3x_4x_1^c}$
exactly when $v_1v_2v_3v_4$ is a tetrahedron.
Another example from Example \ref{eg1} is that $F(1)_{x'_1x'_2x_3}$
consists of all triples $P \in G_{123}$ not containing $v_1$
such that $P_3 = P \cap V_3$ is a neighbour of $v_1$.
We can write this as
$F(1)_{x'_1x'_2x_3} = F(0)^{H+x_1}(v_1^c)_{x'_1x'_2x_3} \sm v_1$,
as by definition $v'_1v_2v_3v_1^c \in F(0)^{H+x_1}_{x'_1x'_2x_3x_1^c}$
exactly when $v'_1v_2v_3 \in G_{123}$ and $v_1v_3 \in G_{13}$.
\end{eg}

Note that the plus complex in Definition \ref{def-plus}(v) depends on $H$,
but we suppress this in the notation, as $H$ will always be clear from the context.
We have also suppressed $x$ from the notation in $R^+$ (etc).
One should note that the definition of $F^+$ is rather different than $G^+$ and $H^+$.
To clarify this definition, we note that $F^+_{H^c}=F^c$, and also that $F^+_{I^c}=F^c_{I^c}$, 
so one could also write $F^{I+x}=G^+[F \cup F^c_{I^c}]$.%
\COMMENT{
Refs comments. Actually, why do we bother with $F^+$?
Maybe I needed it with a former definition of $F^c$ that was not a complex?
}
To justify the definition of the plus complex, 
we note that $F^c$ and $F^+$ are $H^+$-coloured complexes in $G^+$,
as $F^c_{(Sx^c)^\le}$ is the copied version of the complex $F_{Sx^\le}$ for any $S \in H$.
Then $F \cup F^+_{I^c}$ is an $H^+$-coloured complex in $G^+$,
with $(F \cup F^+)_S = F_S$ for $S \in H$
and $(F \cup F^+_{I^c})_{S^c} = F_S^c$ for $S \in I$.
Regarding $G^+$ as an $H^+$-coloured complex in $G^+$,
the plus complex $F^{I+x}$ is a well-defined $H^+$-coloured complex in $G^+$.
By equation (\ref{eq:res}) we have
\[ F^{I+x}_S = F_S \mbox{ for } S \in H \mbox{ and }
 F^{I+x}_{S^c} = F_S^c \mbox{ for } S \in I.\]
As noted in Remark \ref{rem:nhood}, some care must be taken to avoid ambiguity
when defining neighbourhoods. We adopt the following convention:
\[F^{I+x}(v^c)_S = F^{I+x}_{Sx^c}(v^c) \mbox{ for } S \in H.\]

Note that we set $I=H$ in Example \ref{eg2}, and indeed this is the typical application
of this definition. The reason for allowing general $I$ is for
proving the analogue of Lemma \ref{3-x-to-v} in the general case.
Next we prove a lemma which confirms that the plus complex
does describe the update rule in general,
when we map $x$ to $\phi(x)=y$ at time $t$.

\begin{lemma}\label{plus-update}
If $x \notin S \in H$ then $F(t)_{S^\le} = F(t-1)^{H+x}(y^c)_{S^\le} \sm y$.
\COMMENT{
1. Didn't use $Sx$-subscript, but this is well-defined.
2. Allow r-partite S?}
\end{lemma}

\nib{Proof.} By Definition \ref{def-update} we have
$F_{S^\le}(t)=F_{S^\le}(t-1)[F_{S.x^\le}(t-1)(y)] \sm y$.
Thus $P \in F_S(t)$ exactly when $P \in F_S(t-1)$, $y \notin P$
and $P_{S'} y \in F_{S'x}(t-1)$ for all $S' \sub S$ with $S' \in H(x)$.
Since $P_{S'} y \in F_{S'x}(t-1)$ $\Lra$ $P_{S'} y^c \in F_{S'x}(t-1)^c$,
Definition \ref{def-plus} gives $P \in F_S(t)$ exactly when $y \notin P$
and $Py^c \in F(t-1)^{H+x}_{Sx^c}$. \qed

We also note that the plus complex can describe the construction $G^{I_v}$
in a similar manner. For the following identity we could take $H=R$,
considered as an $R$-indexed complex with exactly
one set $S \in H_S$ for every $S \in R$.
Actually, the identity makes sense for any $R$-indexed complex $H$,
when we interpret each $R$-indexed complex as an $H$-coloured complex
as in Definition \ref{def-coloured},
i.e.\ $G^{I_v}$ is the $H$-coloured complex $F$ in $G^{I_v}$
with $F_S = G^{I_v}_S$ for all $S \in H$, and similarly for $G^{I+v}(v^c)$.
\COMMENT{
Formerly: Problem with $G^{I_v}$ def... what is $G(v)_S$?
Def $G^+$ ok? Construction implicitly has one new multi-index with sets $Pv^c$
for each old multi-index with set $P$, so no ambiguity. But does it
amalgamate all possibilities, and is this okay?
}

\begin{lemma}\label{pre-to-plus}
$G^{I+x}(v^c)_S=G^{I_v}_S$ for $S \in H$ and $v \in G_x$.
\end{lemma}

\nib{Proof.} We have $P \in G^{I_v}_S$ exactly when $P \in G_S$
and $P'v \in G$ for all $P' \sub P$ with $i^*(P'v) \in I$.
Since $P'v \in G$ $\Lra$ $P'v^c \in G^c$ this is
equivalent to $Pv^c \in G^{I+x}_{Sx^c}$. \qed

\subsection{The general blow-up lemma}\label{alg}

Now we come to the general blow-up lemma. First we give a couple of definitions.
Suppose $R$ is a multi-$k$-complex on $[r]$. We write $|R|$\index{$|R|$} for the number
of multi-indices in $R$. For $S \in R$, the {\em degree}\index{degree} of $S$ is the
number of $T \in R$ with $S \sub T$.
\COMMENT{Formerly said $|R|$ counts singletons, but all seems cleaner.}

\begin{theo} \label{blowup} {\bf (Hypergraph blow-up lemma)} \index{hypergraph blow-up lemma}
Suppose that
\begin{itemize}
\item[(i)] $0 \ll 1/n \ll 1/n_R \ll \eps \ll \eps' \ll c \ll d_a \ll \theta \ll d,c',1/D_R,1/D, 1/C, 1/k$,
\COMMENT{
1. need union bound over $V(R)$ at the end so make $n$ smaller
2. formerly $1/D_R \ll d \ll c' \ll ...$ but this seems okay
}
\item[(ii)] $R$ is a multi-$k$-complex on $[r]$ of maximum degree at most $D_R$ with $|R| \le n_R$,
\item[(iii)] $H$ is an $R$-indexed complex on $X = \cup_{i \in R} X_i$
of maximum degree at most $D$, $G$ is an $R$-indexed complex on $V = \cup_{i \in R} V_i$,
$G_S$ is defined whenever $H_S$ is defined, and $n \le |X_i|=|V_i|=|G_i| \le Cn$ for $i \in R$,
\COMMENT{It is convenient to assume that $H$ and $G$ span their vertex sets.}
\item[(iv)] $M \sub G_= = \{S \in G: |S|=k\}$
and $(G,M)$ is $(\eps,\eps',d_a,\theta,d)$-super-regular,
\item[(v)] $\GG$ is an $H$-coloured complex in $G$ with $\GG_x$ defined only when $x \in X_*$,
where $|X_* \cap X_i| \le c|X_i|$ for all $i \in R$,\index{$\GG$}\index{$X_*$}
and for $S \in H$, when defined $\GG_S$ is $\eps'$-regular with $d_S(\GG) > c'd_S(G)$,
\end{itemize}
Then there is a bijection $\phi:X \to V$ with $\phi(X_i)=V_i$ for $i \in R$
such that for $S \in H$ we have $\phi(S) \in G_S$, $\phi(S) \in G_S \sm M_S$ when $|S|=k$
and $\phi(S) \in \GG_S$ when defined.
\end{theo}

We make some comments here to explain the statement of Theorem \ref{blowup}.
An informal statement is that we can embed any bounded degree $R$-indexed complex
in any super-regular marked $R$-indexed complex, even with some restricted positions.
The restricted positions are described by assumption (v): \index{restricted positions}
for some sets $S \in H$ we constrain the embedding to satisfy $\phi(S) \in \GG_S$,
for some $H$-coloured complex $\GG$ that is regular and dense
and is not defined for too many vertices.
Note that we now allow the embedding to use all $|R|$ parts of $V$,
provided that $R$ is of bounded degree. Thus this theorem could be used for
embedding spanning hypergraphs, such as Hamilton cycles.
Even in the graph case, embedding spanning subgraphs is a generalisation \index{spanning}
of the graph blow-up lemma in \cite{KSS}. A blow-up lemma for spanning subgraphs
was previously given by Csaba \cite{C} (see \cite{KKO} for another application).
Restricted positions have arisen naturally
in many applications of the graph blow-up lemma, and will no doubt be similarly
useful in future applications of the hypergraph blow-up lemma.
In particular, a simplified form of the condition (where $\GG_S$ is only
defined when $|S|=1$) is used in a forthcoming work \cite{KKMO} on embedding
loose Hamilton cycles in hypergraphs. We allow a general $H$-coloured complex $\GG$
as the proof is the same, and it would be needed for embedding general Hamilton cycles.%
\footnote{
While \cite{KKMO} and the current paper have been under review, more general results on Hamilton cycles 
have been obtained in \cite{KS} and \cite{KMO} without using the hypergraph blow-up lemma.
However, it seems most likely that more complicated embedding problems will require
the hypergraph blow-up lemma.
}
\COMMENT{
... the reader can interpret `$R$-indexed' as `$r$-partite' without any loss for the application
given in this paper. Also, the `restricted positions' form of the theorem
will not be used in our application, so the reader can ignore these for
a further simplification. We allow a restricted position condition as it arises
quite naturally in applications, and we need it in a forthcoming work \cite{KKMO}.
...It is worth noting that complex-indexed complex formulation
is more general than that in \cite{KSS} even in the case of graphs ($k=2$),
as we replace an assumed bound on $v_R$, the number of parts $X_{i^*}$, with
a bound on $D_R$, the maximum degree of $R$. This can be used to obtain
results in some circumstances where the blow-up lemma from \cite{KSS} does
not suffice, as shown by Csaba \cite{C}, who proved a similar type of result,
which was also applied in \cite{KKO}.
[** This surprised me at first. Should think about Csaba conditions to see if there are
any subtleties I'm missing. His termination conditions basically enforce Hall,
so I should be okay here with my argument. The only unclear thing is balancing
the buffer neighbourhood. Formerly thought I didn't need anything, but actually I need
martingales to preserve vertex neighbourhoods. Maybe queue jumping fixes other
conditions? Do we need to use parts evenly?]
}

We prove Theorem \ref{blowup} with an embedding algorithm
that is very similar to that used for Theorem \ref{3blowup}.
We introduce more parameters with the hierarchy
\begin{gather*}
0 \le 1/n \ll 1/n_R \ll \eps \ll \eps' \ll \eps_{0,0} \ll \cdots \ll \eps_{k^3 D,3} \ll \eps_*
\ll p_0 \ll c \ll \gamma \ll \delta_Q \ll p \ll d_u \ll d_a \\
\ll \theta \ll \theta_0 \ll \theta'_0 \ll \cdots \ll \theta_{k^3 D} \ll \theta'_{k^3 D} \ll \theta_*
\ll \delta'_Q \ll \delta_B \ll  d, c', 1/D, 1/D_R, 1/C, 1/k.
\end{gather*}
The roles of the parameters from Theorem \ref{3blowup} are exactly as before.
Our generalisations to $R$-indexed complexes and restricted positions
have introduced some additional parameters, so one should note how they fit into the hierarchy.
The restricted positions hypothesis has two parameters $c$ and $c'$.\index{$c$}\index{$c'$}
Parameter $c$ controls the number of restricted positions and satisfies $p_0 \ll c \ll \gamma$.
Parameter $c'$ gives a lower bound on the density relative to $G$ of the constraints
and satisfies $c' \gg \delta_B$.
The indexing complex $R$ has two parameters $n_R$ and $D_R$.
Parameter $n_R$ is a bound for $|R|$ and can be very large, provided that $n$ is even larger.
Parameter $D_R$ bounds the maximum degree of $R$ and satisfies $D_R \ll 1/\delta_B$.
\index{$n_R$} \index{$D_R$}
\COMMENT{No need to mention r(?)}

\begin{description}

\item[Initialisation and notation.]
Write $X'_* = X_* \cup \bigcup_{x \in X_*} VN_H(x)$.\index{$X'_*$}
We choose a buffer set $B \subset X$ of vertices at mutual distance at least $9$ in $H$
so that $|B \cap X_i| = \delta_B |X_i|$ for $i \in R$ and $B \cap X'_* = \es$.
Since $n \le |X_i| \le Cn$ for $i \in R$ and $H$ has maximum degree $D$
we can construct $B$ by selecting vertices one-by-one greedily.
Every vertex neighbourhood in $H$ has size less than $kD$,
so there are at most $(kD)^8$ vertices at distance less than $9$ from any vertex of $H$.
Similarly, there are at most $(kD_R)^8$ multi-indices $j \in R$
at distance less than $9$ from any fixed multi-index $i \in R$.
Thus at any stage we have excluded at most $(kD_R)^8(kD)^8\delta_BCn < n/2$
vertices from $X_i$, since $\delta_B \ll 1/D_R$.
Similarly, since $|X_* \cap X_j| < c|X_j|$ for all $j \in R$
we have $|X'_* \cap X_i| < (kD_R)^8(kD)^8cCn < \sqrt{c}n$.
Since $|X_i| \ge n$ we can construct $B$ greedily.
Let $N = \cup_{x \in B} VN_H(x)$ be the set of all vertices that have a neighbour in the buffer.
Then $N$ is disjoint from $X_*$, as we chose $B$ disjoint from $X'_*$.
Also, $|N \cap X_i| < (kD_R)(kD)\delta_B cn < \sqrt{\delta_B}|X_i|$ for any $i \in R$.%

For $S \in H$ we set $F_S(0)=G[\GG]_S$.
We define $L=L(0)$, $q(t)$, $Q(t)$, $j(t)$, $J(t)$, $X_i(t)$, $V_i(t)$
as in the $3$-graph algorithm. We let $X(t) = \cup_{i\in R} X_i(t)$
and $V(t) = \cup_{i\in R} V_i(t)$.
\COMMENT{
We order the vertices in a list $L=L(0)$ that starts with $N$ and ends with $B$.
Within $N$, we arrange that $VN_H(x)$ is consecutive for each $x \in B$.
We denote the vertex of $H$ selected for embedding at time $t$ by $s(t)$.
We denote the queue by $q(t)$ and write $Q(t) = \cup_{u \le t}\ q(u)$.
We denote the vertices jumping the queue by $j(t)$ and write $J(t) = \cup_{u \le t}\ j(u)$.
Initially we set $q(0)=Q(0)=j(0)=J(0)=\es$.
We write $X_i(t) = X_i \sm \{s(u): u \le t\}$
and $V_i(t) = V_i \sm \{\phi(s(u)): u \le t\}$.
We let $X(t) = \cup_{i=1}^r X_i(t)$ and $V(t) = \cup_{i=1}^r V_i(t)$.
We also use the convention that $F_S(t)$ is undefined if $S \notin H$.
}

\item[Iteration.] At time $t$, while there are still some
unembedded non-buffer vertices, we select a vertex to embed $x=s(t)$
according to the same selection rule as for the $3$-graph algorithm.
We choose the image $\phi(x)$ of $x$ uniformly at random among all elements
$y \in F_x(t-1)$ that are `good' (a property defined below).
Note that all expressions at time $t$ are
to be understood with the embedding $\phi(x)=y$, for some unspecified vertex $y$.

\nib{Definitions.}

1. For a vertex $x$ we write $\nu_x(t)$ for the number of\index{$\nu_x(t)$}
elements in $VN_H(x)$ that have been embedded at time $t$.
For a set $S$ we write $\nu_S(t) = \sum_{y \in S} \nu_y(t)$.\index{$\nu_S(t)$}
We also define $\nu'_S(t)$ as follows.\index{$\nu'_S(t)$}
When $|S|=k$ we let $\nu'_S(t) = \nu_S(t)$.
When $|S|<k$ we let $\nu'_S(t) = \nu_S(t)+K$, where $K$ is the
maximum value of $\nu'_{Sx'}(t')$ over vertices $x'$ embedded at time
$t' \le t$ with $S \in H(x')$; if there is no such vertex $x'$
we let $\nu'_S(t) = \nu_S(t)$.

2. For any $r$-partite set $S$ we define $F_S(t)=F_S(t-1)^y$ if $x \in S$
or $F_S(t) = F_{S^\le}(t-1)[F_{S.x^\le}(t-1)(y)]_S \sm y$ if $x \notin S$.
We define an {\em exceptional} set $E_x(t-1) \sub F_x(t-1)$ by saying\index{exceptional}
$y$ is in $F_x(t-1) \sm E_x(t-1)$ if and only if
for every unembedded $\es \ne S \in H(x)$,

$\qquad d_S(F(t))=(1 \pm \eps_{\nu'_S(t),0})d_S(F(t-1))d_{Sx}(F(t-1))$,
and $F_S(t)$ is $\eps_{\nu'_S(t),0}$-regular. \hfill $(*_{\ref{alg}})$

3. We define $E^t$, $M_{E^t,E}(t)$, $D_{x,E}(t-1)$ and $U(x)$ as in
the $3$-graph algorithm, replacing `triple' by `$k$-tuple'.
We obtain the set of {\em good} elements $OK_x(t-1)$ from $F_x(t-1)$
by deleting $E_x(t-1)$ and $D_{x,E}(t-1)$ for every $E \in U(x)$.\index{good}\index{$OK$}
\COMMENT{$G$ `good' and $N$ `normal' are taken...}

We embed $x$ as $\phi(x)=y$ where $y$ is chosen uniformly at random from the
good elements of $F_x(t-1)$. We update $L(t)$, $q(t)$ and $j(t)$ as before,
using the same rule for adding vertices to the queue.
We repeat  until the only unembedded vertices are buffer vertices,
but abort with failure if at any time we have
$|Q(t) \cap X_i| > \delta_Q|X_i|$ for some $i \in R$.
Let $T$ denote the time at which the iterative phase terminates
(whether with success or failure).

\item[Conclusion.] When all non-buffer vertices have been embedded,
we choose a system of distinct representatives among the available
slots $A'_x$ (defined as before) for $x \in X(T)$ to complete the embedding,
ending with success if possible, otherwise aborting with failure.
\end{description}

Similarly to Lemma \ref{itworks},
the algorithm embeds $H$ in $G \sm M$ unless it aborts with failure.
Furthermore, when $\GG_S$ is defined we have $F_S(0)=G[\GG]_S=\GG_S$,
so we ensure that $\phi(S) \in \GG_S$.
Note that any vertex neighbourhood contains at most $(k-1)D$ vertices.
Thus in the selection rule, any element of the queue can cause at most $(k-1)D$
vertices to jump the queue. Note also that when a vertex neighbourhood
jumps the queue, its vertices are immediately embedded at consecutive times
before any other vertices are embedded.
We collect here a few more simple observations on the algorithm.

\begin{lemma}\label{observe}$ $
\item[(i)] For any $i \in R$ and time $t$ we have $|V_i(t)| \ge \delta_B n/2$.
\item[(ii)] For any $i \in R$ and time $t$ we have $|J(t) \cap X_i| < \sqrt{\delta_Q}n$.
\item[(iii)] We have $\nu_x(t) \le (k-1)D$ for any vertex and $\nu'_S(t) \le k^3D$
for any $S \in H$. Thus the $\eps$-subscripts are always defined in $(*_{\ref{alg}})$.
\item[(iv)] For any $z \in VN_H(x)$ we have $\nu_z(t)=\nu_z(t-1)+1$,
so for any $S \in H$ that intersects $VN_H(x)$ we have $\nu_S(t)>\nu_S(t-1)$.
\item[(v)] If $\nu_S(t)>\nu_S(t-1)$ then $\nu'_S(t)>\nu'_S(t-1)$.
\item[(vi)] If $z$ is embedded at time $t'\le t$ and $S\in H(z)$ then
$\nu'_S(t) \ge \nu'_{Sz}(t) > \nu'_{Sz}(t'-1)$.
\end{lemma}

\nib{Proof.} The proofs of (i) and (iii-vi) are similar to those
in Lemma \ref{3observe} so we omit them. For (ii) we have to be more
careful to get a good bound inside each part.
Note that we only obtain a new element $x$ of $J(t) \cap X_i$
when $x$ is a neighbour of some $b \in B$ and some $z$
within distance $4$ of $x$ is queued.
In particular $z$ is within distance $5$ of $x$.
Given $i$, there are at most $(kD_R)^5$ choices for $j=i^*(z) \in R$,
at most $|Q(t) \cap X_j| < \delta_Q Cn$ choices for $z$,
then at most $(kD)^5$ choices for $x$.
Therefore $|J(t) \cap X_i| < (kD_R)^5(kD)^5\delta_Q Cn < \sqrt{\delta_Q}n$. \qed

\subsection{Hypergraph regularity properties}

This subsection contains various properties of hypergraph regularity
analogous to those described earlier for $3$-graphs.
We start with the general counting lemma, analogous to Theorem \ref{3count}.
\COMMENT{need k here, even if 3 used earlier}

\begin{theo} (R\"odl-Schacht \cite{RSc2}, see Theorem 13) \label{count}
\index{dense counting lemma}
Suppose $0 < \eps \ll d, \gamma, 1/r, 1/j, 1/k$, that
$J$ and $G$ are $r$-partite $k$-complexes with vertex sets
$Y = Y_1 \cup \cdots \cup Y_r$ and $V = V_1 \cup \cdots \cup V_r$ respectively,
that $|J|=j$, and $G$ is $\eps$-regular with all densities at least $d$.
Then
$$d(J,G) = \mb{E}_{\phi \in \Phi(Y,V)} \left[ \prod_{A \in J} G_A(\phi(A)) \right]
= (1 \pm \gamma) \prod_{A \in J} d_A(G).$$
\end{theo}

A useful case of Theorem \ref{count} is when $r=k$ and $J=[k]^\le$
consists of all subsets of a $k$-tuple; this gives the following
analogue of Lemma \ref{3absolute}.

\begin{lemma}\label{absolute}
Suppose $0 < \eps \ll \eps' \ll d, 1/k$,
$G$ is a $k$-partite $k$-complex on $V = V_1 \cup \cdots \cup V_k$
with all densities $d_S(G) > d$ and $G$ is $\eps$-regular.
Then $d(G_{[k]}) = (1 \pm \eps')\prod_{S \sub [k]}d_S(G)$.
\end{lemma}

Next we give the analogue of Lemma \ref{3neighbour1}.

\begin{lemma}\label{neighbour1} {\bf (Vertex neighbourhoods)} \index{vertex neighbourhood}
Suppose $G$ is a $k$-partite $k$-complex on $V = V_1 \cup \cdots \cup V_k$
with all densities $d_S(G) > d$ and $0 < \eta_I \ll \eta'_I \ll d, 1/k$
for each $I \sub [k]$. Suppose that each $G_I$ is $\eta_I$-regular.
Then for all but at most $2\sum_{I \sub [k-1]} \eta'_{Ik} |G_k|$ vertices $v \in G_k$,
for every $\es \ne I \sub [k-1]$, $G(v)_I$ is $(\eta'_I+\eta'_{Ik})$-regular
with $d_I(G(v)) = (1 \pm \eta'_I\pm \eta'_{Ik})d_I(G)d_{Ik}(G)$.
\end{lemma}

\nib{Proof.} The argument is similar to that in Lemmas \ref{2neighbour} and \ref{3neighbour1}.
We show the following statement by induction on $|\mc{C}|$:
for any subcomplex $\mc{C}$ of $[k-1]^\le$,
for all but at most $2\sum_{I \in \mc{C}} \eta'_{Ik} |G_k|$ vertices $v \in G_k$,
for every $I \in \mc{C}$, $G(v)_I$ is $(\eta'_I+\eta'_{Ik})$-regular
with $d_I(G(v)) = (1 \pm \eta'_I\pm \eta'_{Ik})d_I(G)d_{Ik}(G)$.
The base case is $\mc{C}=\es$, or less trivially any $\mc{C}$ with
$|I| \le 2$ for all $I \in \mc{C}$, by Lemmas \ref{2neighbour} and \ref{3neighbour1}.

For the induction step, fix any maximal element $I$ of $\mc{C}$.
By induction hypothesis, for all but at most $2\sum_{I' \in \mc{C}\sm I} \eta'_{I'k} |G_k|$
vertices $v \in G_k$, for every $I' \in \mc{C}\sm I$, $G(v)_{I'}$ is $(\eta'_{I'}+\eta'_{I'k})$-regular
with $d_{I'}(G(v)) = (1 \pm \eta'_{I'}\pm \eta'_{I'k})d_{I'}(G)d_{I'k}(G)$.
Let $G'_k$ be the set of such vertices. It suffices to show the claim that
all but at most $2\eta'_{Ik}|G_k|$ vertices $v \in G'_k$ have the following property:
if $J^v$ is a subcomplex of $G(v)_{I^<}$ with $|(J^v)_I^*| > (\eta'_I+\eta'_{Ik}) |G(v)_I^*|$
then $|G[J^v]_I| = (1 \pm \eta'_I/2)d_I(G)|(J^v)_I^*|$
and $|G(v)[J^v]_I| = (1 \pm \eta'_{Ik}/2)d_{Ik}(G)|G[J^v]_I|$.

Suppose for a contradiction that this claim is false.
Let $\gamma = \max_{I' \sub I} \eta'_{I'}$. By Theorem \ref{count}, for any $v \in G'_k$ we have
\[|G(v)_I^*| \prod_{i \in I} |V_i|^{-1} = d(G(v)_I^*)
> (1-\gamma)\prod_{I' \subn I} (1-\eta'_{I'}-\eta'_{I'k})d_{I'}(G)d_{I'k}(G)
> \frac{1}{2} d^{2^k}.\]
Then for any $J^v \sub G(v)_{I^<}$ with $|(J^v)_I^*| > (\eta'_I+\eta'_{Ik}) |G(v)_I^*|$
we have  $|(J^v)_I^*| > \eta_I\prod_{i \in I} |V_i| \ge \eta_I|G_I^*|$, so
$|G[J^v]_I| = (d_I(G) \pm \eta_I)|(J^v)_I^*| = (1 \pm \eta'_I/2)d_I(G)|(J^v)_I^*|$, 
since $G_I$ is $\eta_I$-regular.
So without loss of generality, we can assume that we have vertices
$v_1,\cdots,v_t \in G'_k$ with $t>\eta'_{Ik}|G_k|$,
and subcomplexes $J^{v_i} \sub G(v_i)_{I^<}$ with $|(J^{v_i})_I^*| > (\eta'_I+\eta'_{Ik}) |G(v_i)_I^*|$
such that $|G(v_i)[J^{v_i}]_I| < (1 - \eta'_{I'k}/2)d_{Ik}(G)|G[J^{v_i}]_I|$ for $1 \le i \le t$.
Define complexes $A^i = G[J^{v_i}]_I \cup \{v_iS: S \in J^{v_i}\}$ and $A = \cup_{i=1}^t A_i$.

We have $|A_{Ik}^*| = \sum_{i=1}^t |(A^i)_{Ik}^*| = \sum_{i=1}^t |G[J^{v_i}]_I|$.
Now $|G[J^{v_i}]_I| > (d_I(G) - \eta_I)|(J^{v_i})_I^*|$,
$d_I(G)>d$, $|(J^{v_i})_I^*| > (\eta'_I+\eta'_{Ik}) |G(v_i)_I^*|
> \eta'_{Ik} \cdot \frac{1}{2} d^{2^k} \prod_{i \in I} |V_i|$ and $t>\eta'_{Ik}|G_k|$,
so
\[|A_{Ik}^*| > \eta'_{Ik}|G_k| \cdot (d-\eta_I) \cdot \eta'_{Ik} \cdot \frac{1}{2} d^{2^k}\prod_{i \in I} |V_i|
> \eta_{Ik} \prod_{i \in Ik} |V_i| \ge \eta_{Ik}|G_{Ik}^*|.\]
Since $G_{Ik}$ is $\eta_{Ik}$-regular we have
$d_{Ik}(G[A])_{Ik} = d_{Ik}(G) \pm \eta_{Ik}$.
Therefore $|G \cap A_{Ik}^*| > (d_{Ik}(G)-\eta_{Ik})|A_{Ik}^*|=(d_{Ik}(G)-\eta_{Ik})\sum_{i=1}^t |G[J^{v_i}]_I|$.
But we also have
\[|G \cap A_{Ik}^*| = \sum_{i=1}^t |G(v_i)[J^{v_i}]_I|
< \sum_{i=1}^t (1 - \eta'_{I'k}/2)d_{Ik}(G)|G[J^{v_i}]_I| < (d_{Ik}(G)-\eta_{Ik})\sum_{i=1}^t |G[J^{v_i}]_I|,\]
contradiction. This proves the claim, and so completes the induction. \qed

We apply Lemma \ref{neighbour1} in the next lemma showing that
arbitrary neighbourhoods are typically regular.

\begin{lemma}\label{neighbour-set} {\bf (Set neighbourhoods)} \index{set neighbourhoods}
Suppose $0 < \eps \ll \eps' \ll d, 1/k$ and $G$ is an $\eps$-regular
$k$-partite $k$-complex on $V = V_1 \cup \cdots \cup V_k$
with all densities $d_S(G) > d$. Then for any $A \sub [k]$
and for all but at most $\eps' |G_A|$ sets $P \in G_A$,
for any $\es \ne I \sub [k] \sm A$, $G(P)_I$ is $\eps'$-regular
with $d_I(G(P)) =  (1 \pm \eps')\prod_{A' \sub A} d_{A'I}(G)$
(and $d_\es(G(P))=1$ as usual).
\end{lemma}

\nib{Proof.}
For convenient notation suppose that $A=[k']$ for some $k'<k$.
Introduce additional constants with the hierarchy
$\eps \ll \eps_1 \ll \eps'_1 \ll \cdots \ll \eps_{k'} \ll \eps'_{k'} \ll \eps'$.
We prove inductively for $1 \le t \le k'$ that for all but at most $\eps_t |G_{[t]}|$
sets $S \in G_{[t]}$, for any $\es \ne I \sub [k] \sm [t]$,
$G(S)$ is $\eps_t$-regular with  $d_I(G(S)) = (1 \pm \eps_t)\prod_{T \sub [t]} d_{TI}(G)$.
The base case $t=1$ is immediate from Lemma \ref{neighbour1} with $\eps'$ replaced by $\eps_1$.
For the induction step, consider $S \in G_{[t]}$ such that for any $\es \ne I \sub [k] \sm [t]$,
$G(S)$ is $\eps_t$-regular with  $d_I(G(S)) = (1 \pm \eps_t)\prod_{T \sub [t]} d_{TI}(G)$.
We have $d_I(G(S)) > \frac{1}{2}d^{2^t}$, so we can apply Lemma \ref{neighbour1} to $G(S)$
with $\eta_I=\eps_t$, $\eta'_I=2^{-k}\eps'_t$ and $d$ replaced by $\frac{1}{2}d^{2^t}$.
Then for all but at most $\eps'_t|G(S)_{t+1}|$ vertices $v \in G(S)_{t+1}$,
for every $\es \ne I \sub [k] \sm [t+1]$, $G(S)(v)_I=G(Sv)_I$ is $\eps'_t$-regular
with $d_I(G(Sv)) = (1 \pm \eps'_t)d_I(G(S))d_{Ik}(G(S))
= (1 \pm 2\eps'_t)\prod_{T \sub [t+1]} d_{TI}(G)$.
% Since $\eps'_t \ll \eps_{t+1}$ we have the densities required for the induction step.
Also, since $|G_{[t+1]}| > \frac{1}{2}d^{2^{t+1}} \prod_{i=1}^{t+1} |V_i|$ by Theorem \ref{count},
the number of pairs $(S,v)$ for which this fails is at most
$\eps_t |G_{[t]}| \cdot |V_{t+1}| + |G_{[t]}| \cdot \eps'_t|G(S)_{t+1}| < \eps_{t+1} |G_{[t+1]}|$. \qed

We omit the proofs of the next three lemmas, as they are almost identical to those of
the corresponding Lemmas \ref{3restrict}, \ref{3restrict'} and \ref{3average},
replacing $123$ by $[k]$ and using Theorem \ref{count} instead of the triangle-counting lemma.

\begin{lemma}\label{restrict} {\bf (Regular restriction)} \index{regular restriction}
Suppose $G$ is a $k$-partite $k$-complex on $V = V_1 \cup \cdots \cup V_k$
with all densities $d_S(G) > d$, $G_{[k]}$ is $\eps$-regular, where $0 < \eps \ll d, 1/k$,
and $J \sub G$ is a $(k-1)$-complex with $|J^*_{[k]}|>\sqrt{\eps}|G^*_{[k]}|$.
Then $G[J]_{[k]}$ is $\sqrt{\eps}$-regular and $d_{[k]}(G[J])=(1\pm\sqrt{\eps})d_{[k]}(G)$.
\end{lemma}

\begin{lemma}\label{restrict'}
Suppose $G$ is a $k$-partite $k$-complex on $V = V_1 \cup \cdots \cup V_k$
with all densities $d_S(G) > d$ and $G_{[k]}$ is $\eps$-regular, where $0 < \eps \ll d, 1/k$.
Suppose also that $J \sub G$ is a $(k-1)$-complex,
and when defined,
$d_I(J)>d$ and $J_I$ is $\eta$-regular, where $0 < \eta \ll d$.
Then $G[J]_{[k]}$ is $\sqrt{\eps}$-regular and $d_{[k]}(G[J])=(1\pm\sqrt{\eps})d_{[k]}(G)$.
\end{lemma}

\begin{lemma}\label{average}
Suppose $G$ is a $k$-partite $k$-complex on $V = V_1 \cup \cdots \cup V_k$
with all densities $d_S(G) > d$ and $G$ is $\eps$-regular,
where $0 < \eps \ll \eps' \ll d, 1/k$.
Then for any $A \sub [k]$ and for all but at most $\eps' |G_A|$ sets $P \in G_A$
we have $|G(P)_{[k] \sm A}| = (1 \pm \eps')|G_{[k]}|/|G_A|$.
\end{lemma}

Note that we will not need an analogue of the technical Lemma \ref{3technical}.

\subsection{Analysis of the algorithm}

We start the analysis of the algorithm by showing that most free vertices are good.
First we record some properties of the initial sets $F_S(0)$,
taking into account the restricted positions. \index{restricted positions}

\begin{lemma}\label{res-props}
$F_S(0)$ is $\eps'$-regular with $d_S(F(0))>c'd_S(G)$
and $|F_S(0)| > (c')^{2^{|S|}}|G_S|$.
\end{lemma}

\nib{Proof.} By definition we have $F_S(0)=G[\GG]_S$.
Condition (i) of Definition \ref{def-super} tells us that $G_S$ is $\eps$-regular
with $d_S(G) \ge d_a$. Hypothesis (v) of Theorem \ref{blowup}
says that when defined $\GG_{S'}$ is $\eps'$-regular with $d_{S'}(\GG) > c'd_{S'}(G)$ for $S' \sub S$.
If $\GG_S$ is defined then $F_S(0)=\GG_S$ is $\eps'$-regular with $d_S(F(0))=d_S(\GG)>c'd_S(G)$.
Otherwise, by Lemma \ref{restrict'}, $\GG_S$ is $\sqrt{\eps}$-regular
with $d_S(F(0))=(1 \pm \eps)d_S(G)>c'd_S(G)$.
The estimate $|F_S(0)| > (c')^{2^{|S|}}|G_S|$
follows by applying Theorem \ref{count} to $F(0)_{S^\le}$ and $G_{S^\le}$.
Note that $d_{\es}(F(0))=d_{\es}(G)=1$, so one of the $c'$ factors
compensates for the error terms in Theorem \ref{count}. \qed

Our next lemma handles the definitions for regularity and density in the algorithm.

\begin{lemma}\label{exceptional} \index{exceptional}
The exceptional set $E_x(t-1)$ defined by $(*_{\ref{alg}})$ satisfies $|E_x(t-1)| < \eps_*|F_x(t-1)|$,
and $F_S(t)$ is $\eps_{\nu'_S(t),1}$-regular with $d_S(F(t)) \ge d_u$ for every $S \in H$.
\COMMENT{Do we use 0 in second eps co-ordinate? Yes, condition * in alg.}
\end{lemma}

\nib{Proof.} We argue by induction on $t$. At time $t=0$ the first statement is vacuous
and the second follows from Lemma \ref{res-props}, since $d_S(G) \ge d_a$ for $S \in H$.
Now suppose $t \ge 1$ and $\es \ne S \in H$ is unembedded, so $x \notin S$.
We consider various cases for $S$ to establish the bound on the exceptional set
and the regularity property, postponing the density bound until later in the proof.

We start with the case when $S \in H(x)$.
By induction $F_{S'}(t-1)$ is $\eps_{\nu'_{S'}(t-1),1}$-regular
and $d_{S'}(F(t-1)) \ge d_u$ for every $S' \sub Sx$.
Write $\nu = \max_{S' \sub S} \nu'_{S'x}(t-1)$ and
$\nu^* = \max\{\nu'_S(t-1),\nu'_{Sx}(t-1)\}$.
By Lemma \ref{neighbour1}, for all but at most $\eps_{\nu,2}|F_x(t-1)|$ vertices $y \in F_x(t-1)$,
$F_S(t)=F_{Sx}(t-1)(y)$ is $\eps_{\nu^*,2}$-regular with
$d_S(F(t))=(1\pm\eps_{\nu^*,2})d_S(F(t-1))d_{Sx}(F(t-1))$.
Also, we claim that $\nu'_S(t)>\nu^*$. This holds by Lemma \ref{observe}:
(iv) gives $\nu_S(t)>\nu_S(t-1)$,
(v) gives $\nu'_S(t)>\nu'_S(t-1)$, and
(vi) gives $\nu'_S(t)>\nu'_{Sx}(t-1)$.
Thus we have $(*_{\ref{alg}})$ for such $y$.
We also note for future reference that $d_S(F(t))>d_u^2/2$.
In the argument so far we have excluded at most $\eps_{k^3 D,2}|F_x(t-1)|$
vertices $y \in F_x(t-1)$ for each of at most $2^{k-1} D$ sets $S \in H(x)$;
this gives the required bound on $E_x(t-1)$.
We also have the required regularity property of $F_S(t)$,
but for now we postpone showing the density bounds.

Next we consider the case when $S \in H$ and $S \notin H(x)$.
We show by induction on $|S|$ that $F_S(t)$ is $\eps_{\nu'_S(t),1}$-regular
with $d_S(F(t)) > d_u/2$, and moreover $F_S(t)$ is $\eps_{\nu'_S(t-1),2}$-regular
when $S$ intersects $VN_H(x)$. Note that if $S$ intersects $VN_H(x)$
then we have $\nu'_S(t)>\nu'_S(t-1)$ by Lemma \ref{observe}(v).
For the base case when $S=\{v\}$ has size $1$
we have $|F_v(t)|=|F_v(t-1)\sm y| \ge |F_v(t-1)|-1$,
so $d_v(F(t)) \ge d_v(F(t-1))-1/n > d_u/2$.
For the induction step, suppose that $|S| \ge 2$.
Recall that Lemma \ref{build-update} gives $F_S(t) = F_{S^\le}(t-1)[F_{S^<}(t)]_S$.
Since $F_S(t-1)$ is $\eps_{\nu'_S(t-1),1}$-regular, by Lemma \ref{restrict},
$F_S(t)$ is $\eps_{\nu'_S(t-1),2}$-regular with
$d_S(F(t)) = (1 \pm \eps_{\nu'_S(t-1),2})d_S(F(t-1))$.
This gives the required property in the case when $S$ intersects $VN_H(x)$.
Next suppose that $S$ and $VN_H(x)$ are disjoint.
Let $t'$ be the most recent time at which we embedded a vertex $x'$ with a neighbour in $S$.
Then by Lemma \ref{not-local}, $F_{S^\le}(t)$ is obtained from $F_{S^\le}(t')$ just by deleting
all sets containing any vertices that are embedded between time $t'+1$ and $t$.
Equivalently, $F_S(t) = F_S(t')[((F_z(t): z \in S),\{\es\})]$.
Now $F_S(t')$ is $\eps_{\nu'_S(t'-1),2}$-regular with $d_S(F(t'))>d_u$
and $\nu'_S(t) \ge \nu'_S(t') > \nu'_S(t'-1)$,
so $F_S(t)$ is $\eps_{\nu'_S(t),1}$-regular with $d_S(F(t))>d_u/2$ by Lemma \ref{restrict}.

Now we have established the bound on $E_x(t-1)$ and the regularity properties,
so it remains to show the density bounds.
First we consider any unembedded $S \in H$ with $|S|=k$.
Then $F_S(t) = F(0)[F_{S^<}(t)]_S$, similarly to (\ref{refwants}) in Lemma \ref{3exceptional}.
By Lemma \ref{res-props}, $F(0)_S$ is $\eps'$-regular with $d_S(F(0))>c'd$.
Also, we showed above for every $S' \subn S$ that
$F_{S'}(t)$ is $\eps_{\nu'_{S'}(t),1}$-regular with $d_{S'}(F(t)) > d_u/2$.
Now Lemma \ref{restrict'} shows that $F_S(t)$ is $\sqrt{\eps'}$-regular
with $d_S(F(t)) = (1\pm \sqrt{\eps'})d_S(F(0)) > c'd/2$.
In particular we have $d_S(F(t))>d_u$, 
though we also use the bound $c'd/2$ below.

Next we show for $k-1 \ge |S| \ge 2$ that $d_S(F(t)) > 4d_a^{2D^{2(k-|S|)}}$.
We argue by induction on $t$ and reverse induction on $|S|$,
i.e.\ we assume that the bound holds for larger sets than $S$.
When $|S|=k$ we have already proved $d_S(F(t)) > c'd/2 > 4d_a^2$.
Let $t' \le t$ be the most recent time at which we embedded a vertex $x'$
with $S \in H(x')$, or let $t'=0$ if there is no such vertex $x'$.
Note that we may have $t'=t$ if $S \in H(x)$.
Let $J(t)$ be the $1$-complex $((F_z(t):z \in S),\{\es\})$.
As in Lemma \ref{3exceptional} we have $F_S(t)=F(t')[J(t)]_S$,
so Lemma \ref{restrict} gives $d_S(F(t))>d_S(F(t'))/2$.
Also, $(*_{\ref{alg}})$ gives $d_S(F(t'))=(1\pm \eps_{\nu'_S(t'),0})d_S(F(t'-1))d_{Sx'}(F(t'-1))$,
where $d_{Sx'}(F(t'-1)) > 4d_a^{2D^{2(k-|S|-1)}}$ by induction.
Thus the $S$-density starts with $d_S(F(0)) > c'd_a$,
loses a factor no worse than $1/2$ before we embed some $x'$ with $S \in H(x')$,
then loses a factor no worse than $d_a^{2D^{2(k-|S|-1)}}$ on at most $D$ occasions
when we embed some $x'$ with $S \in H(x')$. This gives
$d_S(F(t)) > c'd_a/2 \cdot d_a^{2D^{2(k-|S|-1)+1}} > 4d_a^{2D^{2(k-|S|)}}$.
In particular we have $d_S(F(t))>d_u$,
though we also use the bound $4d_a^{2D^{2(k-|S|)}}$ below.

Finally we consider any unembedded vertex $z$.
Suppose that the current $z$-regime started at some time $t_z\le t$.
If $t_z>0$ then we embedded some neighbour $w=s(t_z)$ of $z$ at time $t_z$.
By $(*_{\ref{alg}})$ and the above bound for pair densities we have
$d_z(F(t_z))>d_{wz}(F(t_z-1))d_z(F(t_z-1))/2 > d_a^{2D^{2(k-2)}} d_z(F(t_z-1))$.
Now we consider cases according to whether $z$ is in the list $L(t-1)$,
the queue $q(t-1)$ or the queue jumpers $j(t-1)$.
Suppose first that $z \in L(t-1)$.
Then the rule for updating the queue in the algorithm
gives $|F_z(t)| \ge \delta'_Q |F_z(t_z)|$.
Next suppose that $z \in j(t-1) \cup q(t-1)$,
and $z$ first joined $j(t') \cup q(t')$ at some time $t'<t$.
Since $z$ did not join the queue at time $t'-1$ we have
$|F_z(t'-1)| \ge \delta'_Q |F_z(t_z)|$.
Also, between times $t'$ and $t$ we only embed vertices that are in the queue
or jumping the queue, or otherwise we would have embedded $z$ before $x$.
Now $|Q(t) \cap X_i| \le \delta_Q Cn$ for $i \in R$,
or otherwise we abort the algorithm,
and $|J(t) \cap X_z| < \sqrt{\delta_Q}n$ by Lemma \ref{observe}(ii),
so we embed at most $2\sqrt{\delta_Q} |V_z|$ vertices in $V_z$
between times $t'$ and $t$. Thus we have catalogued all possible ways
in which the number of vertices free for $z$ can decrease.
It may decrease by a factor of $d_a^{2D^{2(k-2)}}$ when a new $z$-regime starts,
and by a factor $\delta'_Q$ during a $z$-regime before $z$ joins the queue.
Also, if $z$ joins the queue or jumps the queue it may decrease by
at most $2\sqrt{\delta_Q} |V_z|$ in absolute size.
Since $z$ has at most $2D$ neighbours, and $|F_z(0)|>c'|V_z|$, we have
$|F_z(t)| \ge (d_a^{2D^{2(k-2)}}\delta'_Q)^{kD}c'|V_z| - 2\sqrt{\delta_Q} |V_z| > d_u|V_z|$. \qed

From now on it will often suffice and be more convenient to use a crude upper bound
of $\eps_*$ for any epsilon parameter. The estimates in Lemma \ref{3crude}
hold in general (we can replace $12D$ by $k^3 D$ in (vi)).
We also need some similar properties for plus complexes.
In the following statement $H^+$ is to be understood as
in Definition \ref{def-plus} but with $x$ replaced by $z$.\index{plus complex}
\COMMENT{Ref noted $H^z$ was `relic' from earlier version!}

\begin{lemma}\label{plus-props}
Suppose $z \in X$ and $S' \sub S \in H^+$ are unembedded
and $I$ is a subcomplex of $H$. % $\{S': S' \sub S, S' \in H(z)\}$.
Then\COMMENT{`unembedded' technically not defined for $z^c$...}
\begin{itemize}
\item[(i)]
$F(t)^{I+z}_{S^\le}$ is $\eps_{k^3D,3}$-regular. %Sz^c
\item[(ii)]
If $S \in H$ then $d_{S}(F(t)^{I+z})=d_{S}(F(t))$.

If $S \in H^+ \sm H$ then
$d_{S}(F(t)^{I+z})$ is $d_{T}(F(t))$ if $S=T^c$, $T \in I$
or $1$ otherwise.
\item[(iii)]
For all but at most $\eps_*|F(t)^{I+z}_{S'}|$ sets
$P \in F(t)^{I+z}_{S'}$ we have

$|F(t)^{I+z}_S(P)| = (1\pm\eps_*) |F(t)^{I+z}_{S}|/|F(t)^{I+z}_{S'}|$.
\item[(iv)]
$d(F(t)^{I+z}_S)=(1\pm\eps_*)\prod_{T\sub S}d_T(F(t)^{I+z})$.
\item[(v)]$ $

\vspace{-1cm}

\[\frac{|F(t)^{I+z}_S|}{|F(t)^{I+z}_{S'}||F(t)^{I+z}_{S\sm S'}|}
=\frac{d(F(t)^{I+z}_S)}{d(F(t)^{I+z}_{S'})d(F(t)^{I+z}_{S\sm S'})}
= (1\pm 4\eps_*)\prod_{T: T\sub S,T \nsub S',T \nsub S\sm S'}d_T(F(t)^{I+z}).\]
\item[(vi)] Statements (iii-v) hold replacing $F(t)^{I+z}_{S^\le}(t)$ by
$F(t)^{I+z}_{S^\le}(t)[\GG]$ for any $\eps_{k^3 D,3}$-regular subcomplex $\GG$
of $F(t)^{I+z}_{S^\le}(t)$, such that $d_T(\GG) \ge \eps_*$  when defined.
\end{itemize}
\end{lemma}

\nib{Proof.} 
If $S \in H$ then $F(t)^{I+z}_{S^\le}=F(t)_{S^\le}$.
Now suppose $S \in H^+\sm H$, say $S=Tz^c$ with $T \in H$.
Then by Definition \ref{def-plus},
$F(t)^{I+z}_{S^\le} = G^+[F(t)\cup F(t)^+_{I^c}]_{S^\le}$
consists of all $Py^c$ with $P \in F(t)_T$ and
$P_{S'}y^c \in F_{S'z}(t)^c$ when defined for all $S'z \in I$.
Thus we can write
$F(t)^{I+z}_{S^\le} = K(V^+)[F(t)\cup F(t)^+_{I^c}]_{S^\le}$.
Note that $K(V^+)_{S'}$ is $\eps$-regular with $d_{S'}(K(V^+))=1$ for any $S'\in H^+$.
Thus (i) and (ii) follow by regular restriction. The other statements
of the Lemma can be proved as in Lemma \ref{3crude}. \qed

Our next lemma concerns the definitions for marked edges in the algorithm.

\begin{lemma}\label{marked}$ $ \index{marked}
\begin{itemize}
\item[(i)] For every $k$-tuple $E \in H$ we have
$|M_{E^t,E}(t)| < \theta'_{\nu'_{E^t}(t)} |F_{E^t}(t)|$, and in fact

$|M_{E^t,E}(t)| \le \theta_{\nu'_{E^t}(t)} |F_{E^t}(t)|$ for $E \in U(x)$.
\item[(ii)] For every $x$ and $k$-tuple $E \in U(x)$ we have
$|D_{x,E}(t-1)| < \theta_{\nu'_{E^t}(t)} |F_x(t-1)|$.
\end{itemize}
\end{lemma}

\nib{Proof.}
Throughout we use the notation $\ov{E}=E^{t-1}$, $\nu=\nu'_{\ov{E}}(t-1)$, $\nu^*=\nu'_{E^t}(t)$.

(i) To verify the bound for $t=0$ we use our assumption that $(G,M)$ is
$(\eps,\eps',d_a,\theta,d)$-super-regular. We take $I=(\{\es\})$,
when for any $v$ we have $G^{I_v}=G$ by Definition \ref{def-preplus}.
Then condition (iii) in Definition \ref{def-super} gives $|M_E| \le \theta |G_E|$.
Since $E^0=E$ we have $M_{E,E}(0)=M_E \cap F_E(0)$,
where $|F_E(0)| > (c')^{2^k}|G_E|$ by Lemma \ref{res-props}.
Since $\theta \ll c'$ we have $|M_{E,E}(0)| \le |M_E|
\le \theta |G_E| < \theta(c')^{-2^k}|F_E(0)| < \theta_0 |F_E(0)|$.
Now suppose $t>0$. When $E \in U(x)$ we have $|M_{E^t,E}(t)| \le \theta_{\nu^*} |F_{E^t}(t)|$
by definition, since the algorithm chooses $y=\phi(x) \notin D_{x,E}(t-1)$.
Now suppose $E \notin U(x)$, and let $t'<t$ be the most recent time
at which we embedded a vertex $x'$ with $E \in U(x')$. Then $E^{t'}=E^t$, $\nu'_{E^t}(t')=\nu^*$,
and $|M_{E^{t'},E}(t')| \le \theta_{\nu'_{E^t}(t')} |F_{E^{t'}}(t')|$ by the previous case.
For any $z \in E^t$, we can bound $|F_z(t)|$ using the same argument
as that used at the end of the proof of Lemma \ref{exceptional}.
We do not embed any neighbour of $z$ between time $t'+1$ and $t$,
so the size of the free set for $z$ can only decrease by
a factor of $\delta'_Q$ and an absolute term of $2\sqrt{\delta_Q}n$.
Since $d_z(F(t')) \ge d_u \gg \delta_Q$ we have
$|F_z(t)| \ge \delta'_Q |F_z(t')| - 2\sqrt{\delta_Q}n \ge \frac{1}{2}\delta'_Q |F_z(t')|$.
By Lemma \ref{build-update}, for every $\es \ne S \sub E^t$,
$F_S(t)$ is obtained from $F_S(t')$ by restricting to the $1$-complex $((F_z(t):z \in S),\{\es\})$.
If $|S|\ge 2$ then Lemma \ref{restrict} gives $d_S(F(t)) = (1 \pm \eps_*)d_S(F(t'))$.
Now $d(F_{E^t}(t)) = (1 \pm \eps_*) \prod_{S \sub E^t} d_S(F(t))$ by Lemma \ref{absolute}, so
$$\frac{|F_{E^t}(t)|}{|F_{E^t}(t')|} = (1 \pm 2^{k+1}\eps_*) \prod_{S \sub E^t} \frac{d_S(F(t))}{d_S(F(t'))}
= (1 \pm 2^{k+2}\eps_*) \prod_{z \in E^t} \frac{d_z(F(t))}{d_z(F(t'))} > (\delta'_Q/2)^k/2.$$
Therefore $|M_{E^t,E}(t)| \le |M_{E^t,E}(t')| \le \theta_{\nu^*} |F_{E^t}(t')|
< 2 (\delta'_Q/2)^{-k} \theta_{\nu^*} |F_{E^t}(t)| < \theta'_{\nu^*} |F_{E^t}(t)|$.

(ii) The argument when $x \in E$ is identical to that in Lemma \ref{3marked},
so we just consider the case when $x \notin E$. Then $E^t=E^{t-1}=\ov{E}$.
Note that when $E \in U(x)$ we have $\ov{E} \cap VN_H(x) \ne \es$, so $\nu^*>\nu$.
By Lemma \ref{plus-update} we have $F(t)_{\ov{E}}=F(t-1)^{H+x}(y^c)_{\ov{E}}\sm y
= (F(t-1)^{H+x}_{\ov{E}x^c}\sm y)(y^c)$.
Since $M_{\ov{E},E}(t) = M_{\ov{E},E}(t-1) \cap F_{\ov{E}}(t)$ we have
$$D_{x,E}(t-1) = \left\{ y \in F_x(t-1):
\frac{|M_{\ov{E},E}(t-1) \cap (F(t-1)^{H+x}_{\ov{E}x^c}\sm y)(y^c)|}{
|(F(t-1)^{H+x}_{\ov{E}x^c}\sm y)(y^c)|} > \theta_{\nu^*} \right\}.$$
Then
\begin{align*}
\Sa & := \sum_{y \in D_{x,E}(t-1)} |M_{\ov{E},E}(t-1) \cap (F(t-1)^{H+x}_{\ov{E}x^c}\sm y)(y^c)| \\
& > \theta_{\nu^*} \sum_{y \in D_{x,E}(t-1)\sm E_x(t-1)} |(F(t-1)^{H+x}_{\ov{E}x^c}\sm y)(y^c)| \\
& >  (1-2\eps_*)\theta_{\nu^*} (|D_{x,E}(t-1)|-\eps_*|F_x(t-1)|)\
|F(t-1)^{H+x}_{\ov{E}x^c}|/|F_{x^c}(t-1)| \\
& = (1-2\eps_*)\theta_{\nu^*} (|D_{x,E}(t-1)|/|F_x(t-1)|-\eps_*)
|F(t-1)^{H+x}_{\ov{E}x^c}|.
\end{align*}
Here we used $|F_{x^c}(t-1)|=|F_x(t-1)|$,
and in the second inequality we applied Lemma \ref{plus-props},%
\COMMENT{Technically should remove another exceptional set... leave it for now...}
with a factor $1-2\eps_*$ rather than $1-\eps_*$
to account for the error from deleting $y$ (which has a lower order of magnitude).
We also have $\Sa \le \sum_{y \in F_x(t-1)} |M_{\ov{E},E}(t-1) \cap F(t-1)^{H+x}_{\ov{E}x^c}(y^c)|$.
This counts all pairs $(y,P)$ with $P \in M_{\ov{E},E}(t-1)$, $y \in F_x(t-1)$
and $Py^c \in F(t-1)^{H+x}_{\ov{E}x^c}$, so we can rewrite it as
$\Sa \le \sum_{P \in M_{\ov{E},E}(t-1)} |F(t-1)^{H+x}_{\ov{E}x^c}(P)|$.
By Lemma \ref{plus-props} we have% 
\COMMENT{Was worried about ambiguity, but this is fine.}
\[|F(t-1)^{H+x}_{\ov{E}x^c}(P)| =
(1 \pm \eps_*)\frac{|F(t-1)^{H+x}_{\ov{E}x^c}|}{|F(t-1)^{H+x}_{\ov{E}}|}\]
for all but at most $\eps_*|F(t-1)^{H+x}_{\ov{E}}|$ sets $P \in F(t-1)^{H+x}_{\ov{E}}$.
Since $F(t-1)^{H+x}_{\ov{E}} = F_{\ov{E}}(t-1)$ we have
$$\Sa \le |M_{\ov{E},E}(t-1)|(1+\eps_*)\frac{|F(t-1)^{H+x}_{\ov{E}x^c}|}{|F_{\ov{E}}(t-1)|}
+ \eps_*|F_{\ov{E}}(t-1)||F_x(t-1)|.$$
Combining this with the lower bound on $\Sa$  we obtain
$$(1-2\eps_*)\theta_{\nu^*}\left(\frac{|D_{x,E}(t-1)|}{|F_x(t-1)|}-\eps_*\right)
< (1+\eps_*)\frac{|M_{\ov{E},E}(t-1)|}{|F_{\ov{E}}(t-1)|}
+ \eps_*\frac{|F_{\ov{E}}(t-1)||F_x(t-1)|}{|F(t-1)^{H+x}_{\ov{E}x^c}|}.$$
Now $|M_{\ov{E},E}(t-1)| < \theta'_{\nu} |F_{\ov{E}}(t-1)|$ by (i)
and
$\frac{|F_{\ov{E}}(t-1)||F_x(t-1)|}{|F(t-1)^{H+x}_{\ov{E}x^c}|} \le d_u^{-2^k} \ll \eps_*^{-1}$,
by Lemma \ref{plus-props}, so
\begin{equation}
\frac{|D_{x,E}(t-1)|}{|F_x(t-1)|} <
\frac{(1+\eps_*)\theta'_{\nu}+\sqrt{\eps_*}}{(1-2\eps_*)\theta_{\nu^*}}
+ \eps_* < \theta_{\nu^*}.  \tag*{$\Box$}
\end{equation}

The following corollary is now immediate from Lemmas \ref{exceptional} and \ref{marked}.
Recall that $OK_x(t-1)$ is obtained from $F_x(t-1)$ by deleting $E_x(t-1)$
and $D_{x,E}(t-1)$ for $E \in U(x)$, and note that since $H$ has maximum degree
at most $D$ we have $|U(x)| \le (k-1)D^2$.

\begin{coro}\label{ok} $|OK_x(t-1)| > (1-\theta_*)|F_x(t-1)|$.\index{good}\index{$OK$}
\end{coro}

Next we consider the initial phase of the algorithm, \index{initial phase}
during which we embed the neighbourhood $N$ of the buffer $B$.
We give three lemmas that are analogous to those used for the $3$-graph blow-up lemma.
First we recall the key properties of the selection rule during the initial phase.
Since $H$ has maximum degree $D$ we have $|VN_H(x)| \le (k-1)D$ for all $x$.
We embed all vertex neighbourhoods $VN_H(x)$, $x \in B$ at consecutive times,
and before $x$ or any other vertices at distance at most $4$ from $x$. \index{local}
Then Lemma \ref{not-local} implies that if we start embedding $VN_H(x)$ just after some time $T_0$
then $F_z(T_0) = F_z(0) \cap V_z(T_0)$ consists of all vertices in $F_z(0)$ that have not
yet been used by the embedding, for every $z$ at distance at most $3$ from $x$.
Recall that $F_z(0)=G[\GG]_z$ is either a set of restricted positions $\GG_z$ with
$|\GG_z|>c'|G_z|=c'|V_z|$ or is $G_z=V_z$ if $\GG_z$ is undefined. \index{restricted positions}
We chose $B$ disjoint from $X'_* = X_* \cup \bigcup_{x \in X_*} VN_H(x)$,
so $B \cup N$ is disjoint from $X_*$. Thus for $z \in VN_H(x)\cup \{x\}$ we have $F_z(T_0) = V_z(T_0)$.
We also recall that $|B \cap V_z| = \delta_B |V_z|$, $|N \cap V_z| < \sqrt{\delta_B}|V_z|$,
$|Q(T_0) \cap V_z| \le \delta_Q|V_z|$ and $|J(T_0) \cap V_z| \le \sqrt{\delta_Q}|V_z|$
by Lemma \ref{3observe}(ii). Since $\delta_Q \ll \delta_B \ll c'$,
for any $z$ at distance at most $3$ from $x$ we have
\COMMENT{
bottleneck is $c' \gg \theta$ for marked edges, but can go lower
with marking conditions in restricted positions...
}
\begin{equation}\label{eq:fresh}
|F_z(T_0)| = |F_z(0) \cap V_z(T_0)| > (1-\delta_B^{1/3})|F_z(0)|.
\end{equation}

Our first lemma is analogous to Lemma \ref{3uniform}.
We omit the proof, which is almost identical to that for $3$-graphs.
The only modifications are to replace $2D$ by $(k-1)D$, $12D$ by $k^3D$,
and $\sum_{\ell=1}^r$ by a sum over at most $D_R$ neighbours $\ell \in R$ of $i$;
the estimates are still valid as $\delta_B \ll 1/D_R$.
\COMMENT{
Formerly thought optional sampling was necessary to deal with $n_R$,
but still ok for $n \gg n_R$.
}

\begin{lemma}\label{uniform}
With high probability, for every $S\in R$ with $|S|=2$ lying over some
$i,j \in R$, and vertex $v \in G_i$ with $|G_S(v)| \ge d_u|V_j|$, we have
$|G_S(v) \cap V_j(T_I)| > (1-\delta_B^{1/3})|G_S(v)|$.
\end{lemma}

Next we fix a vertex $x \in B$ and write $VN_H(x) = \{z_1,\cdots,z_g\}$,
with vertices listed in the order that they are embedded.
We let $T_j$ be the time at which $z_j$ is embedded.
By the selection rule, $VN_H(x)$ jumps the queue and is embedded
at consecutive times: $T_{j+1}=T_j+1$ for $1 \le j \le g-1$.
For convenience we also define $T_0=T_1-1$.
Note that no vertex of $VN_H(x)$ lies in $X_x$.
Our second lemma shows that for that any $W \sub V_x$ that is not too small, the
probability that $W$ does not contain a vertex available for $x$ is quite small.\index{$W$}

\begin{lemma} \label{initial}
For any $W \sub V_x$ with $|W|> \eps_*|V_x|$,
conditional on any embedding of the vertices $\{s(u): u < T_1\}$
that does not use any vertex of $W$,
we have $\mb{P}[A_x \cap W = \es] < \theta_*$.
\end{lemma}

\nib{Proof.}
The proof is very similar to that of Lemma \ref{3initial},
so we will just describe the necessary modifications.
We note that since $B \cup N$ is disjoint from $X_*$
we do not need to consider restricted positions in this proof.
Suppose $1 \le j \le g$ and that we are considering the embedding of $z_j$.
We interpret quantities at time $T_j$ with the embedding $\phi(z_j)=y$,
for some as yet unspecified $y \in F_{z_j}(T_j-1)$.
We define $W_j$, $[W_j]$, $E^W_{z_j}(T_j-1)$, $D^W_{z_j,E}(T_j-1)$
and the events $A_{i,j}$ as before (replacing triples with $k$-tuples).
The proofs of Claims A, B, C, E and the conclusion of the proof
are almost identical to before.
We need to modify various absolute constants to take account of
the dependence on $k$, e.g.\ changing $20$ to $2^{k+2}$ in Claim A,
$12$ to $k^3$ in Claim E, and $g \le 2D$ to $g \le (k-1)D$.
Also, when we apply equation (\ref{eq:fresh}) instead of (\ref{eq:3fresh})
we will replace $2\sqrt{\delta_B}$ by $\delta_B^{1/3}$.

To complete the proof of the lemma it remains to establish Claim D.
This requires more substantial modifications, similar to those in Lemma \ref{marked},
so we will give more details here. Suppose that $A_{1,j-1}$ and $A_{2,j-1}$ hold
and $E$ is a $k$-tuple containing $x$. As before we write $\ov{E}=E^{T_j-1}$,
$\nu=\nu'_{\ov{E}}(T_j-1)$, $\nu^*=\nu'_{E^{T_j}}(T_j)$
and $B_{z_j} = E_{z_j}(T_j-1) \cup E^W_{z_j}(T_j-1)$.
Again we have $|B_{z_j}|<2\eps_* |F_{z_j}(T_j-1)|$.
We are required to prove that $|D^W_{z_j,E}(T_j-1)| < \theta_{\nu^*} |F_{z_j}(T_j-1)|$.
The proof of Case D.1 when $z_j \in E$ is exactly as before,
so we just consider the case $z_j \notin E$. In Lemma \ref{3initial}
we divided this into Cases D.2 and D.3, but here we will give a unified argument.

Suppose that $z_j \notin E$. Then $E^{T_j}=E^{T_j-1}=\ov{E}$.
Also $x \in \ov{E} \cap VN_H(z_j)$, so $\nu^*>\nu$.
By Lemma \ref{plus-update} we have
$F(T_j)_{\ov{E}}=F(T_j-1)^{H+z_j}(y^c)_{\ov{E}}\sm y
= (F(T_j-1)^{H+z_j}_{\ov{E}z_j^c}\sm y)(y^c)$.
Now $W_j=W_{j-1} \cap F_{xz_j}(T_j-1)(y)
= W_{j-1} \cap F(T_j-1)^{H+z_j}_{xz_j^c}(y^c)$, so
\[F_{\ov{E}}(T_j)[W_j] = (F(T_j-1)^{H+z_j}_{\ov{E}z_j^c}\sm y)(y^c)[W_j]
= (F(T_j-1)^{H+z_j}_{\ov{E}z_j^c}[W_{j-1}]\sm y)(y^c).\]
Since $M_{\ov{E},E}(T_j) = M_{\ov{E},E}(T_j-1) \cap F_{\ov{E}}(T_j)$ we have
$$D^W_{z_j,E}(T_j-1) = \left\{ y \in F_{z_j}(T_j-1):
\frac{|M_{\ov{E},E}(T_j-1) \cap (F(T_j-1)^{H+z_j}_{\ov{E}z_j^c}[W_{j-1}]\sm y)(y^c)|}{
|(F(T_j-1)^{H+z_j}_{\ov{E}z_j^c}[W_{j-1}]\sm y)(y^c)|} > \theta_{\nu^*} \right\}.$$
Then
\begin{align*}
\Sa & := \sum_{y \in D^W_{z_j,E}(T_j-1)}
|M_{\ov{E},E}(T_j-1) \cap (F(T_j-1)^{H+z_j}_{\ov{E}z_j^c}[W_{j-1}]\sm y)(y^c)| \\
& > \theta_{\nu^*} \sum_{y \in D^W_{z_j,E}(T_j-1)\sm B'_{z_j}} |(F(T_j-1)^{H+z_j}_{\ov{E}z_j^c}[W_{j-1}]\sm y)(y^c)| \\
& >  (1-2\eps_*)\theta_{\nu^*} (|D^W_{z_j,E}(T_j-1)|-2\eps_*|F_{z_j}(T_j-1)|)\
|F(T_j-1)^{H+z_j}_{\ov{E}z_j^c}[W_{j-1}]|/|F_{z_j^c}(T_j-1)| \\
& = (1-2\eps_*)\theta_{\nu^*} (|D^W_{z_j,E}(T_j-1)|/|F_{z_j}(T_j-1)|-2\eps_*)
|F(T_j-1)^{H+z_j}_{\ov{E}z_j^c}[W_{j-1}]|.
\end{align*}
Here we used $|F_{z_j^c}(T_j-1)|=|F_{z_j}(T_j-1)|$ and Lemma \ref{plus-props}
with $\GG=(W_{j-1},\{\es\})$, as usual denoting the exceptional set by $B'_{z_j}$;
the factor $1-2\eps_*$ rather than $1-\eps_*$
accounts for the error from deleting $y$ (which has a lower order of magnitude).
We also have $\Sa \le \sum_{y \in F_{z_j}(T_j-1)}
|M_{\ov{E},E}(T_j-1) \cap F(T_j-1)^{H+z_j}_{\ov{E}z_j^c}[W_{j-1}](y^c)|$.
This counts all pairs $(y,P)$ with $P \in M_{\ov{E},E}(T_j-1)[W_{j-1}]$, $y \in F_{z_j}(T_j-1)$
and $Py^c \in F(T_j-1)^{H+z_j}_{\ov{E}z_j^c}[W_{j-1}]$, so we can rewrite it as
$\Sa \le \sum_{P \in M_{\ov{E},E}(T_j-1)[W_{j-1}]} |F(T_j-1)^{H+z_j}_{\ov{E}z_j^c}[W_{j-1}](P)|$.
By Lemma \ref{plus-props}, we have \[|F(T_j-1)^{H+z_j}_{\ov{E}z_j^c}[W_{j-1}](P)| =
(1 \pm \eps_*)\frac{|F(T_j-1)^{H+z_j}_{\ov{E}z_j^c}[W_{j-1}]|}{|F(T_j-1)^{H+z_j}_{\ov{E}}[W_{j-1}]|}\]
for all but at most $\eps_*|F(T_j-1)^{H+z_j}_{\ov{E}}[W_{j-1}]|$
sets $P \in F(T_j-1)^{H+z_j}_{\ov{E}}[W_{j-1}]$.

Since $F(T_j-1)^{H+z_j}_{\ov{E}} = F_{\ov{E}}(T_j-1)$ we have
$$\Sa \le |M_{\ov{E},E}(T_j-1)[W_{j-1}]|(1+\eps_*)
\frac{|F(T_j-1)^{H+z_j}_{\ov{E}z_j^c}[W_{j-1}]|}{|F_{\ov{E}}(T_j-1)[W_{j-1}]|}
+ \eps_*|F_{\ov{E}}(T_j-1)[W_{j-1}]||F_{z_j}(T_j-1)|.$$
Combining this with the lower bound on $\Sa$  we obtain
\begin{align*}
& (1-2\eps_*)\theta_{\nu^*}\left(\frac{|D^W_{z_j,E}(T_j-1)|}{|F_{z_j}(T_j-1)|}-2\eps_*\right) \\
& < (1+\eps_*)\frac{|M_{\ov{E},E}(T_j-1)[W_{j-1}]|}{|F_{\ov{E}}(T_j-1)[W_{j-1}]|}
+ \eps_*\frac{|F_{\ov{E}}(T_j-1)[W_{j-1}]||F_{z_j}(T_j-1)|}{|F(T_j-1)^{H+z_j}_{\ov{E}z_j^c}[W_{j-1}]|}.
\end{align*}
Now $|M_{\ov{E},E}(T_j-1)[W_{j-1}]| < \theta_{\nu} |F_{\ov{E}}(T_j-1)|$ by $A_{2,j-1}$ and
\[\frac{|F_{\ov{E}}(T_j-1)[W_{j-1}]||F_{z_j}(T_j-1)|}{|F(T_j-1)^{H+z_j}_{\ov{E}z_j^c}[W_{j-1}]|}
 \le d_u^{-2^k} \ll \eps_*^{-1},\]
by Lemma \ref{plus-props}, so
\begin{equation}
\frac{|D^W_{z_j,E}(T_j-1)|}{|F_{z_j}(T_j-1)|} <
\frac{(1+\eps_*)\theta_{\nu}+\sqrt{\eps_*}}{(1-2\eps_*)\theta_{\nu^*}}
+ 2\eps_* < \theta_{\nu^*}. \tag*{$\Box$}
\end{equation}

Our final lemma for the initial phase is similar to the previous one,
but instead of asking for a set $W$ of vertices to contain
an available vertex for $x$, we ask for some particular vertex $v$
to be available for $x$.

\begin{lemma} \label{x-to-v}
For any $v \in V_x$, conditional on any embedding of the vertices $\{s(u): u < T_1\}$
that does not use $v$, with probability at least $p$ we have
$\phi(H(x)) \sub (G \sm M)(v)$, so $v \in A_x$.
\end{lemma}

\nib{Proof.}
We follow the proof of Lemma \ref{3-x-to-v}, indicating the necessary modifications.
We note again that since $B \cup N$ is disjoint from $X_*$,
restricted positions have no effect on any $z \in VN_H(x)\cup \{x\}$.
However, we need to consider all vertices within distance $3$ of $x$,
so some of these may have restricted positions.
The bound in equation (\ref{eq:fresh}) will be adequate to deal with these.
By Remark \ref{rem:nhood} we also need to clarify the meaning of neighbourhood
constructions, which are potentially ambiguous:
$F(T_j)(v)_S$ is $F(T_j)_{Sx}(v)$ when $Sx \in H$ 
or undefined when $Sx \notin H$.
\COMMENT{Need to say this: earlier comment was just for plus nhoods.}

For $z \in VN_H(x)$ we define $\alpha_z$ as before.
We also define $\alpha_S=1$ for $S\in H$ with $|S|\ge 2$.\index{$\alpha_S$}
As before, we define $\nu''_S(t)$ similarly to $\nu'_S(t)$,
replacing `embedded' with `allocated'.
Suppose $1 \le j \le g$ and that we are considering the embedding of $z_j$.
We interpret quantities at time $T_j$ with the embedding $\phi(z_j)=y$,
for some as yet unspecified $y \in F_{z_j}(T_j-1)$.
We define $E^v_{z_j}(T_j-1)$ as before, except that the condition
for $|S|=2$ now applies whenever $|S| \ge 2$.
We define $Y$, $H'$, $F(T_j)^{Z*v}_{S^\le}$ and $D^{Z*v}_{z_j,E}(T_j-1)$ as before.
Properties (i-iv) of $Z$ hold as before.
We define the events $A_{i,j}$ as before.

Recall that we used the notation $Z \sub Y$, $Z'=Zz_j$,
$I = \{S \sub Z: S \in H(x)\}$, $I' = \{S \sub Z': S \in H(x)\}$.
Here we also define $J = \{S \sub Zx: S \in H\}$ and
$J' = \{S \sub Z'x: S \in H\}$.
Using the plus complex notation we can write
\begin{equation}\label{eq:+}\tag{$+_{\ref{x-to-v}}$}
F(T_j)^{Z*v}_{S^\le} = F(T_j)^{J+v}(v^c)_{S^\le}.
\end{equation}
To see this, we need to show that
$P \in F(T_j)^{Z*v}_S$ $\Lra$ $Pv^c \in F(T_j)^{J+v}_{Sx^c}$.
Recall that $F(T_j)^{Z*v}_S$ consists of all sets $P \in F_S(T_j)$
such that $P_{S'} v \in F_{S' x}(T_j)$ for all $S' \sub S$ with $S' \in I$.
Also, from Definition \ref{def-plus}(v) we have $Pv^c \in F(T_j)^{J+v}_{Sx^c}$
if and only if $P \in F(T_j)_S$ (this is the restriction from $F$)
and $P_{S'}v^c \in (F(T_j)^+_{J^c})_{S'x^c}=F(T_j)^c_{S'x^c}$ for every 
$S' \sub S$ with $S'x \in J$. Since $S' \in I$ $\Lra$ $S'x \in J$,
this is equivalent to the condition for $F(T_j)^{Z*v}_S$, as required.

The proof of Claim A is similar to before. Instead of the bound
$|F_z(T_0)| > (1-2\sqrt{\delta_B})|V_z|$ from equation (\ref{eq:3fresh}) we use
$|F_z(T_0)| > (1-\delta_B^{1/3})|V_z|$ from equation (\ref{eq:fresh}).
We again have $|F_z(T_0) \cap G(v)_z| > (1-\delta_B^{1/3})|G(v)_z|$
for $z \in VN_H(x)$ by Lemma \ref{uniform}, using the fact that $z \notin X_*$.
Again, $A_{4,0}$ holds by definition and Lemma \ref{uniform}
implies that $A_{1,0}$ holds with high probability.
The arguments for $A_{2,0}$ and $A_{3,0}$ are as before,
modifying the absolute constants to take account of their dependence on $k$.
(In the $A_{2,0}$ argument we have $F(0)^{Z*v}_E = G^{J_v}_E$,
where we let $J$ also denote the submulticomplex $\{i^*(S): S \in J\}$ of $R$.)
The proofs of Claims B and E and the conclusion of the lemma are also similar to before.
As usual we replace triple by $k$-tuple, $2D$ by $(k-1)D$ and $12D$ by $k^3D$.
Also, in Claim E we previously estimated $2D^2$ choices for $E$ then $8$ choices for $Z$;
here we estimate $(k-1)D^2$ choices for $E$ then $2^k$ choices for $Z$.
To complete the proof of the lemma it remains to establish Claims C and D.
These require more substantial modifications, so we will give the details here.

We start with Claim C. Suppose that $A_{1,j-1}$ and $A_{3,j-1}$ hold.
We are required to prove that
$|E^v_{z_j}(T_j-1) \sm E_{z_j}(T_j-1)| < \eps_* |F_{xz_j}(T_j-1)(v)|$.
For any $S \in H$ we write $\nu''_S=\nu''_S(T_j-1)$ and $\nu^*_S=\nu''_S(T_j)$.
Consider any unembedded $\es \ne S \in H(x) \cap H(z_j)$.
Note that $\nu^*_S > \max\{\nu''_S,\nu''_{Sx},\nu''_{Sz_j}\}$.
Applying the definitions, it suffices to show that for all but at most
$\eps_{k^3D,3} |F_{xz_j}(T_j-1)(v)|$ vertices $y \in F_{xz_j}(T_j-1)(v) \sm E_{z_j}(T_j-1)$,
$F_{Sx}(T_j)(v)$ is $\eps_{\nu^*_S,1}$-regular with
$d_S(F(T_j)(v))=(1\pm\eps_{\nu^*_S,1})d_{Sx}(F(T_j))d_S(F(T_j))\alpha_S$.
We claim that
\[F_{Sx}(T_j)(v) = F(T_j-1)^{H+z_j}(vy^c)_S.\]
To see this we apply Definition \ref{def-plus},
which tells us that for $P \in F(T_j-1)_S$,
$v_0 \in F(T_j-1)_x$ and $y_0 \in F(T_j-1)_{z_j}$
we have $Pv_0y_0^c \in F(T_j-1)^{H+z_j}_{Sxz_j^c}$
exactly when $Pv_0 \in F(T_j-1)_{Sx}$ and
$(Pv_0)_{S'}y_0 \in F(T_j-1)_{S'z_j}$ for all $S' \sub Sx$, $S' \in H(z_j)$.
Therefore $F(T_j-1)^{H+z_j}(vy^c)_S$ consists of all $P \in F(T_j-1)_S$
such that $Pv \in F(T_j-1)_{Sx}$ and
$(Pv)_{S'} \in F(T_j-1)_{S'z_j}(y)$ for all $S' \sub Sx$, $S' \in H(z_j)$.
By Definition \ref{def-update} this is equivalent to
$Pv \in F(T_j)_{Sx}$, i.e.\ $P \in F_{Sx}(T_j)(v)$ as claimed.
(We do not need to delete $y$ as $S$ and $x$ are in $H(z_j)$.)
\COMMENT{Formerly:
$F(T_j-1)^{H(z_j)+z_j}_{Sxz_j^c}$ is either $F(T_j-1)_{Sxz_j}^c$ if $Sxz_j \in H$
or obtained from $K(V)_{Sxz_j^c}$ by regular restriction if $Sxz_j \notin H$.
In the first case it is $\eps_{\nu''_{Sxz_j},1}$-regular by Lemma \ref{exceptional}
and in the second case it is $\eps$-regular by Lemma \ref{restrict'}.
}

By Definition \ref{def-plus}, $F(T_j-1)^{H+z_j}_{Sxz_j^c}$ is
$F(T_j-1)_{Sxz_j}^c$ if $Sxz_j \in H$
or consists of all $Py^c$ with $P \in F(T_j-1)_{Sx}$, $y \in F(T_j-1)_{z_j}$
and $P_{S'}y \in F(T_j-1)_{S'z_j}$ for $S'\sub Sx$ with $S' \in H(z_j)$.
Thus $F(T_j-1)^{H+z_j}_{Sxz_j^c}(v)$ is $F(T_j-1)_{Sxz_j}^c(v)$ if $Sxz_j \in H$
or $F(T_j-1)[F(T_j-1)_{Sxz_j^\le}(v)]_{Sz_j}^c$ if $Sxz_j \notin H$.
Either way we can see that is $\eps_{\nu''_{Sz_j},2}$-regular:
in the first case we write $F(T_j-1)_{Sxz_j}(v)=F(T_j-1)^{Sz_j*v}_{Sz_j}$ and use $A_{3,j-1}$;
in the second case we use Lemma \ref{restrict'} and the fact that $F(T_j-1)_{Sz_j}^c$
is a copy of $F(T_j-1)_{Sz_j}$,
which is $\eps_{\nu''_{Sz_j},1}$-regular by Lemma \ref{exceptional}.
Also by $A_{3,j-1}$,
$d_{Sz_j^c}(F(T_j-1)^{H+z_j}(v))$
is $(1\pm\eps_{\nu''_{Sz_j},2}) d_{Sz_j}(F(T_j-1))d_{Sxz_j}(F(T_j-1))$ if $Sxz_j \in H$
or $(1\pm\eps_{\nu''_{Sz_j},2}) d_{Sz_j}(F(T_j-1))$ if $Sxz_j \notin H$.\COMMENT{Details?}
Similarly, $F(T_j-1)^{H+z_j}_{Sx}(v)$ is $\eps_{\nu''_S,2}$-regular with
$d_S(F(T_j-1)^{H+z_j}(v)) = (1\pm\eps_{\nu''_{S},2})d_S(F(T_j-1))d_{Sx}(F(T_j-1))\alpha_S$.
\COMMENT{1. More plus props? 2. Apply $A_{3,j-1}$ with $Z=Sz_j$.}

Since $x \in H(z_j)$ we have $F(T_j-1)^{H+z_j}(v)_{z_j^c}=F_{xz_j}(T_j-1)^c(v)$.
By Lemma \ref{neighbour1}, for all but at most $\sum_{S' \sub Sxz_j} \eps_{\nu''_{S'},3} |F_{xz_j}(T_j-1)(v)|$
vertices $y \in F_{xz_j}(T_j-1)(v)$, writing $\eta = \eps_{\nu''_{S},3} + \eps_{\nu''_{Sz_j},3}$,
$F_{Sx}(T_j)(v)=F(T_j-1)^{H+z_j}(v)(y^c)_S$ is $\eta$-regular
with\COMMENT{also more details here for 3-graphs...}
\begin{align*}
& d_S(F(T_j)(v))= d_S(F(T_j-1)^{H+z_j}(vy^c)) \\
& = (1\pm \eta)d_S(F(T_j-1)^{H+z_j}(v))d_{Sz_j}(F(T_j-1)^{H+z_j}(v)). \\
& = (1\pm 3\eta)d_S(F(T_j-1))d_{Sx}(F(T_j-1))\alpha_S \cdot
d_{Sz_j}(F(T_j-1))d_{Sxz_j}(F(T_j-1))^{1_{Sxz_j \in H}}.
\end{align*}
This gives the required regularity property for $F_{Sx}(T_j)(v)$. Next, $(*_{\ref{alg}})$ gives
\begin{align*}
&d_S(F(T_j))=(1 \pm \eps_{\nu^*_S,0})d_S(F(T_j-1))d_{Sz_j}(F(T_j-1)) \mbox{ and} \\
& d_{Sx}(F(T_j))=(1\pm \eps_{\nu^*_{Sx},0})d_{Sx}(F(T_j-1))d_{Sxz_j}(F(T_j-1))^{1_{Sxz_j \in H}},
\end{align*}
so since $\nu^*_S > \max\{\nu''_S,\nu''_{Sx},\nu''_{Sz_j}\}$
we have $d_S(F(T_j)(v)) = (1\pm\eps_{\nu^*_S,1})d_{Sx}(F(T_j))d_S(F(T_j))\alpha_S$.
This proves Claim C.

It remains to prove Claim D.
Suppose that $A_{1,j-1}$, $A_{2,j-1}$ and $A_{3,j-1}$ hold,
$E \in U(z_j)$, $Z \sub E$ and $Z' = Z \cup z_j$.
As before we write $\ov{E}=E^{T_j-1}$,
$\nu=\nu''_{\ov{E}}(T_j-1)$, $\nu^*=\nu''_{E^{T_j}}(T_j)$,
$B_{z_j} = E_{z_j}(T_j-1) \cup E^v_{z_j}(T_j-1)$.
Since $E \in U(z_j)$ we have $\nu^*>\nu$.
We are required to prove that $|D^{Z'*v}_{z_j,E}(T_j-1)| < \theta_{\nu^*} |F_{xz_j}(T_j-1)(v)|$.
The proof of Case D.1 when $z_j \in E$ is exactly as before,
so we just consider the case $z_j \notin E$.

Suppose that $z_j \notin E$. Then $E^{T_j}=E^{T_j-1}=\ov{E}$.
By Lemma \ref{plus-update} we have 
\[F(T_j)_{\ov{E}}=F(T_j-1)^{H+z_j}(y^c)_{\ov{E}}\sm y
= (F(T_j-1)^{H+z_j}_{\ov{E}z_j^c}\sm y)(y^c).\]
We claim that
\COMMENT{** double plus complex, both index r+1, interact?}
\begin{equation}\label{eq:d}\tag{$\dagger_{\ref{x-to-v}}$}
F(T_j)^{Z*v}_{\ov{E}^\le} = F(T_j-1)^{H+z_j}(y^c)^{J+v}(v^c)_{\ov{E}^\le} \sm y
= F(T_j-1)^{J'+v}(v^c)^{H+z_j}(y^c)_{\ov{E}^\le} \sm y.
\end{equation}
For the first equality we use ($+_{\ref{x-to-v}}$) to get
$F(T_j)^{Z*v}_{\ov{E}^\le}=F(T_j)^{J+v}(v^c)_{\ov{E}^\le}$
and substitute $F(T_j) = F(T_j-1)^{H+z_j}(y^c)\sm y$
from Lemma \ref{plus-update}.%
\COMMENT{Don't specify subscript so as to include $x$.}
For the second equality we apply Definition \ref{def-plus} as follows.
Suppose $S \sub \ov{E}$. We have $P \in F(T_j-1)^{H+z_j}(y^c)^{J+v}(v^c)_S$
exactly when $P \in F(T_j-1)^{H+z_j}(y^c)_S$
and $P_{S'} v \in F(T_j-1)^{H+z_j}(y^c)_{S'x}$ for all $S' \sub S$ with $S'x \in J$.
Equivalently, $P \in F(T_j-1)_S$, $P_U y \in F(T_j-1)_{Uz_j}$ for $U \sub S$ with $U \in H(z_j)$,
$P_{S'} v \in F(T_j-1)_{S'x}$ and $(Pv)_{S''} y \in F(T_j-1)_{S''z_j}$
for $S' \sub S$ with $S'x \in J$ and $S'' \sub S'x$ with $S'' \in H(z_j)$.
Note that it is equivalent to assume $x \in S''$, as otherwise the
$S''$ condition is covered by the $U$ condition.
Writing $W=S'' \sm x$ and using $J=\{A \sub Zx: A \in H\}$ we have 

$P \in F(T_j-1)^{H+z_j}(y^c)^{J+v}(v^c)_S$ if and only if %$\Lra$
$P \in F(T_j-1)_S$, $P_U y \in F(T_j-1)_{Uz_j}$ for $U \sub S$ with $U \in H(z_j)$,
$P_{S'} v \in F(T_j-1)_{S'x}$ for $S' \sub S \cap Z$, $S' \in H(x)$,
and $P_Wvy \in F(T_j-1)_{Wxz_j}$ for $W \sub S \cap Z$ with $W \in H(xz_j)$.

On the other hand, we have $P \in F(T_j-1)^{J'+v}(v^c)^{H+z_j}(y^c)_S$
exactly when $P \in F(T_j-1)^{J'+v}(v^c)_S$
and $P_U y \in F(T_j-1)^{J'+v}(v^c)_{Uz_j}$ for all $U \sub S$ with $U \in H(z_j)$.
Equivalently, $P \in F(T_j-1)_S$, $P_{S'}v \in F(T_j-1)_{S'x}$ for $S' \sub S$ with $S'x \in J'$,
$P_U y \in F(T_j-1)_{Uz_j}$ and $(Py)_{U'}v \in F(T_j-1)_{U'x}$
for $U \sub S$ with $U \in H(z_j)$ and $U' \sub Uz_j$ with $U'x \in J'$.
Note that it is equivalent to assume $z_j \in U'$, as otherwise the
$U'$ condition is covered by the $S'$ condition.
Writing $W=U' \sm z_j$ and using $J'=\{A \sub Z'x: A \in H\}$ we have

$P \in F(T_j-1)^{J'+v}(v^c)^{H+z_j}(y^c)_S$ if and only if
$P \in F(T_j-1)_S$, $P_{S'}v \in F(T_j-1)_{S'x}$ for $S' \sub S \cap Z' = S \cap Z$ with $S' \in H(x)$,
$P_U y \in F(T_j-1)_{Uz_j}$ for $U \sub S$ with $U \in H(z_j)$,
and $P_W yv \in F(T_j-1)_{Wz_jx}$ for $W \sub S \cap Z$ with $W \in H(xz_j)$.

This proves (\ref{eq:d}). Now, since $M_{\ov{E},E}(T_j) = M_{\ov{E},E}(T_j-1) \cap F_{\ov{E}}(T_j)$
we have $$ D^{Z'*v}_{z_j,E}(T_j-1) = \left\{ y \in F_{z_j}(T_j-1):
\frac{|M_{\ov{E},E}(T_j-1) \cap (F(T_j-1)^{J'+v}(v^c)^{H+z_j}_{\ov{E}z_j^c} \sm y)(y^c)|}{
|(F(T_j-1)^{J'+v}(v^c)^{H+z_j}_{\ov{E}z_j^c} \sm y)(y^c)|} > \theta_{\nu^*} \right\}. $$
Writing $B'_{z_j}$ for the set of vertices
$y \in F(T_j-1)^{J'+v}(v^c)^{H+z_j}_{z_j} = F(T_j-1)_{xz_j}(v)$ for which we do not have
\[|F(T_j-1)^{J'+v}(v^c)^{H+z_j}_{\ov{E}z_j^c}(y^c)|=
(1\pm\eps_*)|F(T_j-1)^{J'+v}(v^c)^{H+z_j}_{\ov{E}z_j^c}|/|F(T_j-1)_{xz_j}(v)|,\]
we have $|B'_{z_j}| < \eps_*|F(T_j-1)_{xz_j}(v)|$ by Lemma \ref{average}. 
Here we use the fact that the `double plus' complex is $\eps_{k^3D,3}$-regular; 
this proof of this is similar to that of Lemma \ref{plus-props}(i):
$F(T_j-1)^{J'+v}$ is $\eps_{k^3D,1}$-regular,
$F(T_j-1)^{J'+v}(v^c)$ is $\eps_{k^3D,2}$-regular by Lemma \ref{neighbour1},
$F(T_j-1)^{J'+v}(v^c)^{H+z_j}$ is $\eps_{k^3D,3}$-regular.
Then%
\COMMENT{Formerly: Since $|B_{z_j}| < 2\eps_*|F(T_j-1)_{z_j}| < \sqrt{\eps_*}|F(T_j-1)_{xz_j}(v)|$ we have}
\begin{align*}
\Sa & := \sum_{y \in D^{Z'*v}_{z_j,E}(T_j-1)} |M_{\ov{E},E}(T_j-1)
\cap (F(T_j-1)^{J'+v}(v^c)^{H+z_j}_{\ov{E}z_j^c} \sm y)(y^c)| \\
& > \theta_{\nu^*} \sum_{y \in D^{Z'*v}_{z_j,E}(T_j-1)\sm B'_{z_j}}
|(F(T_j-1)^{J'+v}(v^c)^{H+z_j}_{\ov{E}z_j^c} \sm y)(y^c)| \\
& >  (1-2\eps_*)\theta_{\nu^*} (|D^{Z'*v}_{z_j,E}(T_j-1)|-\eps_*|F(T_j-1)_{xz_j}(v)|)\
\frac{|F(T_j-1)^{J'+v}(v^c)^{H+z_j}_{\ov{E}z_j^c}|}{|F(T_j-1)_{xz_j}(v)|} \\
& = (1-2\eps_*)\theta_{\nu^*} (|D^{Z'*v}_{z_j,E}(T_j-1)|/|F(T_j-1)_{xz_j}(v)|-\eps_*)
|F(T_j-1)^{J'+v}(v^c)^{H+z_j}_{\ov{E}z_j^c}|.
\end{align*}
We also have
\[\Sa \le \sum_{y \in F_{xz_j}(T_j-1)(v)}
|M_{\ov{E},E}(T_j-1) \cap F(T_j-1)^{J'+v}(v^c)^{H+z_j}_{\ov{E}z_j^c}(y^c)|.\]
This counts all pairs $(y,P)$ with
$P \in M_{\ov{E},E}(T_j-1) \cap F(T_j-1)^{J'+v}(v^c)^{H+z_j}_{\ov{E}}$,
$y \in F_{xz_j}(T_j-1)(v)$ and $Py^c \in F(T_j-1)^{J'+v}(v^c)^{H+z_j}_{\ov{E}z_j^c}$,
so we can rewrite it as
\[\Sa \le \sum_{P \in M_{\ov{E},E}(T_j-1)\cap F(T_j-1)^{J'+v}(v^c)^{H+z_j}_{\ov{E}}}
|F(T_j-1)^{J'+v}(v^c)^{H+z_j}_{\ov{E}z_j^c}(P)|.\]
For all but at most $\eps_* |F(T_j-1)^{J'+v}(v^c)^{H+z_j}_{\ov{E}}|$
sets $P \in F(T_j-1)^{J'+v}(v^c)^{H+z_j}_{\ov{E}}$, we have
\[|F(T_j-1)^{J'+v}(v^c)^{H+z_j}_{\ov{E}z_j^c}(P)|=(1\pm\eps_*)
\frac{|F(T_j-1)^{J'+v}(v^c)^{H+z_j}_{\ov{E}z_j^c}|}{|F(T_j-1)^{J'+v}(v^c)^{H+z_j}_{\ov{E}}|}.\]
(Recall that the `double plus' complex is $\eps_{k^3D,3}$-regular and use Lemma \ref{average}.)
Therefore
\begin{align*}
\Sa \le & \ |M_{\ov{E},E}(T_j-1) \cap F(T_j-1)^{J'+v}(v^c)^{H+z_j}_{\ov{E}}| \cdot (1+\eps_*)
\frac{|F(T_j-1)^{J'+v}(v^c)^{H+z_j}_{\ov{E}z_j^c}|}{|F(T_j-1)^{J'+v}(v^c)^{H+z_j}_{\ov{E}}|} \\
& +\ \eps_* |F(T_j-1)^{J'+v}(v^c)^{H+z_j}_{\ov{E}}| |F(T_j-1)_{xz_j}(v)|.
\end{align*}
Combining this with the lower bound on $\Sa$ we have
\begin{align*}
& (1-2\eps_*)\theta_{\nu^*} (|D^{Z'*v}_{z_j,E}(T_j-1)|/|F(T_j-1)_{xz_j}(v)|-\eps_*) \\
& \le (1+\eps_*) \frac{|M_{\ov{E},E}(T_j-1) \cap F(T_j-1)^{J'+v}(v^c)^{H+z_j}_{\ov{E}}|}{
|F(T_j-1)^{J'+v}(v^c)^{H+z_j}_{\ov{E}}|} \\
& \ +\ \eps_* \frac{|F(T_j-1)^{J'+v}(v^c)^{H+z_j}_{\ov{E}}||F(T_j-1)_{xz_j}(v)|}{
|F(T_j-1)^{J'+v}(v^c)^{H+z_j}_{\ov{E}z_j^c}|}.
\end{align*}
Since $z_j \notin \ov{E}$, Definition \ref{def-plus} and equation (\ref{eq:+}) give
$F(T_j-1)^{J'+v}(v^c)^{H+z_j}_{\ov{E}}=F(T_j-1)^{J'+v}(v^c)_{\ov{E}}
=F(T_j-1)^{J+v}(v^c)_{\ov{E}}=F(T_j-1)^{Z*v}_{\ov{E}}$,
so
\[ |M_{\ov{E},E}(T_j-1) \cap F(T_j-1)^{J'+v}(v^c)^{H+z_j}_{\ov{E}}|
< \theta'_{\nu}|F(T_j-1)^{J'+v}(v^c)^{H+z_j}_{\ov{E}}| \]
by $A_{2,j-1}$. We also have
\[ \frac{|F(T_j-1)^{J'+v}(v^c)^{H+z_j}_{\ov{E}}||F(T_j-1)_{xz_j}(v)|}{
|F(T_j-1)^{J'+v}(v^c)^{H+z_j}_{\ov{E}z_j^c}|} \le d_u^{-2^k} \ll \eps_*^{-1},\]
similarly to Lemma \ref{plus-props}(v) 
%but with $F(t)$ replaced by $F(T_j-1)^{J'+v}(v^c)$ (if (vi) or double plus without...
(the statement is only for $F(t)$, but the estimate
for the densities is valid for any $\eps_{k^3D,3}$-regular complex). 
Therefore\COMMENT{more details?}
\begin{equation}
\frac{|D^{Z'*v}_{z_j,E}(T_j-1)|}{|F_{xz_j}(T_j-1)(v)|} < \
\frac{(1+\eps_*)\theta'_{\nu}+\sqrt{\eps_*}}{(1-2\eps_*)\theta_{\nu^*}}
+ \eps_* < \theta_{\nu^*}. \tag*{$\Box$}
\end{equation}

The analysis for the conclusion of the algorithm is very similar to that
for $3$-graphs, with the usual modifications to absolute constants
to account for their dependence on $k$. The only important difference is to
take account of restricted positions. Lemma \ref{3main} (the `main lemma') holds, 
provided that we assume that the set $Y$ is disjoint from the set $X_*$ of vertices
with restricted positions. We applied Lemma \ref{3main} in the proof of Theorem \ref{3final}
to show that it is very unlikely that the iteration phase aborts with failure. 
This required an estimate for the probability that a given set $Y \sub X_i$ of 
size $\delta_Q|X_i|$ is contained in $Q(T)$ (the vertices that have ever been queued).
Since $|X_* \cap X_i| \le c|X_i|$ and $c \ll \delta_Q$, we can apply the
same argument to $Y\sm X_*$, which has size at least $\frac{1}{2}\delta_Q|X_i|$.
The remainder of the proof of Theorem \ref{3final} is checking Hall's condition
for the sets $\{A'_z: z\in S\}$, where $S \sub X_i(T) \sub B$.
No changes are required here, as we chose the buffer $B$ to be disjoint from $X_*$.
This completes the proof of Theorem \ref{blowup}.

\subsection{Obtaining super-regularity and robust universality}

We conclude with some lemmas that will be useful when applying the blow-up lemma.
We start with the analogue of Lemma \ref{3del}, showing that one can delete
a small number of vertices to enforce super-regularity.
\COMMENT{Say `at most' for max degrees of $H$ or $R$.}

\begin{lemma} \label{del} Suppose
$0 < \eps_0 \ll \eps \ll \eps' \ll d_a \ll \theta \ll 1/D_R, d, 1/k$,
we have a multi-$k$-complex $R$ on $[r]$ with maximum degree at most $D_R$,
and $(G,M)$ is an $R$-indexed marked complex on $V = \cup_{i \in R} V_i$
such that when defined $G_S$ is $\eps_0$-regular, $|M_S| \le \theta|G_S|$,
$d_S(G) \ge d_a$ when $|S|\ge 2$ and $d_S(G) \ge d$ when $|S|=k$.
Then we can delete at most $2\theta^{1/3}|G_i|$ vertices from each $G_i$, $i \in R$
to obtain an $(\eps,\eps',d_a/2,2\sqrt{\theta},d/2)$-super-regular marked complex $(G^\sharp,M^\sharp)$.
\end{lemma}

\nib{Proof.} The proof is similar to that of Lemma \ref{3del}, so we will just sketch the necessary
modifications. Similarly to before, for any $i,S$ such that $i \in S$, $|S|=k$ and
$G_S$ is defined we let $Y_{i,S\sm i}$ be the set of vertices $v \in G_i$ for which
$|M_S(v)|>\theta|G_S(v)|$. We also let $Z_{i,S\sm i}$ be the set of vertices $v \in G_i$
such that we do not have $|G_S(v)|=(1\pm\eps)|G_S|/|G_i|$ and $G_{S^\le}(v)$
is $\eps$-regular with $d_{S'\sm i}(G_{S^\le}(v))= (1\pm \eps)d_{S'\sm i}(G)d_{S'}(G)$
for $i \subn S'\sub S$.
As before we have $|Z_{i,S\sm i}|<\eps|G_i|$ and $|Y_{i,S\sm i}|<2\sqrt{\theta}$.
Next, consider any $k$-tuple $S$ containing at least one neighbour of $i$ in $R$
such that $G_S$ is defined, and any subcomplex $I$ of $Si^\le$ such that $G_{S'}$ is defined
for all $S'\in I$. We let $Y^I_{i,S}$ be the set of vertices $v \in G_i$
for which $|(M \cap G^{I_v})_S| > \sqrt{\theta}|G^{I_v}_S|$.
Note that we only need to consider $S$ containing at least one neighbour of $i$ in $R$,
as otherwise we have $G^{I_v}_S=G_S$, and $|M_S| \le \theta|G_S|$ by assumption.
We let $Z_{i,S}$ be the union of all $Z_{i,S'}$ with $i \notin S' \sub S$, $|S'|=k-1$.
Recall from Lemma \ref{pre-to-plus} that $G^{I_v}_S = G^{I+i}(v^c)_S$.
If $v \notin Z_{i,S}$ then $G^{I+i}_{Si^c}$ is $\sqrt{\eps}$-regular by regular restriction
and $|G^{I+i}(v^c)_S| = (1\pm\eps')|G^{I+i}_{Si^c}|/|G_i|$ by Lemma \ref{average}.%
\COMMENT{Formerly: (the calculations in). Why does ref want explanation of regularity?}
Then
$\Sa = \sum_{v \in Y^I_{i,S}} |(M \cap G^{I_v})_S| = \sum_{v \in Y^I_{i,S}} |(M \cap G^{I+i}(v^c))_S|$
satisfies\COMMENT{$G^{I+i}_{Si^c}$ is well-defined.}
\[\Sa > \sqrt{\theta} \sum_{v \in Y^I_{i,S}\sm Z_{i,S}} |(M \cap G^{I+i}(v^c))_S|
> \sqrt{\theta}(|Y^I_{i,S}|-k\eps|G_i|)(1-\eps')|G^{I+i}_{Si^c}|/|G_i|.\]
We also have $\Sa \le \sum_{v \in G_i}  |(M \cap G^{I+i}(v^c))_S|$, which
counts all pairs $(v,P)$ with $P \in M_S$, $v \in G_i$ and $Pv^c \in G^{I+i}$.
By Lemma \ref{average} we have $|G^{I+i}(P)_{i^c}|=(1\pm\eps')|G^{I+i}_{Si^c}|/|G^{I+i}_S|$
for all but at most $\eps'|G^{I+i}_S|$ sets $P \in G^{I+i}_S$.
Since $G^{I+i}_S=G_S$ we have
\[\Sa \le \sum_{P \in M_S} |G^{I+i}(P)_{i^c}| \le |M_S|(1+\eps')|G^{I+i}_{Si^c}|/|G_S|
+ \eps'|G_S||G_i|.\]
Combining this with the lower bound on $\Sa$ and using $|M_S|\le \theta|G_S|$
we obtain $|Y^I_{i,S}|/|G_i| < \frac{(1+\eps')\theta + \sqrt{\eps'}}{(1-\eps')\sqrt{\theta}} + k\eps
< 2\sqrt{\theta}$. Let $Y_i$ be the union of all such sets $Y_{i,S\sm i}$ and $Y^I_{i,S}$.
Since $\theta \ll 1/D_R$ we have $|Y_i| < \theta^{1/3}|G_i|$ as in Step 1 of Lemma \ref{3del}.
We define $Z'_{i,j}$, $Z_i$ and obtain $|Z_i|<\sqrt{\eps}|G_i|$ as in Step 2 of Lemma \ref{3del}.
Now we delete $Y_i \cup Z_i$ from $G_i$ for every $i \in R$; as in Step 3 of Lemma \ref{3del}
this gives an $(\eps,\eps',d_a/2,2\sqrt{\theta},d/2)$-super-regular marked complex $(G^\sharp,M^\sharp)$. \qed

The next lemma is analogous to Lemma \ref{3super-restrict};
we omit its very similar proof.

\begin{lemma}\label{super-restrict}{\bf (Super-regular restriction)}\index{super-regular restriction}
Suppose $0< \eps \ll \eps' \ll \eps'' \ll d_a \ll \theta \ll d, d', 1/k$,
we have a multi-$k$-complex $R$,
and $(G,M)$ is a $(\eps,\eps',d_a,\theta,d)$-super-regular
marked $R$-indexed complex on $V = \cup_{i \in R} V_i$ with $G_i=V_i$ for $i \in R$.
Suppose also that we have $V'_i \sub V_i$ for $i \in R$,
write $V' = \cup_{i \in R} V'_i$, $G'=G[V']$, $M'=M[V']$,
and that $|V'_i|\ge d'|V_i|$ and $|G'_S(v) \cap V'_i| \ge d'|G_S(v)|$
whenever $S\in R$ with $|S|=2$ lies over $i,j\in R$ and $v\in G_j$. 
Then $(G',M')$ is $(\eps',\eps'',d_a/2,\sqrt{\theta},d/2)$-super-regular.
\end{lemma}

More generally, the same proof shows that super-regularity is preserved on restriction
to a dense regular subcomplex $\GG$, provided that the singleton parts of $\GG$ have large \index{$\GG$}
intersection with every vertex neighbourhood. More precisely, suppose $(G,M)$ is as
in Lemma \ref{super-restrict} and $\GG$ is an $\eps'$-regular subcomplex of $G$
with $|\GG_S|\ge d'|G_S|$ when defined and $|G_S(v) \cap \GG_i| \ge d'|G_S(v)|$
whenever $S\in R$ with $|S|=2$ lies over $i,j\in R$, $v\in G_j$ and $\GG_i$ is defined.
Then $(G',M')$ is $(\eps',\eps'',d'd_a,\sqrt{\theta},d'd)$-super-regular.
Next we will reformulate the blow-up lemma in a more convenient
`black box' form. The following definition of robustly universal \index{black box}
is more general than that used for $3$-graphs, in that it allows
for restricted positions.% 
\COMMENT{Formerly: There is an extra parameter $c$ for these
conditions, so there is no conflict with the previous definition.}

\begin{defn}\label{def-robust} {\bf (Robustly universal)} \index{robustly universal}
Suppose $R$ is a multi-$k$-complex $R$ and
$J$ is an $R$-indexed complex on $Y = \cup_{i \in R} Y_i$ with $J_i=Y_i$ for $i \in R$.
We say that $J$ is {\em $(c^\sharp,c)$-robustly $D$-universal} if whenever
\begin{itemize}
\item[(i)] $Y'_i \sub Y_i$ with $|Y'_i| \ge c^\sharp|Y_i|$ such that $Y'=\cup_{i\in R} Y'_i$,
$J'=J[Y']$ satisfy $|J'_S(v)| \ge c^\sharp|J_S(v)|$ whenever
$|S|=k$,  $J_S$ is defined, $i\in S$, $v\in J'_i$,
\item[(ii)] $H'$ is an $R$-indexed complex on $X'=\cup_{i\in R}X'_i$
of maximum degree at most $D$ with $|X'_i|=|Y'_i|$ for $i \in R$,
\item[(iii)] $X_* \sub X'$ with $|X_* \cap X'_i| \le c|X_i|$ for all $i \in R$,
and $\GG_x \sub Y'_x$ with $|\GG_x| \ge c^\sharp|Y'_x|$ for $x \in X_*$, \index{$\GG$}
\end{itemize}
then there is a bijection $\phi:X' \to V'$ with $\phi(X'_i)=V'_i$ for $i \in R$
such that $\phi(S) \in J_S$ for $S \in H'$ and $\phi(x) \in \GG_x$ for $x \in X_*$.
\end{defn}

More generally, one can allow restrictions to regular subcomplexes in both 
conditions (i) and (iii) of Definition \ref{def-robust},
but for simplicity we will not formulate the definition here.
As before, one can delete a small number of vertices from a regular complex
with a small number of marked $k$-tuples to obtain a robustly universal complex.
As for Theorem \ref{3robust}, the proof is immediate from Lemma \ref{super-restrict},
Definition \ref{def-robust} and Theorem \ref{blowup}.%
\COMMENT{
(i) $\GG' \sub J$, $J'=J[\GG]$, change $Y'_i$ to $\GG'_i$;
(iii) $\GG$ is $H'$-col in $J'$, $\GG_x$ only def for $x \in X_*$;
need another param for regularity: `$\eps$-regularly' [?]
}

\begin{theo}\label{robust}
Suppose
$0 < 1/n \ll 1/n_R \ll \eps \ll c \ll d^\sharp \ll d_a \ll \theta \ll c^\sharp, d, 1/k, 1/D_R, 1/D$,
we have a multi-$k$-complex $R$ on $[r]$ with maximum degree at most $D_R$ and $|R| \le n_R$,
$G$ is an $\eps$-regular $R$-indexed complex on $V = \cup_{i \in R} V_i$
with $n\le |V_i|=|G_i|\le Cn$ for $i\in R$,
$d_S(G) \ge d_a$ when $|S|\ge 2$ and $d_S(G) \ge d$ when $|S|=k$,
and $M \sub G_=$ with $|M_S| \le \theta|G_S|$ when defined.
Then we can delete at most $2\theta^{1/3}|G_i|$ vertices from $G_i$ for $i \in R$
to obtain $G^\sharp$ and $M^\sharp$ so that
\begin{itemize}
\item[(i)] $d(G^\sharp_S)>d^\sharp$ and $|G^\sharp_S(v)|>d^\sharp|G^\sharp_S|/|G^\sharp_i|$
whenever $|S|=k$, $G^\sharp_S$ is defined, $i\in S$, $v\in G^\sharp_i$, and
\item[(ii)] $G^\sharp\sm M^\sharp$ is $(c^\sharp,c)$-robustly $D$-universal.
\end{itemize}
\end{theo}

Finally, we mention that one can allow much smaller densities in the restricted
positions, provided that one makes an additional assumption to control the marking
edges. We can replace $c'$ by $d_a$ in condition (v) of Theorem \ref{blowup},
provided that we add the following additional assumptions:
\begin{itemize}
\item[(v.1)] $|M_S(v)| \le \theta|\GG_S(v)|$ when $|S|=k$, $\GG_S$ is defined, $i \in S$ and $v \in \GG_i$,
\item[(v.2)] $|(M \cap G[\GG]^{I_v})_S| \le \theta|G[\GG]^{I_v}_S|$ for any submulticomplex $I$ of $R$,
when $|S|=k$, $v \in G_i$ and $S \cap VN_R(i) \ne \es$.
\end{itemize}
Note that these conditions ensure that the marked edges are controlled in $G[\GG]$
exactly as in conditions (ii) and (iii) of super-regularity, so the proof goes through as before.
In this general form there is no simplification to be gained by reformulating
the statement in a black box form. We suppressed this refined form in the statement
of Theorem \ref{blowup} to avoid overburdening the reader with technicalities,
but we note that it may be needed in some applications. Indeed, one may well have to generate restricted
positions using neighbourhood complexes in $G$, and then $c'$ will be of the order of the densities in $G$.

\subsection{Concluding remarks}

The theory of regularity and super-regularity for hypergraphs is
considerably more involved than that for graphs.
As explained in Section \ref{3reg}, these technicalities cannot be avoided,
but the black box reformulation in Lemma \ref{robust} should make the
hypergraph blow-up lemma more convenient for future applications.
The graph blow-up lemma has had many applications in modern graph theory,
so it is natural to look for hypergraph generalisations of these results.
However, many such applications build on basic results for graphs for which the
hypergraph analogue is unknown.
For example, in our application in Section \ref{apps} we only needed a matching,
and were able to rely on Kahn's matching theorem, which is already quite a difficult result.
Thus one may expect it will take longer for the hypergraph blow-up lemma to achieve its
full potential.

Another question for future research is to obtain an algorithmic version of our theorem,
along the lines of the algorithmic graph blow-up lemma in \cite{KSS4}. In
applications this could be combined with an algorithmic version of hypergraph regularity
given by \cite{DHNR}. A rather different direction of research would be along the
lines of the `infinitary' versions of hypergraph regularity theory, whether
probabilistic \cite{Au, T1}, analytic \cite{LS}, model theoretic \cite{ES} or algebraic \cite{FLS,R}.
It is natural to ask whether the blow-up lemma has an interpretation
in any of these frameworks.

Further refinements could include estimating the number of embeddings,
rather than just proving the existence of a single embedding as in this paper.
Here it may be helpful to note that one can combine Lemma \ref{x-to-v} with martingale
estimates to show that with high probability there will be at least $\frac{1}{2}p|B \cap X_v|$
vertices $x \in B \cap X_v$ such that $\phi(H(x)) \sub (G \sm M)(v)$.
We also note that small improvements to the tail decay of our martingales
may be obtained from the Optional Sampling Theorem (see e.g.\ \cite{GS} p.\ 462).
One could also try to obtain (nearly) perfect edge-decompositions of super-regular complexes
into copies of a given bounded degree hypergraph. For example, one could ask for hypergraph
generalisations of a result of Frieze and Krivelevich \cite{FK} that one can cover almost all edges
of an $\eps$-regular graph by edge-disjoint Hamilton cycles.%
\footnote{Note added during revisions: Frieze and Krivelevich now have a preprint on this topic.}

\nib{Acknowledgements.}
I thank Daniela K\"uhn, Deryk Osthus, Richard Mycroft and Mathias Schacht
for drawing my attention to a significant oversight in an earlier version of this paper.
That version contained a weaker version of the blow-up lemma, which assumed
a dense environment that cannot be guaranteed in many applications. In particular
it was not adequate for the application that I give in this paper.
I also thank the anonymous referees, Daniela K\"uhn, Deryk Osthus and Richard Mycroft
for carefully reading and commenting extensively on various versions of this paper,
and Mathias Schacht for explaining \cite{RSc1, RSc2} to me.
\COMMENT{
to do:
1. robust universality:
(i) allow lower $c^\sharp$ with marking condition (`strongly'?)
(ii) allow complex restriction, both for res pos and in def conditions
2. more pictures?
3. motivate generalities (res pos, indexed complexes)
4. remarks on reduced complex
5. notation index
}

\printindex

\end{document}